\documentclass [12pt] {article}
\title {Quantum one-cocycles for knots}
\author{Thomas Fiedler}
\usepackage{amsfonts, amssymb, latexsym, graphicx}
\begin{document}
\newtheorem{proposition}{Proposition}
\newtheorem{theorem}{Theorem}
\newtheorem{lemma}{Lemma}
\newtheorem{corollary}{Corollary}
\newtheorem{example}{Example}
\newtheorem{remark}{Remark}
\newtheorem{definition}{Definition}
\newtheorem{question}{Question}
\newtheorem{conjecture}{Conjecture}
\newtheorem{observation}{Observation}
\maketitle

\begin{center}
{\em pour S\'everine}
\end{center}

\begin{abstract}
We give a method to construct  non symmetric solutions of a global tetrahedron equation from solutions of the Yang-Baxter equation. The solution in the HOMFLYPT case  gives rise to the first combinatorial quantum 1-cocycle  which represents a  non trivial cohomology class in the topological moduli space of long knots. We conjecture that the quotient of its values on Hatchers loop and on the rotation around the long axis  is related to the simplicial volume of the knot complement in the 3-sphere and we prove this for the figure eight knot.

Surprisingly, the formula for the solution in the HOMFLYPT case of the positive global tetrahedron equation gives also a solution in the case of the 2-variable Kauffman invariant. But there is also a second solution giving rise to yet another non trivial quantum 1-cocycle.

\footnote{2000 {\em Mathematics Subject Classification\/}: 57M25 {\em Keywords\/}:
quantum 1-cocycles, simplicial volume, invariants of string links, global tetrahedron equation, cube equations.}

\end{abstract}  

\tableofcontents
\section{Introduction}
\subsection{Results}
There are two important directions in 3-dimensional topology: quantum invariants for knots from irreducible representations of quantum groups on one hand and the geometry of 3-manifolds via geometrisation on the other hand. It seems that the main  connections  between them are the volume conjecture, which generalizes Kashaevs conjecture for hyperbolic knots (see \cite{Mu} and \cite{Ka}), the AJ-conjecture (see \cite{G} and \cite{Le}) and the relations, via certain incompressible surfaces in a link complement, of the extremal coefficients of the n-colored Jones polynomial with the geometric structure of the link complement (see \cite{FKP}, \cite{FKP2}, \cite{FKP3}).

In this paper we propose a potential new connection. Quantum knot invariants can be naturally seen as combinatorial 0-cocycles in the {\em topological moduli space of smooth knots} (i.e. the space of all smooth knots isotopic to a given one). Our main idea is to go one dimension higher. We construct non trivial combinatorial 1-cocycles in the corresponding topological moduli spaces of long knots and of tangles without closed components.

The connection to geometry is based on a result of Hatcher. It is well known that the classification problem for knots is equivalent to the classification problem of long knots, i.e. a smoothly embedded arcs in 3-space which go to infinity outside a compact set  as a straight line. Let $K$ be a framed long knot and let $M_K$ be the topological moduli space of the unframed knot $K$. There are two natural loops in $M_K$: the rotation of $K$ around the long axis, which we call $rot(K)$. (This loop is often called {\em Gramain's cycle} because it has appeared first in \cite{Gr}.) Notice that this loop does not depend on the framing of $K$. The other loop, which we call {\em Hatchers loop} $hat(K)$ (compare \cite{H}) is defined as follows: one puts a pearl (i.e. a small 3-ball $B$ with a frame at the origin, where one vector is tangent to the oriented knot) on the (closure of the) framed knot $K$ in the 3-sphere. The part of $K$ in $S^3 \setminus B$ is a long knot. Pushing $B$ once along the knot following the framing induces Hatchers loop in $M_K$. Changing the framing of $K$ by +1 adds the homology class  $[rot(K)]$ to $[hat(K)]$. The following theorem is an immediate consequence of a result of Hatcher (see \cite{H}).

\begin {theorem}\textbf{(Hatcher)}

Let $K$ be a  long knot with trivial framing (standard near infinity) and which is not a satellite (i.e. there is no incompressible torus in its complement in the 3-sphere).
Then  $hat(K)$ represents a non zero integer multiple of $rot(K)$ in $H_1(M_K; \mathbb{Z})$ if and only if $K$ is a non trivial torus knot.

\end {theorem}

Hatcher has proven that for the standard diagram of a long $(p,q)$-torus knot $K$ the loops $rot(K)$ and $hat(K)$ are homotopic.
The standard diagram of a positive torus knot has blackboard framing $p(q-1)$ (we assume $p>q$) and consequently for $K$ with trivial framing  $[hat(K)]=(1-p(q-1)) [rot(K)]$ in $H_1(M_K; \mathbb{Z})$. Hence, $[hat(K)]$ is a non zero multiple of $[rot(K)]$ if $K$ is a non trivial torus knot.

The hard direction of this theorem is of course to prove that for a hyperbolic knot these two loops are linearly independent. This follows from the minimal models for the topological moduli spaces of hyperbolic and of torus knots which were constructed by Hatcher (see \cite{Bu} and \cite{BC} for the case of satellites).
Hatchers construction is mainly based on very deep results in 3-dimensional topology: the Smith conjecture, the Smale conjecture, the Linearization conjecture (which is a consequence of the spherical case in Perelmann's work) and the result of Gordon-Luecke that a knot in the 3-sphere is determined by its complement.

In this paper we construct a combinatorial 1-cocycle for $M_K$ which is  based on the HOMFLYPT invariant, see Theorem 4 in Section 11. It is called $\bar R^{(1)}$  ("R" stands for Reidemeister).

\begin{figure}
\centering
\includegraphics{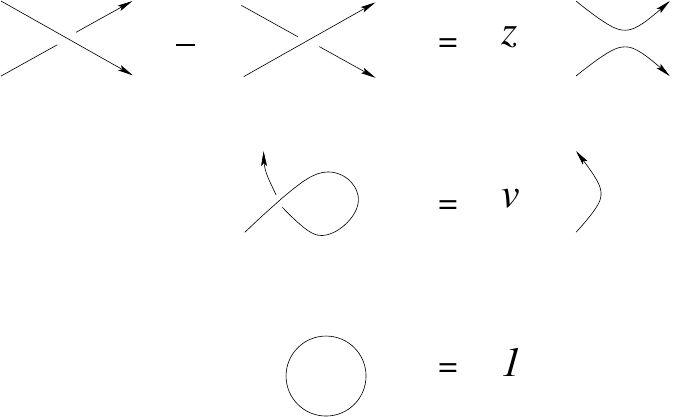}
\caption{\label{HOMFLYPT}  skein relations for the HOMFLY polynomial}  
\end{figure}

The HOMFLYPT polynomial is defined by the skein relations shown in Fig.~\ref{HOMFLYPT}. It is normalized as usual by $v^{-w(D)}$, where $w(D)$ is the writhe of the diagram $D$ of the knot $K$.

\begin{conjecture}
Let $K$ be a long knot, let $v_2(K)$ be its Vassiliev invariant of degree 2 and let $P_K$ be its HOMFLYPT polynomial. Let $\delta$ denote as usual the HOMFLYPT polynomial of the trivial 2-component link. Then

$\bar R^{(1)}(rot(K))=-\delta v_2(K) P_K $. 

\end{conjecture}

We prove this conjecture for the trivial knot, the trefoil and the figure eight knot (Section 9, Example 3 but compare also Remark 14 in Section 11).

On the other hand let $K$ be the figure eight knot with trivial framing. A calculation by hand gives  $\bar R^{(1)}(hat(K)) = 0$.
Consequently, we have re-proven that for the figure eight knot $[hat(K)]$ is not a non zero multiple of $[rot(K)]$ (and hence in particular the figure eight knot  is  hyperbolic) by using only quantum topology instead of deep results in 3-dimensional topology!

\begin{conjecture}

Let $K$ be a long knot with trivial framing and with non trivial Vassiliev invariant $v_2(K)$. Then  
$\bar R^{(1)}(hat(K))$ is a non zero integer multiple of  $\bar R^{(1)}(rot(K))$ if and only if $K$ is a torus knot. 

\end{conjecture}

The conjecture is true for the trefoils and for the figure eight knot (Section 9, Example 5).
It could be possible that $\bar R^{(1)}(hat(K))/ \bar R^{(1)}(rot(K))$ contains information about the hyperbolic volume in the case of a hyperbolic knot $K$. Unfortunately, there is not yet a computer program in order to calculate lots of examples and to formulate a more precise conjecture.

We construct also another 1-cocycle, called $R^{(1)}_{reg}$, which is well defined only for regular isotopies (i.e. Reidemeister moves of type I are not allowed, or in other words the isotopy preserves the blackboard framing). The corresponding topological moduli space is called $M_K^{reg}$.

In fact, our 1-cocycle $\bar R^{(1)}$ is the result of gluing together two 1-cocycles which are of a completely different nature. First we construct a quantum 1-cocycle $R^{(1)}$ (i.e. its construction uses skein relations of quantum invariants) which is well defined only in the complement of certain strata of codimension two in $M_K$ which correspond to diagrams which have in the projection to the plane a branch which passes transversally through a cusp. Luckily, the value of $R^{(1)}$ on the meridians of these strata is just an integer multiple of $\delta P_K$.
Next we construct an integer valued finite type Gauss diagram 1-cocycle $V^{(1)}$ (i.e. it is constructed by Gauss diagram formulas of finite degrees) in the complement of exactly the same strata of codimension two. We define then the completion $\bar R^{(1)}=R^{(1)}-\delta P_K V^{(1)}$ which turns out to be a 1-cocycle in the whole topological moduli space $M_K$.

All quantum invariants verify {\em skein relations}  and hence give rise to {\em skein modules}, i.e. formal linear combinations over some ring of coefficients of tangles $T$ (with the same oriented boundary $\partial T$) modulo the skein relations. The multiplication of tangles is defined by their composition (if possible). The main property of quantum invariants is their {\em multiplicativity}.
We denote the HOMFLYPT skein module associated to the oriented boundary of a tangle $T$ by $S(\partial T)$.

We will consider only a special class of tangles.
A {\em n-string link T} are n ordered oriented properly and smoothly embedded arcs in the 3-ball. They are considered up to ambient isotopy which is the identity on the boundary of the 3-ball. They generalize long knots. We chose an abstract closure of the string link to a circle.

 {\em So, a tangle is for us a knot with only a part of it embedded in 3-space.}

 It is well known that if the complement of the string link $T$ does not contain an incompressible torus then the topological moduli space $M_T$ is a contractible space. This is always the case for braids (see e.g. \cite{B}).

Hence there aren't any non trivial cohomology classes in this case. However our 1-cocycles $\bar R^{(1)}$ and $R^{(1)}_{reg}$ have a remarkable property, which we call the {\em scan-property}. We fix an orthogonal projection of the 3-ball into a disc, such that the string link $T$ is represented by a generic diagram. We fix an arbitrary abstract closure $\sigma$ of $T$ to an oriented circle, we fix a {\em point at infinity} in $\partial T$. Let $R^{(1)}$ be a combinatorial 1-cocycle in $M_T^{reg}$ or in $M_T$ (i.e. we sum up contributions from Reidemeister moves).

\begin{definition}
The 1-cocycle $R^{(1)}$ has the {\em scan-property} if the contribution of each Reidemeister move $t$ doesn't change when a branch of $T$ has moved under the Reidemeister move $t$ to the other side of it.
\end{definition}

Let us add a small positive curl to an arbitrary component of $T$ near to the boundary $\partial T$ of $T$. 

\begin{figure}
\centering
\includegraphics{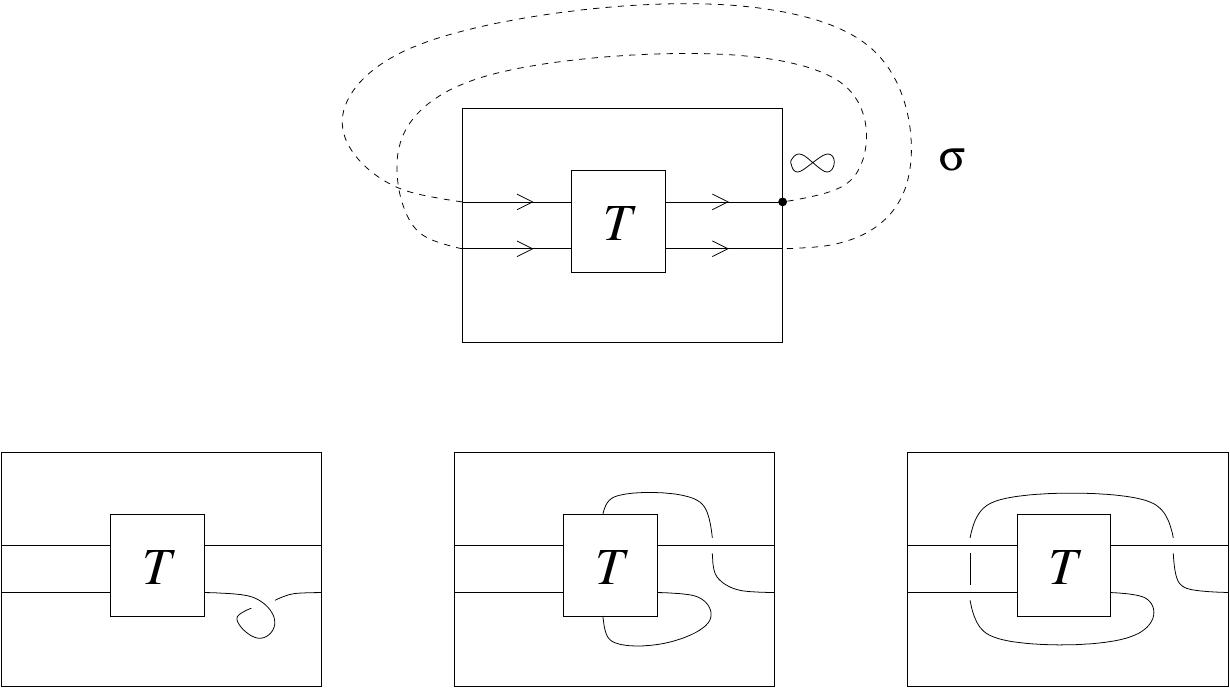}
\caption{\label{scan}  scan-arc for a tangle $T$}  
\end{figure}

\begin{definition}
 The {\em scan-arc scan(T)} in $M^{reg}_T$ is the regular isotopy which makes the small curl big under the rest of $T$ up to being near to the whole boundary of the disc, compare Fig.~\ref{scan} (it is convenient to replace the ball by the cube). 
\end{definition}

Each point in $\partial T$ comes with a sign as a component of the boundary of the oriented tangle $T$. We consider all cyclically ordered (by $\sigma \cup  T$) subsets A in $\partial T$ with alternating signs and which contain the point at infinity. The 1-cocycle $R^{(1)}_{reg}$ will be {\em graded} by the sets $A$, in contrast to the 1-cocycle $\bar R^{(1)}$ which can not be graded.

\begin{theorem}

$\bar R^{(1)}(scan(T))$ and $R^{(1)}_{reg}(A)(scan(T))$ with values in the HOMFLYPT skein module of $\partial T$ do not depend on the chosen regular diagram  of $T$   (but they depend on the chosen closure, on the chosen component for the curl, of the choice of $\infty$ in $\partial T$ and of the grading $A$ in the case of $R^{(1)}_{reg}(A)(scan(T))$)
and they are  consequently  invariants of regular isotopy for the string link $T$.

\end{theorem}

This means that the evaluation of the 1-cocycles on some canonical arc (instead of a loop) are already  invariants. This is a consequence of the scan-property. In other words, the 1-cocycles $\bar R^{(1)}$ and $R^{(1)}_{reg}$ can be useful even when they represent the trivial cohomology class. Notice that the scan-property is a priori a property of the cocycle and not of the corresponding cohomology class.
(There  are "dual" 1-cocycles  which have the scan-property for small curls which become big over everything instead of under everything. It could well be that they represent the same cohomolgy class as $\bar R^{(1)}$ respectively $ R^{(1)}_{reg}$.)

Let us come back to the case of long knots. Let $flip$ denote a rotation in 3-space around a vertical axis in the plan (i.e. it interchanges the two ends of a long knot) followed by an orientation reversing. It is easy to see that a long knot $K$ is isotopic to the long knot $flip(K)$ if and only if the closure of $K$ is an invertible knot in the 3-sphere (see e.g. \cite{F1}, \cite{DK}).

If the long oriented knot is framed, i.e. a trivialization (standard at infinity) of its normal bundle is chosen, then we can replace the knot $K$ by n parallel oriented copies with respect to this framing. The result is a n-string link which is called the non-twisted n-cable $Cab_n(K)$ if the framing was zero. The rotation around the long axis is no longer possible because the end points have to stay fixed. However, the loop $hat(Cab_n(K))$ is still well defined. Indeed, we put a pearl B on the braid closure of $Cab_n(K)$ in $S^3$ in such a way that in its interior we have n almost straight parallel lines and we push it once along the closed $Cab_n(K)$. The result is a loop in $M_{Cab_n(K)}$. We fix now an abstract closure of $Cab_n(K)$ to a circle. Consequently, $\bar R^{(1)}(hat(Cab_n(K)))$ is well defined.  It takes its values in a skein module which has basis elements which are not invariant under $flip$ starting from $n=3$ (e.g. the permutation 3-braids $\sigma_1 \sigma_2 \not= \sigma_2 \sigma_1$). It would be very interesting to test if 
$\bar R^{(1)}(hat(Cab_3(K)))$ can detect the non-invertibility of a knot $K$. The same question is interesting for $\bar R^{(1)}(scan(Cab_3(K)))$ and for $R^{(1)}_{reg}(scan(Cab_3(K)))$ too (where calculations are much simpler because of the grading).
Most of the ingredients in our construction are not invariant under {\em flip} and hence at least apriori the 1-cocycles $\bar R^{(1)}$ and $R^{(1)}_{reg}$ are not invariant under {\em flip}. Unfortunately (for me), however calculations have to be done with a computer program.

(Notice that there are finite type invariants which can sometimes distinguish the orientation of 2-string links, see \cite{DK}.)

{\em The most important properties of our 1-cocycle $\bar R^{(1)}$ are the facts  that it represents a non trivial cohomology class and that it has the scan-property. Hence it gives a completely new sort of calculable  invariants for tangles.}

\begin{figure}
\centering
\includegraphics{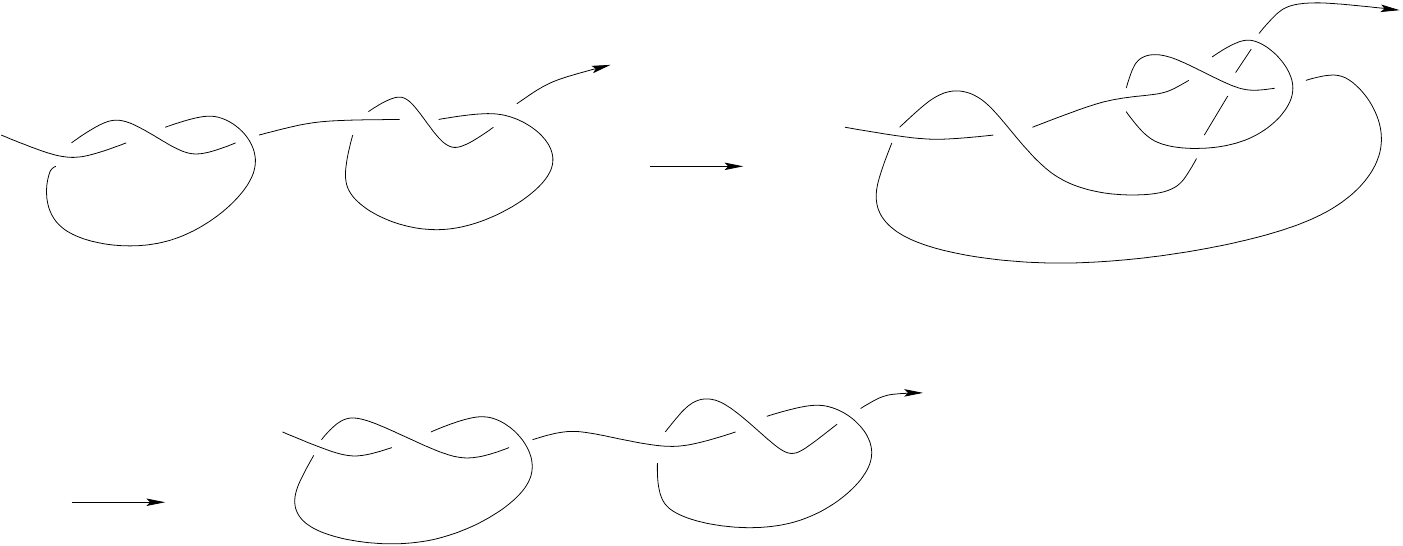}
\caption{\label{drag}  dragging a trefoil through another trefoil}  
\end{figure}

Our 1-cocycle $\bar R^{(1)}$ allows sometimes to answer (in the negative) a difficult question: are two given complicated  isotopies of long knots homologous in $M_K$? We don't even need to find out first whether these are isotopies of the same knot.
For example, $rot(K)$ and $hat(K)$ are not homologous for the figure eight knot.
Let us take the product of two long standard diagrams of the right trefoil. We consider the dragging of one through the other, denoted by $drag 3_1^+$, see Fig.~\ref{drag}.  A calculation by hand gives
 $\bar R^{(1)}(drag 3_1^+)=-3\delta (P_{3_1^+})^2$ (Section 9, Example 4).  

Let $D$ be the standard diagram of the knot  $3_1^+$. Our calculation of  $\bar R^{(1)}(rot (3_1^+))$ implies easily that  $\bar R^{(1)}(rot (3_1^+ \sharp 3_1^+)) =-2\delta (P_{3_1^+})^2$. Consequently, we have re-proven that $drag D$ and $rot (3_1^+ \sharp 3_1^+)$ are not homologous for the connected sum of two right trefoils (compare \cite{H} and \cite{Bu}).  Notice that we could have altered these loops by lots of unnecessary Reidemeister moves. {\em There isn't known any algorithm to detect the homology class of a loop in $M_K$.}

More generally, let $drag(T,T')$ be the loop which consists of dragging the framed tangle $T$ through the framed tangle  $T'$ (if possible) followed by dragging $T'$ through $T$. For example $drag(Cab_n(K), Cab_n(K'))$ is well defined for each couple of framed long knots $K$ and $K'$. This is still a regular loop and we can calculate $\bar R^{(1)}$, $R^{(1)}_{reg}$ and  $V^{(1)}$ for it. These are invariants of the couple $(K,K')$ for each abstract closure of $Cab_n(K) \sharp Cab_n(K')$ to a circle and for each choice of a point at infinity in $\partial Cab_n(K) \sharp Cab_n(K')$. The cocycles $\bar R^{(1)}$ and $R^{(1)}_{reg}$ take their values in the skein module $S(\partial Cab_n(K) \sharp Cab_n(K'))$ and $V^{(1)}$ is just an integer. Are these new knot invariants?
\vspace{0,2 cm}

Easily calculable finite type 1-cocycles of all degrees for connected closed braids were found in the unpublished paper \cite{F2}. It was the starting point of our present work. As well known quantum invariants can be decomposed into finite type invariants (see \cite{BN}). Are there non trivial quantum 1-cocycles which can perhaps be decomposed into finite type 1-cocycles?

\begin{remark}
 Our example for "$ hat / rot$" gives  0  for the figure eight knot . On the other hand $\bar R^{(1)}(hat(K)) / \bar R^{(1)}(rot(K))$ is a non zero constant for the trefoils (which have simplicial volume 0, see e.g. \cite{Mu}). 
Is the "0" for the figure eight perhaps related to the fact that the figure eight is the knot with the smallest hyperbolic volume amongst all hyperbolic knots (see \cite{Th})? 
\end{remark}

\begin{remark}
In all our examples the value of $\bar R^{(1)}$ on a loop contains the HOMFLYPT invariant of the knot as a factor. We don't know whether this happens in general. Let e.g. $K$ be a satellite and let $V$ be the non trivially embedded solid torus in the 3-sphere which contains $K$. The knot $K$ is not isotopic to the core $C$ of the solid torus $V$. We can assume that $K$ is tangential to $C$ at some point and we take this point as the point at infinity in the 3-sphere. The rotation of $V$ around its core $C$ induces now a loop for the long knot $K$. We call this loop $rot_C(K)$. It would be very interesting to calculate $\bar R^{(1)}(rot_C(K))$ in examples.

However,  the value of $\bar R^{(1)}$ on scan-arcs instead of loops does in general not contain the HOMFLYPT invariant as a factor (compare Examples 1 and 2 in Section 8)! 
\end{remark}

\begin{question}
In particular, does $\bar R^{(1)}(hat(Cab_2(K)))$ contain $P_{Cab_2(K)}$ in the Hecke algebra $H_2(z,v)$ as a factor for each choice of a point at infinity in $\partial Cab_2(K)$?
\end{question}

Let us conclude this part of the introduction by pointing out some essential differences between invariants which come from our quantum 1-cocycles and the usual quantum invariants (coming from 0-cocycles). Most of the usual quantum invariants (if not all) can be extended to links and arbitrary tangles (see \cite{Tu}), virtual knots (see \cite{K2}) and singular knots (see \cite{KV}). This is {\em not} the case for the invariants coming from 1-cocycles!
We use in an essential way that the string link has an abstract closure to a circle. So it is not clear at all how to generalize it to links. The loops $rot$ and $hat$ would always contain forbidden moves for a long virtual knot. The loop $rot$ is still defined for a long singular knot, but the loop $hat$ is no longer defined.

Very simple examples show that the invariants from 1-cocycles are no longer multiplicative at least for string links with more than one component (Section 8, Example 1). This seems to be the price to pay for invariants which make connections to geometry.

\subsection{Method}

\begin{figure}
\centering
\includegraphics{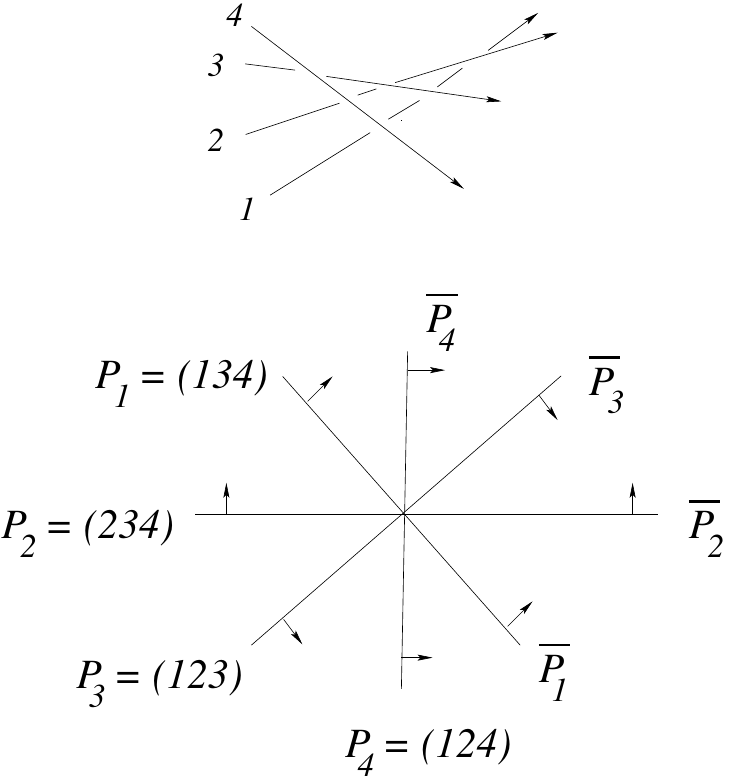}
\caption{\label{quadcros}  intersection of a normal $2$-disc of a positive quadruple crossing with the strata of triple crossings}  
\end{figure}

The construction of $\bar R^{(1)}$ is based on a completely new solution of some tetrahedron equation. The classical tetrahedron equation is a 3-dimensional generalization of the Yang-Baxter equation. It is a local equation which is fundamental for studying integrable models in $2+1$-dimensional mathematical physics. This equation has many solutions and the first one was found by Zamolodchikov (see \cite{Kor}, \cite{Ku} \cite{KaKo}). We look at the tetrahedron equation from the point of view of singularities of projections of lines. We fix an orthogonal projection $pr: I^3 \rightarrow I^2$. Consider four oriented straight lines in the cube $I^3$ which form a braid and for which the intersection of their projection into the square $I^2$ consists of a single point. We call this an {\em ordinary quadruple crossing}. After a generic perturbation of the four lines we will see now exactly six ordinary crossings. We assume that all six crossings are positive and we call the corresponding quadruple crossing a positive quadruple crossing. Quadruple crossings form smooth strata of codimension 2 in  the topological moduli space of lines in 3-space which is equipped with a fixed projection $pr$. Each generic point in such a stratum is adjacent to exactly eight smooth strata of codimension 1. Each of them corresponds to configurations of lines which have exactly one ordinary triple crossing besides the remaining ordinary crossings. We number the lines from 1 to 4 from the lowest to the highest (with respect to the projection $pr$).
The eight strata of triple crossings glue pairwise together to form four smooth strata which intersect pairwise transversally in the stratum of the quadruple crossing (see e.g. \cite{FK}). The strata of triple crossings are determined by the names of the three lines which give the triple crossing. For shorter writing we give them names from $P_1$ to $P_4$ and $\bar P_1$ to $\bar P_4$ for the corresponding stratum on the other side of the quadruple crossing.  We show the intersection of a normal 2-disc of the stratum of codimension 2 of a positive quadruple crossing with the strata of codimension 1 in Fig.~\ref{quadcros}. We give in the figure also the coorientation of the strata of codimension 1, which will be defined below in Definition 1. (We could interpret the six ordinary crossings as the edges of a tetrahedron and the four triple crossings likewise as the vertices's or the 2-faces of the tetrahedron.) 

\begin{figure}
\centering
\includegraphics{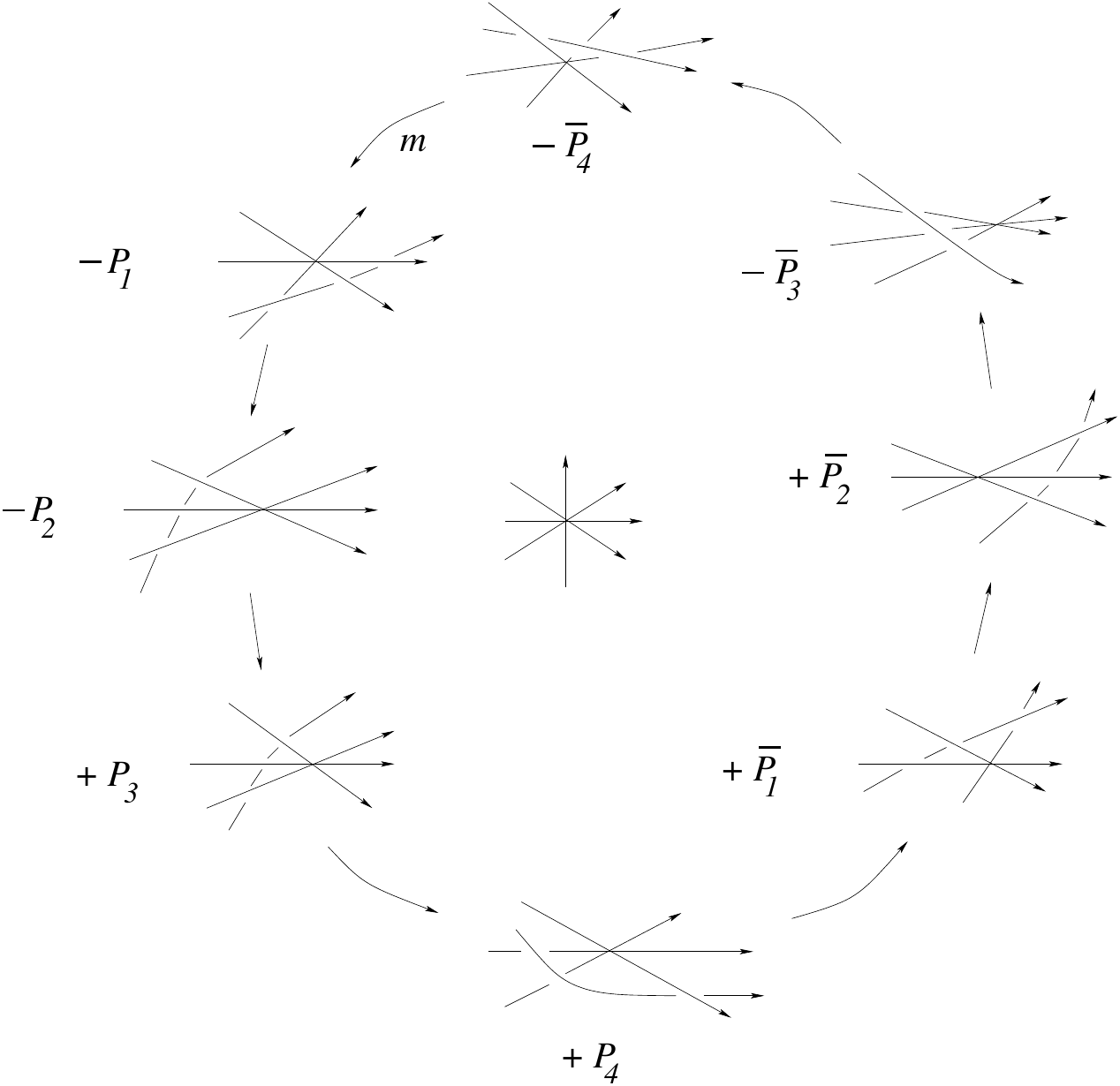}
\caption{\label{unfoldquad}  unfolding of a positive quadruple crossing}  
\end{figure}

Let us consider a small circle in the normal 2-disc and which goes once around the stratum of the quadruple crossing. We call it a {\em meridian} $m$ of the quadruple crossing. Going along the meridian $m$ we see ordinary diagrams of 4-braids and exactly eight diagrams with an ordinary triple crossing. We show this in Fig.~\ref{unfoldquad}. (For simplicity we have drawn the triple crossings as triple points, but the branches do never intersect.) For the classical tetrahedron equation one associates to each stratum
 $P_i$  some operator (or some R-matrix) which depends  only on the names of the three lines and to each stratum $\bar P_i$ the inverse operator. The tetrahedron equation says now that if we go along the meridian $m$ then the product of these operators is equal to the identity. Notice, that in the literature (see e.g. \cite{KaKo}) one considers planar configurations of lines. But this is of course equivalent to our situation because all crossings are positive and hence the lift of the lines into 3-space is determined by the planar picture. Moreover, each move of the lines in the plane which preserves the transversality lifts to an isotopy of the lines in 3-space.

However, the solutions of the classical tetrahedron equation are not well adapted in order to construct 1-cocycles for moduli spaces of knots. There is no natural way to give names to the three branches of a triple crossing in an arbitrary knot isotopy besides in the case of braids. But it is not hard to see that in the case of braids Markov moves  would make big trouble (see e.g. \cite{B} for the definition of Markov moves). As well known a Markov move leads only to a normalization factor in the construction of 0-cocycles (see e.g. \cite{Tu}). However, the place in the diagram and the moment in the isotopy of a Markov move become important  in the construction of 1-cocycles (as already indicated by the non-triviality of  Markov's theorem). To overcome this difficulty we search for a solution of the tetrahedron equation which does not use names of the branches.  

{\em The idea is to associate to a triple crossing $p$ of a diagram of a string link $T$, a global element in the skein module $S(\partial T)$ instead of a local operator.}

\begin{figure}
\centering
\includegraphics{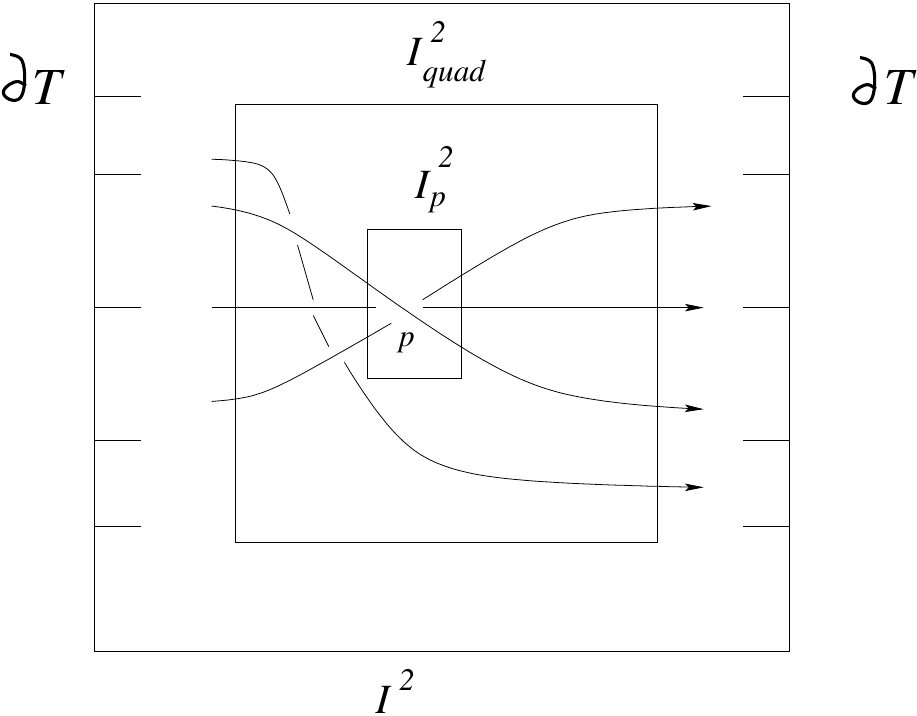}
\caption{\label{increascubes} increasing cubes}  
\end{figure}

 Let us consider a diagram of a string link with a positive quadruple crossing $quad$ and let $p$ be e.g. the positive triple crossing associated to the adjacent stratum $P_1$. We consider three cubes: $I^2_p \times I \subset I^2_{quad} \times I \subset I^2 \times I$. Here $I^2_p$ is a small square which contains the triple point $pr(p)$ and $I^2_{quad}$ is a slightly bigger square which contains the quadruple point $pr(quad)$. We illustrate this in Fig.~\ref{increascubes}. Let us replace the triple crossing $p$ in $I^2_p \times I$ by some standard element, say $L(p)$, in the skein module $S( 3-braids)$.
Gluing to $L(p)$ the rest of the tangle $T$ outside of $I^2_p \times I$ gives an element, say $T(p)$, in the skein module $S(\partial T)$. A triple crossing in an oriented isotopy corresponds to a Reidemeister move of type III (see e.g. \cite{BZ} for the definition of Reidemeister moves). We use the standard notations for braid groups (see e.g. \cite{B}).

\begin {definition}
The Reidemeister move $\sigma_1 \sigma_2 \sigma_1 \rightarrow \sigma_2 \sigma_1 \sigma_2$ has $sign= +1$ and its inverse has $sign =-1$.

\end{definition}

Let $m$ be the oriented meridian of the positive quadruple crossing. Our tetrahedron equation in $S(\partial T)$ is now \vspace{0,2 cm}

(1)      $\sum_{p \in m} sign(p) T(p)=0$    \vspace{0,2 cm}                                                                                            

where the sum is over all eight triple crossings in $m$.

It follows immediately from the definition of a skein module and from the fact that the tangle $T$ is the same for all eight triple crossings outside of $I^2_{quad} \times I$ that the equation (1) reduces to the local tetrahedron equation (2) in the skein module  $S(4-braids)$: \vspace{0,2 cm}

(2)      $\sum_{p \in m} sign(p) L(p)=0$    \vspace{0,2 cm}

(Of course we have to add here to each tangle in $L(p)$ the corresponding fourth strand in  $I^2_{quad} \times I$, but we denote it still by $L(p)$.)  

\begin{proposition}

The equation (2) has a unique solution (up to a common factor), which is $L(p)= \sigma_2(p)-\sigma_1(p)$.

\end{proposition}

We call this solution the {\em solution of constant weight} (the name will become clear later).

{\em Proof.} We prove only the existence. The unicity is less important and we left the proof to the reader. In general, for each permutation 3-braid $\beta$ we will write $\beta(p)$ for the element in the skein module, which is obtained by replacing the triple crossing by $\beta$.
One has to use that the skein module of 3-braids is generated by the permutation braids $1$, $\sigma_1$, $\sigma_2$, $\sigma_1 \sigma_2$,  $\sigma_2 \sigma_1$ and $\sigma_1 \sigma_2 \sigma_1$. The latter permutation braid is not interesting, because it is always just the original tangle. We call  $\sigma_1$, $\sigma_2$ the {\em partial smoothings with one crossing} and $\sigma_1 \sigma_2$, $\sigma_2 \sigma_1$ {\em partial smoothings with two crossings}. The Equation (2) reduces now to a finite number of equations which have a single solution.

The signs of the triple crossings in the meridian are shown in Fig.~\ref{unfoldquad} too.

First of all we observe that for every basis element of the skein module of 3-braids the contribution of the stratum $P_2$ cancels out with that of the stratum $\bar P_2$. The same is true for $P_3$ and $\bar P_3$. This comes from the fact that the signs are opposite but the resulting tangles are isotopic (there is just one branch which moves over or under everything else). Hence, we are left only with the strata $P_1$, $P_4$, $\bar P_1$, $\bar P_4$.

The substitution of the solution of constant weight into equation (2) gives now:

$-\sigma_1(\sigma_3-\sigma_2)\sigma_1\sigma_2+\sigma_2\sigma_1(\sigma_3-\sigma_2)\sigma_1
+\sigma_2\sigma_3(\sigma_2-\sigma_1)\sigma_3-\sigma_3(\sigma_2-\sigma_1)\sigma_3\sigma_2
=(\sigma_2\sigma_1+z\sigma_1\sigma_2\sigma_1-\sigma_3\sigma_2-z\sigma_1\sigma_3\sigma_2)
-(\sigma_1\sigma_2+z\sigma_1\sigma_2\sigma_1-\sigma_2\sigma_3-z\sigma_2\sigma_1\sigma_3)
-(\sigma_2\sigma_1+z\sigma_2\sigma_1\sigma_3-\sigma_3\sigma_2-z\sigma_3\sigma_2\sigma_3)
+(\sigma_1\sigma_2+z\sigma_1\sigma_3\sigma_2-\sigma_2\sigma_3-z\sigma_2\sigma_3\sigma_2)=0$

$\Box$

In order to construct a 1-cocycle from a solution of the tetrahedron equation we have to solve now in addition the {\em cube equations}. This will be explained later.

This was the point which we have  reached in 2008, see \cite{BNF} and also \cite{F3} for an earlier unsuccessful attempt (the one dimensional cohomology classes are all trivial) using state sum models.
Karen Chu (personal communication) has observed that the 1-cocycle which is obtained from the solution of constant weight satisfies the skein relations of the HOMFLYPT polynomial, when it is evaluated on $scan K$ for a long knot $K$. Consequently, this 1-cocycle is not interesting and in fact it represents the trivial cohomology class.

Notice, that from the six basis elements of $S(3-braids)$ we have used only two. Moreover, we haven't yet used at all the point at infinity of a long knot. But it is of crucial importance to use the point at infinity. This follows from another result of Hatcher. He has proven (see \cite{H}) that the fundamental group of the topological moduli space $M_{K \subset S^3}$ is finite for a prime knot $K$ in the 3-sphere (instead of long knots). Consequently, there can't be any non trivial 1-dimensional cohomology classes with values in a torsion free ring in this case. Hence, in order to define a non trivial 1-cocycle for prime knots it is essential to consider long knots instead of knots in the 3-sphere!

 It was only in spring 2012 that we have finally understood the complicated combinatorics which is needed in order to find a solution of the tetrahedron equation which uses the point at infinity and which uses exactly one of the only  non symmetric (under {\em flip}) partial smoothings $\sigma_1\sigma_2$ and $\sigma_2\sigma_1$.

Obviously, the contribution of $\sigma_2\sigma_1$ from $-P_2$ cancels out with that from $\bar P_2$ as well as the contribution from $P_3$ with that from $-\bar P_3$.
Let us make a table of the contributions to $S(4-braids)$ of the partial smoothing $\sigma_2\sigma_1$ from the remaining four adjacent strata to a positive quadruple crossing together with the contribution $zL(p)$ of $P_3$. It is convenient to give them as planar pictures in Fig.~\ref{tablesmooth}, using the fact that all crossings are positive.
Looking at Fig.~\ref{tablesmooth} we make an astonishing discovery: The contribution of $\sigma_2\sigma_1$ cancels out in $P_4$ with that from $-\bar P_4$. Its contribution in $-P_1$ and $\bar P_1$ would cancel out with the contribution of $zL(p)$ in $P_3$ and $-\bar P_3$ if the latter triple crossings would count with some coefficient $W$ such that $W(P_3)=W(\bar P_3)-1$. (Of course we have a symmetric result for $\sigma_1\sigma_2$ by interchanging $P_4$ with $P_1$ and $P_3$ with $P_2$.)

So, we are looking for an integer $W(p)$ which is defined for each diagram of a long knot with a positive triple crossing $p$ and which has the following properties:

\begin{figure}
\centering
\includegraphics{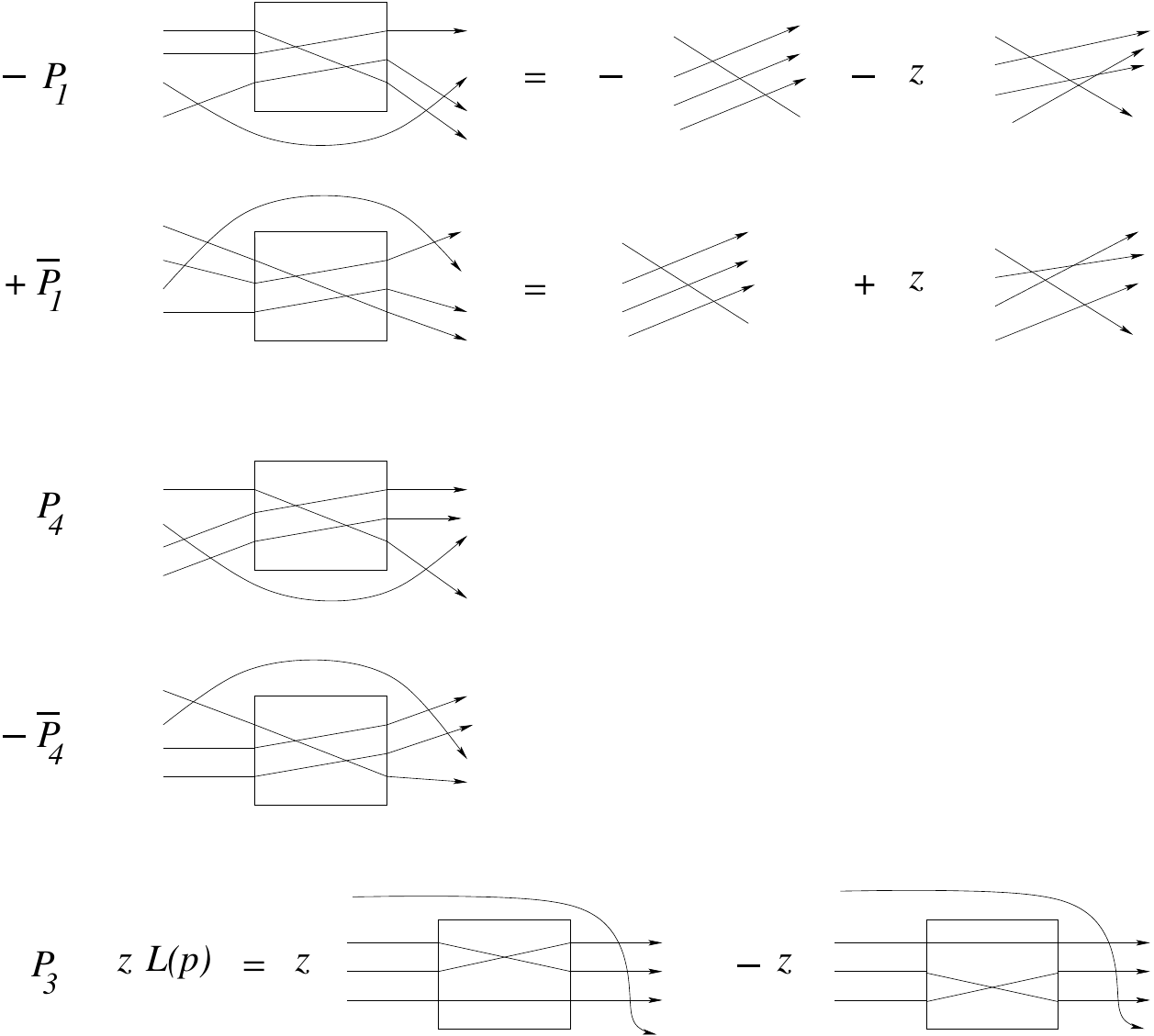}
\caption{\label{tablesmooth}  the remaining partial smoothings $\sigma_2\sigma_1$ for the meridian of a positive quadruple crossing}  
\end{figure}

\begin{itemize}
\item $W(p)$ is an {\em isotopy invariant for the Reidemeister move $p$}, i.e. $W(p)$ is invariant under any isotopy of the rest of the tangle outside of $I^2_p \times I$. Notice that $I^2_p \times I$ is not a small cube around the triple crossing, but its hight is the hight of the whole cube. In particular, no branch is allowed to move over or under the triple crossing.

\item $W(p)$ depends non trivially on the point at infinity

\item $W(P_3)= W(\bar P_3)-1$ if and only if the triple crossing in $P_1$ (and hence in $\bar P_1$) is of a certain global type

\item $W(P_2)=W(\bar P_2)$ (this implies in fact the scan-property)

\item $W(P_1)=W(\bar P_1)=W(P_4)=W(\bar P_4)$

\end{itemize}

Assume that we have found such a $W(p)$.
We define then  \vspace{0,2 cm}

\begin{figure}
\centering
\includegraphics{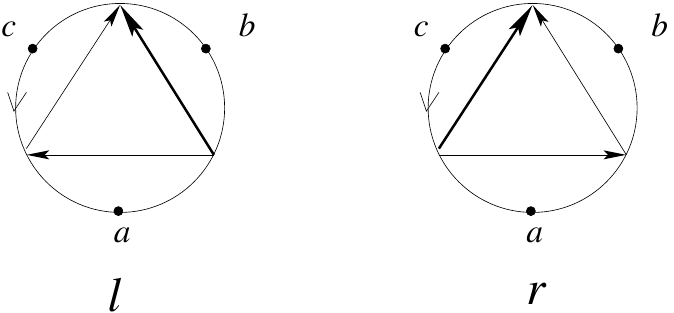}
\caption{\label{globtricross}  the six global types of triple crossings}  
\end{figure}

$R(m)=\sum_{p \in m} sign(p) \sigma_2\sigma_1(p) + \sum_{p \in m}sing(p)zW(p)(\sigma_2(p)-\sigma_1(p))$  \vspace{0,2 cm}

where the first sum is only over certain global types of triple crossings (compare Fig.~\ref{globtricross}).
It follows from the above considerations that $R(m)=0$ and consequently the 1-cochain $R$ it is a new solution of the tetrahedron equation (2).

It turns out that $W(p)$ does exist. We call it a {\em weight for the partial smoothing}. It is a finite type invariant of degree 1 for long knots relative to the triple crossing $p$ and it is convenient to define it as a Gauss diagram invariant (for the definition of Gauss diagram invariants see \cite{PV} and also \cite{F1}). However, it is only invariant under regular isotopy and the corresponding 1-cocycle behaves uncontrollable under Reidemeister  moves of type I.  

So, let us recall the situation:

\begin{itemize}

\item There is a solution of the tetrahedron equation of {\em constant} weight, i.e. the partial smoothing $\sigma_2\sigma_1(p)$ enters with coefficient 0 and $L(p)$ enters with a coefficient of degree 0 (non-zero constant)

\item There is a solution of the tetrahedron equation of {\em linear} weight, i.e. $\sigma_2\sigma_1(p)$ enters with a coefficient of degree 0 (non-zero constant) and $L(p)$ enters with a coefficient of degree 1

\end{itemize}

The solution with linear weight leads to our 1-cocycle $R^{(1)}_{reg}$. The grading, mentioned in the introduction,  comes from the surprising fact that {\em the new contributing to $W(p)$ crossing in $\bar P_3$ with respect to $ P_3$ is always  identical to the crossing between the highest and the middle branch in the triple crossings $P_1$ and $\bar P_1$}. This grading is very useful for $R^{(1)}_{reg}(scan(T))$.

But these two solutions do not lead to a 1-cocycle which represents a non trivial cohomology class in the special case of long knots, and moreover the grading is trivial in this case. It turns out that we have to go one step further.

\begin{itemize}

\item There is a solution of the tetrahedron equation of {\em quadratic} weight, i.e. $\sigma_2\sigma_1(p)$ enters with a coefficient of degree 1 and $L(p)$ enters with a coefficient of degree 2
\end{itemize}

Notice that the weights $W(p)$ will depend on the whole long knot $K$ and not only on the part of it in $I^2_{quad} \times I$. Moreover, they depend crucially on the point at infinity and they can not be defined for compact knots in $S^3$.
Therefore we have to consider {\em twenty four} different positive tetrahedron equations, corresponding to the six different abstract closures of the four lines in $I^2_{quad} \times I$ to a circle and to the four different choices of the point at infinity in each of the six cases. $R^{(1)}$ and $R^{(1)}_{reg}$ are common solutions of all twenty four tetrahedron equations and we call them therefore  {\em solutions of a global positive tetrahedron equation}.

The existence of these solutions looks like a miracle. Over the years we have constructed some interesting 1-cocycles, notably one where virtual tangles appear in the target and with coefficients which are polynomials in five(!) variables (will be published later). However, it seems likely that it represents the trivial cohomology class for long knots because it doesn't break the symmetry. (Let us consider for example the loop $rot(K)$. Instead of moving the knot $K$ we keep it fixed but we rotate the plan to which we project it. Then we see each triple crossing exactly twice, for one direction and for the opposite direction and their contributions enter with opposite signs.)

In order to define a 1-cocycle which represents a non trivial cohomology class for long knots one has to break the symmetry in a very strong way: {\em the partial smoothings $\sigma_2\sigma_1(p)$ and  $\sigma_2(p)-\sigma_1(p)$ enter both with  coefficients which depend on the whole long knot $K$ but the partial smoothing $\sigma_1\sigma_2(p)$ doesn't enter at all.} For solutions with this property all twenty four tetrahedron equations become independent (there aren't any symmetries and the system becomes very over determined)! Their solutions give rise (after in addition the solving of the corresponding cube equations) to our 1-cocycles $R^{(1)}_{reg}$ and $R^{(1)}$. Moreover, $R^{(1)}$ vanishes for the meridians 
 of a degenerate cusp and of the transverse intersections of arbitrary strata of codimension 1. It follows from the list of singularities of codimension two for projections of knots into the plan that $R^{(1)}$ is a 1-cocycle for all isotopies in the complement of the strata corresponding to a cusp with a transverse branch. We have succeeded to complete $R^{(1)}$ with a finite type Gauss diagram 1-cocycle $V^{(1)}$ of degree 3 (which has also the scan-property) in order to make it zero on the meridians of these remaining strata of codimension two. The result is our 1-cocycle $\bar R^{(1)}$ (compare the important Remark 11 in Section 10).

\begin{remark}
We could make both $R^{(1)}$ and $V^{(1)}$ to 1-cocycles in the whole $M_T$ by "symmetrizing" them (i.e. by considering linear combinations with "dual" 1-cocycles obtained by applying $flip$ to the configurations or by using "heads" instead of "foots" of arrows, compare Sections 4 and 5). But it seems likely that each of them would represent the trivial cohomology class.
\end{remark}

\begin{figure}
\centering
\includegraphics{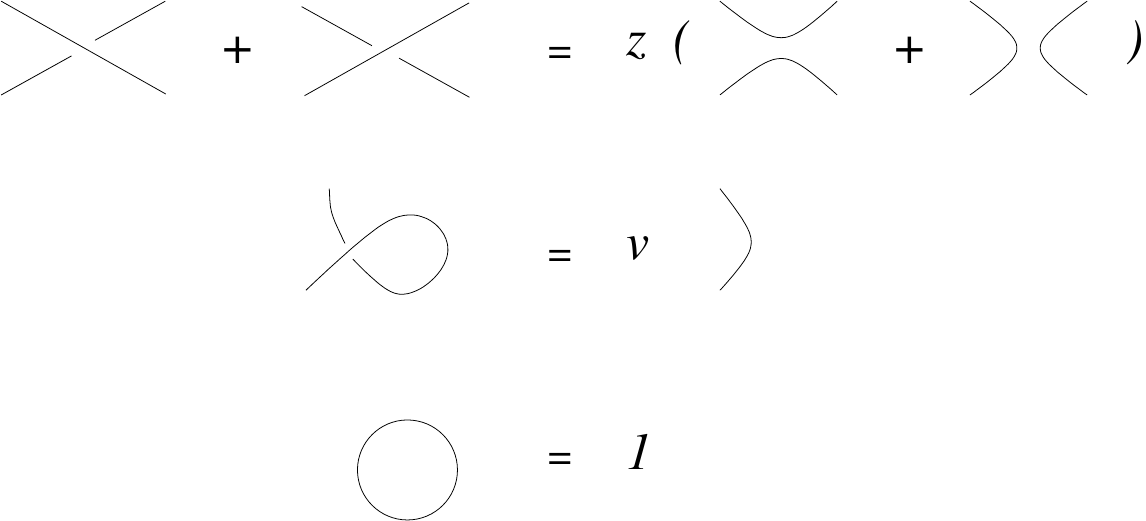}
\caption{\label{basisKauf} skein relations for the Kauffman invariant}  
\end{figure}

We have carried out our construction also for the 2-variable Kauffman polynomial $F$ instead of the HOMFLYPT polynomial. 
The skein relations for the 2-variable Kauffman invariant are shown in Fig.~\ref{basisKauf} (see \cite{K}).

\begin{remark}
Notice that the solution of the global positive tetrahedron equation in the case of the HOMFLYPT polynomial  can be expressed by using only partial smoothings with {\em two} crossings: \vspace{0,2 cm}

$R(m)=\sum_{p \in m} sign(p) \sigma_2\sigma_1(p) + \sum_{p \in m}sing(p)W(p)(\sigma^2_2(p)-\sigma^2_1(p))$  \vspace{0,2 cm}

where the first sum is only over certain global types of triple crossings.
Amazingly, exactly the same formula is also a solution in the case of the Kauffman polynomial!
\end{remark}

However, we haven't solved the cube equations in this case. They are very complicated. 

\begin{figure}
\centering
\includegraphics{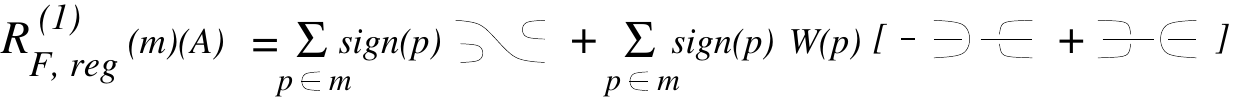}
\caption{\label{solkauf}  surprising solution of the positive tetrahedron equation in Kauffman's case}  
\end{figure}

{\em Surprisingly, there is a second solution of the global positive tetrahedron equation in Kauffman's case (which does not exist in the HOMFLYPT case). We give it in Fig.~\ref{solkauf}. The weights and the global types of triple crossings are exactly the same as in the previous solutions.}

We solve the cube equations for this solution with coefficients in $\mathbb{Z}/2\mathbb{Z}$.
The result are two other quantum 1-cocycles, called $R^{(1)}_{F, reg}$ and $R^{(1)}_F$ which also have the scan-property. We can complete $R^{(1)}_F$ with exactly the same 1-cocycle $V^{(1)}$ as in the HOMFLYPT case: 

$\bar R^{(1)}_F=R^{(1)}_F+F_T V^{(1)}$. 

A calculation by hand gives now

$\bar R^{(1)}_F(rot(3^+_1))=F_{3^+_1}$ mod 2.

Consequently, $\bar R^{(1)}_F$ represents also a non trivial cohomology class in the topological moduli spaces $M_T$ and it has the scan-property too.

It is well known that the colored HOMFLYPT polynomial, i.e. the knot invariant associated to representations of $sl_N$, can be expressed as a linear combination of usual HOMFLYPT polynomials applied to cables, see e.g. \cite{L}. Hence we define the corresponding 1-cocycles  just as linear combinations of the 1-cocycles $\bar R^{(1)}$ applied to the cables $Cab_n(K)$ of long knots $K$. (We define a substitute for $rot(Cab_n(K))$ in Remark 13 in Section 11.)

We haven't done it in the case of the $g_2$-Kuperberg invariant, which is extremely complicated (see e.g. \cite{BS}).

Let us finish this long introduction with the vague formulation of a challenging problem.
\begin{question}
 Let $\Delta^2=(\sigma_1 \sigma_2...\sigma_{n-1})^n$ be Garsides braid (as a canonical braid word) which is the generator of the center of the braid group $B_n$ (see e.g. \cite{B}). Let $\beta \Delta^2 \rightarrow \Delta^2 \beta$ be the isotopy which pushes $\Delta^2$ through $\beta$. This is a well defined homotopy class of an arc in $M_{ \beta \Delta^2}=M_{ \Delta^2 \beta}$. We consider the standard closure of $\beta$ in the solid torus. We assume that it is a knot and we chose a point at infinity in $\partial \beta$. 

{\em Does there exist a higher representation theory $Rep$ for homotopy classes of one parameter families of n-braids with a marked point in their boundary (the entries of the matrices associated to the applications of the braid relations should be elements of the non commutative Hecke algebra $H_n$) and such that 

$R^{(1)}(\beta \Delta^2 \rightarrow \Delta^2 \beta)=
\\log(det(Rep(\beta \Delta^2 \rightarrow \Delta^2 \beta)))=tr(log(Rep(\beta \Delta^2 \rightarrow \Delta^2 \beta)))$?} 
\end{question}

However, we hope that perhaps some day our new solution of a tetrahedron equation using knot theory will become interesting in mathematical physics too (usually it is the other way round).

{\em Acknowledgments}

{\em I am grateful to Alan Hatcher for patiently answering my questions about topological moduli spaces, to Christian Blanchet for his encouragements over many years and to Ryan Budney who has pointed out an error in the application of Hatcher's loop in an earlier version. Finally, let me mention that without S\'everine, who has created all the figures, this paper wouldn't exist.}

\section{ The topological moduli space of diagrams of string links and its stratification}

We work in the smooth category.

We fix an orthogonal  projection $pr: I^3 \rightarrow I^2$. Let $T$ be a string link. $T$ is called {\em regular} if $pr:T \rightarrow I^2$ is an immersion. Let $T$ be a regular string link and let $M_T$ be the (infinite dimensional) space of all string links isotopic to $T$. We denote by $M_T^{reg}$ the subspace of all string links which are regularly isotopic to $T$, i.e. all string links in the isotopy are regular. The infinite dimensional spaces $M_T$ and $M_T^{reg}$ have a natural {\em stratification with respect to $pr$}:

$M_T=\Sigma^{(0)}\cup\Sigma^{(1)}\cup\Sigma^{(2)}\cup\Sigma^{(3)}\cup\Sigma^{(4)}...$

Here, $\Sigma^{(i)}$ denotes the union of all strata of codimension i.
%%%%%%%%%%%%%%%
The strata of codimension 0 correspond to the usual generic {\em diagrams of tangles}, i.e. all singularities in the projection are ordinary double points. So, our {\em discriminant} is the complement of $\Sigma^{(0)}$ in $M_T$. Notice that this discriminant is very different from Vassiliev's discriminant of singular knots \cite{V1}.

 The three types of strata of codimension 1 correspond to the {\em Reidemeister moves}, i.e. non generic diagrams which have exactly one ordinary triple point, denoted by $\Sigma^{(1)}_{tri}$,  or one ordinary self-tangency, denoted by $\Sigma^{(1)}_{tan}$, or one ordinary cusp, denoted by $\Sigma^{(1)}_{cusp}$, in the projection $pr$. We call the triple point together with the under-over information (i.e. its embedded resolution in $I^3$ given by the tangle T) a {\em triple crossing}.

\begin{figure}
\centering
\includegraphics{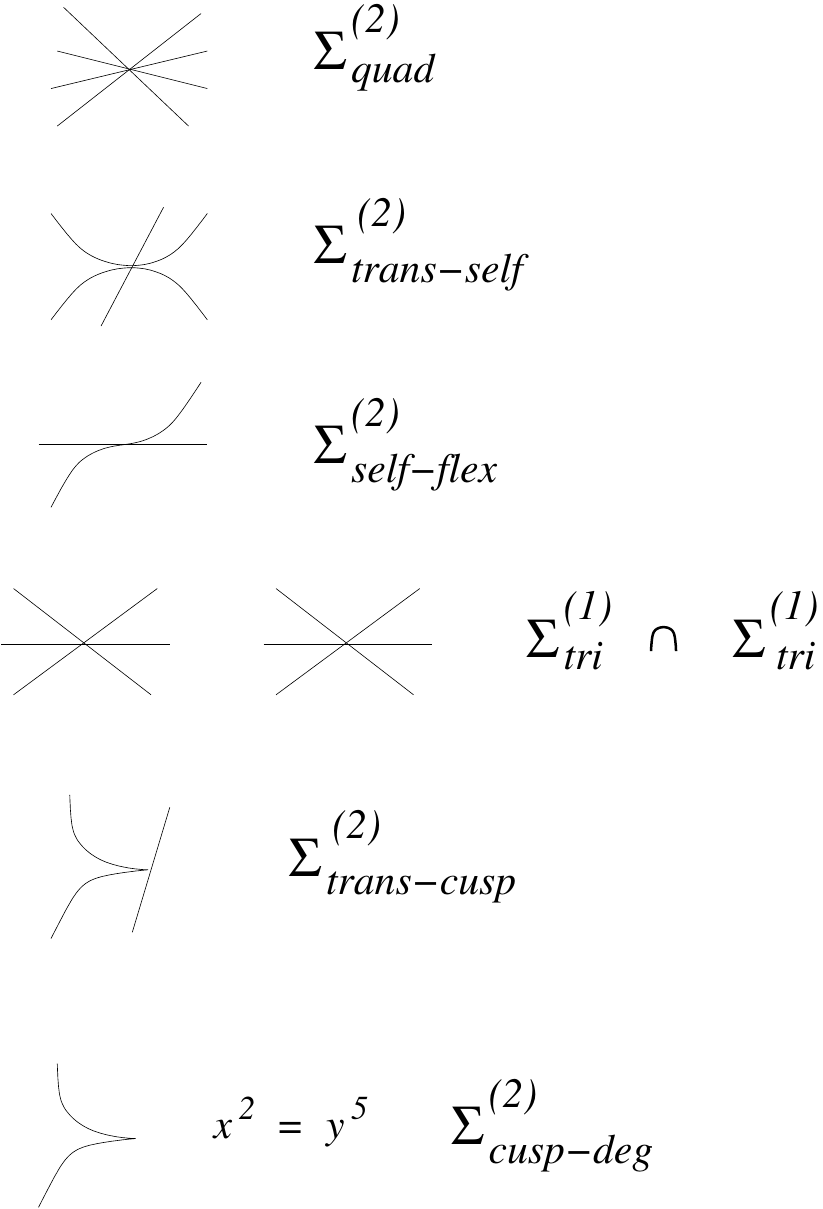}
\caption{\label{sing}  the strata of codimension $2$ of the discriminant}  
\end{figure}

There are exactly six types of strata of codimension 2. They correspond to non generic diagrams which have exactly either 
\begin{itemize}
\item one ordinary quadruple point, denoted by $\Sigma^{(2)}_{quad}$
\item one ordinary self-tangency with a transverse branch passing through the tangent point, denoted by  
$\Sigma^{(2)}_{trans-self}$
\item one self-tangency in an ordinary flex ($x=y^3$), denoted by  $\Sigma^{(2)}_{self-flex}$
\item two singularities of codimension 1 in disjoint small discs (this corresponds to the transverse intersection of two strata from 
$\Sigma^{(1)}$)
\item one ordinary cusp with a transverse branch passing through the cusp, denoted by  $\Sigma^{(2)}_{trans-cusp}$
\item one degenerate cusp, locally given by $x^2=y^5$, denoted by  $\Sigma^{(2)}_{cusp-deg}$

We show these strata in Fig.~\ref{sing}.

(Compare \cite{FK}, which contains all the necessary preparations from singularity theory.)

\end{itemize}

The strata which we need from $\Sigma^{(3)}$ will be described later.

The list above of singularities of codimension two can be seen as an analogue for one parameter families of diagrams of the Reidemeister theorem (which gives the list of singularities of codimension one).

{\em Reidemeisters Theorem} \cite{BZ} (sometimes its refinement {\em Markovs Theorem} \cite{B}) was the main tool for almost 90 years to construct calculable combinatorial knot invariants, notably quantum knot invariants. The latter can be often upgraded to a TQFT, see \cite{BHMV}. The application of Reidemeisters Theorem  has culminated in the {\em categorification} of quantum invariants  initiated by Khovanov (see \cite{Kh}). (Other approaches which do not use projections and Reidemeisters Theorem are knot invariants which arrive from the geometric structure of the knot complement, going back to  \cite{Th}, invariants of 3-manifolds obtained by Dehn surgery of the knot or by coverings ramified in the knot e.g. \cite{BZ}, finite type invariants via the Kontsevich integral \cite{Kon}, the Heegard-Floer knot homology \cite{OS}, invariants via configuration spaces \cite{BuCo}, \cite{PV2}, \cite{DTh}, \cite{Oh}, \cite{Vol} and probably still other approaches.) The procedure in the combinatorial approach is always the following:

\begin{itemize}
\item (1) associating something calculable (e.g. an element in some algebra) to a point in $\Sigma^{(0)}$
\item (2) proving the invariance under passing strata from $\Sigma^{(1)}$
\item (3) simplifying this proof by using strata from $\Sigma^{(2)}$
\end{itemize}

A generic isotopy of string links is an arc in $M_T$ and the intersection with $\Sigma^{(i)}$ is empty for $i>1$.

Invariants which are constructed in this way are consequently 0-cocycles for $M_T$ with values in some algebra.

The main tool for (1) are {\ quantum invariants}, i.e. elements $I$ of a {\em skein module} coming from irreducible representations of quantum groups $U_q(g)$ (here $U_q(g)$ is the quantum enveloping algebra of a complex simple Lie algebra $g$ at a nonzero complex number q which is not a root of unity) or more generally of ribbon Hopf algebras ( see e.g. \cite{J}, \cite{Tu}, \cite{L}).

(2) reduces essentially to the proof of invariance of $I$ under passing a stratum of $\Sigma^{(1)}$ which corresponds to a non generic diagram which has exactly one ordinary positive triple crossing. In the case of quantum invariants this reduces to prove that quantum invariants are solutions of the {\em Yang-Baxter equation}.

(3) is not explicitly mentioned (or only in a sloppy way) in many papers. We will see later that there are exactly two global (without taking into account the point at infinity) and eight local types of triple crossings. The Yang-Baxter equation does not depend on the global type of the triple crossing. Hence we stay with exactly eight different Yang-Baxter equations. Different local types of triple crossings come together in 
points of $\Sigma^{(2)}_{trans-self}$. We make a graph $\Gamma$ in the following way: the vertices's correspond to the different local types of triple crossings. We connect two vertices's by an edge if and only if the corresponding strata of triple crossings are adjacent to a stratum of $\Sigma^{(2)}_{trans-self}$. One easily sees that the resulting graph is the 1-skeleton of the 3-dimensional cube $I^3$ (compare Section 6).
In particular, it is connected. Studying the normal discs to $\Sigma^{(2)}_{trans-self}$ in $M_T$ one shows that if $I$ is invariant under passing all $\Sigma^{(1)}_{tan}$ and just {\em one} local type of a stratum $\Sigma^{(1)}_{tri}$ then it is invariant under passing all types of triple crossings because $\Gamma$ is connected. Hence the use of certain strata of $\Sigma^{(2)}$ allows to reduce the eight Yang-Baxter equations to a single one.

\begin{figure}
\centering
\includegraphics{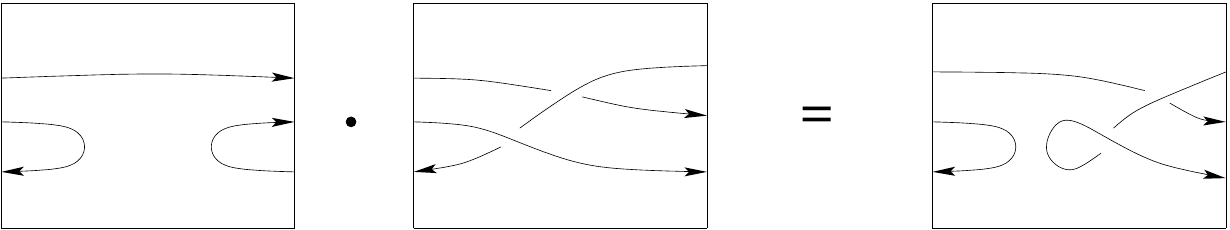}
\caption{\label{comptangle}  multiplication of tangles}  
\end{figure}

We call the resulting 0-cocycles {\em quantum 0-cocycles}. This are the usual well known quantum invariants.
All quantum invariants verify {\em skein relations} (see Fig.~ \ref{HOMFLYPT} and  Fig.~ \ref{basisKauf}) and hence give rise to {\em skein modules}, i.e. formal linear combinations over the corresponding ring of coefficients of tangles (with the same oriented boundary) modulo the skein relations. Skein modules were introduced by Turaev \cite{Tu2} and Przytycki \cite{Pr}. The {\em multiplication of tangles} is defined by their composition (if possible), see Fig.~ \ref{comptangle}. The main property of quantum invariants (and which follows directly from their construction) is their {\em multiplicativity}:
let $T_1$ and $T_2$ be tangles which can be composed and let $I(T_i)$ be the corresponding quantum invariants in the skein modules. Then $I(T_1T_2)=I(T_1)I(T_2)$, where we have to decompose, by using the skein relations again, the product of the generators of the skein modules for 
$\partial T_1$ with those for $\partial T_2$ into a linear combination of generators for the skein module of $\partial (T_1T_2)$.
In particular, if the n-tangle $T_2$ consists of n trivial arcs with just a local knot $K$  tied into one of the arcs then  $I(T_1T_2)=I(T_1)I(K)$. Consequently, a Reidemeister move of type I (i.e. a connected sum with a small diagram of the unknot which has just one crossing) changes $I(T_1)$ only  by some standard factor.

The most famous quantum invariants are the HOMFLYPT polynomial $P$, the 2-variable Kauffman polynomial $F$ and Kuperberg's  polynomial $I^{g_2}$ and their generalizations for tangles.  In the case of Kuperberg's invariant we have to allow planar trivalent graphs \cite{BS}.

\section{The strategy of the construction and notations}

The construction of the 1-cocycles is very long and difficult. But once we have them their applications will be much simpler.

{\em The main idea consists in using strata $\Sigma^{(i)}$ of the discriminant  of codimensions $i=0,1,2,3$ in order to construct quantum 1-cocycles  in a similar way as the strata of codimensions 0,1 and 2 were used in order to construct quantum 0-cocycles.}

\begin{figure}
\centering
\includegraphics{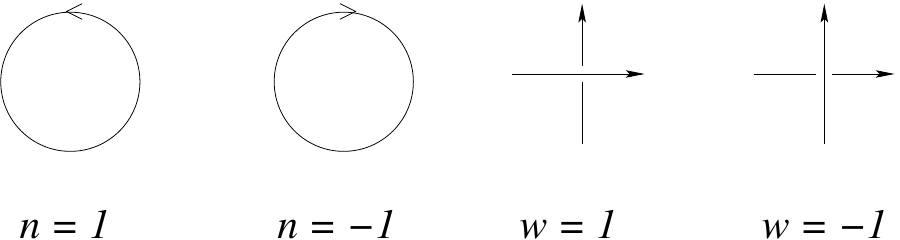}
\caption{\label{conv} Whitney index and writhe}  
\end{figure}

First of all, it turns out that starting from 1-cocycles we have to replace in general isotopy by regular isotopy. This comes from the fact that the invariants which we will construct are no longer multiplicative and hence the moment (in a family of diagrams) and the place (in the diagram) of a Reidemeister move of type I become important. Fortunately, this is not restrictive in order to distinguish tangles. It is well known that two (ordered oriented) tangles are regularly isotopic if and only if they are isotopic and the corresponding components share the same {\em Whitney index} $n(T)$ and the same {\em writhe} $w(T)$ (see e.g. \cite{F1}). Here as usual, the Whitney index is the algebraic number of oriented circles in the plane after smoothing all double points of a component with respect to the orientation and the writhe is the algebraic number of all crossings of a component. We use the standard conventions shown in Fig.~ \ref{conv}. Hence, by adding if necessary small curls to components we can restrict ourself to regular isotopy.

We choose a quantum 0-cocycle e.g. the regular HOMFLYPT invariant in the regular HOMFLYPT skein module $S(\partial T)$. We fix moreover an abstract closure $\sigma$ of $T$ such that it becomes an oriented circle (just a part of the circle is embedded as $T$ in $I^3$) and we fix a {\em point at infinity} in $\partial T$.

\begin{itemize}
\item (1) A generic arc $s$ in $M^{reg}_T$ intersects $\Sigma^{(1)}_{tri}$ and $\Sigma^{(1)}_{tan}$ transversally in a finite number of points.  To each such intersection point we associate an element in the HOMFLYPT skein module $S(\partial T)$. 
We sum up with signs the contributions over all intersection points in $s$ and obtain an element, say $R(s)$, in  $S(\partial T)$.

\item (2) We use now the strata from $\Sigma^{(2)}$ to show the invariance of $R(s)$ under all generic homotopies of $s$ in $M_T^{reg}$ which fix the endpoints of $s$. The contributions to $R(s)$ come from discrete points in $s$ and hence they survive under generic Morse modifications of $s$. It follows that $R(s)$ is a 1-cocycle for $M_T^{reg}$.

\item (3) Again, we use some strata from $\Sigma^{(3)}$ in order to simplify the proof that $R(s)$ is a 1-cocycle.

\end{itemize}

In (1) we do something new for triple crossings and for self-tangencies, but we treat the usual crossings as previously in the construction of 0-cocycles. This is forced by those strata in $\Sigma^{(2)}$ which consist of the transverse intersections of two strata in $\Sigma^{(1)}$. 

In (2) we have to study normal 2-discs for the strata in $\Sigma^{(2)}$. For each type of stratum in $\Sigma^{(2)}$ we have to show that $R(m)=0$ for the boundary $m$ of the corresponding normal 2-disc. We call $m$ a {\em meridian}. $\Sigma^{(2)}_{quad}$ is by fare the hardest case and we call the corresponding equation $R(m)=0$ the {\em tetrahedron equation}. Indeed, as already mentioned in the Introduction the quadruple crossing deforms into four different triple crossings (the vertices's of a tetrahedron) and the six involved ordinary crossings form the edges of the tetrahedron. (This is an analogue of the Yang-Baxter or triangle equation in the case of 0-cocycles.) It turns out that without considering the point at infinity there are 6 global types (corresponding to the {\em Gauss diagrams} but without the writhe of the quadruple crossings, see the end of the section) and 48 local types of quadruple crossings. At each point of $\Sigma^{(2)}_{quad}$ there are adjacent exactly eight strata of $\Sigma^{(1)}_{tri}$ (compare Fig.~ \ref{quadcros}). 

The edges of the graph $\Gamma = skl_1(I^3)$, which was constructed in the last section, correspond to the types of strata in $\Sigma^{(2)}_{trans-self}$. The solution of the tetrahedron equation tells us what is the contribution to $R$ of a positive triple crossing (i.e. all three involved crossings are positive) . The meridians of the strata from $\Sigma^{(2)}_{trans-self}$ give equations which allow us to determine the contributions of all other types of triple crossings as well as the contributions of self-tangencies. However, each loop in $\Gamma$ could give an additional equation. Evidently, it suffices to consider the loops which are the boundaries of the 2-faces from $skl_2(I^3)$. In fact, they do give additional equations (see Sections 6 and 7) and force us to consider smoothings of self-tangencies too. We call all the equations which come from the meridians of $\Sigma^{(2)}_{trans-self}$ and from the loops in $\Gamma = skl_1(I^3)$ the {\em cube equations}. (Notice that a loop in $\Gamma$ is more general than a loop in $M_T$. For a loop in $\Gamma$ we come back to the same local type of a triple crossing but not necessarily to the same whole diagram of $T$.)

\begin{figure}
\centering
\includegraphics{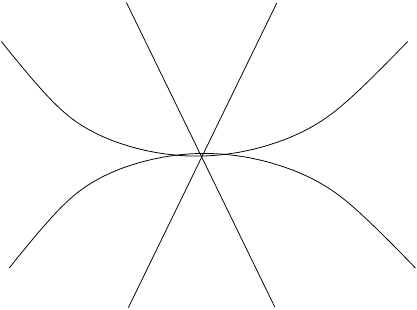}
\caption{\label{degquad}  a quadruple crossing with two tangential branches}  
\end{figure}

\begin{figure}
\centering
\includegraphics{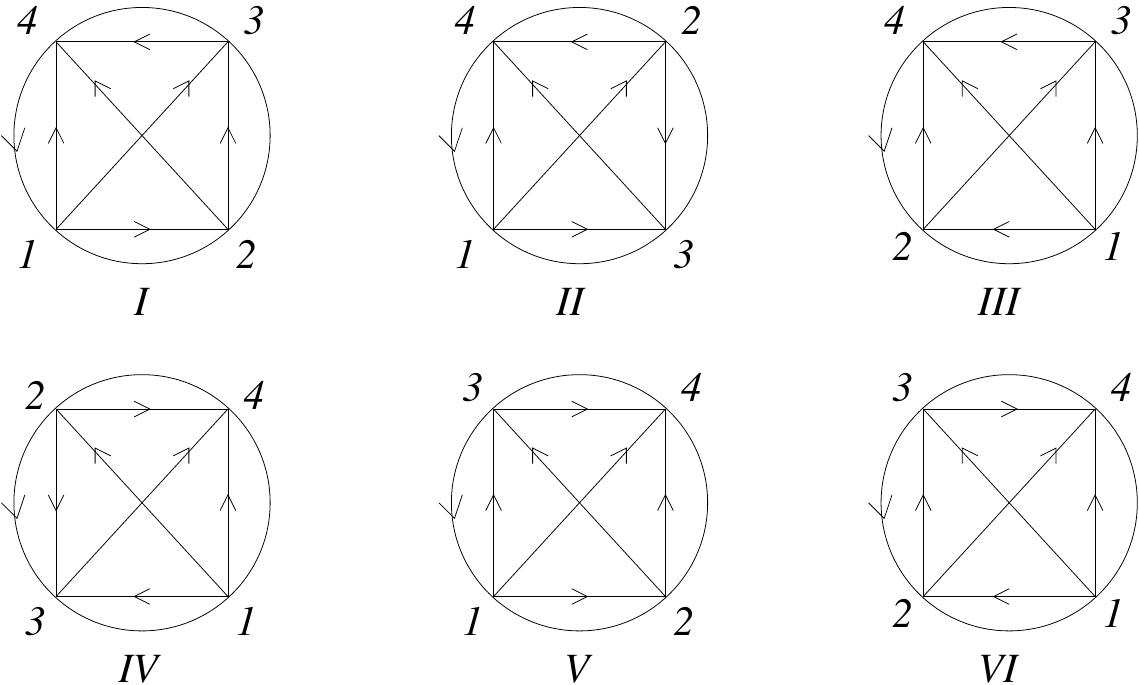}
\caption{\label{globquad} global types of quadruple crossings}  
\end{figure}

In (3) we consider the strata  from $\Sigma^{(3)}$ which correspond to diagrams which have a degenerate quadruple crossing where exactly two branches have an ordinary self-tangency in the quadruple point, denoted by $\Sigma^{(3)}_{trans-trans-self}$ (see Fig.~\ref{degquad}). (These strata are the analogue of the strata $\Sigma^{(2)}_{trans-self}$ in the construction of 0-cocycles.) Again, we form a graph with the local types of quadruple crossings as vertices's and the adjacent strata of $\Sigma^{(3)}_{trans-trans-self}$ as edges. One easily sees that the resulting graph $\Gamma$ has exactly 48 vertices's and that it is again connected. We don't need to study the unfolding of $\Sigma^{(3)}_{trans-trans-self}$ in detail. It is clear that each meridional 2-sphere for $\Sigma^{(3)}_{trans-trans-self}$ intersects $\Sigma^{(2)}$  transversally in a finite number of points, namely exactly in two strata from $\Sigma^{(2)}_{quad}$ and in lots of strata from $\Sigma^{(2)}_{trans-self}$ and from 
$\Sigma^{(1)} \cap \Sigma^{(1)}$ and there are no other intersections with strata of codimension 2.  If we know now that $R(m)=0$ for the meridian $m$ of one of the quadruple crossings, that $R(m)=0$ for all meridians of  $\Sigma^{(2)}_{trans-self}$ (i.e. $R$ satisfies the cube equations) and for all meridians of  $\Sigma^{(1)} \cap \Sigma^{(1)}$ then $R(m)=0$ for the other quadruple crossing too. It follows that for each fixed global type (see Fig.~\ref{globquad}) the 48 tetrahedron equations reduce to a single one, which we call the {\em positive tetrahedron equation} (compare the Introduction). This phenomenon is analog to that for the Yang-Baxter or triangle equation, which was explained in the previous section. We need now a solution of the positive tetrahedron equation which depends on a point at infinity but which is a solution independent  of the position of this point at infinity in the knot diagram.

Let $T$ be a generic n-tangle in $I^3$. We assume that a half of $\partial T$ is contained in $0 \times I^2$ and the other half in $1 \times I^2$ and we fix an abstract closure $\sigma$ of $T$ such that it becomes an oriented circle (compare e.g. Example 1 in Section 8).

\begin{figure}
\centering
\includegraphics{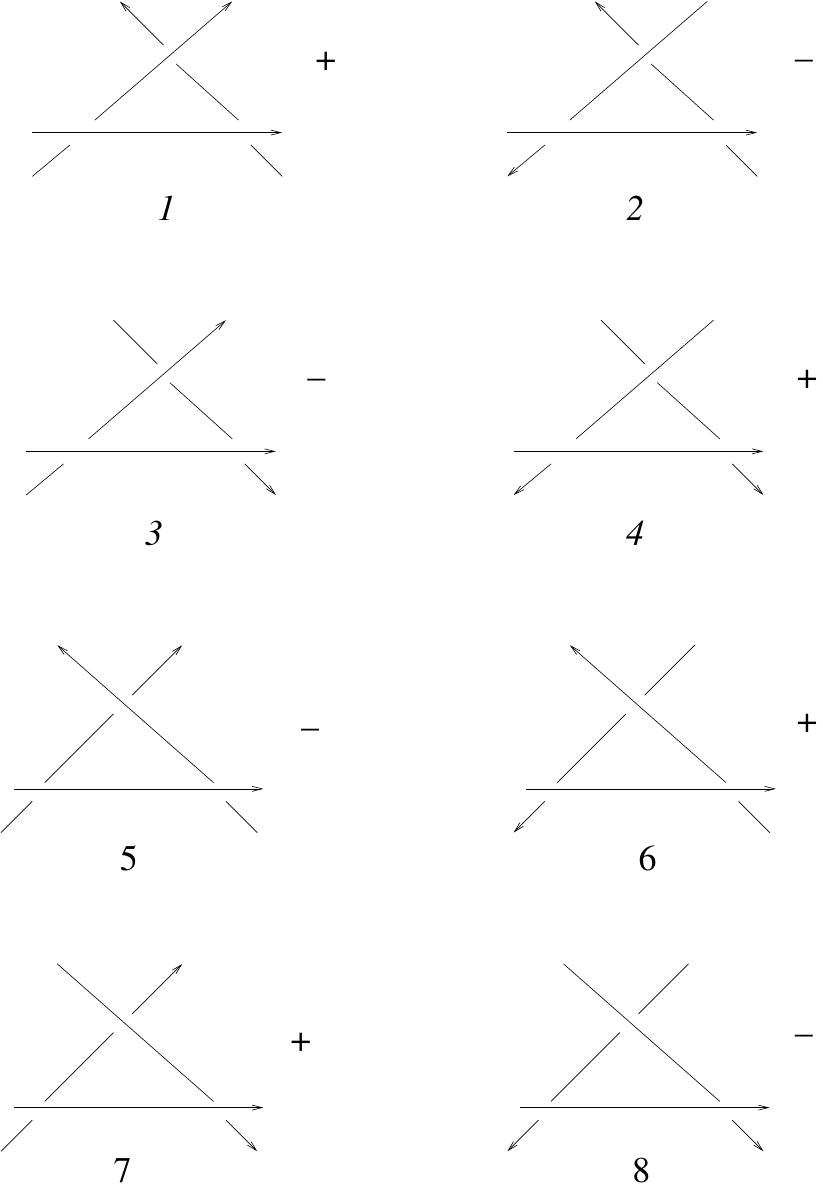}
\caption{\label{loc-trip} the local types of triple crossings}  
\end{figure}

The sign of a positive Reidemeister III move was already defined in the Introduction. The signs of the remaining seven Reidemeister III moves are determined by the meridians of the strata  $\Sigma^{(2)}_{trans-self}$. It is convenient to equip the strata of $\Sigma^{(1)}_{tri}$ with a co-orientation. The sign is then defined as the local intersection index, which is $+1$ if the oriented isotopy intersects $\Sigma^{(1)}$ transversally from the negative side to the positive side. We give the local types of Reidemeister moves of type III together with the co-orientation  in 
Fig.~\ref{loc-trip}. The + or - sign here denotes the side of the complement of $\Sigma^{(1)}_{tri}$ with respect to the co-orientation.
An easy verification shows that they satisfy the cube equations, i.e. the signs are consistent for the boundary of each 2-cell in  $skl_2(I^3)$. The six local types $1, 3, 4, 5, 7, 8$ can appear for braids and are therefore called {\em braid-like}. The remaining two types $2, 6$ are called {\em star-like}. In the star-like case we denote by $mid$ the boundary point which corresponds to the ingoing middle branch.

\begin{remark}
In order to handle complicated combinatorial facts we had to create our own language. We apologize by the reader for the arbitrary kind of notations, here and in other places. But we are used to them for many years and we feel that changing them could  possibly lead to errors.
\end{remark}

\begin{figure}
\centering
\includegraphics{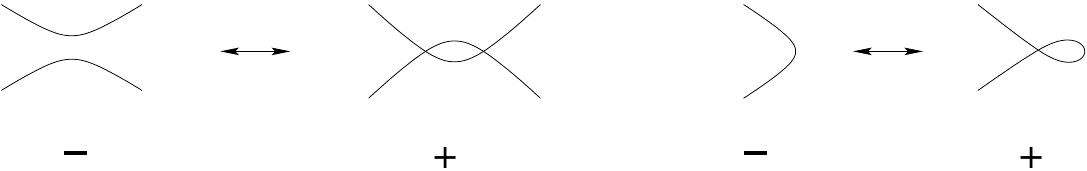}
\caption{\label{signRII}  coorientation for Reidemeister moves of typ II and I}  
\end{figure}

\begin{definition}
The sign of a Reidemeister II move and of a Reidemeister I move is shown in Fig.~\ref{signRII}.
\end{definition}

Notice that these signs depend only on the underlying planar curve in contrast to the signs for triple crossings.

To each Reidemeister move of type III corresponds a diagram with a {\em triple 
crossing} $p$: three branches of the tangle (the highest, middle and lowest with respect to the projection $pr$) have a common point in the projection into the 
plane. A small perturbation of the triple crossing leads to an ordinary diagram with three crossings near $pr(p)$.

\begin{figure}
\centering
\includegraphics{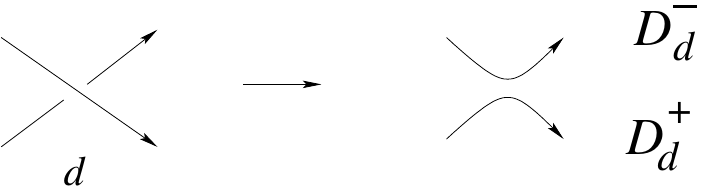}
\caption{\label{splitd}  the two ordered knot diagrams associated to a distinguished crossing}  
\end{figure}

\begin{definition}
 We call the crossing between the 
highest and the lowest branch of the triple crossing $p$ the {\em distinguished crossing} of $p$ and we denote it by $d$. The crossing between the highest branch and the middle branch is denoted by $hm$ and that of the middle branch with the lowest is denoted by $ml$.
 Smoothing the distinguished crossing with respect to the orientation splits $T \cup \sigma$ into two oriented and ordered circles. We call $D^+_d$ the component which goes from the under-cross to the over-cross at $d$ and by $D^-_d$ the remaining component (compare Fig.~\ref{splitd}).

In a Reidemeister move of type II both new crossings are considered as distinguished and denoted by $d$.
\end{definition}

For the definition of the weights $W(p)$ we use Gauss diagram formulas (see \cite{PV} and also \cite{F1}).
A {\em Gauss diagram of\/} $K$ is an oriented circle with oriented chords and a marked point.
There exists an orientation preserving diffeomorphism from the oriented line to oriented the knot $T \cup \sigma$ such that each
chord connects a pair of points which are mapped onto a crossing of
$pr(T)$ and infinity is mapped to the marked point. The chords are oriented from the preimage of the under crossing
to the pre-image of the over crossing (here we use the orientation of
the $I$-factor in $I^2 \times I$). The circle of a Gauss diagram in the plan is always equipped with the counter-clockwise orientation.

A {\em Gauss diagram formula\/} of degree $k$ is an expression assigned to a knot
diagram which is of the following form:
\begin{displaymath}
\sum{\textrm{function( writhes of the crossings)}}
\end{displaymath}
where the sum is taken over all possible choices of $k$ (unordered)
different crossings in the knot diagram such that the chords arising
from these crossings in the knot diagram of $T \cup \sigma$ build a given sub-diagram
with given marking. The marked sub-diagrams are called
{\em configurations\/}.
If the function is the product of the writhes of the crossings in the configuration, then we will denote the
sum shortly by the configuration itself.

A Gauss diagram formula which is invariant under regular isotopies of $T \cup \sigma$ is called a {\em Gauss diagram invariant}.

Let us consider the {\em Gauss diagram of a triple crossing $p$}, but without taking into account the writhe of the crossings. In the oriented circle $T \cup \sigma$ we connect the preimages of the triple point $pr(p)$ by arrows which go from the under-cross to the over-cross and we obtain a triangle. The distinguished crossing $d$ is always drawn by a thicker arrow.

There are exactly six different {\em global types} of triple crossings with a point at infinity. We give names to them and show them in Fig.~\ref{globtricross}. (Here "r" indicates that the crossing between the middle and the lowest branch goes to the right and "l" indicates that it goes to the left.) Notice that the involution {\em flip} exchanges the types $r_a$ and $l_a$ as well as the types $r_b$ with $l_c$ and $r_c$ with $l_b$.

The following  figures Fig.~\ref{Iglob}...Fig.~\ref{VI2glob} are derived from Fig.~\ref{unfoldquad} and Fig.~\ref{globquad}. They show the six global types of positive quadruple crossings where we give the names $1, 2, 3, 4$ to the four points at infinity. Each of the six crossings involved in a quadruple crossing gets the name of the two crossing lines. 

For each of the adjacent eight strata of triple crossings we show the Gauss diagrams of the six crossings with the names of the crossings which are not in the triangle (the names of the latter can easily be established from the figure too). These twelve figures are our main instrument in this paper.

We make the following convention for the whole paper: crossings which survive in isotopies (i.e. they are not involved in Reidemeister moves of type I or II) are identified. In particular, they get the same name.

\begin{figure}
\centering
\includegraphics{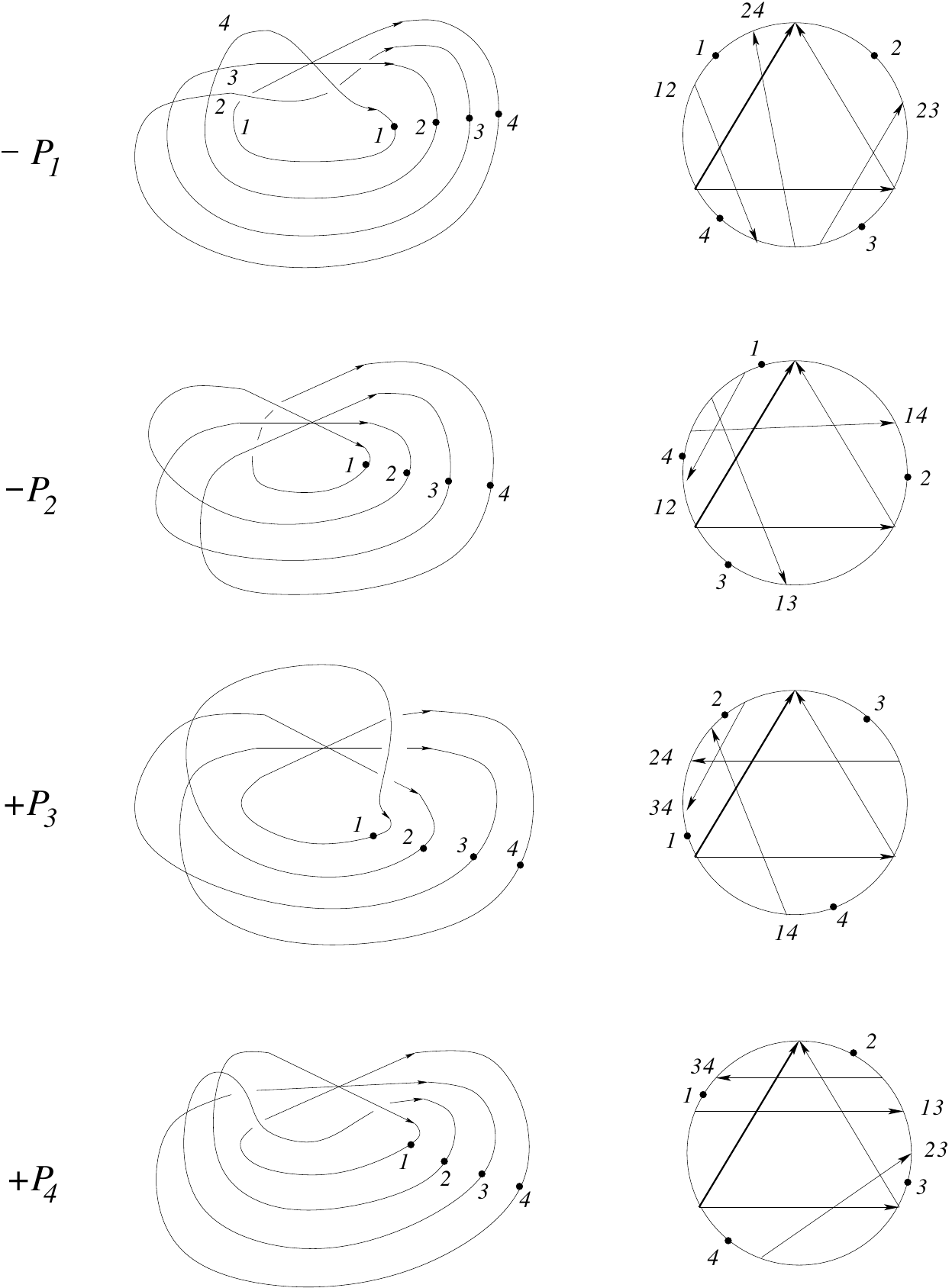}
\caption{\label{Iglob}  first half of the meridian for global type I}  
\end{figure}

\begin{figure}
\centering
\includegraphics{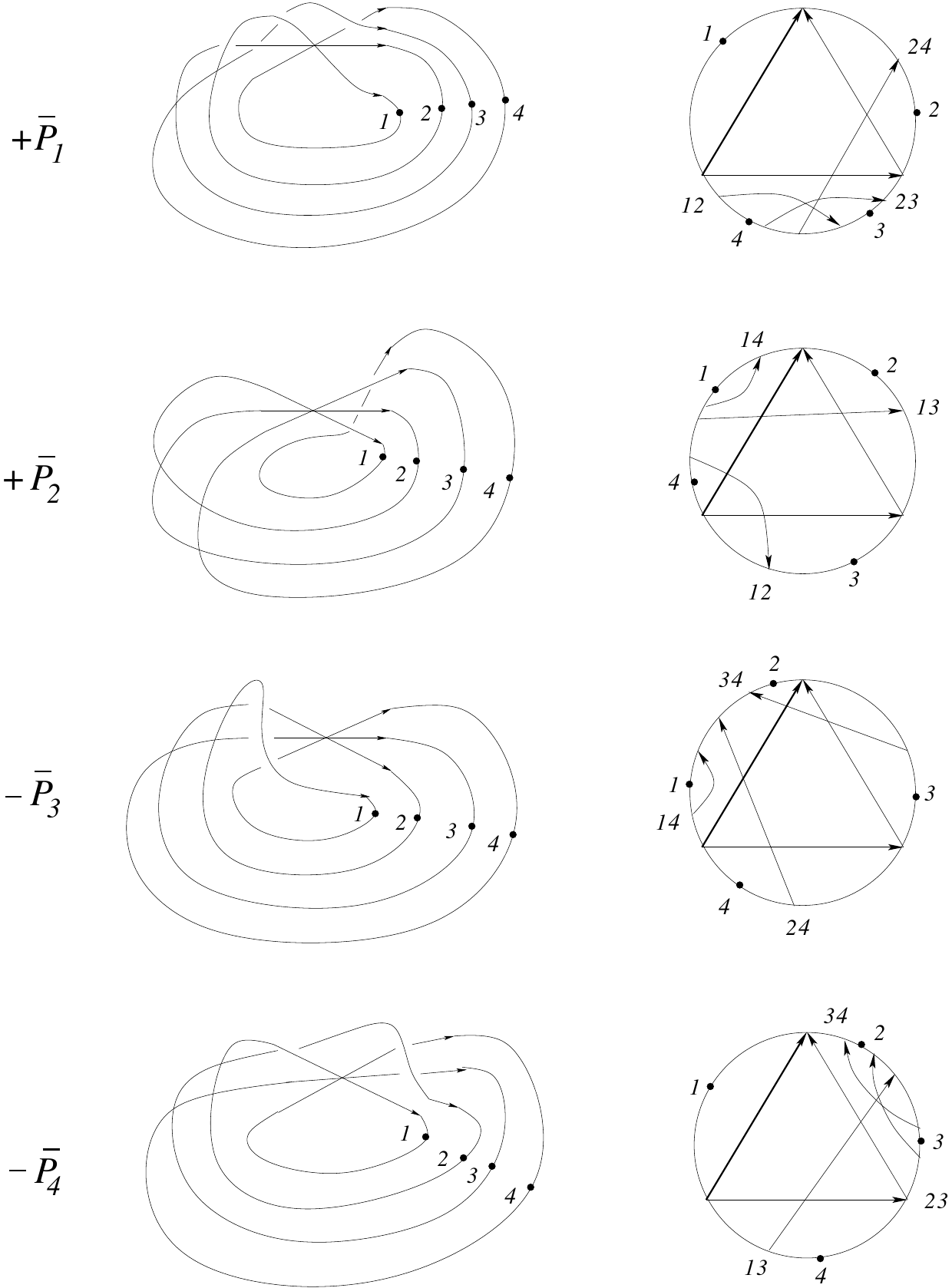}
\caption{\label{I2glob} second half of the meridian for global type I}  
\end{figure}

\begin{figure}
\centering
\includegraphics{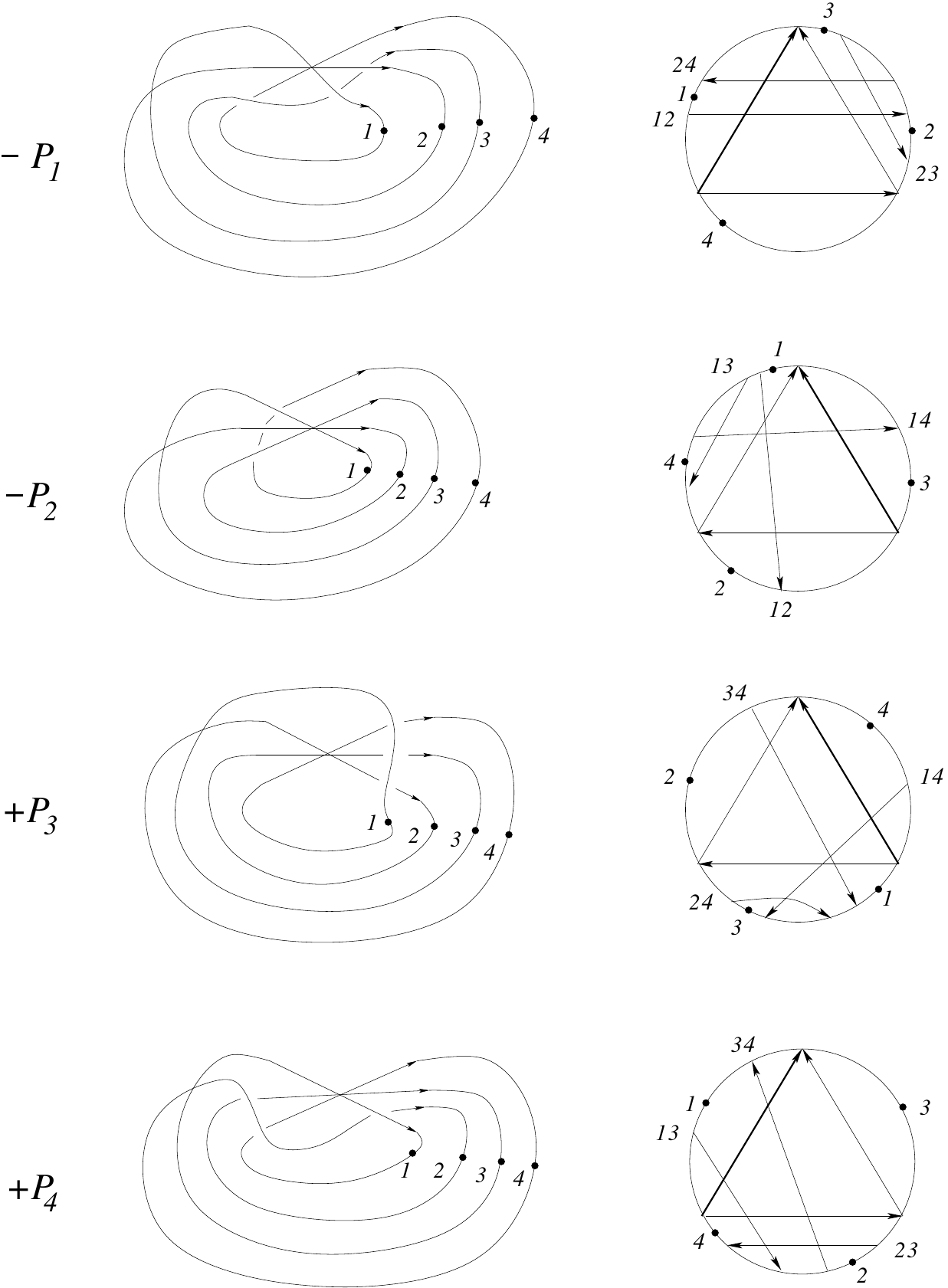}
\caption{\label{IIglob} first half of the meridian for global type II}  
\end{figure}

\begin{figure}
\centering
\includegraphics{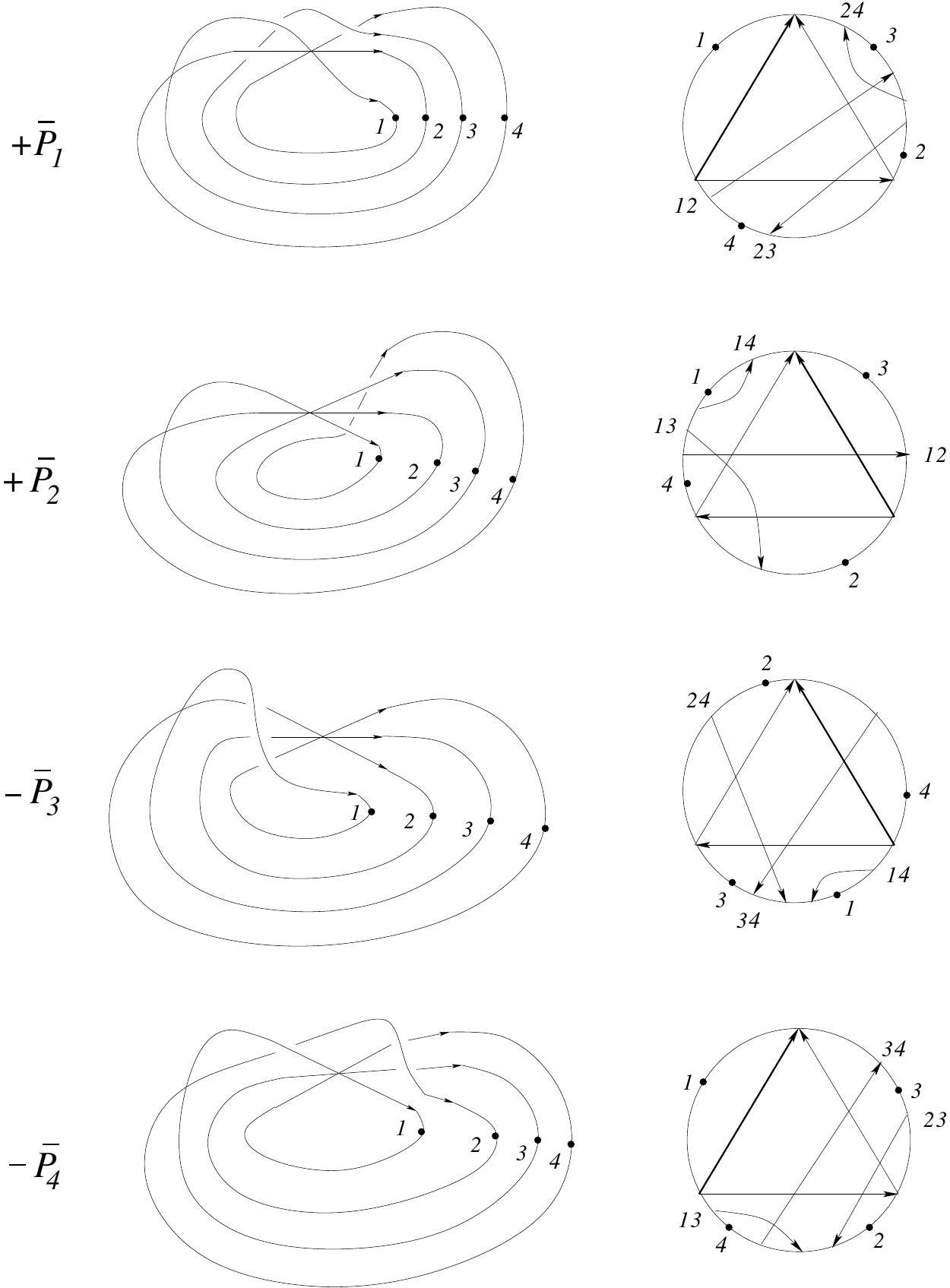}
\caption{\label{II2glob} second half of the meridian for global type II}  
\end{figure}

\begin{figure}
\centering
\includegraphics{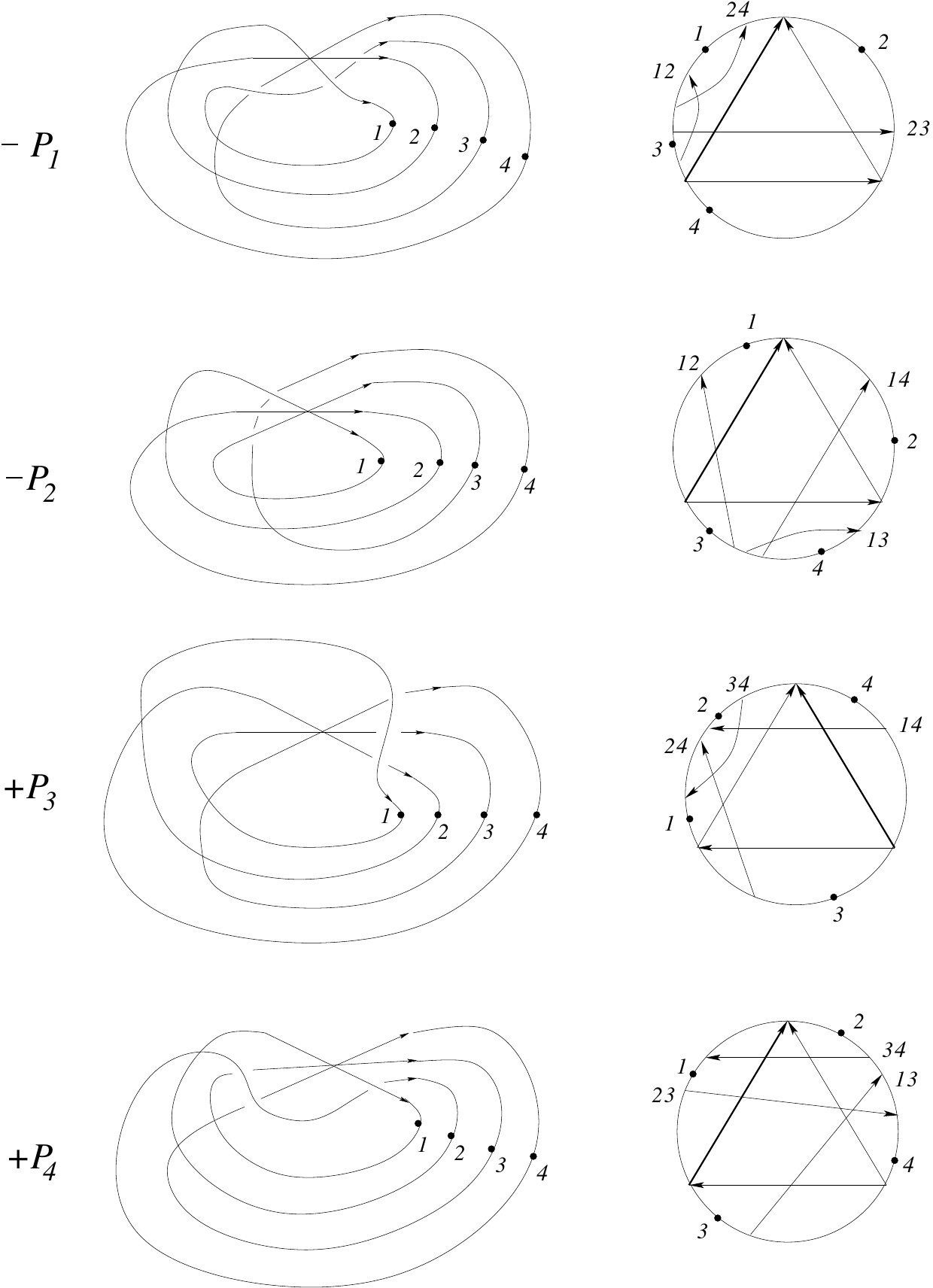}
\caption{\label{IIIglob} first half of the meridian for global type III}  
\end{figure}

\begin{figure}
\centering
\includegraphics{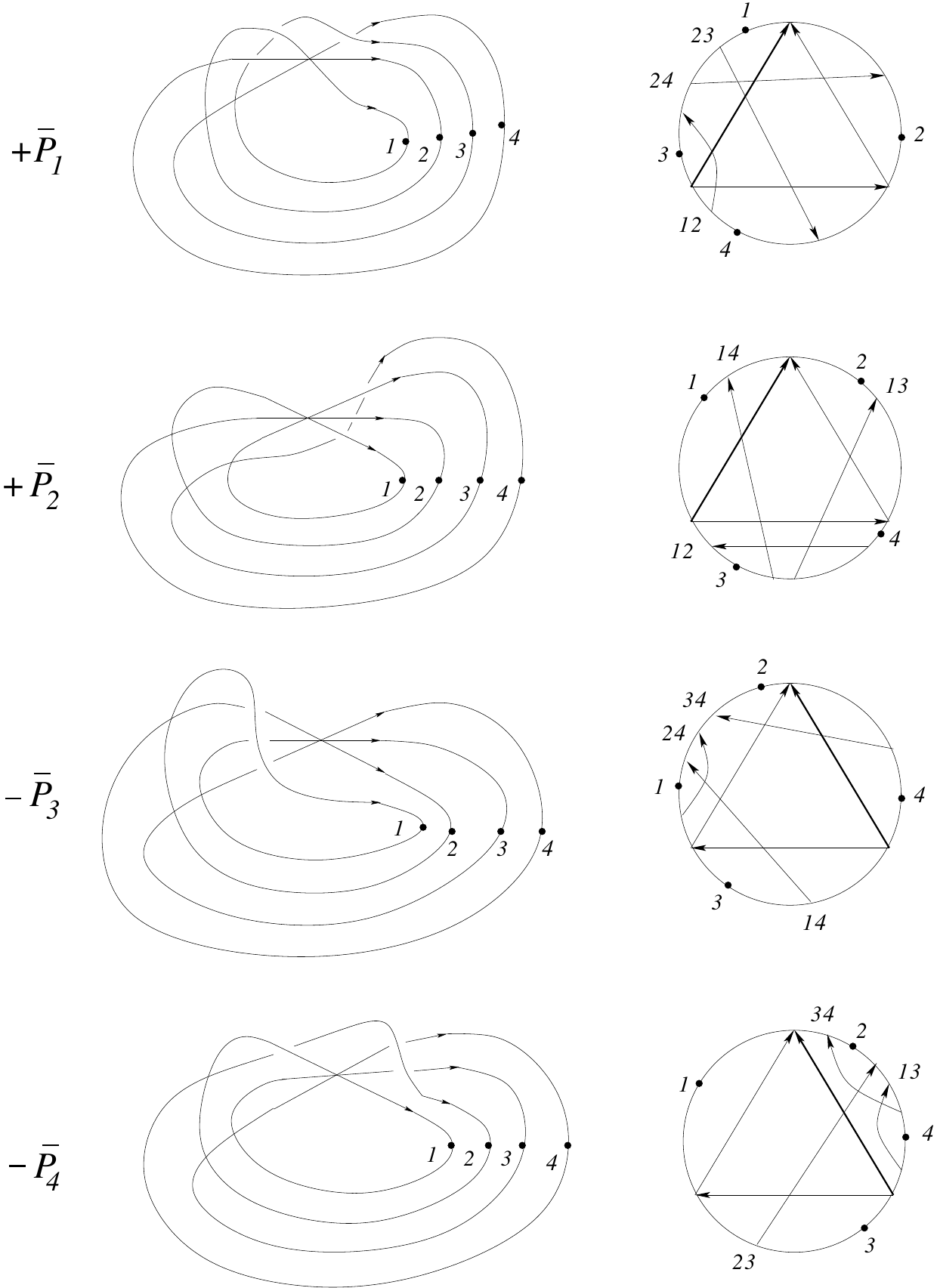}
\caption{\label{III2glob} second half of the meridian for global type III}  
\end{figure}

\begin{figure}
\centering
\includegraphics{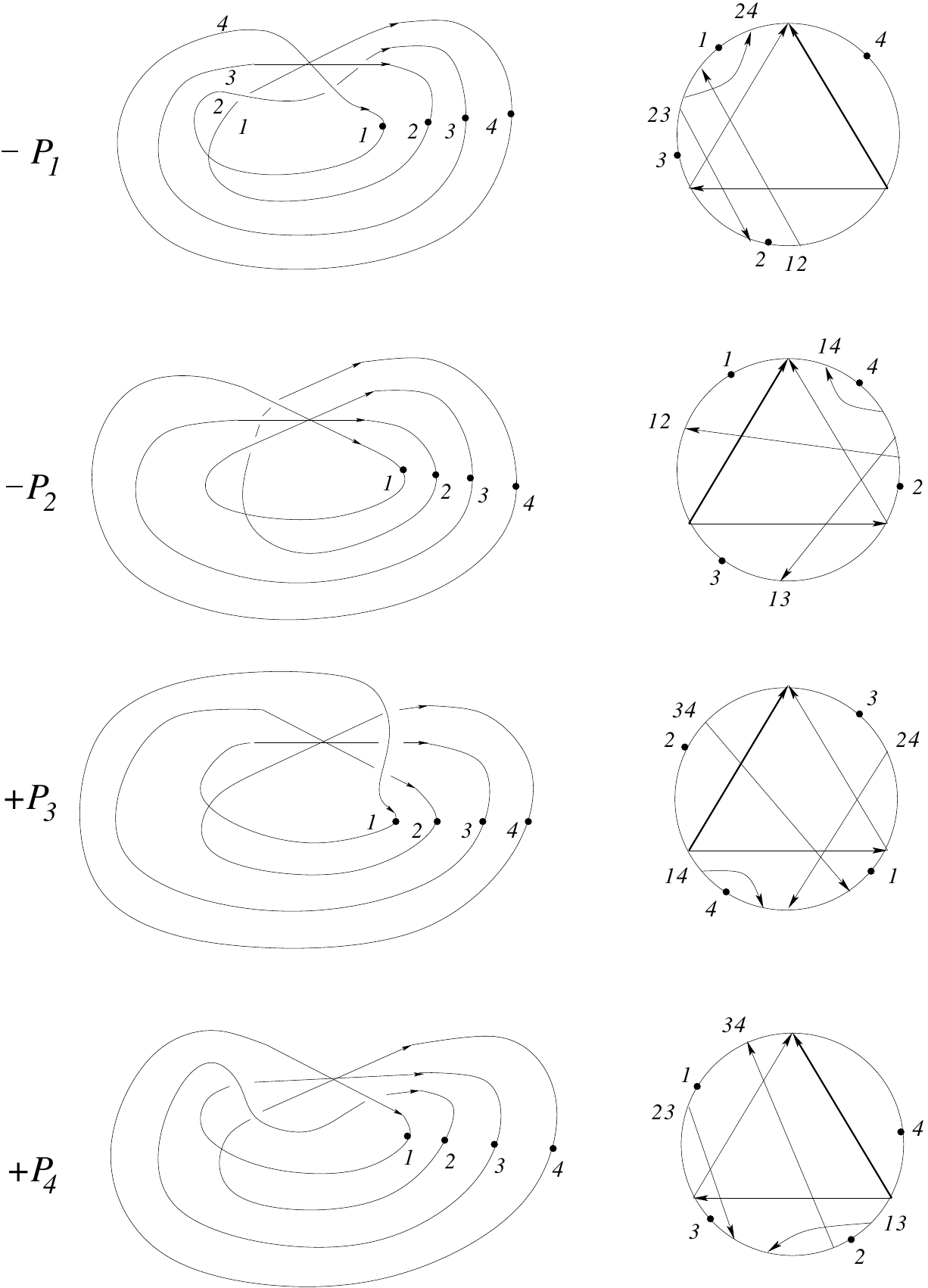}
\caption{\label{IVglob} first half of the meridian for global type IV}  
\end{figure}

\begin{figure}
\centering
\includegraphics{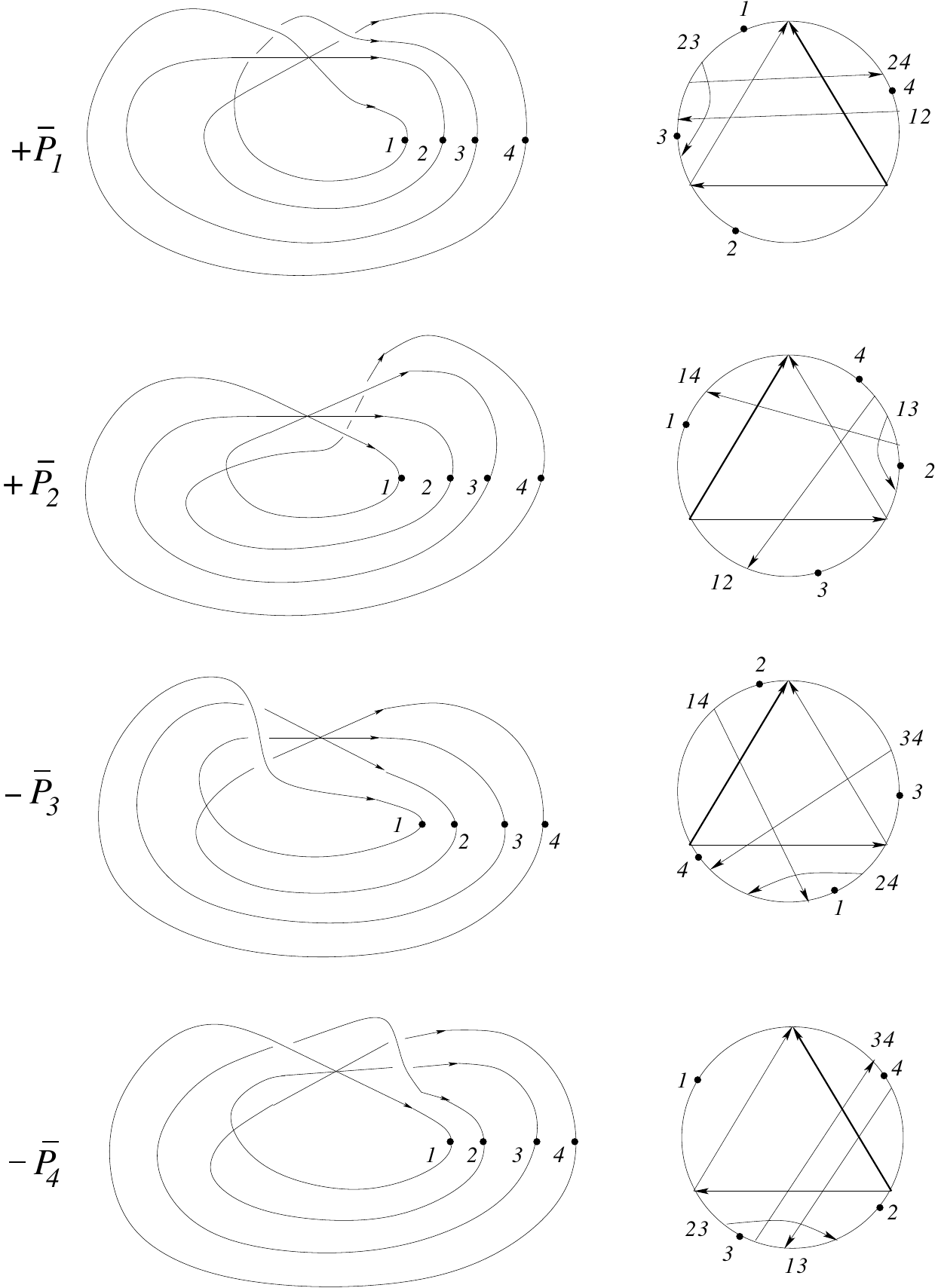}
\caption{\label{IV2glob} second half of the meridian for global type IV}  
\end{figure}

\begin{figure}
\centering
\includegraphics{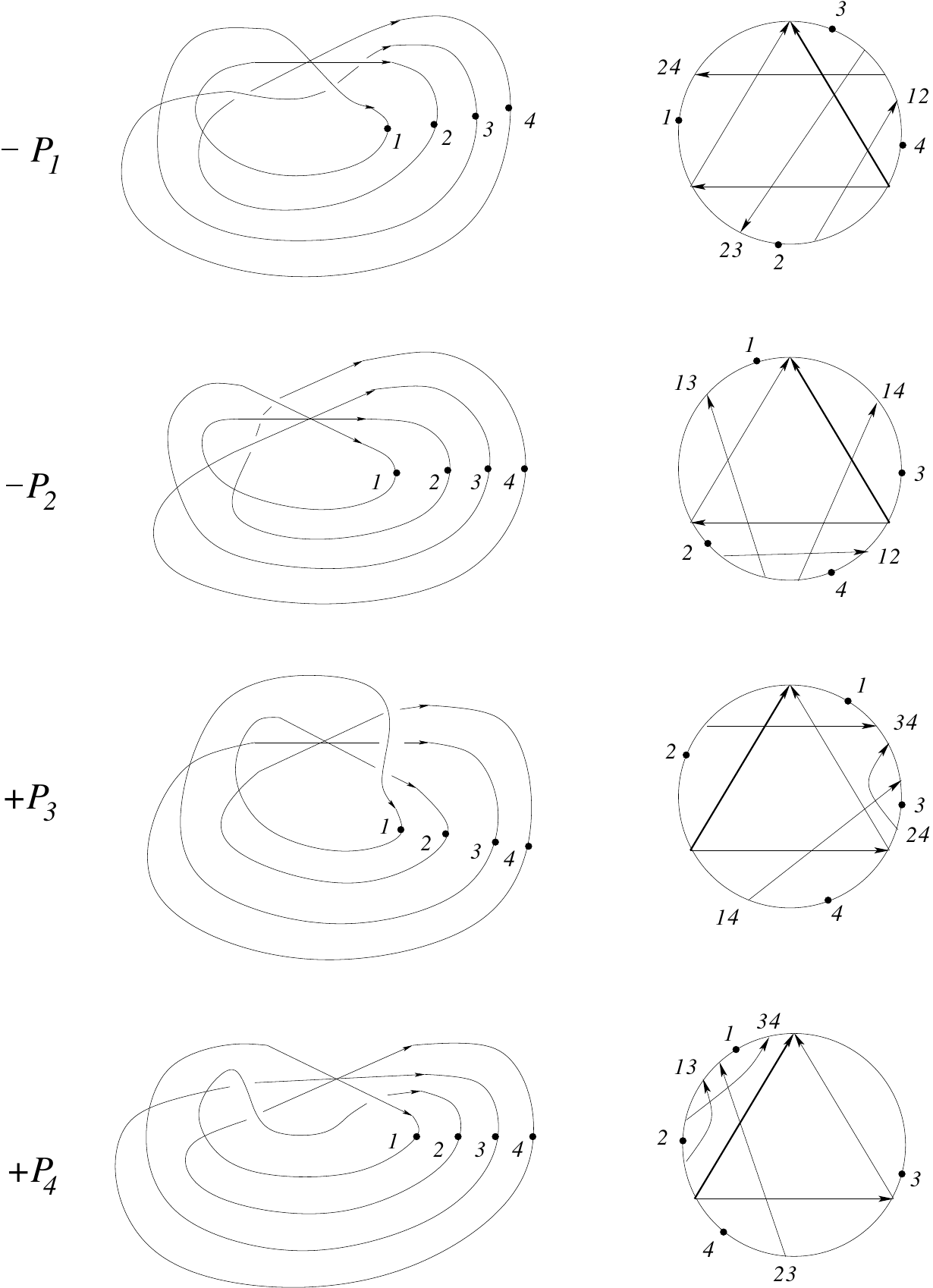}
\caption{\label{Vglob} first half of the meridian for global type V}  
\end{figure}

\begin{figure}
\centering
\includegraphics{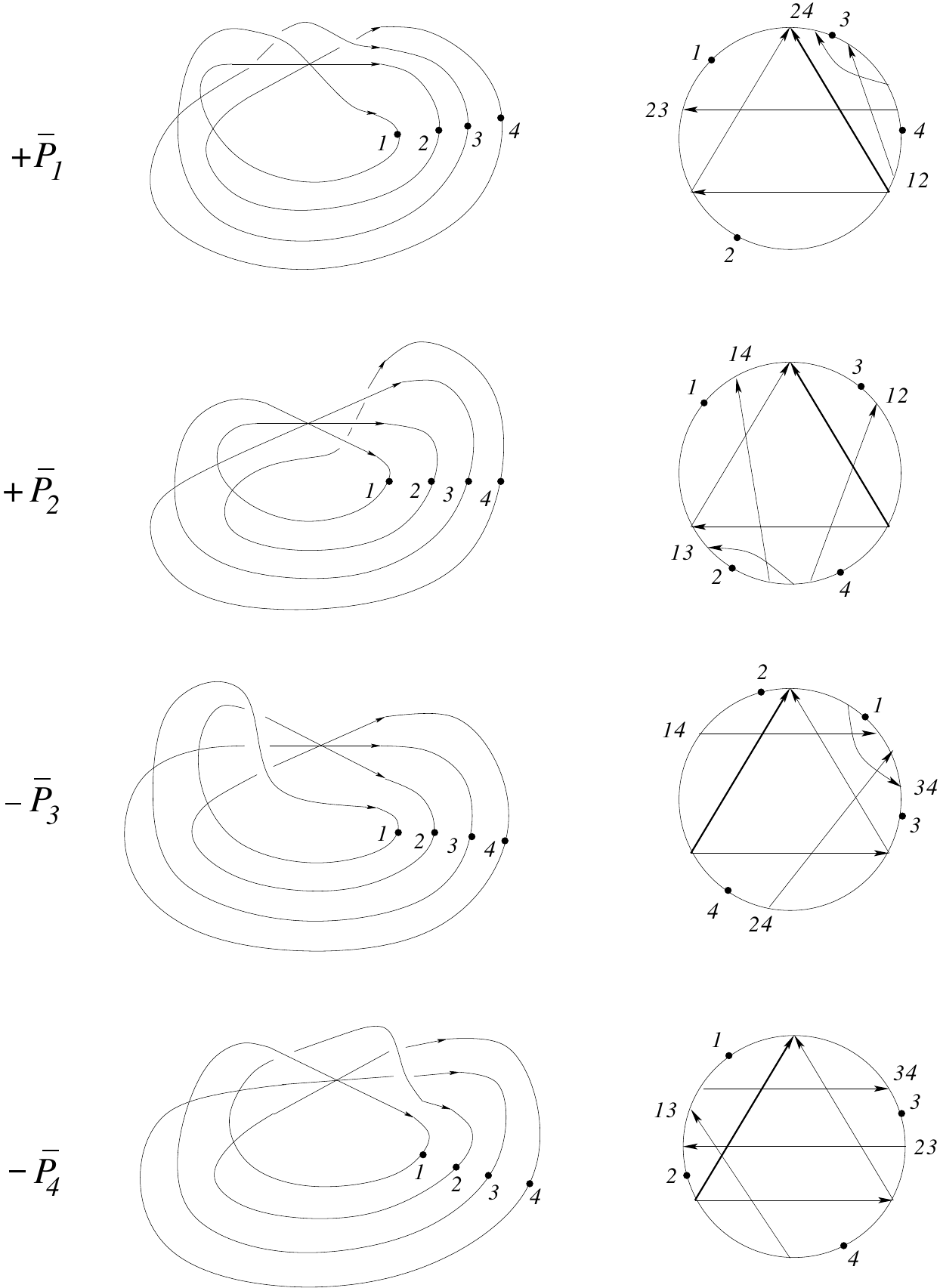}
\caption{\label{V2glob} second half of the meridian for global type V}  
\end{figure}

\begin{figure}
\centering
\includegraphics{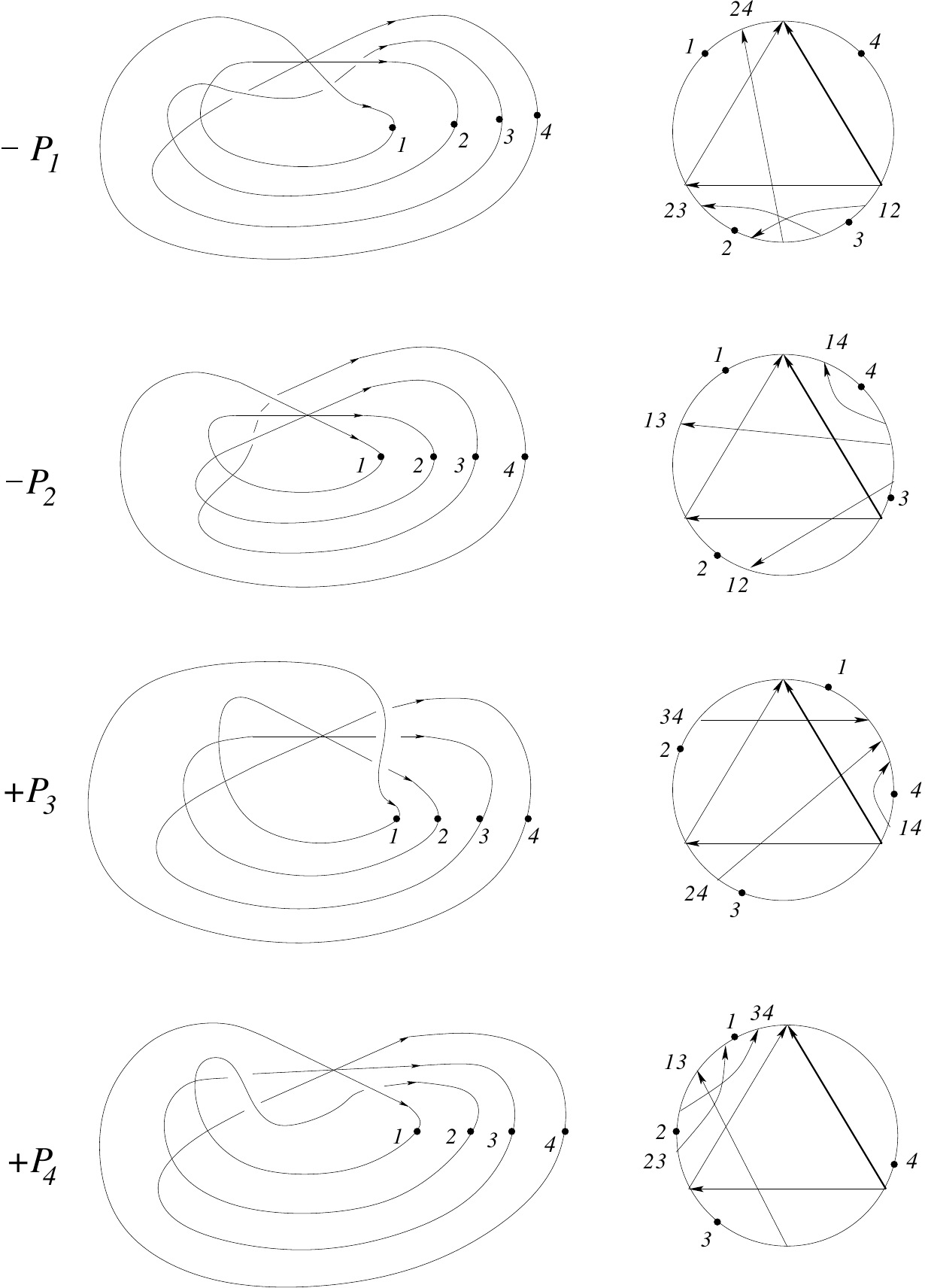}
\caption{\label{VIglob} first half of the meridian for global type VI}  
\end{figure}

\begin{figure}
\centering
\includegraphics{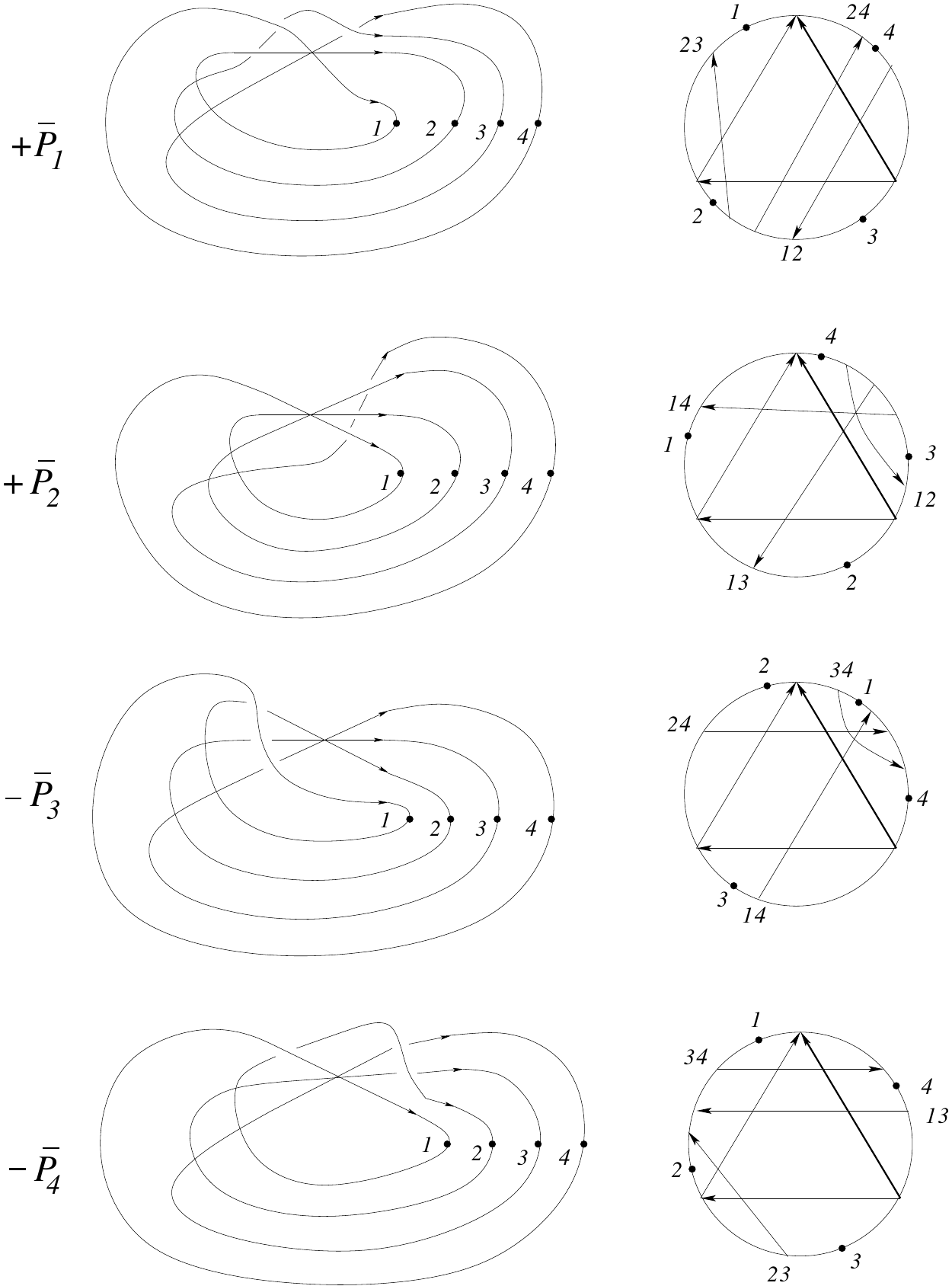}
\caption{\label{VI2glob} second half of the meridian for global type VI}  
\end{figure}

\begin{figure}
\centering
\includegraphics{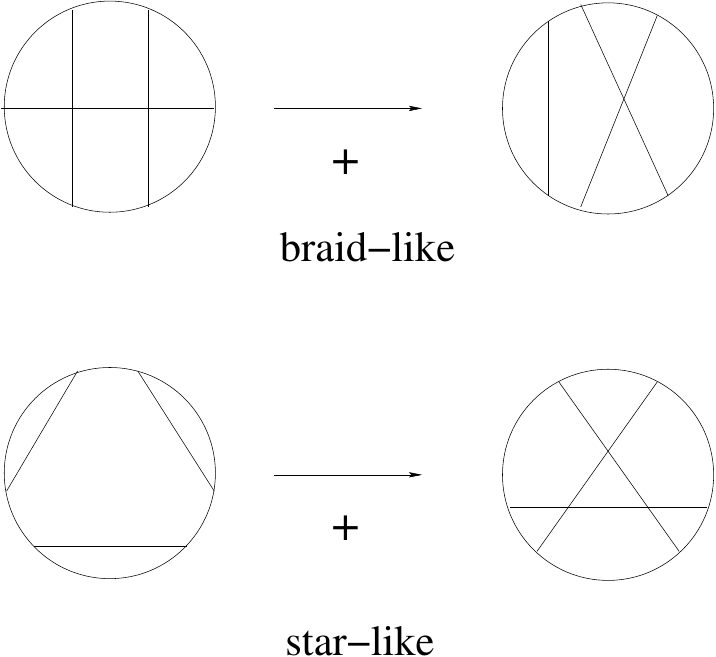}
\caption{\label{homoco}  the global co-orientation}  
\end{figure}

\begin{remark}
Our co-orientation for Reidemeister III moves is purely {\em local}. In \cite{F2} we have used a different co-orientation which is completely determined by the unoriented planar curves. We show it in Fig.~\ref{homoco}. Let us call it the {\em global co-orientation}. To the three crossings of a Reidemeister III move we associate chords in a circle which parametrizes the knot.
Notice that for positive triple crossings of global type $r$ the local and the global co-orientation coincide and that for the global type $l$ they are opposite. Comparing the global types $I$ and $III$ of positive quadruple points (compare Fig.~\ref{globquad}) we see that there is no solution with constant weight of the global positive tetrahedron equation for the global co-orientation. Indeed, for the type $I$ we have the equation $-P_1+P_4+\bar P_1-\bar P_4=0$ and for the type $III$ we have $-P_1-P_4+\bar P_1+\bar P_4=0$. Hence the use of our local co-orientation for triple crossings is essential in our approach!
\end{remark}

\section{Solution with linear weight of the positive global tetrahedron equation}

We denote the point at infinity in the circle $\sigma \cup T$ by $\infty$.

\begin{definition}
An ordinary crossing $q$ in a diagram $T$ with closure $\sigma$ is of {\em type 1} if $\infty \in D^+_q$ and is of {\em type 0} otherwise.
\end{definition}

Let $T(p)$ be a generic diagram with  a triple crossing or a self-tangency $p$ (in other words $T(p)$ is an interior point of $\bar \Sigma^{(1)}$) and let $d$ be the distinguished crossing for $p$. In the case of a self-tangency we identify the two distinguished crossings. {\em We assume that the distinguished crossing $d$ is of type 0.}

\begin{definition}
A crossing $q$ of $T(p)$ is called a {\em f-crossing} if $q$ is of type 1 and  the under cross of $q$ is in the oriented arc from $\infty$ to the over cross of $d$ in $\sigma \cup T$.
\end{definition}

\begin{figure}
\centering
\includegraphics{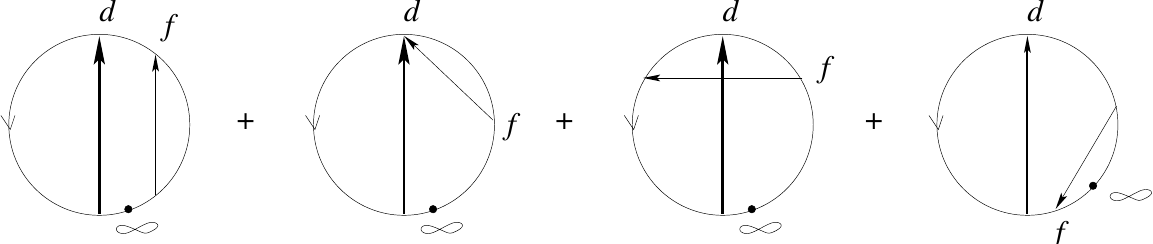}
\caption{\label{foot}  the f-crossings}  
\end{figure}

\begin{figure}
\centering
\includegraphics{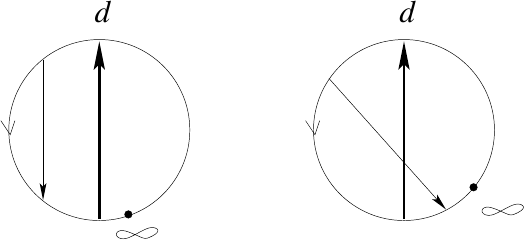}
\caption{\label{notfoot} crossings of type 1 which are not f-crossings}  
\end{figure}

In other words, the "foot" of $q$ in the Gauss diagram is in the sub arc of the circle which goes from $\infty$ to the head of $d$ and the head of $q$ is in the sub arc which goes from the foot of $q$ to $\infty$.   We illustrate this in Fig.~\ref{foot} and Fig.~\ref{notfoot}. (The letter "f" stands for "foot" or "Fu\ss" and not for "fiedler".) Notice that we make essential use of the fact that no crossing can ever move over the point at infinity.

\begin{definition}
The {\em linear weight $W(p)$} is  the Gauss diagram invariant shown in  Fig.~\ref{foot}. In other words, it is the sum of the writhes of all f-crossings with respect to $p$ in a given diagram $T(p)$. Notice that if $p$ is a self-tangency then the degenerate second configuration in  Fig.~\ref{foot} can not appear.
\end{definition}

\begin{definition}
Let $q$ be a f-crossing. The {\em grading} $A=\partial q$ of $q$ is defined by $\partial q = D^+_q \cap \partial T$. 

\end{definition}

The grading $\partial q$ is a cyclically ordered subset of $\partial T$ with alternating signs and which contains $\infty$ (compare the Introduction).

\begin{lemma}
The graded $W(p)$ is an  isotopy invariant for each Reidemeister move  of type II or III, i.e. $W(p)$ is invariant under any isotopy of the rest of the tangle outside of $I^2_p \times I$ (compare the Introduction).
\end{lemma}

In other words, $W(p)$ doesn't change under a homotopy of an arc which passes through $\Sigma^{(1)} \cap \Sigma^{(1)}$ in $M_T^{reg}$.

{\em Proof.} This is obvious. One has only to observe that the two new crossings from a Reidemeister move of type II are either both not f-crossings or are both f-crossings. In the latter case they have the same grading but they have different writhe. 

$\Box$

{\em We denote by $W_A(p)$ the invariant $W(p)$ restricted to all f-crossings with a given grading $A$.}

 It is clear that $W_A(p)$ is an invariant of degree 1 for $(T,p)$ for each fixed grading $A$ (see e.g. \cite{PV} or \cite{F1} for a proof that Gauss diagram invariants are of finite type).

Let us consider a global positive quadruple crossing with a fixed point at infinity (compare the Introduction and Section 3).
The following lemma is of crucial importance.

\begin{figure}
\centering
\includegraphics{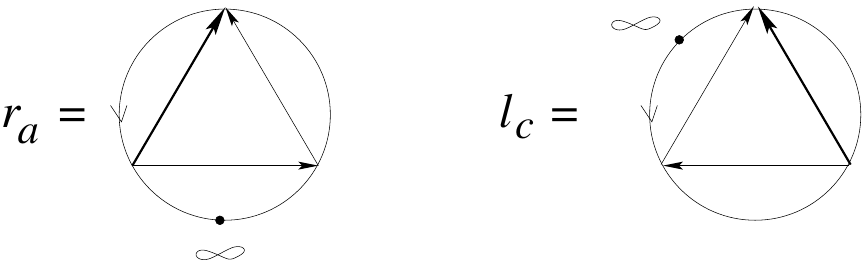}
\caption{\label{fglob} the global types for which a new f-crossing in $\bar P_3$ appears}  
\end{figure}

\begin{lemma}
The f-crossings of the eight adjacent strata of triple crossings have the following properties:
\begin{itemize}

\item the f-crossings in  $P_2$ are identical with those in $\bar P_2$ (with our convention above)

\item the f-crossings in $P_1$, $\bar P_1$, $P_4$ and $\bar P_4$ are all identical

\item the f-crossings in  $P_3$ are either identical with those in $\bar P_3$ or there is exactly one new f-crossing in 
$\bar P_3$ with respect to  $P_3$. In the latter case the new crossing is always exactly the crossing $hm$ in 
$P_1$ and $\bar P_1$

\item the new f-crossing in $\bar P_3$ appears if and only if $P_1$ (and hence also $\bar P_1$) is of one of the two global types shown in Fig.~\ref{fglob}.

\end{itemize}
\end{lemma}

{\em Proof.} We have checked  the assertions of the lemma in all twenty four cases (denoted by the global type of the quadruple crossing together with the point at infinity)  using the figures  Fig.~\ref{Iglob}, ... Fig.~\ref{VI2glob}. Notice that the crossing $hm$ in $P_1$ and $\bar P_1$ is always the crossing $34$ and that $d$ is of type 0 if $\infty$ is on the right side of it in the figures.
We consider just some examples. In particular we show that it is necessary to add the "degenerate" configuration in Fig.~\ref{foot} and we left the rest of the verification to the reader.

The f-crossings are:

case $I_1$. non at all

case $I_2$. non at all 

case $I_3$. In $P_1, \bar P_1$: $34$ (both are the degenerate case). In $P_4, \bar P_4$: $34$ (the third configuration and the first configuration in Fig.~\ref{foot}). In $P_2, \bar P_2$: $34$. In $P_3$: non. In $\bar P_3$: $34$.

 case $I_4$. In $P_1, \bar P_1$, $P_4, \bar P_4$: $34, 24, 23$. In $P_2, \bar P_2$: non. In $P_3$: $23, 24$. In $\bar P_3$: $23, 24, 34$.

case $VI_3$: In $P_2$: $12, 13, 14$. In $\bar P_2$: $12, 13, 14$ ($12$ shows that the fourth configuration in Fig.~\ref{foot} is necessary too). 

case $VI_1$: In $P_3$: non. In $\bar P_3$: $34$ (shows that the second configuration in Fig.~\ref{fglob} is necessary too).

$\Box$

Let $\gamma$ be an oriented generic arc in $M^{reg}_T$ which intersects $\Sigma^{(1)}$ only in positive triple crossings and let $A$ be a fixed element of the gradings.

\begin{definition}
The evaluation of the 1-cochain $R^{(1)}_{reg}$ with grading $A$ on $\gamma$  is defined by \vspace{0,2 cm}

 $R^{(1)}_{reg}(A)(\gamma)=\sum_{p \in \gamma} sign(p) \sigma_2\sigma_1(p) + \sum_{p \in \gamma}sing(p)zW_A(p)(\sigma_2(p)-\sigma_1(p))$     \vspace{0,2 cm}

where the first sum is over all triple crossings of the global types shown in Fig.~\ref{fglob} (i.e. the types $l_c$ and $r_a$) and such that $\partial  (hm)=A$.  The second sum is over all triple crossings $p$ which have a distinguished crossing $d$ of type 0 (i.e. the types $r_a$, $r_b$ and $l_b$).

\end{definition}

(Remember that e.g. $\sigma_2\sigma_1(p)$ denotes the element in the HOMFLYPT skein module $S(\partial T)$ which is obtained by replacing the triple crossing $p$ in the tangle $T$ by the partial smoothing $\sigma_2\sigma_1$, compare the Introduction.)

Notice that a triple crossing of the first type in Fig.~\ref{fglob} could also contribute to the second sum in $R^{(1)}_{reg}(A)(\gamma)$ but a triple crossing of the second type  in Fig.~\ref{fglob} can't because $d$ is of type 1.

\begin{proposition}
Let $m$ be the meridian for a positive quadruple crossing. Then  $R^{(1)}_{reg}(A)(m)=0$ for each grading $A$.
\end{proposition}
{\em Proof.} We have to consider the contributions of the eight strata $P_1$,...$\bar P_4$.

{\em Case 1:  The f-crossings of grading $A$ are identical in $P_3$ and $\bar P_3$.}

It follows from Lemma 2 that either the new f-crossing in $\bar P_3$ is not of grading $A$ or there is no new f-crossing in $\bar P_3$ at all. Moreover, in the first case the crossing $hm$ in $P_1$ is not of grading $A$ and hence $P_1$ and $\bar P_1$ do not contribute to the first sum in  $R^{(1)}_{reg}(A)(m)$. In the second case the triple crossing $P_1$ (and consequently $\bar P_1$ too) is not of the type required in Fig.~\ref{fglob} and hence again does not contribute to the first sum.   

It follows from Lemma 2 that $W_A(P_2)=W_A(\bar P_2)$ as well as $W_A(P_3)=W_A(\bar P_3)$. Moreover, each partial smoothing of $P_2$ (respectively $P_3$) is regularly isotopic to the same partial smoothing of $\bar P_2$ (respectively $\bar P_3$). (This follows from the fact that the fourth branch just moves over or under everything.) But $P_i$ and $\bar P_i$ enter always with different signs. It follows that all contributions of these strata to  $R^{(1)}_{reg}(A)(m)$ cancel out.

The contributions of $P_4$ and $\bar P_4$ to the first sum in  $R^{(1)}_{reg}(A)(m)$ cancel out as was shown in Fig.~\ref{tablesmooth}.

It follows again from Lemma 2 that $P_1$, $\bar P_1$, $P_4$, $\bar P_4$ share always the same $W_A$. Consequently, the contribution of these strata to the second sum in $R^{(1)}_{reg}(A)(m)$ is just a multiple of the solution with constant weight (compare the Introduction) and hence it is 0 by Proposition 1.

{\em Case 2:  The f-crossings of grading $A$ are not identical in $P_3$ and $\bar P_3$.}

In this case the new f-crossing in $\bar P_3$ is of grading $A$. It follows from Lemma 2 that the crossing $hm$ in $P_1$ (and hence in $\bar P_1$ too) is of grading $A$. $W_1(\bar P_3) = W_A(P_3) +1$ because all crossings are positive. Using Fig.~\ref{unfoldquad} we find the contributions of $P_1$, $\bar P_1$, $P_3$, $\bar P_3$ to $R^{(1)}_{reg}(A)(m)$ and the calculation was already given in Fig.~\ref{tablesmooth}.
 The rest of the arguments is the same as in Case 1.

$\Box$

We have proven that  $R^{(1)}_{reg}(A)$ is a solution of the positive global tetrahedron equation.

\begin{remark}
There are of course {\em dual} solutions. One of them  uses the same global types as shown in Fig.~\ref{fglob}, but it replaces triple crossings with distinguished crossing $d$ of type $0$ by those of type $1$. It uses now  $P_2$ and $\bar P_2$ instead of $P_3$ and $\bar P_3$, $P_4$ and $\bar P_4$ instead of $P_1$ and $\bar P_1$ and the crossing $ml$ in $P_4$ instead of the crossing $hm$ in $P_1$. The f-crossings become crossings of type $0$ and such that the "heads" of the arrows are in the oriented arc from $\infty$ to the undercross of $d$. 
 We haven't carried this out in detail. But it would be of course interesting to find out whether the corresponding 1-cocycles are independent of our 1-cocycles $R^{(1)}_{reg}(A)$ and $R^{(1)}$.

Notice that there should be also 1-cocycles which use the global types $l_a$ and $r_b$ instead of the types $l_c$ and $r_a$ (compare Fig.~\ref{globtricross}).
\end{remark}

\section{Solution with quadratic weight of the positive global tetrahedron equation}

\begin{figure}
\centering
\includegraphics{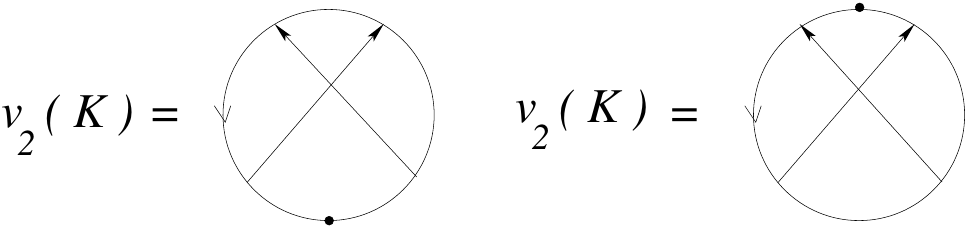}
\caption{\label{PV}  the Polyak-Viro formulas for $v_2(K)$}  
\end{figure}

We start with recalling the Gauss diagram formula of Polyak and Viro for the Vassiliev invariant $v_2(K)$ (see \cite{PV}). In fact, there are two formulas shown in Fig.~\ref{PV}. These formulas are for knots in the 3-sphere and the marked point can be chosen arbitrary on the knot. In particular, in the case of long knots we chose of course $\infty$ as the marked point.

Let $T(p)$ be a generic diagram with  a triple crossing or a self-tangency $p$ (in other words $T(p)$ is an interior point of $\Sigma^{(1)}$) and let $d$ be the distinguished crossing for $p$. We assume that $d$ is of type 0.

\begin{figure}
\centering
\includegraphics{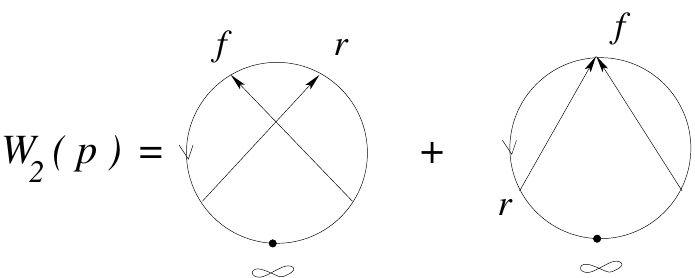}
\caption{\label{W2} the quadratic weight $W_2$}  
\end{figure}

\begin{definition}
The {\em quadratic weight} $W_2(p)$ is defined by the Gauss diagram formula shown in Fig.~\ref{W2}. Here $f$ denotes an f-crossing of $p$ and $r$ denotes an arbitrary crossing with a position as shown in the configurations.
\end{definition}
(Remember that with our conventions the formula in Fig.~\ref{W2} says that $W_2(p)$ is the sum of $w(f)w(r)$ over all couples of crossings $(f,r)$ in $T(p)$ which form one of the  configurations from Fig.~\ref{W2}.)

The second (degenerate) configuration in Fig.~\ref{W2} can of course not appear for a self-tangency $p$. Notice that the crossing $r$ is always of type 0. 

The crossings $r$ in Fig.~\ref{W2} are called the {\em r-crossings}.

\begin{lemma}
$W_2(p)$ is an  isotopy invariant for each Reidemeister move of type I, II or III, i.e. $W_2(p)$ is invariant under any isotopy of the rest of the tangle outside of $I^2_p \times I$ (compare the Introduction).
\end{lemma}
The lemma implies that as in the linear case $W_2(p)$ doesn't change under a homotopy of an arc which passes through $\Sigma^{(1)} \cap \Sigma^{(1)}$, but this time even in $M_T$. However, we have lost the grading.

\begin{figure}
\centering
\includegraphics{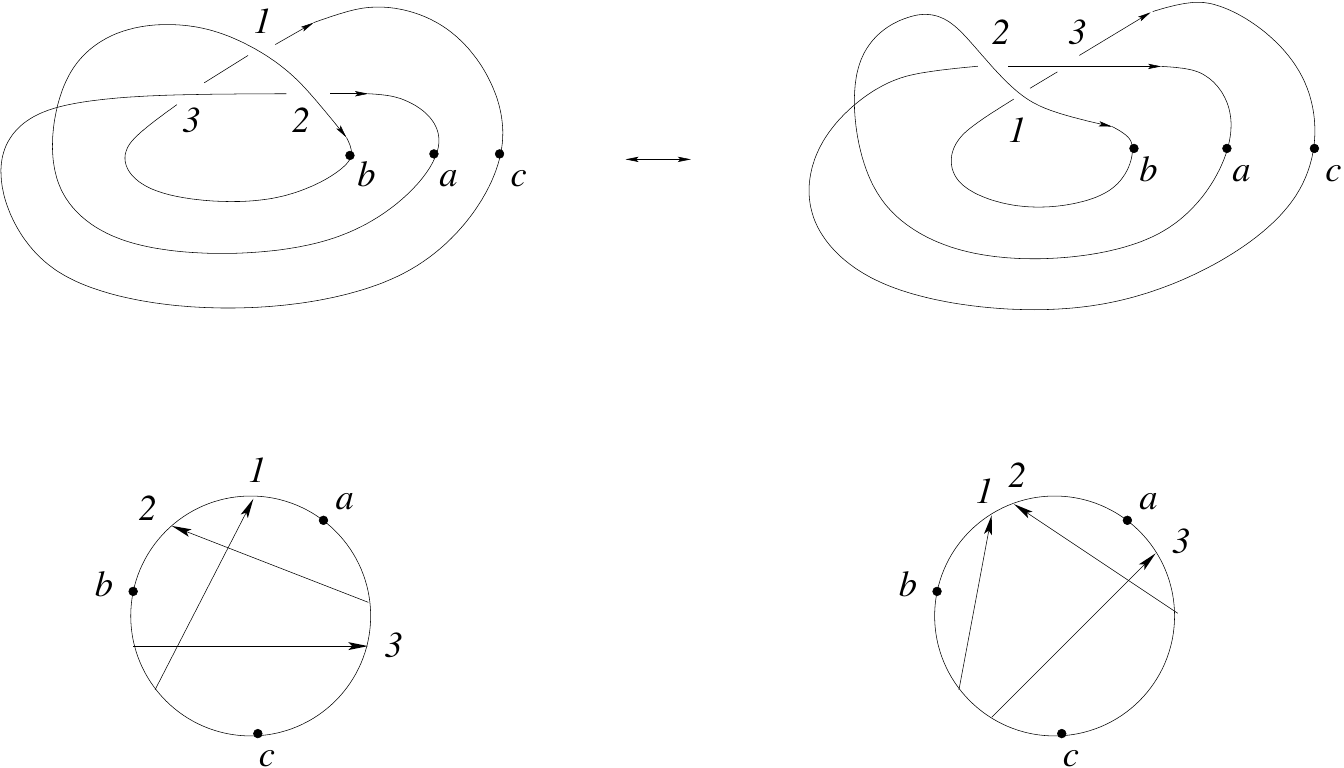}
\caption{\label{III1typer} $RIII$ of type $r$}  
\end{figure}

\begin{figure}
\centering
\includegraphics{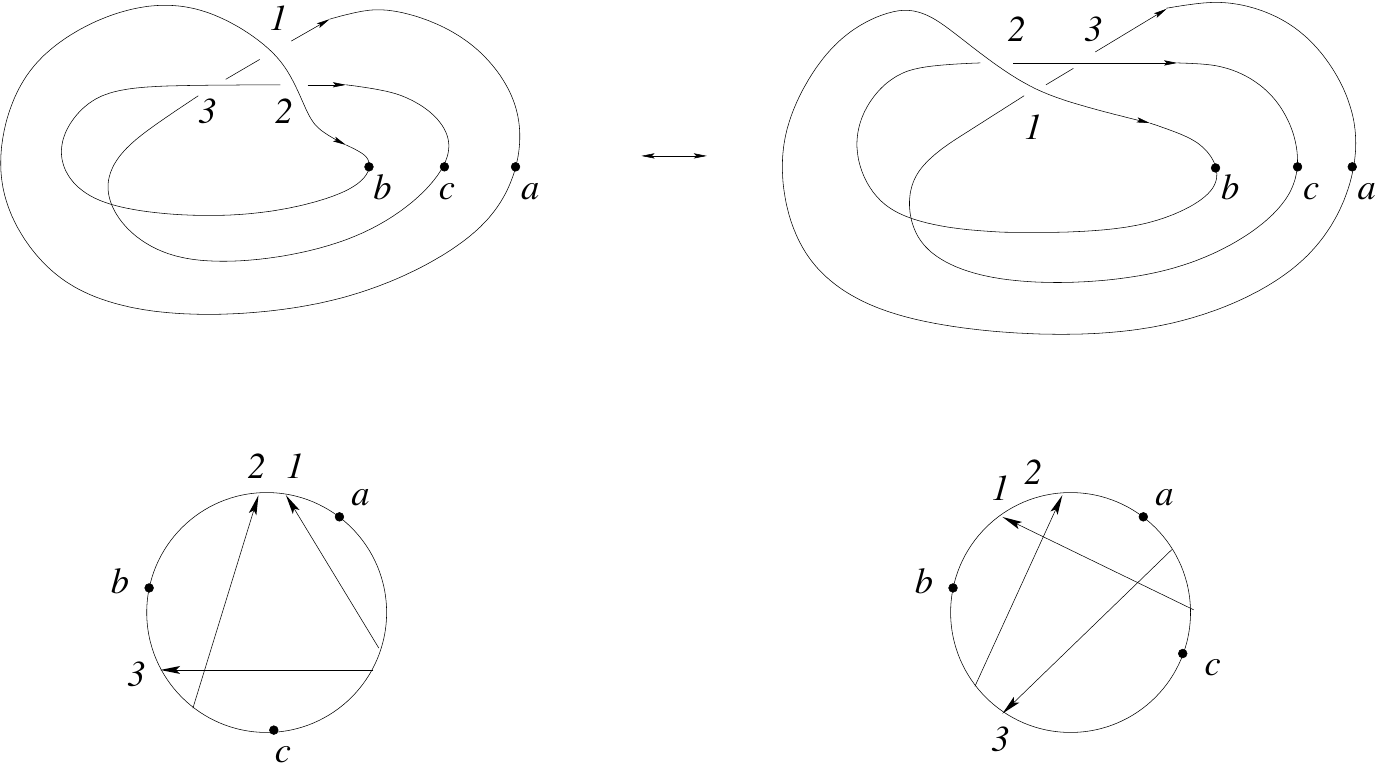}
\caption{\label{III2typel}  $RIII$ of type $l$}  
\end{figure}

{\em Proof.} It is obvious that $W_2(p)$ is invariant under Reidemeister moves of type I and II. The latter comes from the fact that both new crossings are simultaneously crossings of type $f$ or $r$ and that they have different writhe. As was explained in the Introduction the graph $\Gamma$ implies now that it is sufficient to prove the invariance of $W_2(p)$ only under positive Reidemeister moves of type III. There are two global types and for each of them there are three possibilities for the point at infinity. We give names $1, 2, 3$ to the crossings and $a, b, c$ to the points at infinity and we show the six cases in Fig.~\ref{III1typer} and Fig.~\ref{III2typel}. Evidently, we have only to consider the mutual position of the three crossings in the pictures because the contributions with all other crossings do not change.

$r_a$: there is only 3 which could be a f-crossing but it does not contribute with another crossing in the picture  to $W_2(p)$. 

$r_b$: no r-crossing at all.

$r_c$: only 2 could be a f-crossing. In that case it contributes on the left side exactly with the r-crossing 1 and on the right side exactly with the r-crossing 3.

$l_a$: no f-crossing at all.

$l_b$: 1 and 2 could be  f-crossings but non of them contributes with the crossing 3 to $W_2(p)$.

$l_c$: 1 and 3 could be f-crossings. But the foots of 1 and 3 are arbitrary close. Consequently, they can be f-crossings only simultaneously! In that case 3 contributes with 2 on the left side and 1 with 2 on the right side.

Consequently, we have proven that $W_2(p)$ is invariant. But we have lost the grading in $l_c$. The f-crossings which contribute are not the same on the left side and on the right side and we can not guaranty that the crossings 1 and 3 have the same grading.

$\Box$

We want to replace $W(p)$ by $W_2(p)$ for the partial smoothing $z(\sigma_2(p)-\sigma_1(p))$. 

{\em Let us denote by $W_2(p,f)$ the contribution of a given f-crossing $f$ to $W_2(p)$.}

 In the case that the f-crossings in $P_3$ and $\bar P_3$ are not the same we get now \vspace{0,2 cm}

$W_2(\bar P_3) = W_2(P_3) + W_2(p,f)$
\vspace{0,2 cm}

where $f$ is the new f-crossing in $\bar P_3$. But luckily the new f-crossing $f$ is always just the crossing $hm$ in $P_1$ (compare Lemma 2). This leads us to the following definition.

\begin{figure}
\centering
\includegraphics{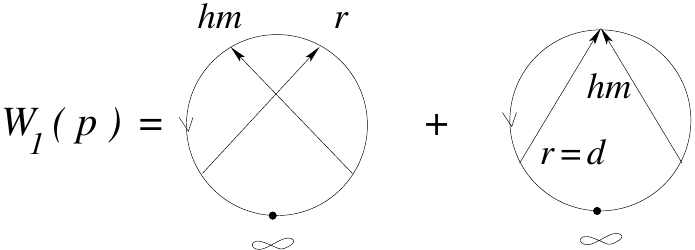}
\caption{\label{W1} the linear weight $W_1$}  
\end{figure}

\begin{definition}
Let $p$ be a triple crossing of one of the types shown in Fig.~\ref{fglob}. The linear weight $W_1(p)$ is defined as the sum of the writhes of the crossings $r$ in $T(p)$ which form one of the configurations shown in Fig.~\ref{W1}.
\end{definition}

The linear weight $W_1(p)$ is defined with respect to the single f-crossing $hm$ (normally we should denote it by $W_1(p,hm)$).
Notice that we do not multiply by $w(hm)$. For the positive tetrahedron equation this wouldn't make any difference. But we will see later that this is forced by the cube equations.

Let $\gamma$ be an oriented generic arc in $M^{reg}_T$ which intersects $\Sigma^{(1)}$ only in positive triple crossings.

\begin{definition}
The evaluation of the 1-cochain $R^{(1)}$ on $\gamma$  is defined by \vspace{0,2 cm}

 $R^{(1)}(\gamma)=\sum_{p \in \gamma} sign(p)W_1(p) \sigma_2\sigma_1(p) + \sum_{p \in \gamma}sing(p)zW_2(p)(\sigma_2(p)-\sigma_1(p))$     \vspace{0,2 cm}

where the first sum is over all triple crossings of the global types shown in Fig.~\ref{fglob}.  The second sum is over all triple crossings $p$ which have a distinguished crossing $d$ of type 0.

\end{definition}

The following lemma is a "quadratic refinement" of Lemma 2.
\begin{lemma}
(1) $W_2(P_1)=W_2(\bar P_1)=W_2(P_4)=W_2(\bar P_4)$

(2) $W_2(P_2)= W_2(\bar P_2)$

(3) Let $P_i$ for some  $i \in\{1, 2\}$ be of one of the types shown in Fig.~\ref{fglob}. Then $W_1(P_i)=W_1(\bar P_i)$.

(4) If the f-crossings in $P_3$ and $\bar P_3$ coincide then $W_2(P_3)= W_2(\bar P_3)$

(5) Let $f$ be the new f-crossing in $\bar P_3$ with respect to $P_3$ and let $f=hm$ be the corresponding crossing in $P_1$. Then $W_2(\bar P_3)-W_2(P_3)=W_2(\bar P_3,f)$ and
$W_2(\bar P_3,f)=W_1(P_1)=W_1(\bar P_1)$.

(6) Let $P_i$ for some  $i \in\{3, 4\}$ be of one of the types shown in Fig.~\ref{fglob}. Then either simultaneously $W_1(P_3)=W_1(\bar P_3)$ and $W_1(P_4)=W_1(\bar P_4)$ or simultaneously $W_1(P_3)= W_1(\bar P_3) +1$ and  $W_1(\bar P_4)= W_1(P_4) +1$.

\end{lemma}
{\em Proof.} 
The proof is by inspection of all f-crossings, all r-crossings and all hm-crossings in all twenty four cases in the figures. We give it below in details. (Those strata which can never contribute are dropped.) It shows in particular  that it is necessary to add the degenerate configurations in Fig.~\ref{W2} and Fig.~\ref{W1}.

$I_1$: nothing at all
\vspace{0,2 cm}

$I_2$: $W_2(P_1)=W_2(\bar P_1)=W_2(P_4)=W_2(\bar P_4)=0$. $W_2(P_2)=W_2(\bar P_2)=0$.
\vspace{0,2 cm}

$I_3$: $W_2(P_1)=W_2(\bar P_1)=W_2(P_4)=W_2(\bar P_4)=2$. $W_2(P_2)=W_2(\bar P_2)=2$.
$W_1(P_1)=W_1(\bar P_1)=2$. $W_1(P_2)=W_1(\bar P_2)=2$. $W_2(P_3)=0$ and $W_2(\bar P_3)=W_2(\bar P_3,f)=2$.
\vspace{0,2 cm}

$I_4$: $W_2(P_1)=W_2(\bar P_1)=W_2(P_4)=W_2(\bar P_4)=4$.  
$W_1(P_1)=W_1(\bar P_1)=1$. $W_1(P_3)=W_1(\bar P_3)=1$. $W_1(P_4)=W_1(\bar P_4)=2$. 
$W_2(\bar P_3)-W_2(P_3)=W_2(\bar P_3,f)=1$.
\vspace{0,2 cm}

$II_1$: $W_1(P_2)=W_1(\bar P_2)=0$.
\vspace{0,2 cm}

$II_2$: $W_2(P_1)=W_2(\bar P_1)=W_2(P_4)=W_2(\bar P_4)=2$. $W_1(P_3)=1$ and $W_1(\bar P_3)=0$. $W_1(P_4)=1$ and $W_1(\bar P_4)=2$. 
\vspace{0,2 cm}

$II_3$: $W_2(P_1)=W_2(\bar P_1)=W_2(P_4)=W_2(\bar P_4)=0$. $W_2(P_2)=W_2(\bar P_2)=0$.
\vspace{0,2 cm}

$II_4$: $W_2(P_1)=W_2(\bar P_1)=W_2(P_4)=W_2(\bar P_4)=3$.  $W_1(P_1)=W_1(\bar P_1)=2$. 
$W_1(P_2)=W_1(\bar P_2)=2$. $W_2(\bar P_3)-W_2(P_3)=W_2(\bar P_3,f)=2$. $W_1(P_4)=W_1(\bar P_4)=1$.
\vspace{0,2 cm}

$III_1$: $W_1(P_3)=W_1(\bar P_3)=W_1(P_4)=W_1(\bar P_4)=0$.
\vspace{0,2 cm}

$III_2$:  $W_2(P_1)=W_2(\bar P_1)=W_2(P_4)=W_2(\bar P_4)=0$. $W_2(P_2)=W_2(\bar P_2)=0$. $W_1(P_3)=W_1(\bar P_3)=0$.
\vspace{0,2 cm}

$III_3$: $W_2(P_2)=W_2(\bar P_2)=3$. $W_1(P_2)=W_1(\bar P_2)=1$. 
\vspace{0,2 cm}

$III_4$: $W_2(P_1)=W_2(\bar P_1)=W_2(P_4)=W_2(\bar P_4)=2$. $W_1(P_1)=W_1(\bar P_1)=2$. $W_2(P_2)=W_2(\bar P_2)=2$. $W_1(P_2)=W_1(\bar P_2)=2$. $W_2(\bar P_3)-W_2(P_3)=W_2(\bar P_3,f)=2$.
\vspace{0,2 cm}

$IV_1$: $W_1(P_1)=W_1(\bar P_1)=0$. $W_2(\bar P_3)-W_2(P_3)=W_2(\bar P_3,f)=0$. $W_1(P_3)=W_1(\bar P_3)=1$.
$W_1(P_4)=W_1(\bar P_4)=0$.\vspace{0,2 cm}

$IV_2$: $W_2(P_2)=W_2(\bar P_2)=3$.\vspace{0,2 cm}

$IV_3$: $W_1(P_1)=W_1(\bar P_1)=1$. $W_1(P_2)=W_1(\bar P_2)=1$. $W_2(P_2)=W_2(\bar P_2)=2$. $W_2(\bar P_3)-W_2(P_3)=W_2(\bar P_3,f)=1$. \vspace{0,2 cm}

$IV_4$: $W_2(P_1)=W_2(\bar P_1)=W_2(P_4)=W_2(\bar P_4)=0$. $W_2(P_2)=W_2(\bar P_2)=0$. $W_1(P_3)=W_1(\bar P_3)=1$. $W_2(P_3)=W_2(\bar P_3)=1$.\vspace{0,2 cm}

$V_1$:  $W_1(P_1)=W_1(\bar P_1)=0$. $W_1(P_2)=W_1(\bar P_2)=0$. $W_2(\bar P_3)-W_2(P_3)=W_2(\bar P_3,f)=0$.
\vspace{0,2 cm}

$V_2$: nothing at all\vspace{0,2 cm}

$V_3$: $W_2(P_1)=W_2(\bar P_1)=W_2(P_4)=W_2(\bar P_4)=0$. $W_2(P_2)=W_2(\bar P_2)=0$. $W_2(P_3)=W_2(\bar P_3)=0$.\vspace{0,2 cm}

$V_4$: $W_2(P_1)=W_2(\bar P_1)=W_2(P_4)=W_2(\bar P_4)=4$.  $W_2(P_3)=W_2(\bar P_3)=4$.  $W_2(P_4)=W_2(\bar P_4)=4$. $W_1(P_3)=3$. $W_1(\bar P_3)=2$. $W_1(P_4)=1$. $W_1(\bar P_4)=2$.\vspace{0,2 cm}

$VI_1$:  $W_1(P_1)=W_1(\bar P_1)=0$. $W_1(P_2)=W_1(\bar P_2)=0$. $W_1(P_4)=W_1(\bar P_4)=0$. $W_2(\bar P_3)-W_2(P_3)=W_2(\bar P_3,f)=0$.\vspace{0,2 cm}

$VI_2$: $W_1(P_3)=1$. $W_1(\bar P_3)=0$. $W_1(P_4)=0$. $W_1(\bar P_4)=1$.\vspace{0,2 cm}

$VI_3$: $W_2(P_2)=W_2(\bar P_2)=4$.\vspace{0,2 cm}

$VI_4$: $W_2(P_1)=W_2(\bar P_1)=W_2(P_4)=W_2(\bar P_4)=0$. $W_2(P_2)=W_2(\bar P_2)=0$.  $W_2(P_3)=W_2(\bar P_3)=0$.\vspace{0,2 cm}

$\Box$
\begin{proposition}
Let $m$ be the meridian for a positive quadruple crossing. Then  $R^{(1)}(m)=0$.
\end{proposition}
{\em Proof.} The points $(1)$ to $(5)$ in Lemma 4 imply that we can repeat word for word (without the grading) the proof of Proposition 2 but with the coefficient $1$ of $\sigma_2 \sigma_1(p)$ replaced by $W_1(p)$ and the coefficient $W(p)$ of $z(\sigma_2(p)-\sigma_1(p))$ replaced  by $W_2(p)$.

\begin{figure}
\centering
\includegraphics{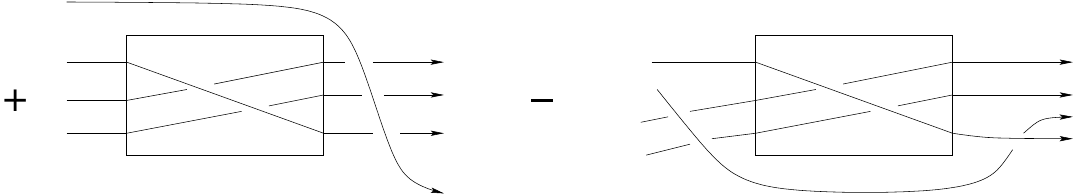}
\caption{\label{P3P4} cancellation in (6) of Lemma 4}  
\end{figure}

The only new feature is point $(6)$ in Lemma 4. If $W_1(P_3)= W_1(\bar P_3) +1$ then $P_3-\bar P_3$ contribute to $R^{(1)}$ the element associated to the tangle on the left in Fig.~\ref{P3P4} and $P_4-\bar P_4$ contribute the element associated to the 
tangle on the right in Fig.~\ref{P3P4}. But these contributions cancel out (even without using the skein relations).

$\Box$

We have proven that  $R^{(1)}$ is a solution of the positive global tetrahedron equation.

\begin{remark}
For the "dual" 1-cocycles we have of course to construct  the weights $W_1$ and $W_2$ from the second formula for $v_2$ of Polyak and Viro (compare Fig.~\ref{PV})  by using the "heads" of the arrows.
\end{remark}

\begin{question}
In the construction of $R^{(1)}$ we have used relative finite type invariants $W_1(p,hm)$ and $W_2(p)$ of degrees 1 and 2. Can one push this further with relative finite type invariants of higher degrees? 

\end{question}
\section{Solution with linear weight of the cube equations}

We start with some observations and notations. The identity $\sigma_2(p)-\sigma_1(p)=\sigma^{-1}_2(p)-\sigma^{-1}_1(p)$ is sometimes useful to simplify calculations. We have given names to the global types of strata of $\Sigma^{(1)}_{tri}$ in Fig.~\ref{globtricross}. {\em The union of all strata of the same type is denoted simply by this name.} The weights and the partial smoothings in the case of positive triple crossings were determined in the last two sections. We have to consider now all other local types of triple crossings together with the self-tangencies. 

\begin{definition}
The weight $W(p)$ is the same as in Definition 8 for all local types of triple crossings $p$ .
\end{definition}

For a triple crossing $p$ of a given local type from $\{1,...,8\}$ we denote by $T_{r_a}(type)(p)=T_{l_c}(type)(p)$ and $T(type)(p)$ the corresponding partial smoothing of $p$ in the skein module. We will determine these partial smoothings from the cube equations.

\begin{definition}
We distinguish four types of self-tangencies, denoted by 

$II_0^+, II_0^-, II_1^+, II_1^-$. Here "0" or "1" is the type of the distinguished crossing $d$ and "+" stands for opposite tangent directions of the two branches and "-" stands for the same tangent direction.
\end{definition}

\begin{figure}
\centering
\includegraphics{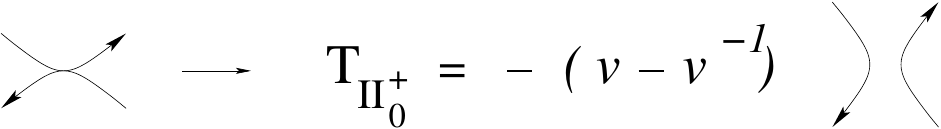}
\caption{\label{sts}  partial smoothing of a self-tangency $II_0^+$}  
\end{figure}

\begin{definition}
The partial smoothing $T_{II^+_0}(p)$ of a self-tangency with opposite tangent direction is defined in Fig.~\ref{sts}.
The partial smoothing $T_{II^-_0}(p)$ of a self-tangency with equal tangent direction is defined as $-zP_T$ (where $P_T$ is the element represented by the original tangle $T$ in the HOMFLYPT skein module $S(\partial T)$).
\end{definition}

Notice that the distinguished crossing $d$ has to be of type 0 in both cases.

\begin{definition}
Let $p$ be a self-tangency. Its contribution to $R^{(1)}_{reg}$ is defined by

 $R^{(1)}_{reg}=sign(p) W(p)T_{II^+_0}(p)$ respectively $sign(p) W(p)T_{II^-_0}(p)$.\vspace{0,2 cm}

Its contribution to $R^{(1)}$ is defined by

 $R^{(1)}=sign(p) W_2(p)T_{II^+_0}(p)$  respectively $sign(p) W_2(p)T_{II^-_0}(p)$.

(Don't forget that in the case of $R^{(1)}$ the HOMFLYPT invariants have to be normalized by $v^{-w(T)}$.)
\end{definition}

There are exactly two local types of self-crossings with opposite tangent direction as well as for equal tangent direction. The corresponding strata from $\Sigma^{(1)}_{tan}$ are adjacent in $\Sigma^{(2)}_{self-flex}$.

\begin{lemma}
Let $m$ be the meridian of $\Sigma^{(2)}_{self-flex}$ (compare Section 2). Then $R^{(1)}_{reg}(m)=0$ and $R^{(1)}(m)=0$.
\end{lemma}

\begin{figure}
\centering
\includegraphics{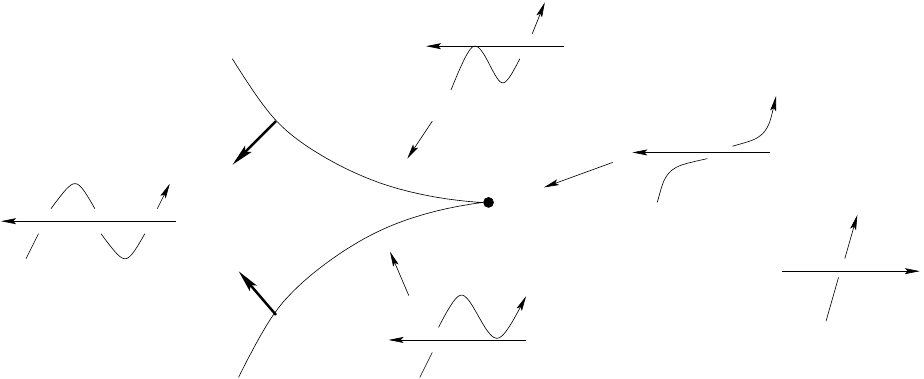}
\caption{\label{unflex}  the meridional disc for a self-tangency in a flex}  
\end{figure}

\begin{figure}
\centering
\includegraphics{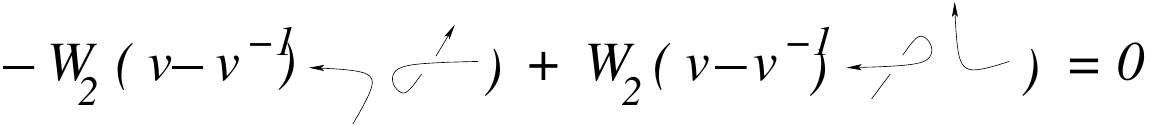}
\caption{\label{calself}  $R^{(1)}$ vanishes on the meridian of a self-tangency in a flex}  
\end{figure}

{\em Proof.} We consider the case with opposite tangent direction. We show the unfolding (i.e. the intersection of a meridional disc with $\Sigma$ in $M_T$) in Fig.~\ref{unflex} (compare \cite{FK}). There are two cases: either all three crossings in the picture are of type 1 or all three are of type 0. In the first case both self-tangencies in the unfolding do not contribute to $R^{(1)}_{reg}$ and to $R^{(1)}$ (because $d$ has to be of type $0$).
In the second case both self-tangencies have different sign but the same weight $W(p)$, because the crossing which is not in $d$ is not a f-crossing. But they have also the same weight $W_2(p)$. Indeed, the crossing which is not in $d$ has for both self-crossings the same writhe and if it is a r-crossing  for a f-crossing  for one of the  self-tangencies then this is the case for the other self-tangency too. It follows that $R^{(1)}_{reg}(m)=0$ and $R^{(1)}(m)=0$ as shown in Fig.~\ref{calself}.
The case with equal tangent direction is analogue (the third crossing can never be a f-crossing).
$\Box$

This lemma implies that we can restrict ourself to consider  in the cube equations only one of the two local types of self-tangencies.

\begin{figure}
\centering
\includegraphics{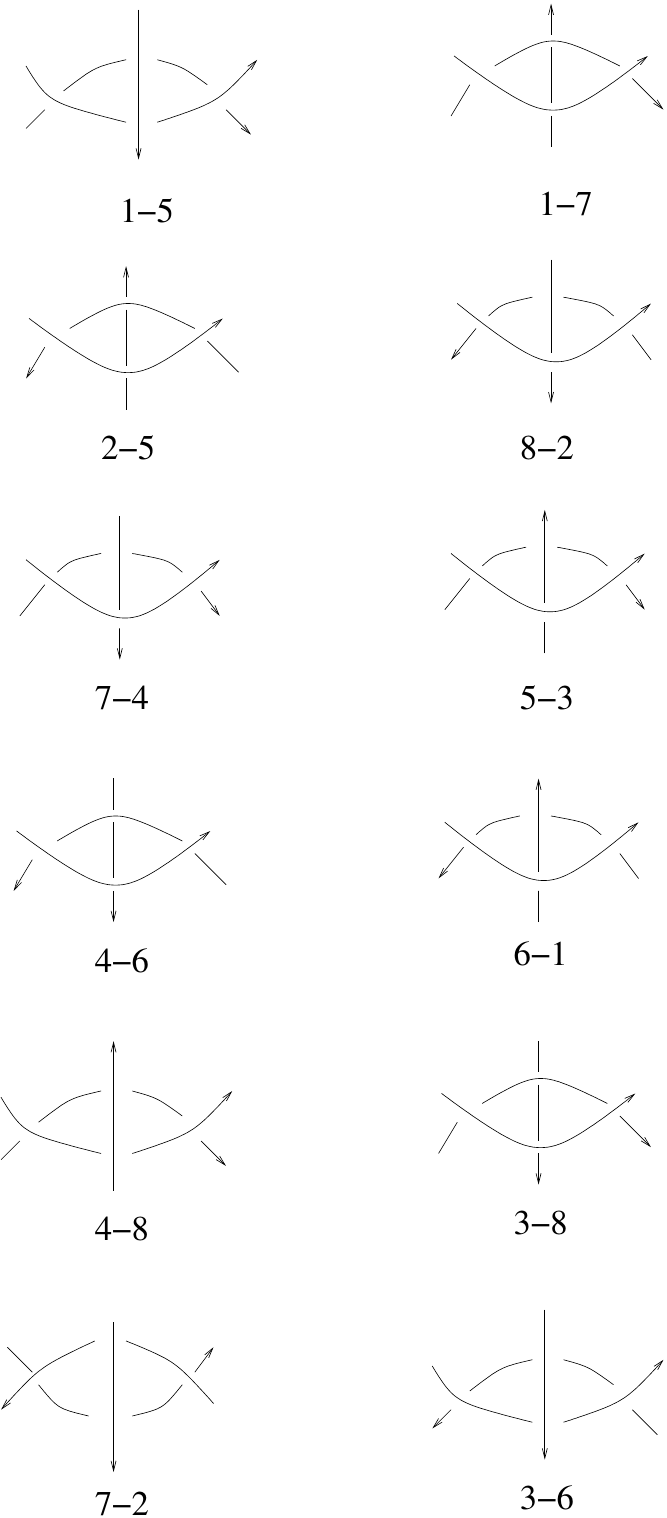}
\caption{\label{edgegam}  the twelve edges of the graph $\Gamma$}  
\end{figure}

\begin{figure}
\centering
\includegraphics{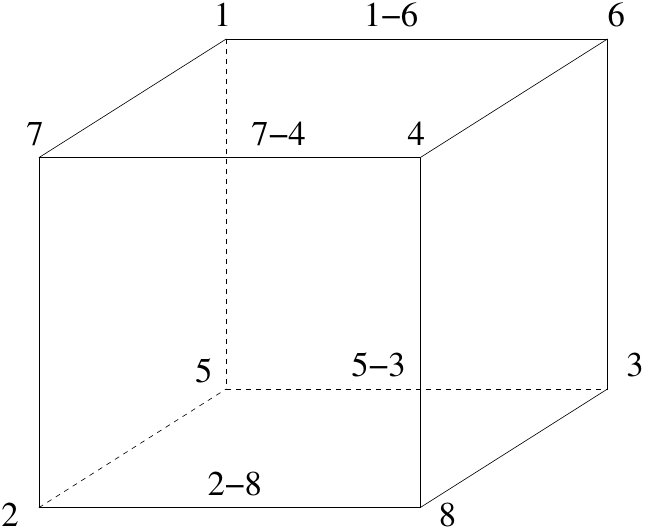}
\caption{\label{gamma}  the graph $\Gamma$}  
\end{figure}

\begin{figure}
\centering
\includegraphics{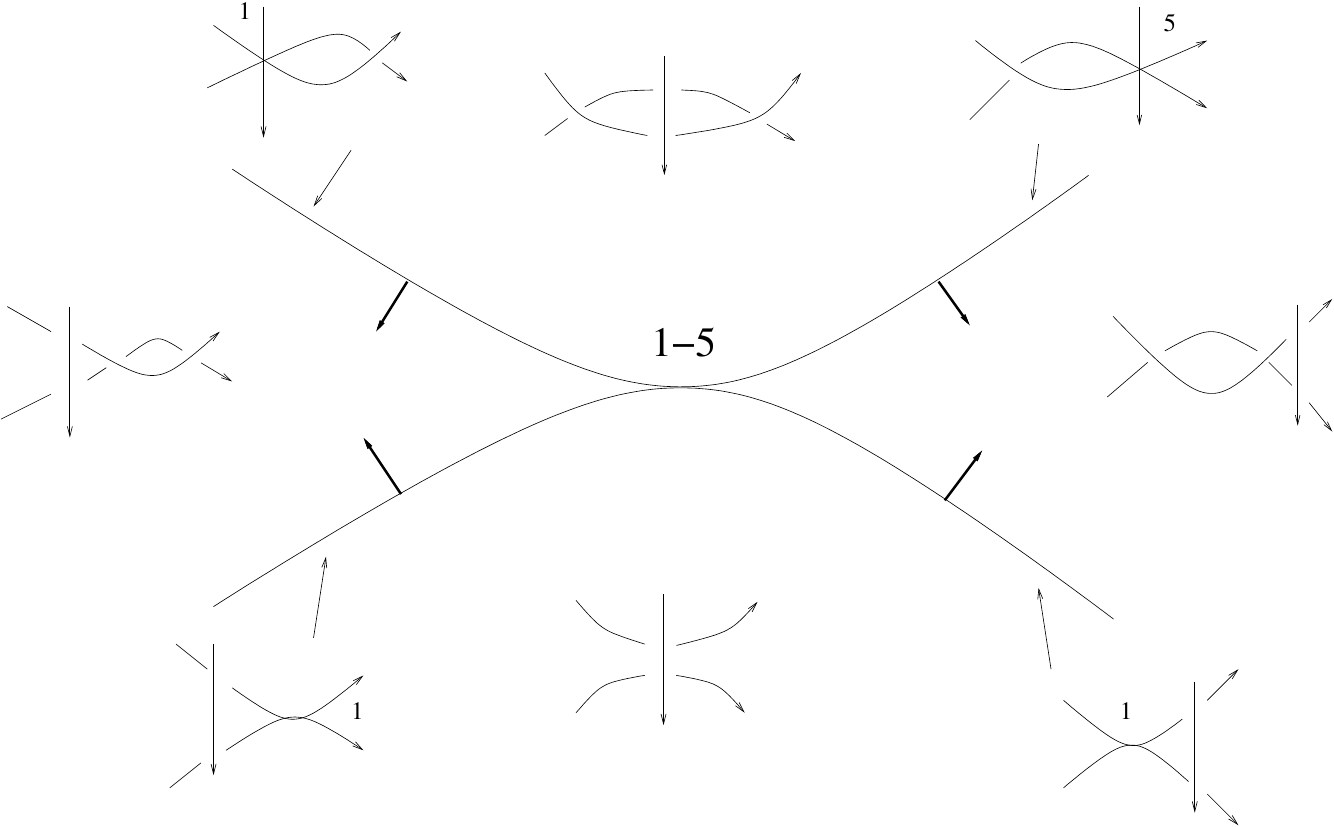}
\caption{\label{unfold1-7}  the unfolding of a self-tangency with a transverse branch}  
\end{figure}

The local types of triple crossings were shown in Fig.~\ref{loc-trip}. The diagrams which correspond to the edges of the graph $\Gamma$ (compare Section 2) are shown in Fig.~\ref{edgegam}. We show the corresponding graph $\Gamma$ now in 
Fig.~\ref{gamma}. The unfolding of e.g. the edge $1-5$ is shown in Fig.~\ref{unfold1-7} (compare \cite{FK}).

\begin{figure}
\centering
\includegraphics{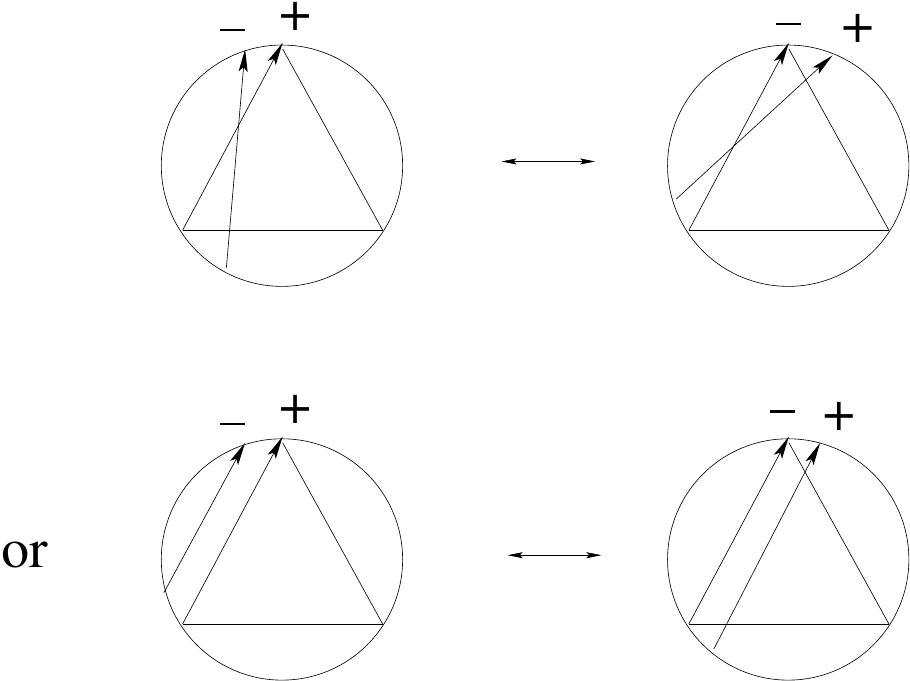}
\caption{\label{triedge} the two vertices of an edge of $\Gamma$}  
\end{figure}

\begin{observation}
The diagrams corresponding to the two vertices's of an edge differ just by the two crossings of the self-tangency  which replace each other in the triangle as shown in Fig.~\ref{triedge}. Consequently, the two vertices of an edge have always the same global type and always different local types and the two crossings which are interchanged are simultaneously f-crossings  and have the same grading $A$. Moreover, the two self-tangencies in the unfolding (of an edge) can  contribute non trivially to $R^{(1)}_{reg}$ if they have different weight $W(p)$ or if their partial smoothings are not isotopic. The latter happens exactly for the edges "1-6" and "2-8" (where the third branch passes between the two branches with opposite orientation of the self-tangency).
\end{observation}

It follows from Observation 1 and Lemma 5 that we have to solve the cube equations for the graph $\Gamma$ exactly six times: one solution for each global type of triple crossings.

\begin{figure}
\centering
\includegraphics{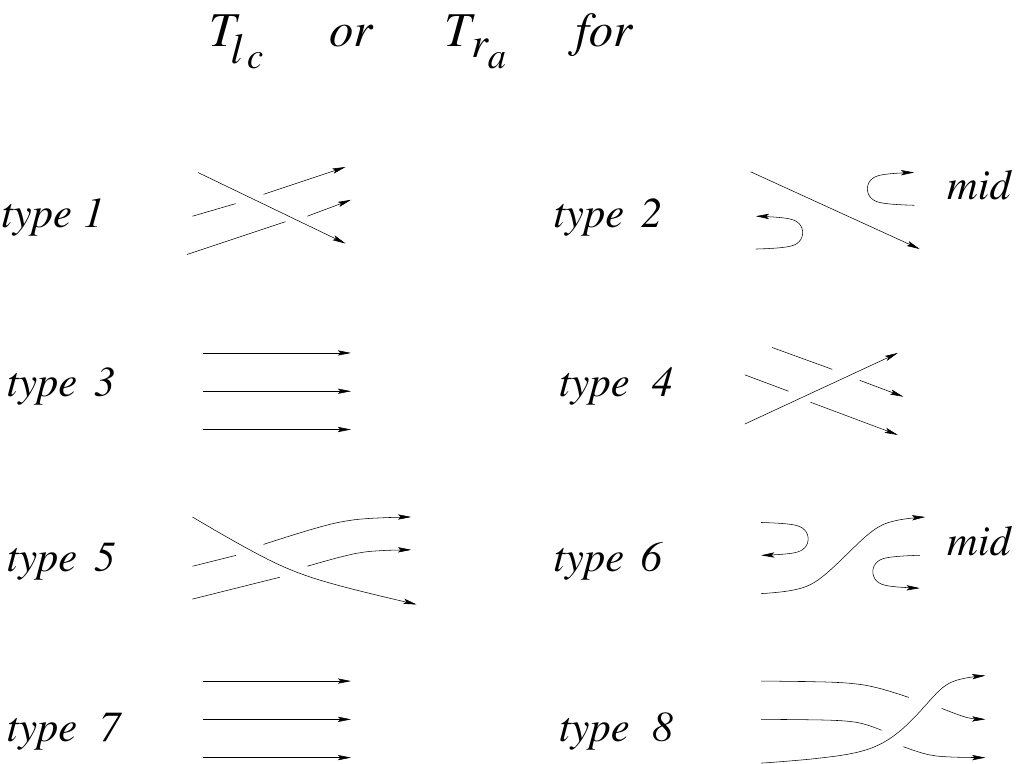}
\caption{\label{typesmooth}  the partial smoothings $T_{l_c}$ and $T_{r_a}$}  
\end{figure}

\begin{figure}
\centering
\includegraphics{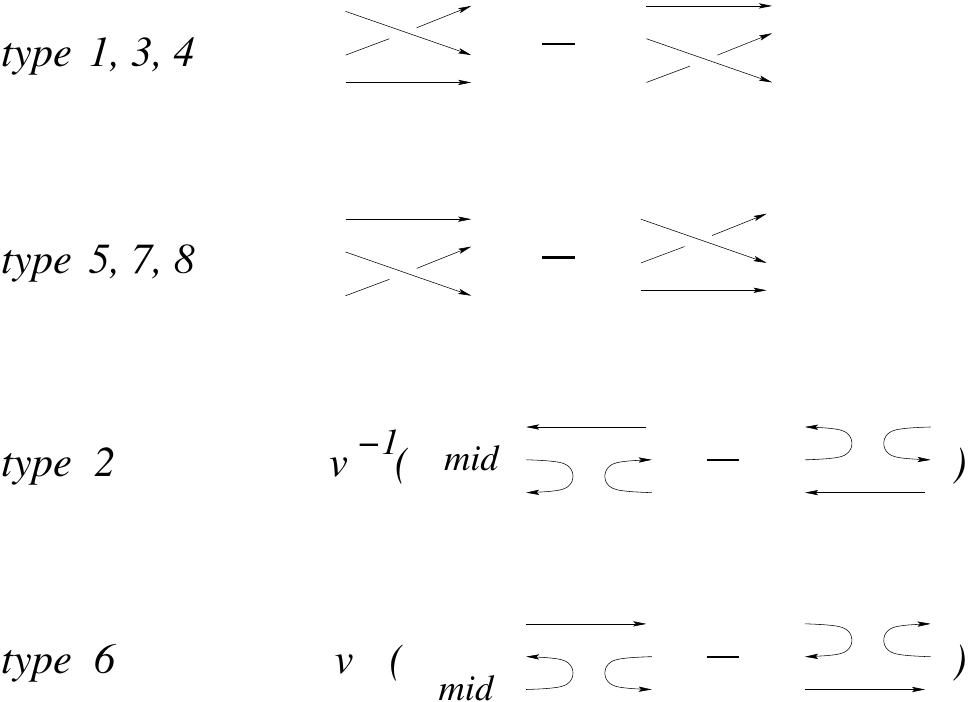}
\caption{\label{smooth} the partial smoothings $T$ (type)}  
\end{figure}

\begin{definition}
The partial smoothings for the local and global types of triple crossings are given in Fig.~\ref{typesmooth} and  Fig.~\ref{smooth} (remember that {\em mid} denotes the ingoing middle branch for a star-like triple crossing). 
\end{definition}

\begin{definition}
Let $s$ be a generic oriented arc in $M^{reg}_T$ (with a fixed abstract closure $T \cup \sigma$ to an oriented circle). Let $A \subset \partial T$ be a given grading.  The 
1-cochain $R^{(1)}_{reg}(A)$ is defined by \vspace{0,2 cm}

$R^{(1)}_{reg}(A)(s ) = \sum_{p \in s \cap l_c} sign(p)T_{l_c} (type)(p) +\sum_{p \in s \cap r_a }sign(p)T_{r_a}(type)(p) +$\vspace{0,2 cm}

 $\sum_{p \in s \cap r_a }sign(p)zW(p)T(type)(p) +$\vspace{0,2 cm}

 $\sum_{p \in s \cap r_b }sign(p)zW(p)T(type)(p) +$\vspace{0,2 cm}

 $\sum_{p \in s \cap l_b }sign(p)zW(p)T(type)(p) +$\vspace{0,2 cm}

 $\sum_{p \in s \cap II^+_0 } sign(p)W(p)T_{II^+_0}(p)+$\vspace{0,2 cm} 

 $\sum_{p \in s \cap II^-_0 } sign(p)W(p)T_{II^-_0}(p)$\vspace{0,2 cm} 

Here all weights $W(p)$ are defined only over the f-crossings $f$ with $\partial f= A$ and in the first two sums (i.e. for $T_{l_c}$ and $T_{r_a}$) we require that $\partial (hm)=A$ for the triple crossings.

Notice that the stratum $r_a$ contributes in two different ways.
\end{definition}

\begin{figure}
\centering
\includegraphics{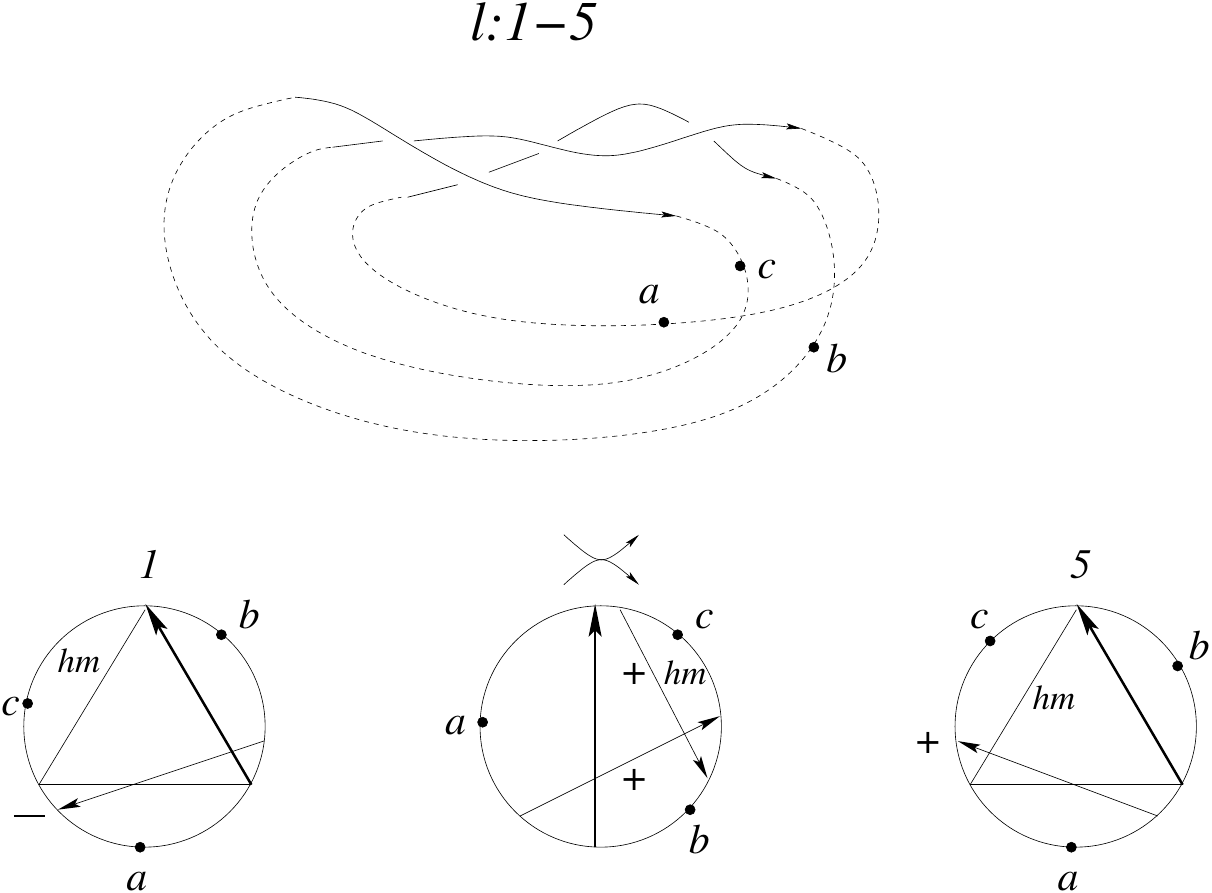}
\caption{\label{l1-5}  l1-5}  
\end{figure}

\begin{figure}
\centering
\includegraphics{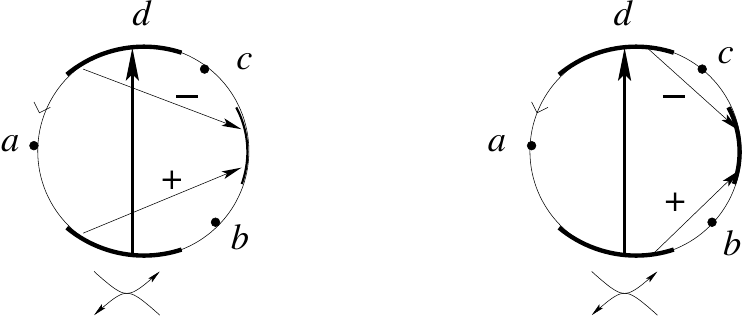}
\caption{\label{self2}  the two self-tangencies for the edge $l7-2$}  
\end{figure}

This turns out to be the solution of the cube equations with linear weight.
\begin{proposition}
Let $m$ be a meridian of $\Sigma^{(2)}_{trans-self}$ or a loop in $\Gamma$ and let $A$ be a given grading. Then $R^{(1)}_{reg}(A)(m)=0$ for the partial smoothings given in the Definitions 17 and 18.
\end{proposition}
The proof will be in Fig.~\ref{l1-5},... Fig.~\ref{self2}.

\begin{figure}
\centering
\includegraphics{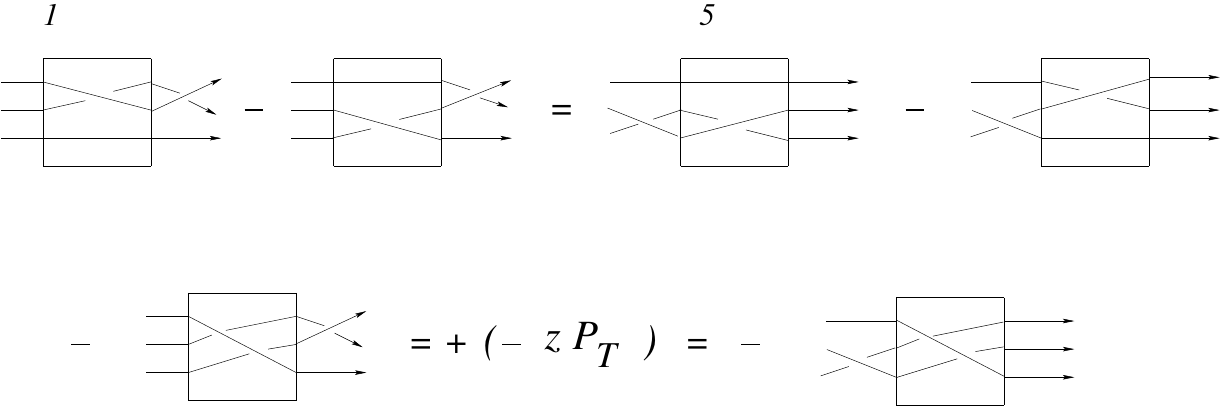}
\caption{\label{l1-5s}  smoothings for l1-5}  
\end{figure}

\begin{figure}
\centering
\includegraphics{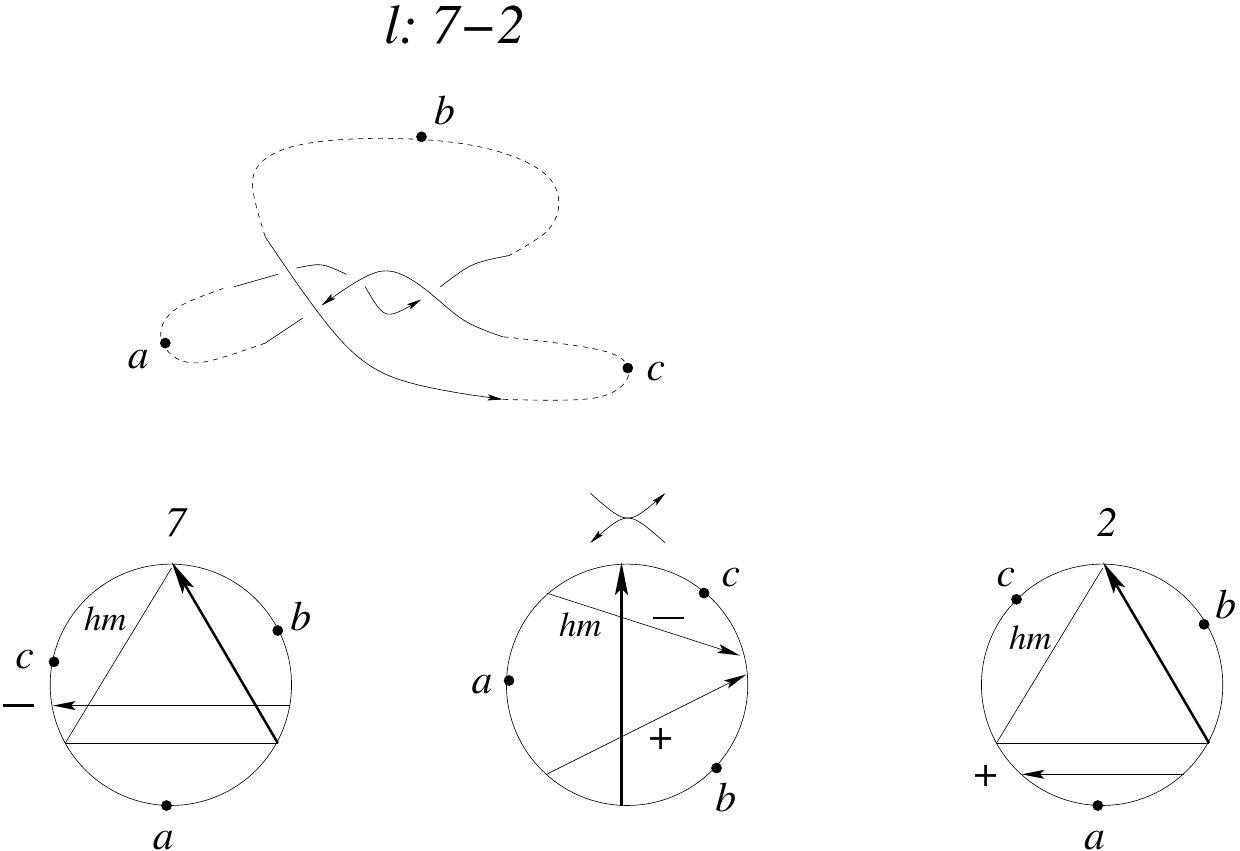}
\caption{\label{l7-2} l7-2}  
\end{figure}

\begin{figure}
\centering
\includegraphics{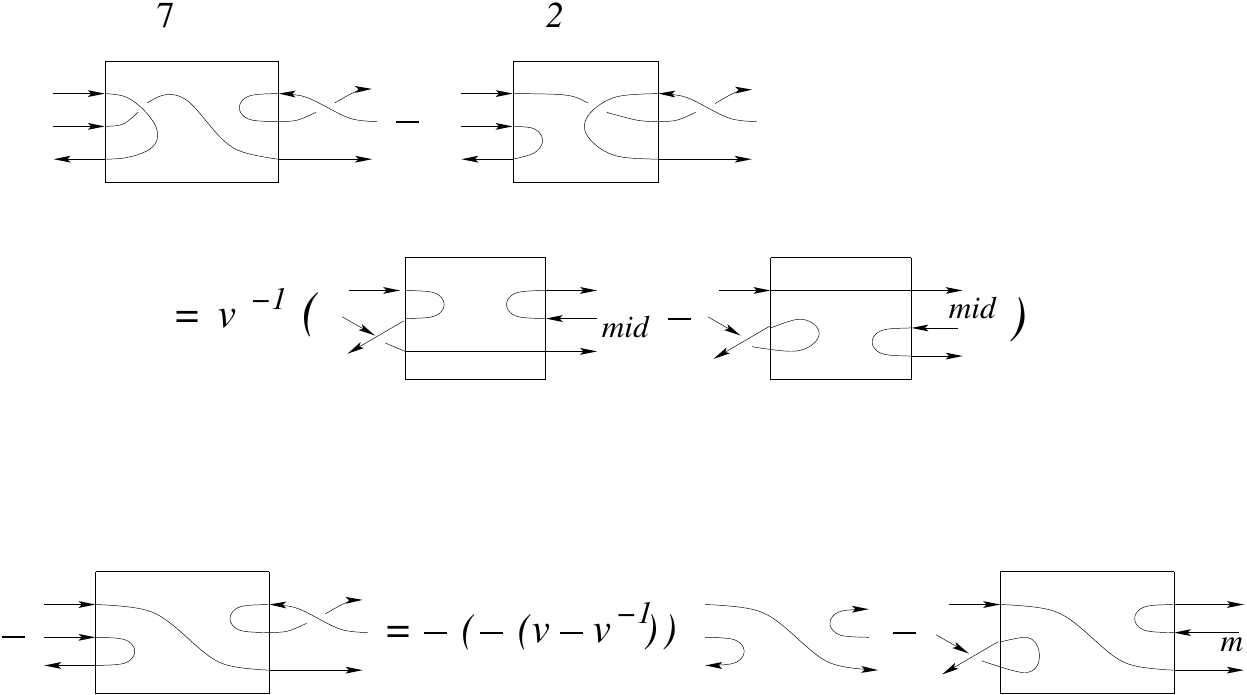}
\caption{\label{l7-2s} smoothings for l7-2}  
\end{figure}

\begin{figure}
\centering
\includegraphics{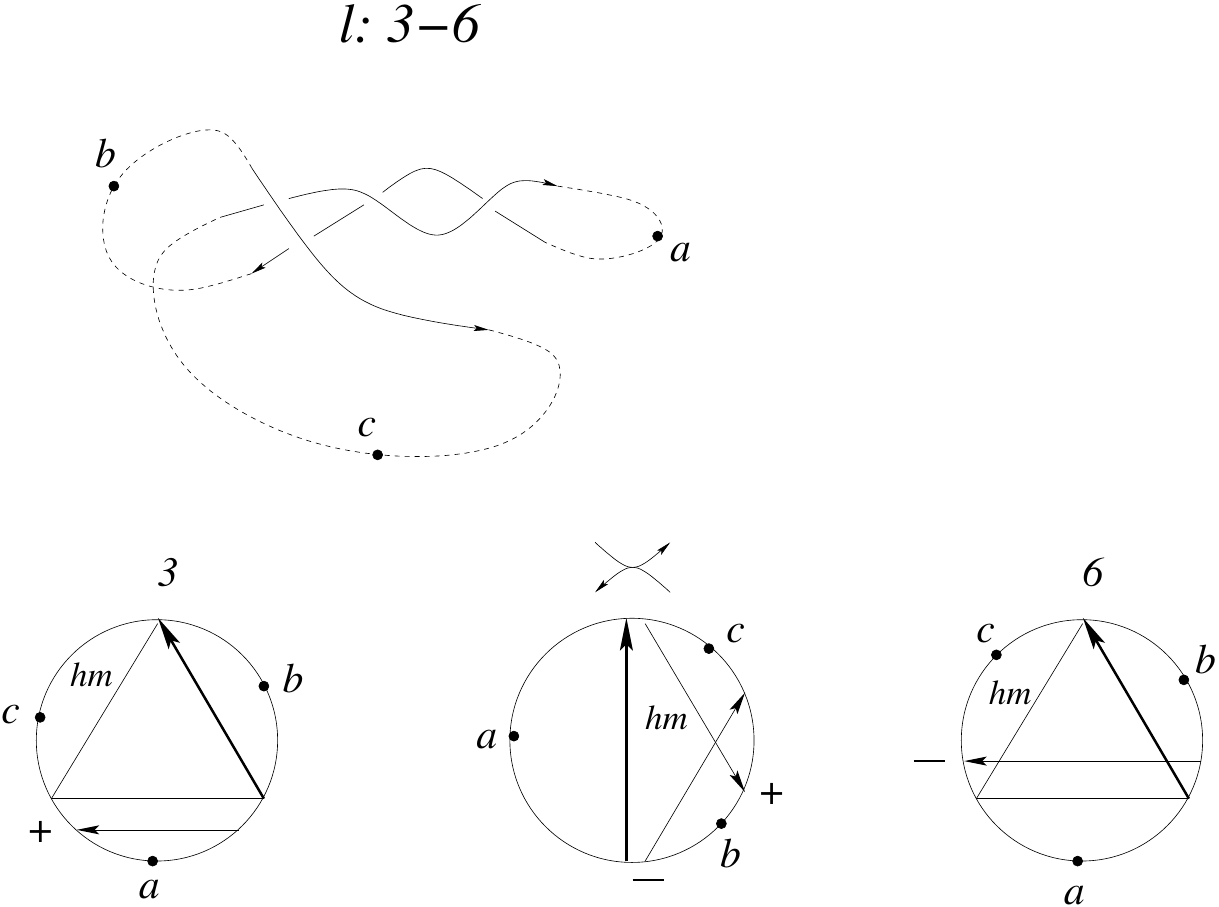}
\caption{\label{l3-6} l3-6}  
\end{figure}

{\em Proof.} 
First of all we observe that the two vertices of an edge share the same weight $W(p)$. The weight could only change if the foot of a f-crossing slides over the head of the crossing $d$. But this is not the case as it was shown in Fig.~\ref{triedge}. Indeed, the foot of the crossing which changes in the triangle can not coincide with the head of $d$ because the latter coincides always with the head of another crossing.

\begin{figure}
\centering
\includegraphics{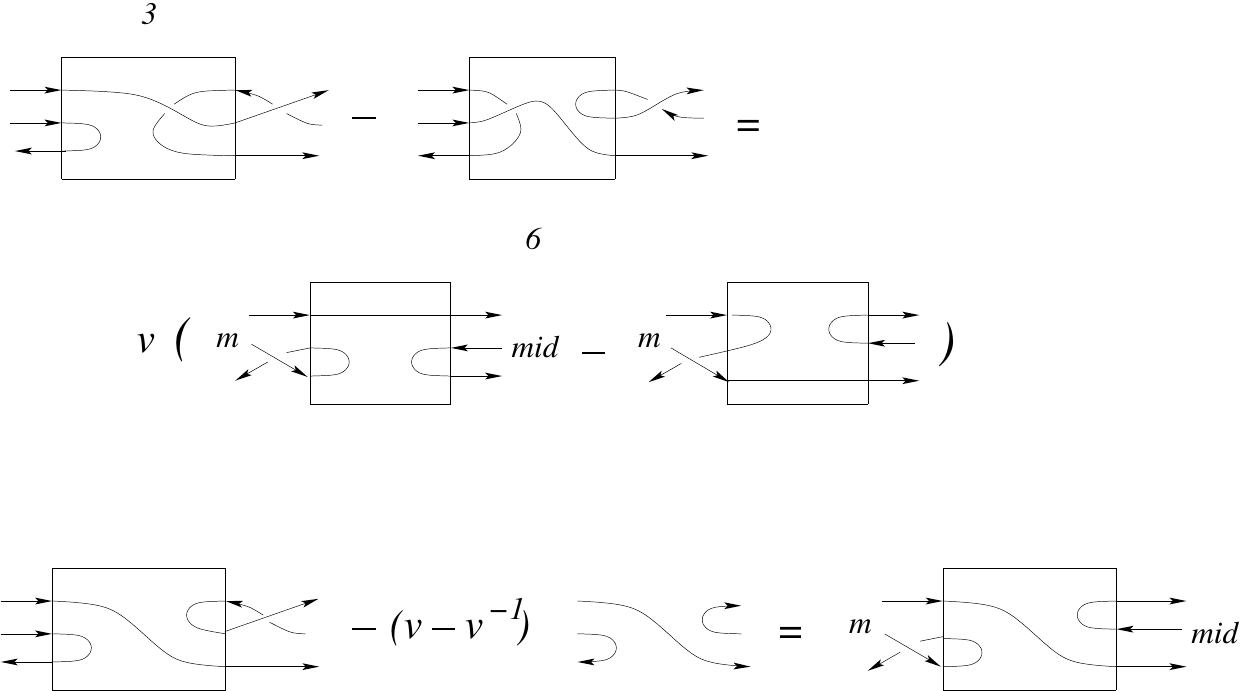}
\caption{\label{l3-6s}  smoothings for l3-6}  
\end{figure}

\begin{figure}
\centering
\includegraphics{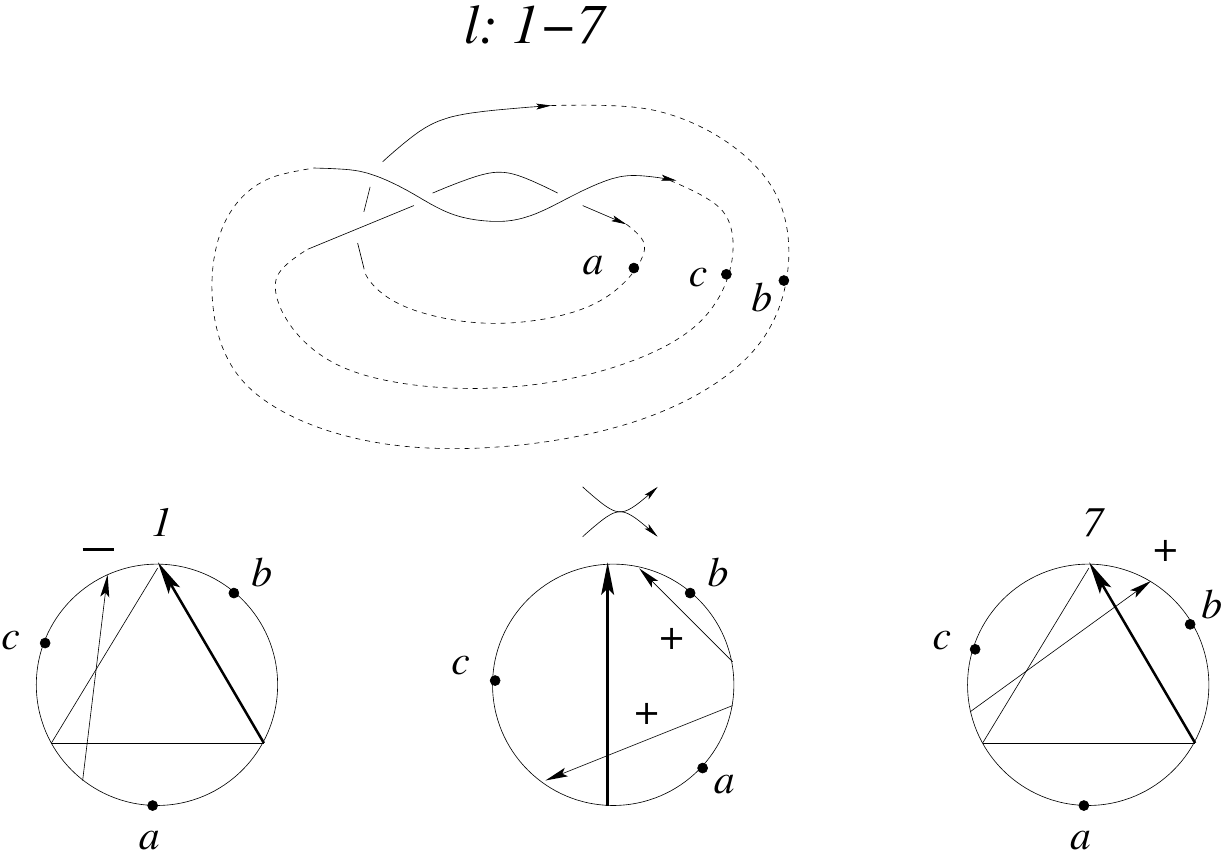}
\caption{\label{l1-7} l1-7 }  
\end{figure}

\begin{figure}
\centering
\includegraphics{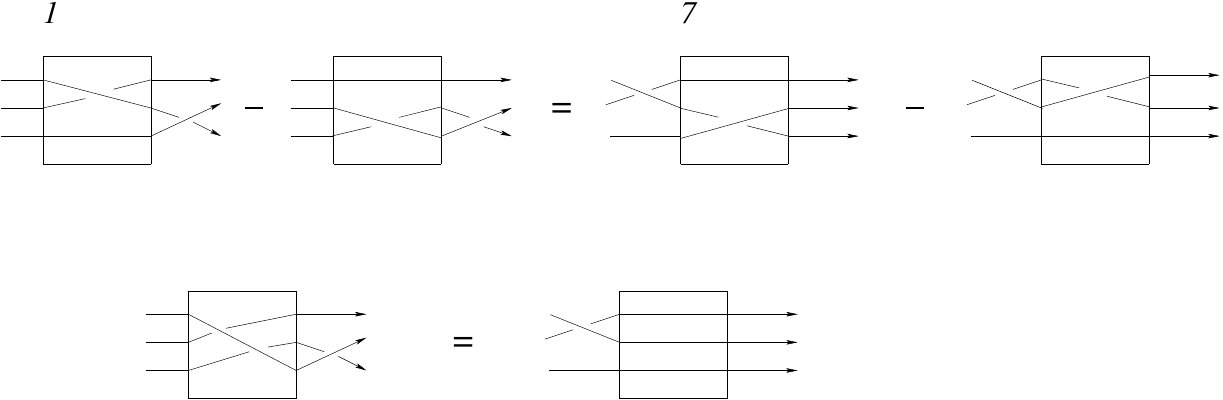}
\caption{\label{l1-7s} smoothings for l1-7}  
\end{figure}

\begin{figure}
\centering
\includegraphics{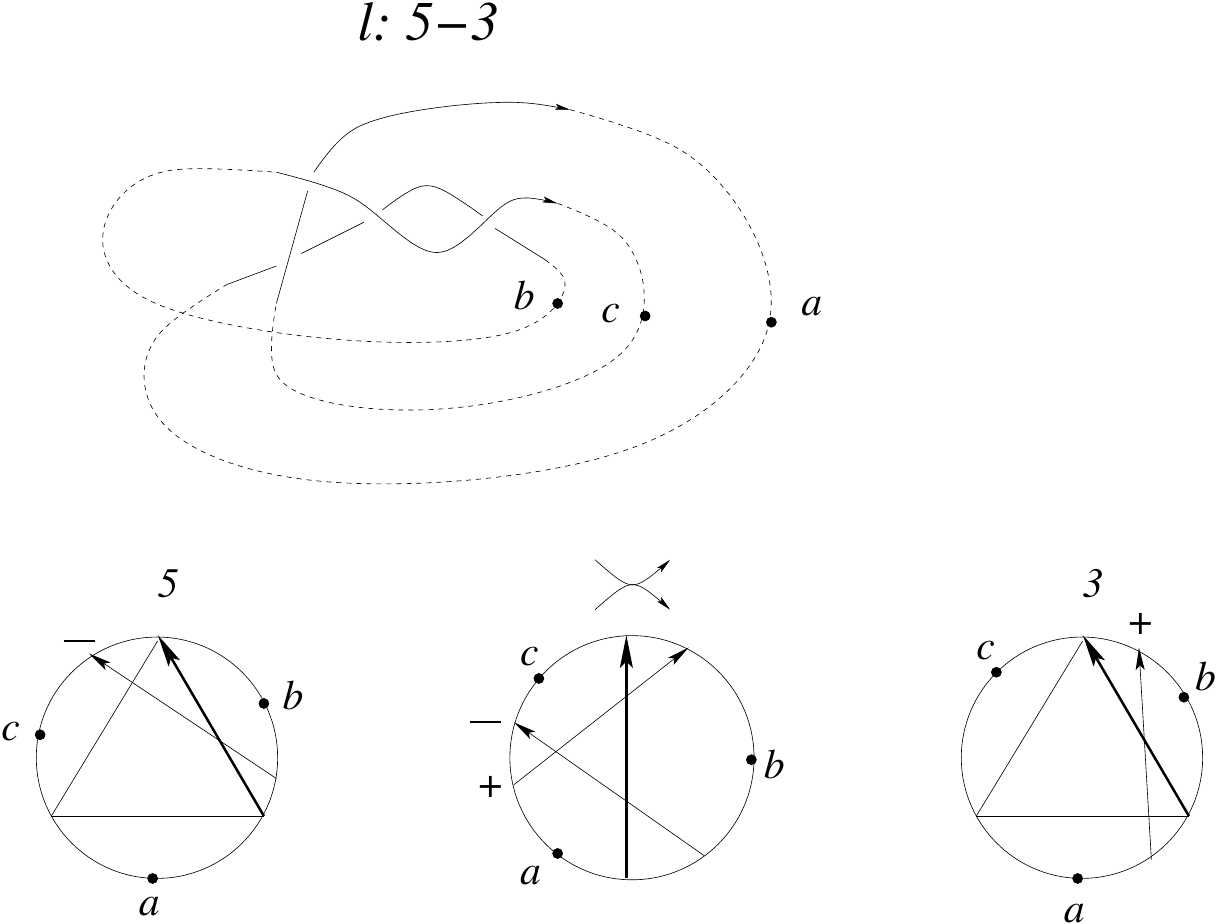}
\caption{\label{l5-3}  l5-3}  
\end{figure}

\begin{figure}
\centering
\includegraphics{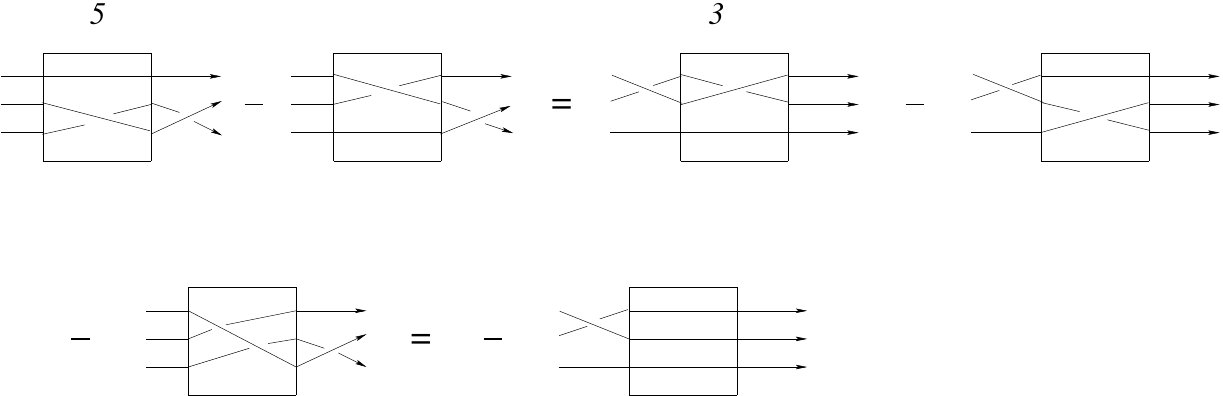}
\caption{\label{l5-3s} smoothings for l5-3}  
\end{figure}

It suffices of course to consider only the four crossings involved in an edge because their position with respect to all other crossings do not change.

 Our strategy is the following: we try to solve the equations separately for the partial smoothings with constant weight of triple crossings and the partial smoothings with linear weight. However they can not be separated completely. It turns out that the two self-tangencies in the unfolding for the edge "1-6" as well as for the edge "2-8" have always the same weight $W(p)$. Therefore they will contribute to the smoothings with linear weight of triple crossings. For the edges "2-7" and "3-6" the smoothings of the two self-tangencies are isotopic. However the weight $W(p)$ changes by 1. Consequently, these self-tangencies have to contribute to the partial smoothings with constant weight of triple crossings. For the edges "2-5" and "4-6" the two self-tangencies have the same weight $W(p)$ and do not contribute.

\begin{figure}
\centering
\includegraphics{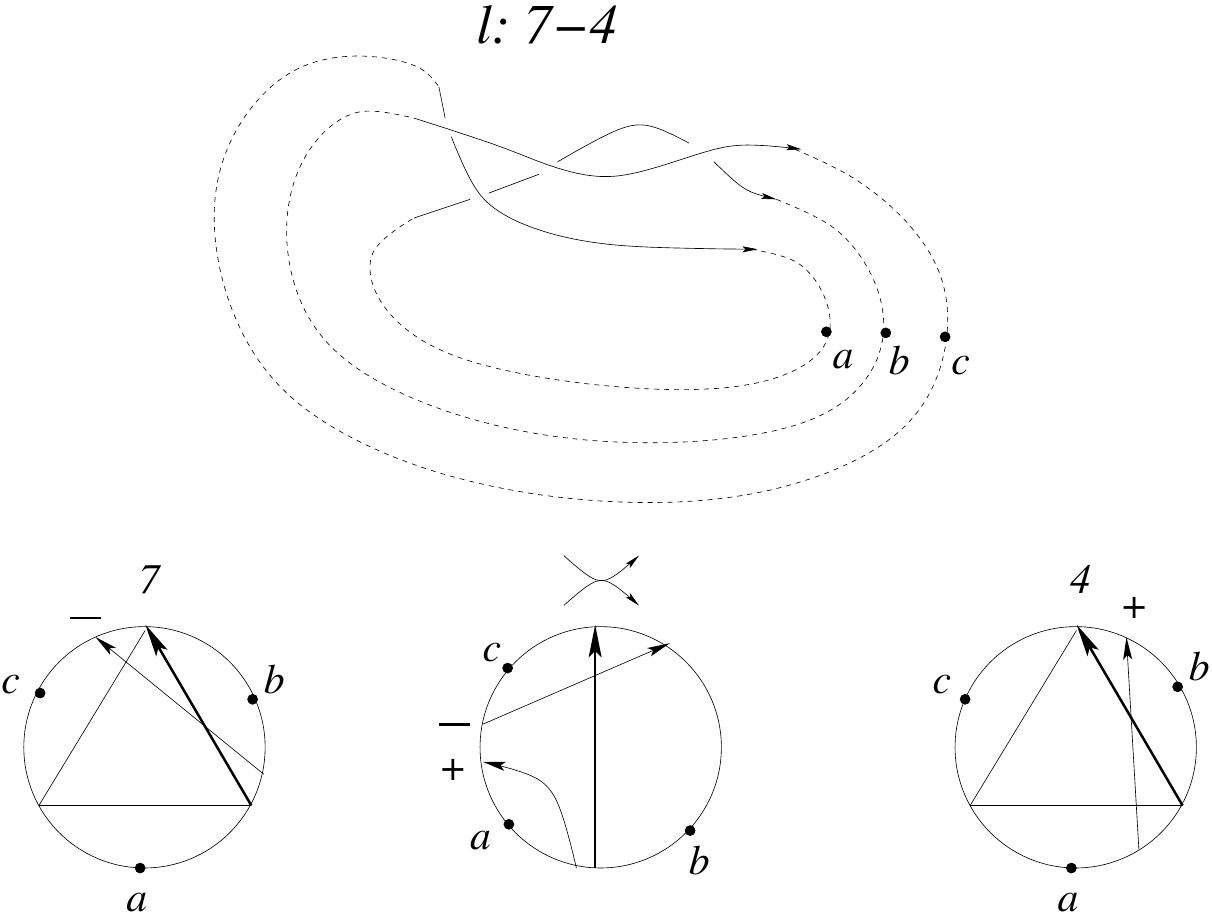}
\caption{\label{l7-4} l7-4}  
\end{figure}

\begin{figure}
\centering
\includegraphics{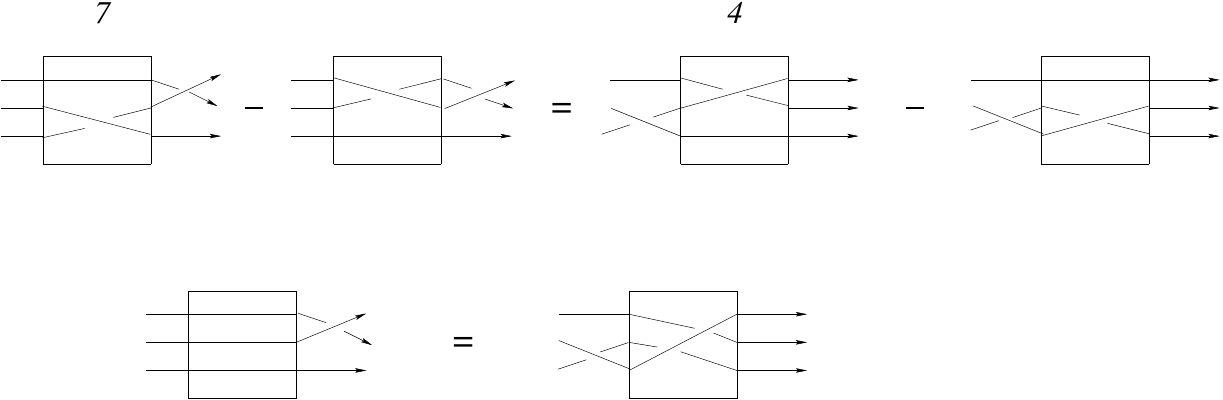}
\caption{\label{l7-4s} smoothings for l7-4}  
\end{figure}

\begin{figure}
\centering
\includegraphics{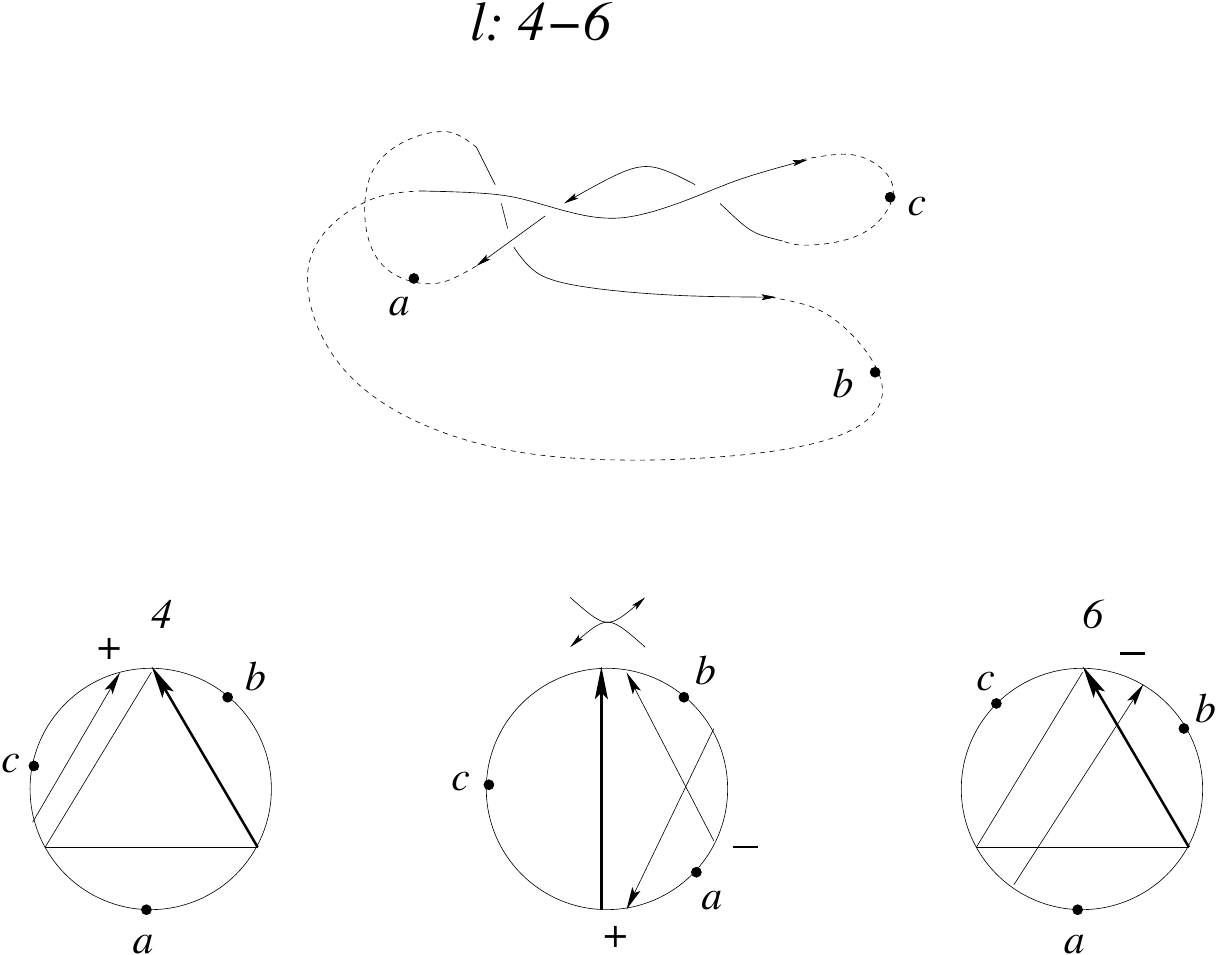}
\caption{\label{l4-6}  l4-6}  
\end{figure}

\begin{figure}
\centering
\includegraphics{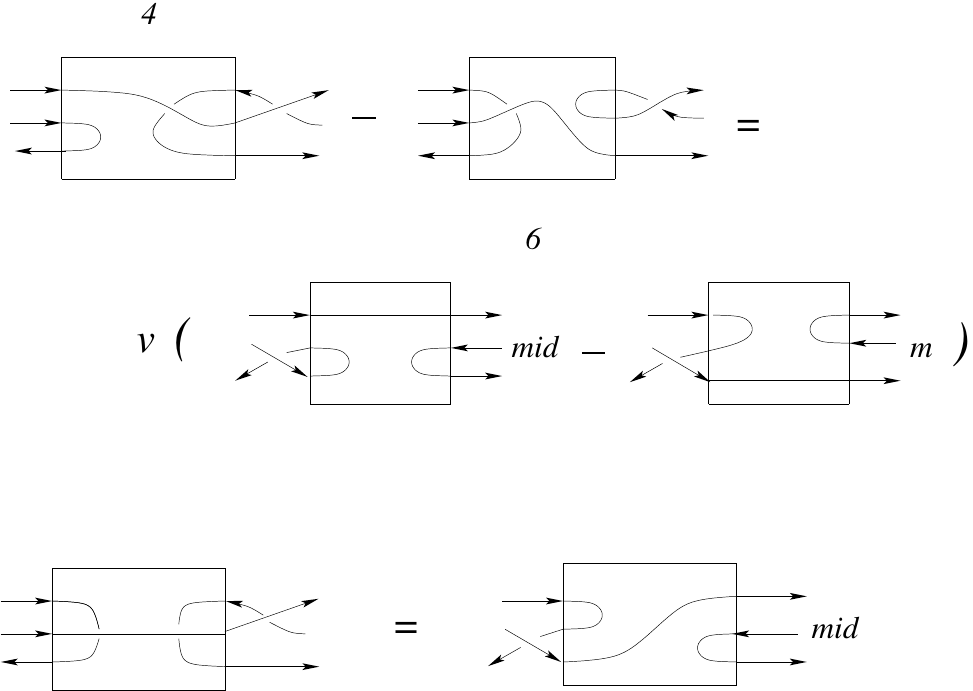}
\caption{\label{l4-6s} smoothings for l4-6}  
\end{figure}

\begin{figure}
\centering
\includegraphics{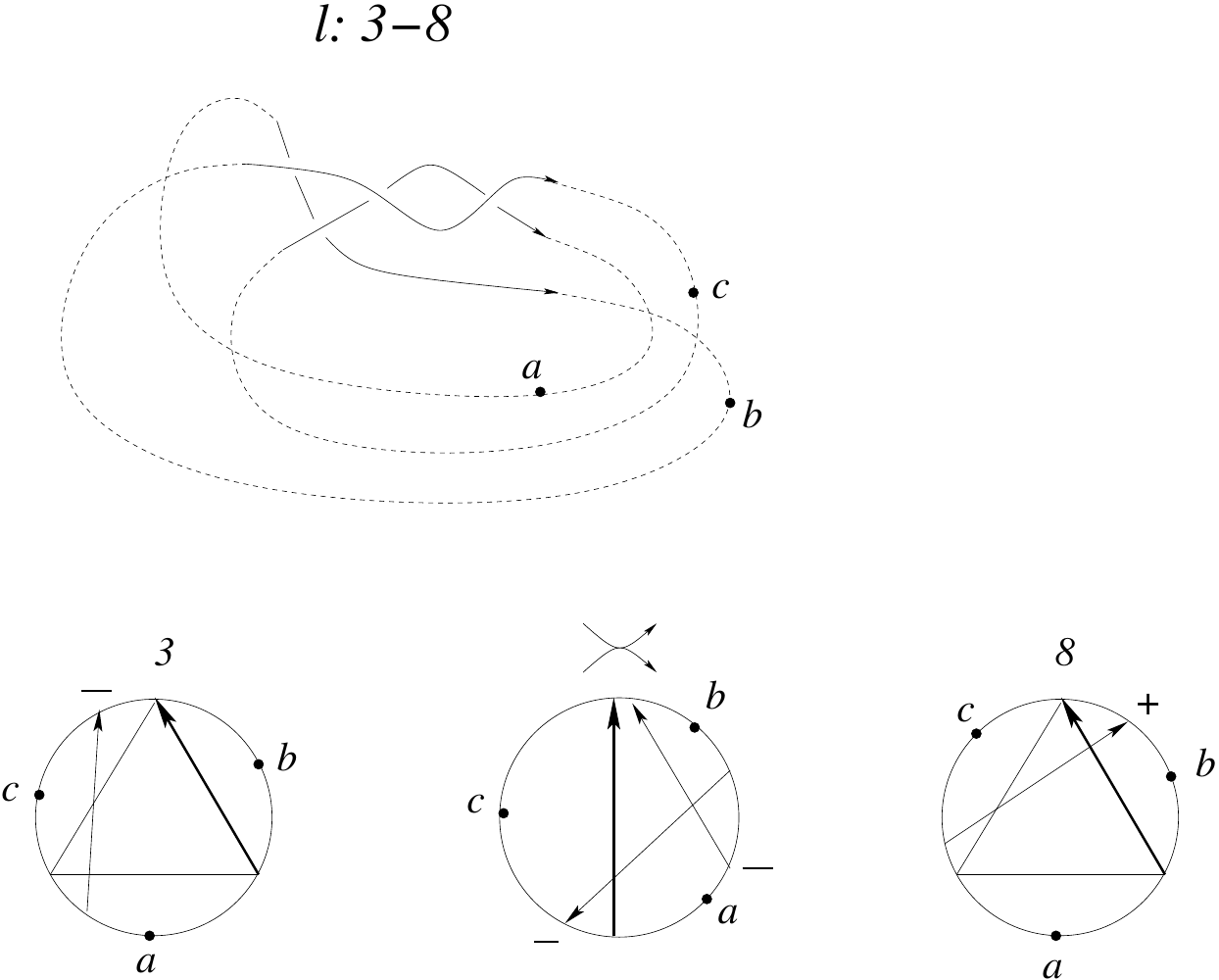}
\caption{\label{l3-8} l3-8}  
\end{figure}

\begin{figure}
\centering
\includegraphics{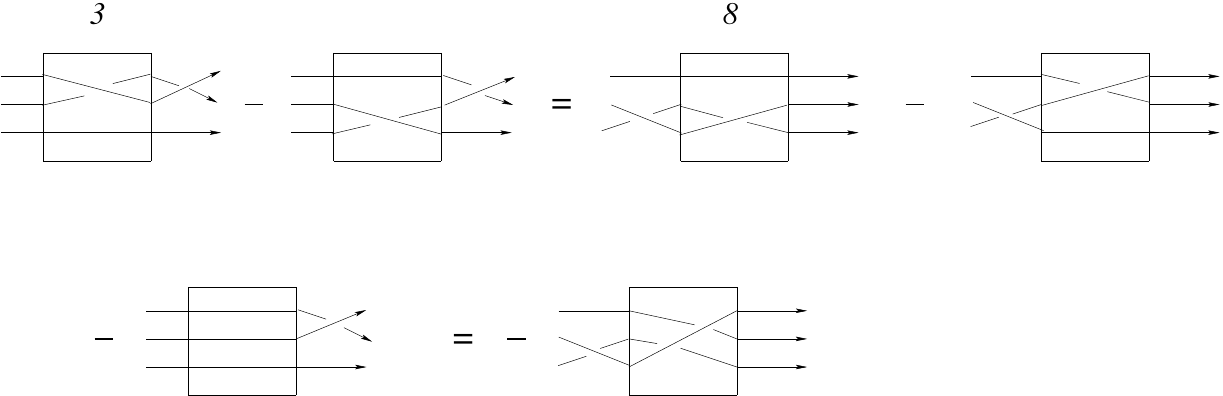}
\caption{\label{l3-8s} smoothings for l3-8}  
\end{figure}

We start with the solutions for triple crossings of global type $l$. In the figures we show the Gauss diagrams of the two triple crossings together with the points at infinity and the Gauss diagram of just one of the two self-tangencies. The Gauss diagram of the second self-tangency is derived from the first one in the following way: the two arrows slide over the arrow $d$ but their mutual position does not change. We show an example in  Fig.~\ref{self2}. Notice that in the thick part of the circle there aren't any other heads or foots of arrows.

To each figure of an edge corresponds a second figure. In the first line on the left we have the smoothing with constant weight which is already known. (Of course we draw only the case when the triple crossing contributes, i.e. global type $l_b$ or $l_c$.) This determines the new smoothing on the right. We put the smoothings in rectangles for better visualizing them. In the second line we give the corresponding smoothings with linear weight (but we drop the common factor $W(p)$). Here we often use the identity mentioned at the beginning of this section.

After establishing the partial smoothings for all eight local types of triple crossings we have to verify that the solution is consistent for the loops in $\Gamma$. We do this for the boundary of the three different types of 2-faces "1-7-2-5", "1-6-4-7" and "1-5-3-6" of the cube. The rest of the 2-faces is completely analogous and we left the verification to the reader.

We proceed then in exactly the same way for the global type $r$. But notice the difference which comes from the fact that we have broken the symmetry: the global type $r_a$ contributes both with constant weight as well as with linear weight.

The rest of the proof is in the figures Fig.~\ref{l1-5},... Fig.~\ref{self2}.
$\Box$

\begin{figure}
\centering
\includegraphics{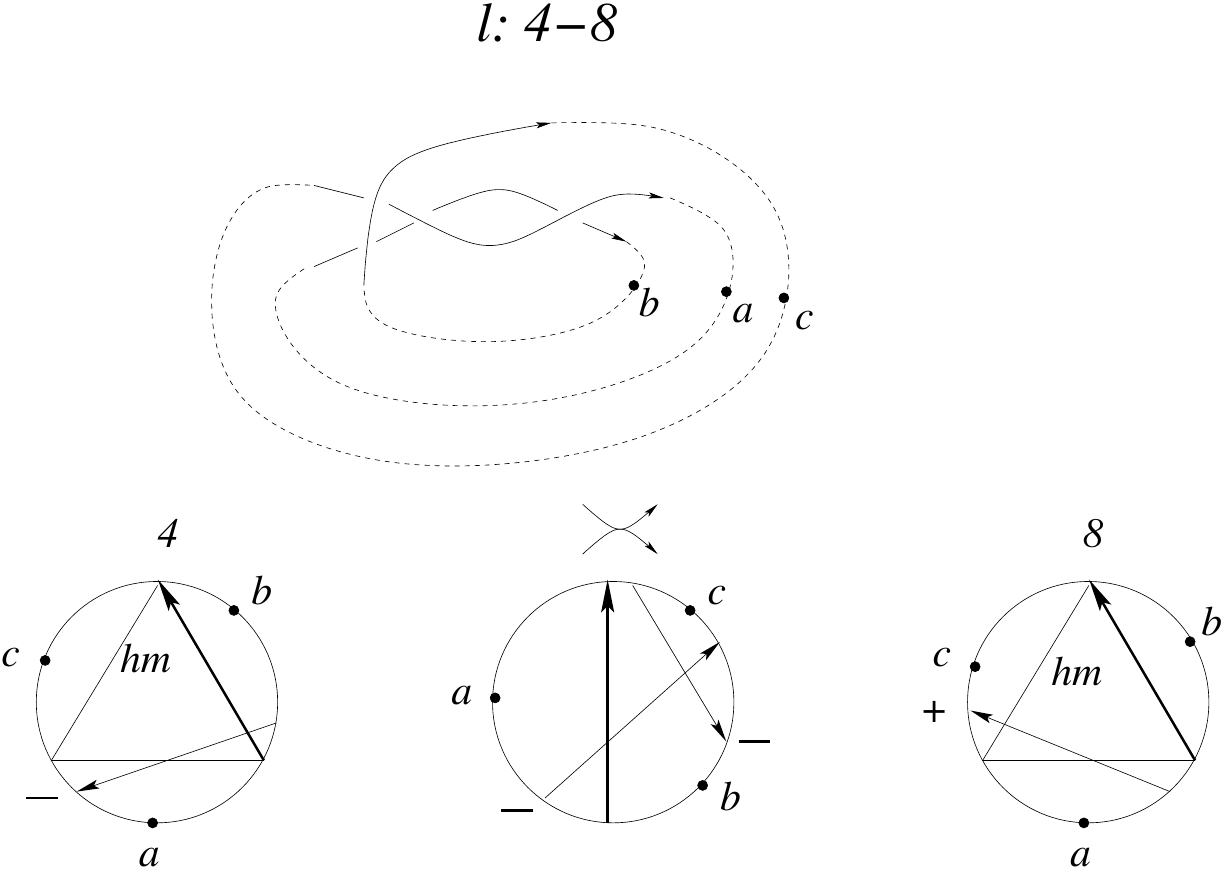}
\caption{\label{l4-8} l4-8}  
\end{figure}

\begin{figure}
\centering
\includegraphics{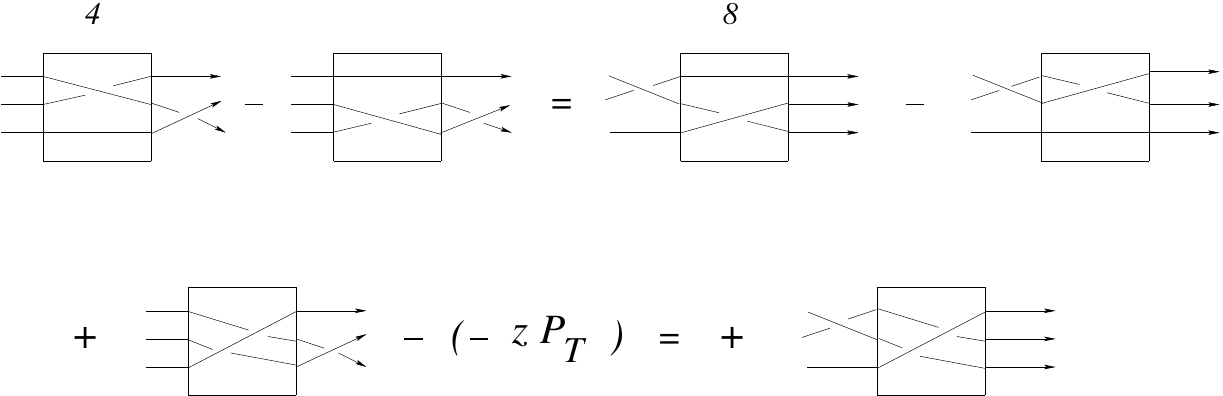}
\caption{\label{l4-8s} smoothings for l4-8}  
\end{figure}

\begin{figure}
\centering
\includegraphics{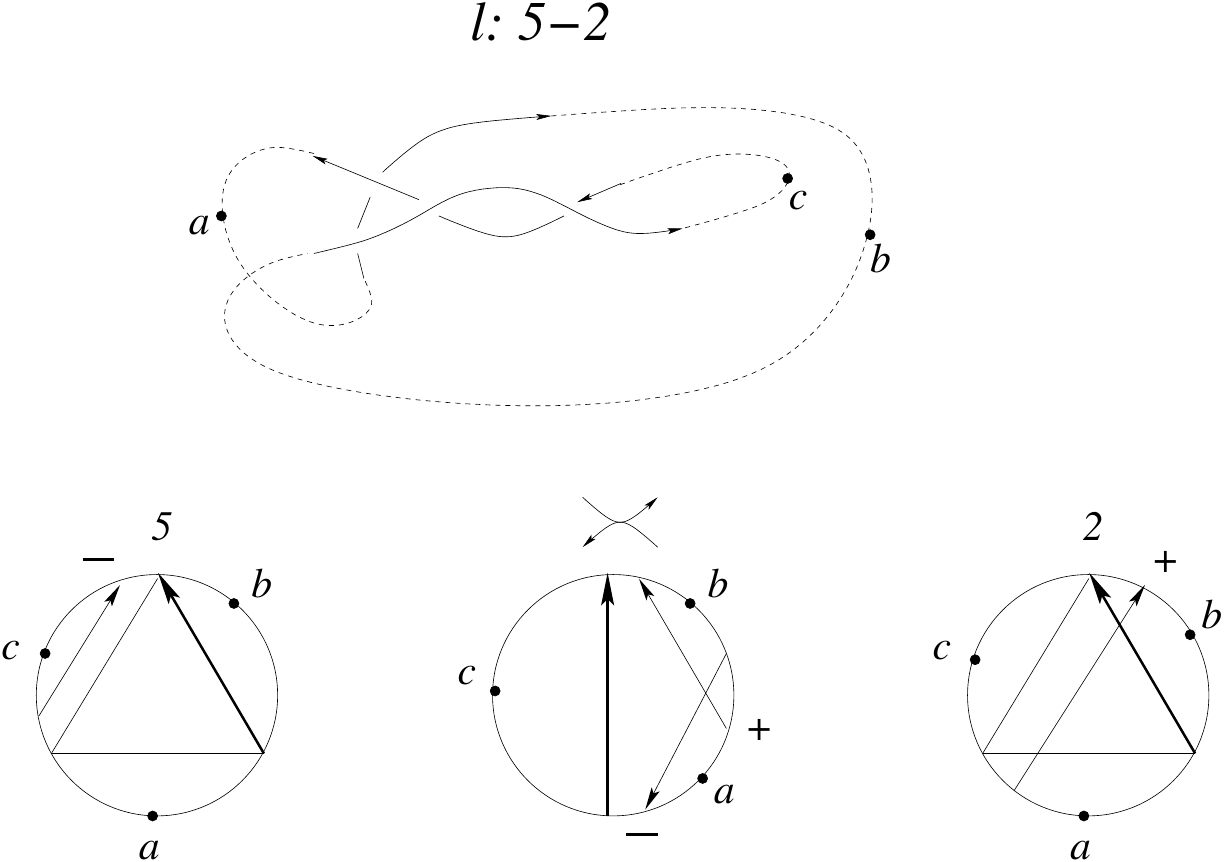}
\caption{\label{l5-2} l5-2}  
\end{figure}

\begin{figure}
\centering
\includegraphics{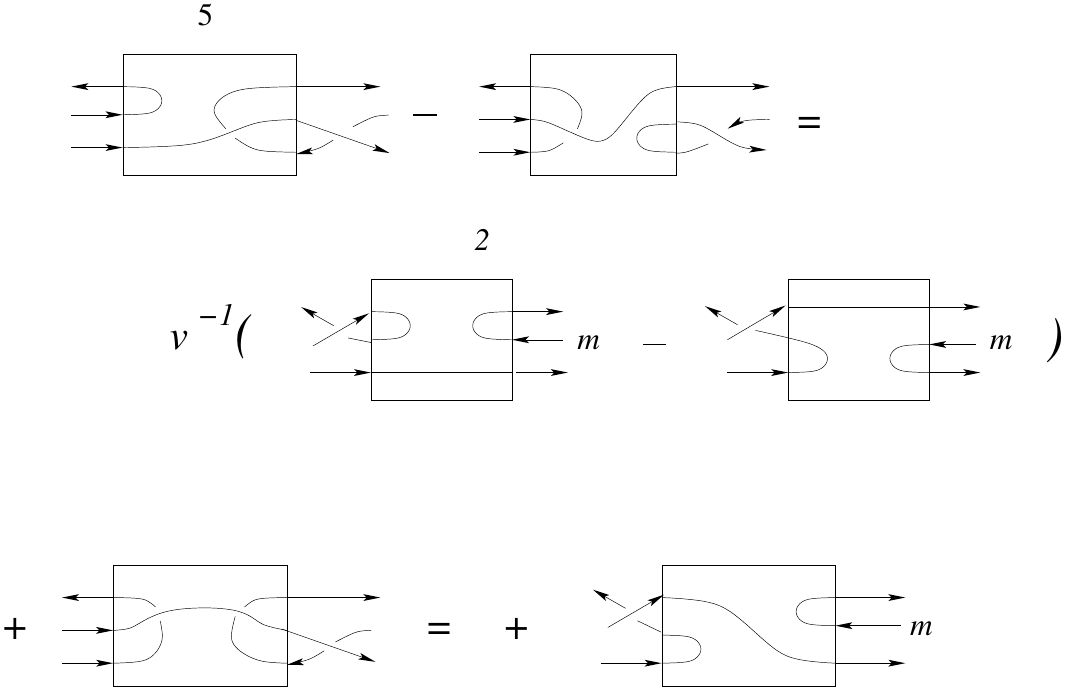}
\caption{\label{l5-2s} smoothings for l5-2}  
\end{figure}

\begin{figure}
\centering
\includegraphics{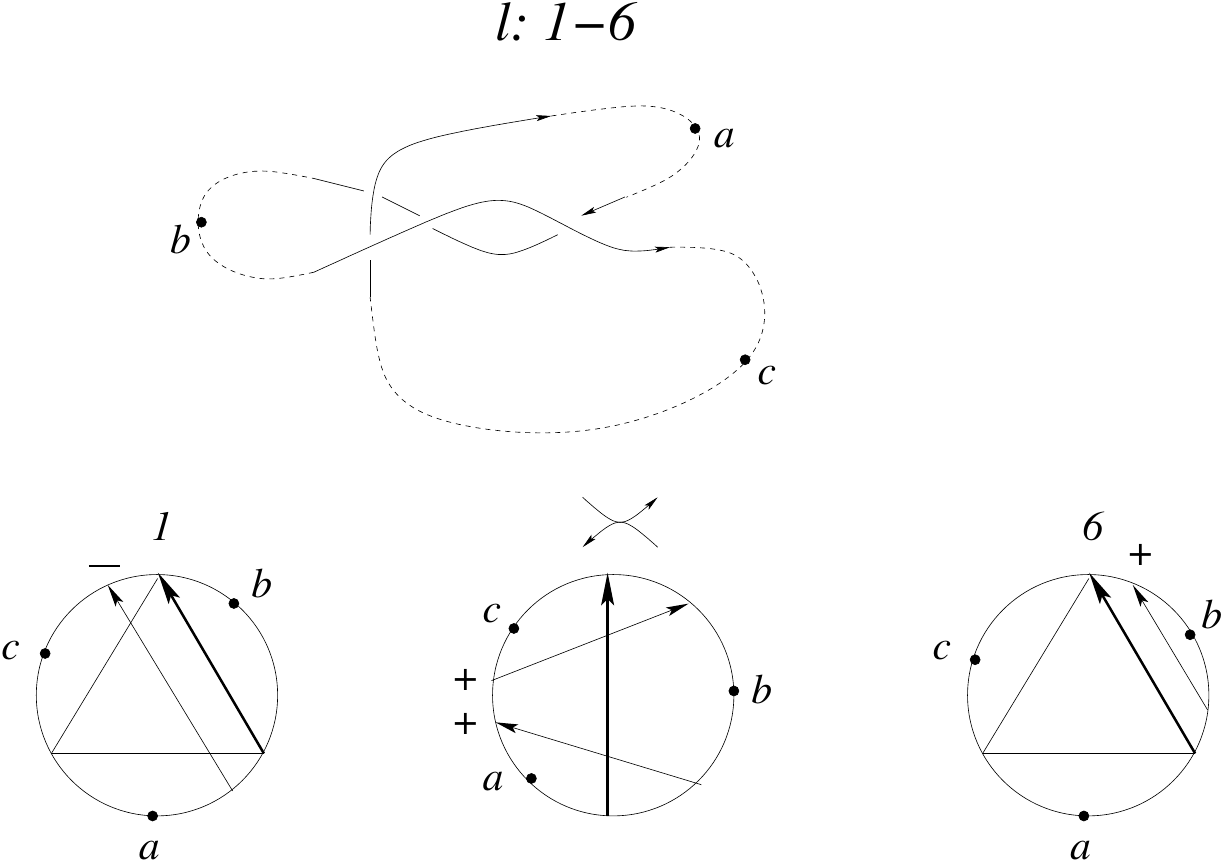}
\caption{\label{l1-6} l1-6}  
\end{figure}

\begin{figure}
\centering
\includegraphics{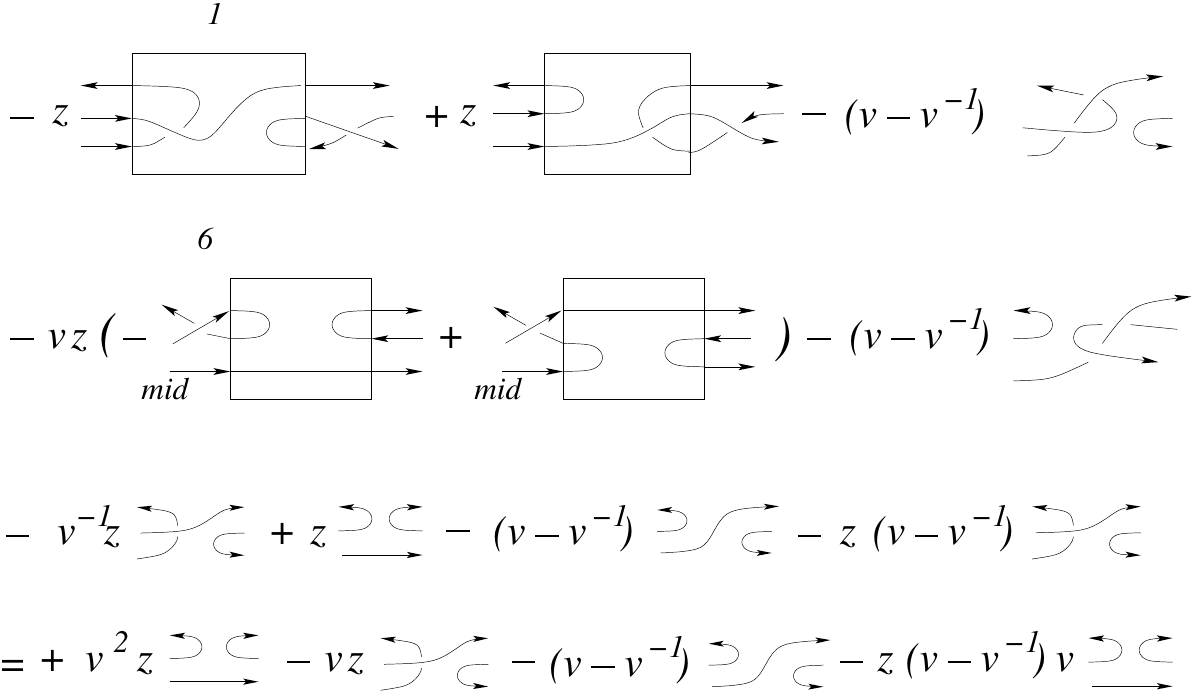}
\caption{\label{l1-6sb} l1-6sb}  
\end{figure}

\begin{figure}
\centering
\includegraphics{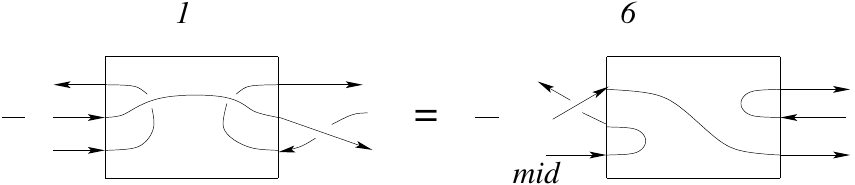}
\caption{\label{l1-6sc} l1-6sc}  
\end{figure}

\begin{figure}
\centering
\includegraphics{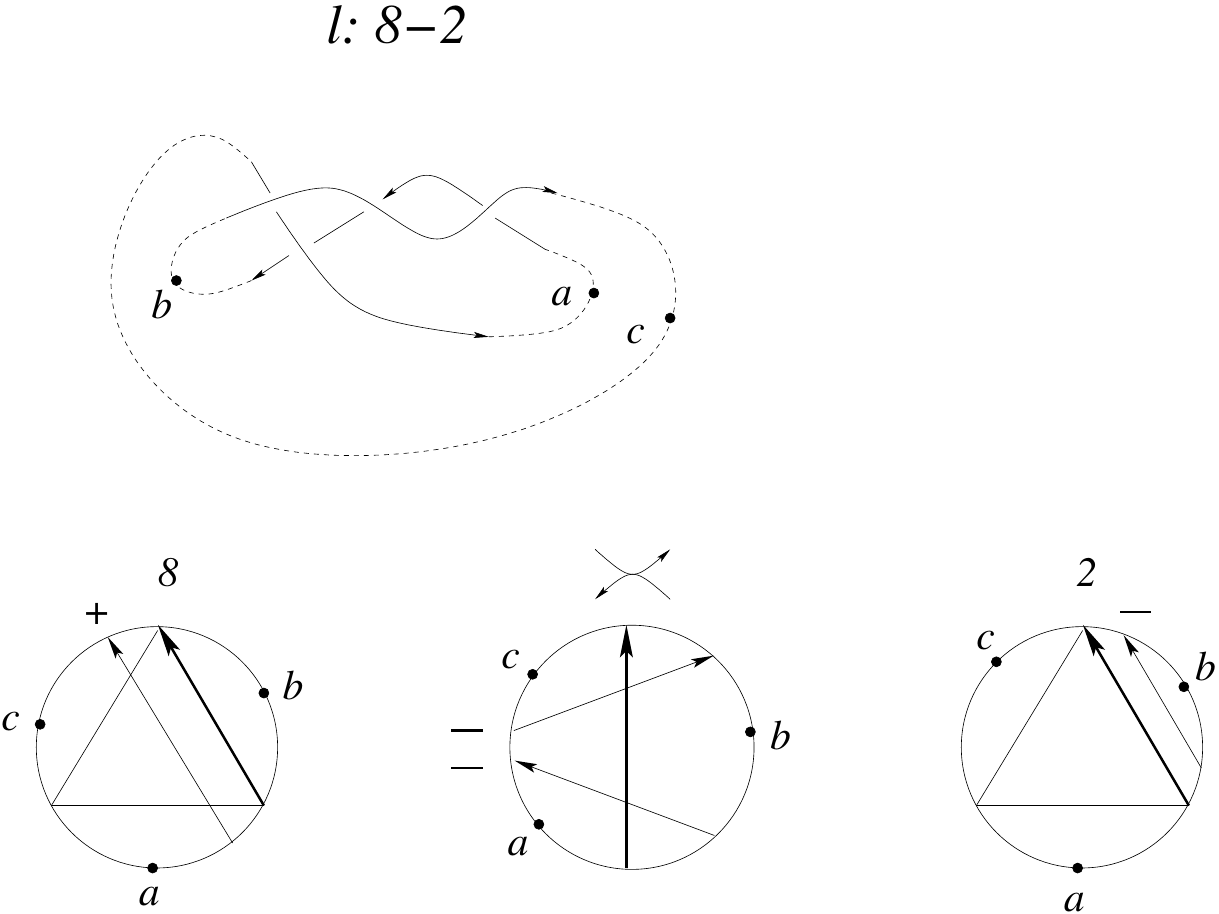}
\caption{\label{l8-2}  l8-2}  
\end{figure}

\begin{figure}
\centering
\includegraphics{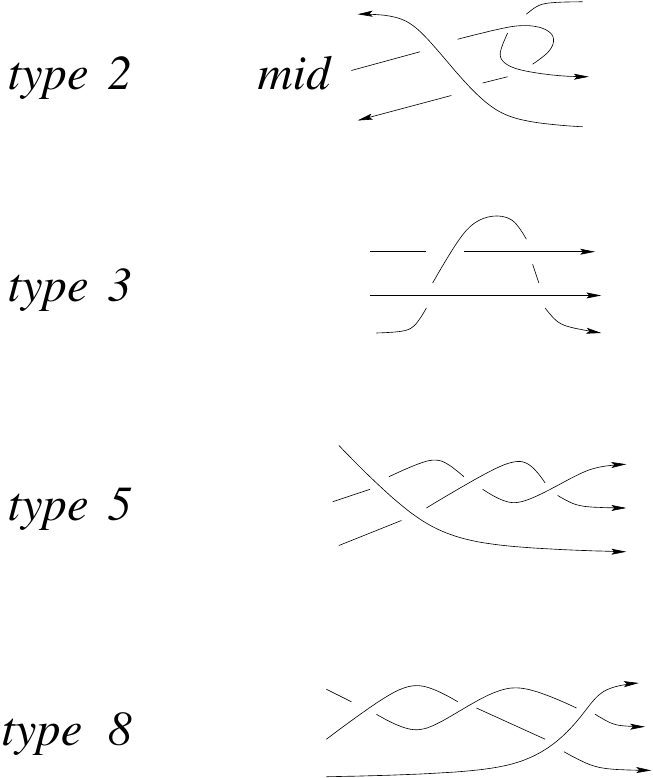}
\caption{\label{altpartsmo} alternative solution for the cube equations without type
$II_0^-$}  
\end{figure}

\begin{remark}
It is not hard to see that one can find a solution $R^{(1)}_{reg}(A)(m)=0$ for the meridians of $\Sigma^{(2)}_{trans-self}$ without using self-tangencies at all. But this solution does not survive for the loops in $\Gamma$. Hence, the contributions of the self-tangencies with opposite tangent direction are necessary.

\begin{figure}
\centering
\includegraphics{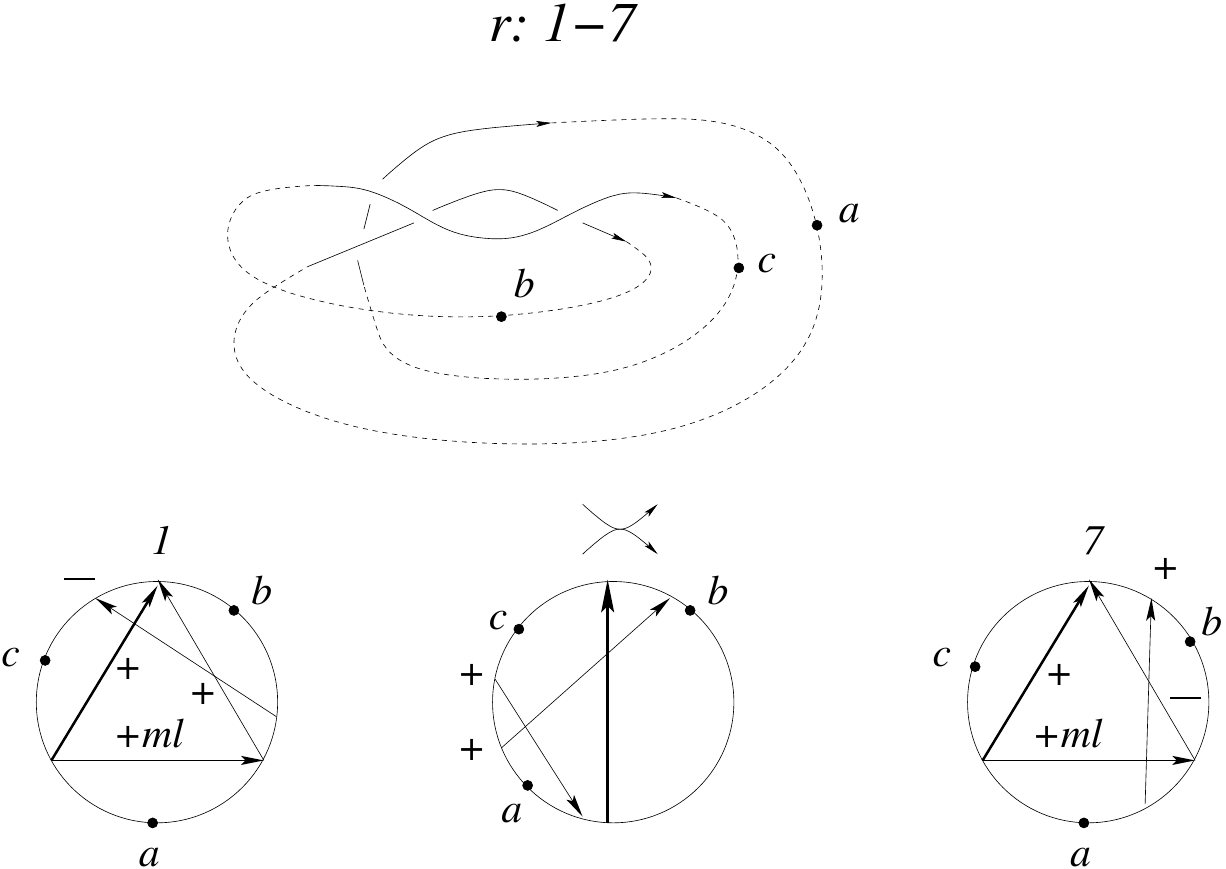}
\caption{\label{r1-7}  r1-7}  
\end{figure}

\begin{figure}
\centering
\includegraphics{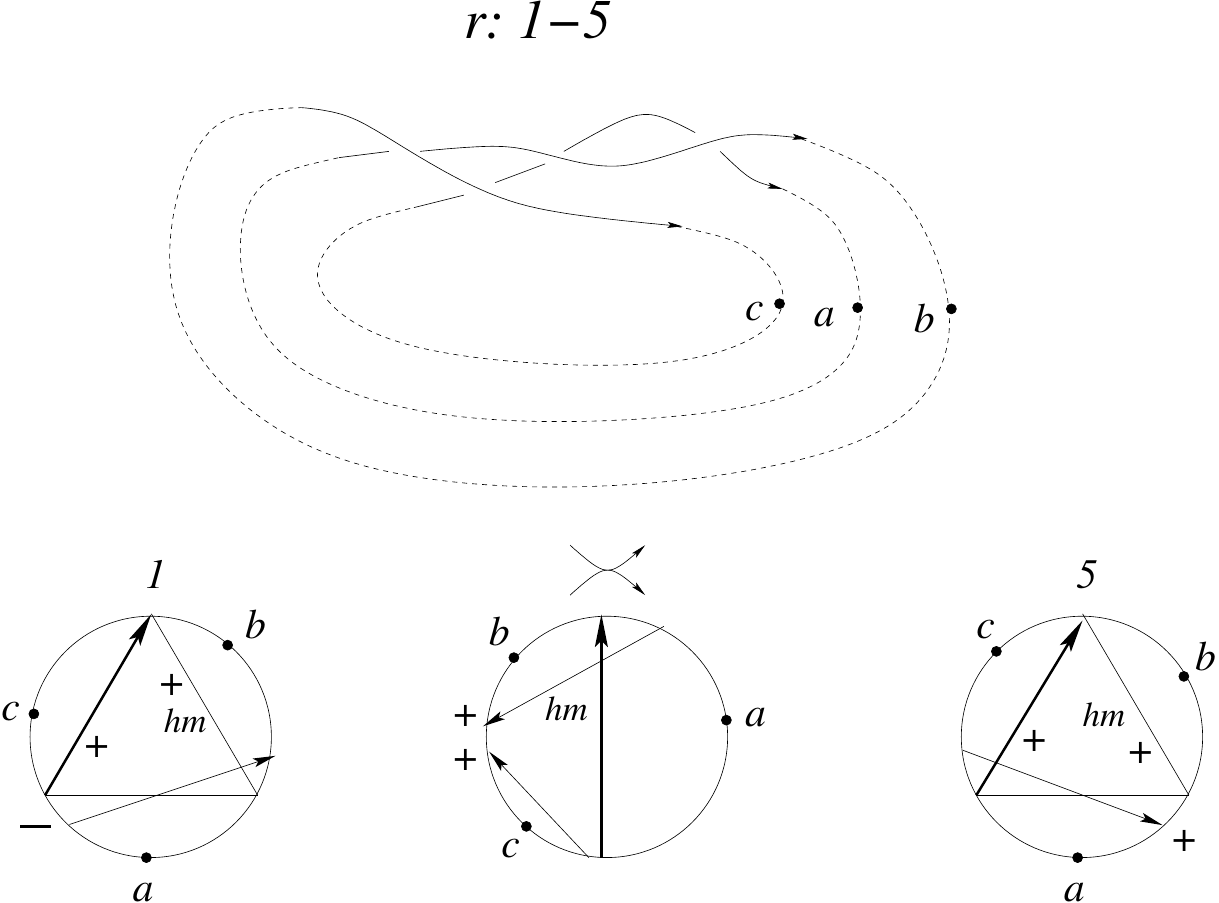}
\caption{\label{r1-5}  r1-5}  
\end{figure}

\begin{figure}
\centering
\includegraphics{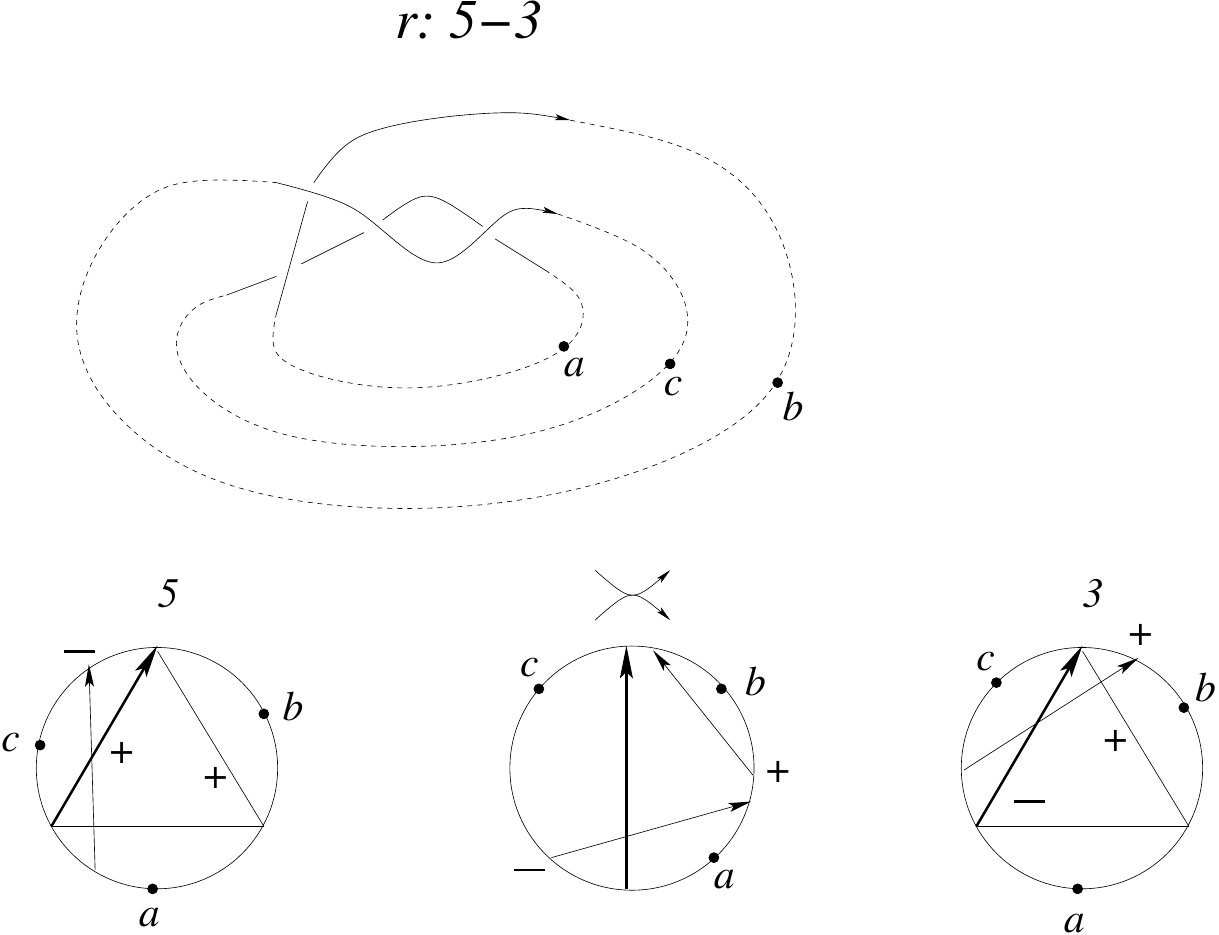}
\caption{\label{r5-3} r5-3}  
\end{figure}

However, our solution of the cube equations is not unique!

\begin{figure}
\centering
\includegraphics{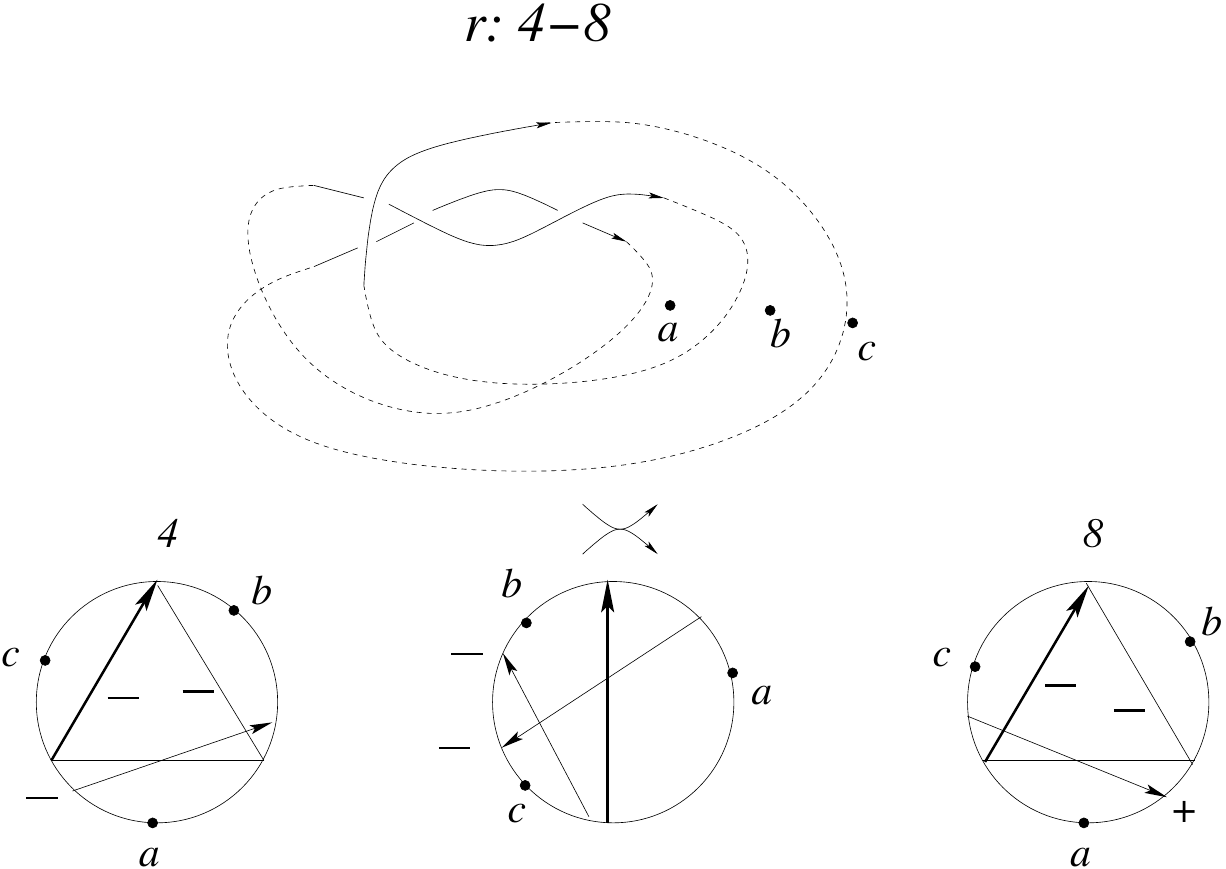}
\caption{\label{r4-8}  r4-8}  
\end{figure}

\begin{figure}
\centering
\includegraphics{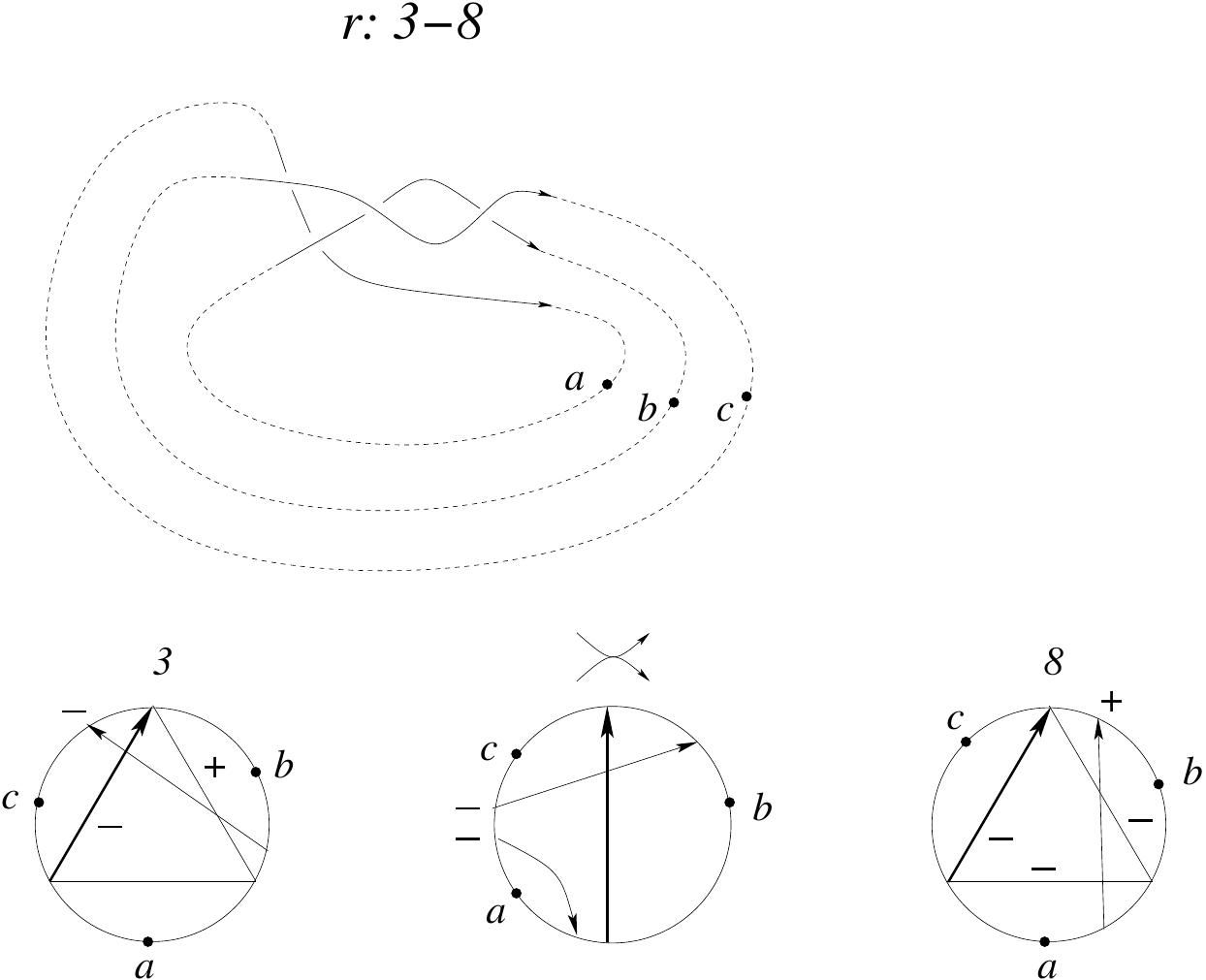}
\caption{\label{r3-8} r3-8}  
\end{figure}

The reader can easily verify (using our figures) that the following is another solution of the cube equations: we do not use self-tangencies of type $II^-_0$ at all. But for the contributions with constant weight of the triple crossings in the boundary of the 2-face "2-5-3-8" we add two additional negative crossings between two branches, as shown in Fig.~\ref{altpartsmo}. Moreover, we add $sign(p)W(p)zP_T$ to the contributions of the self-tangencies of type $II^+_0$.

\begin{figure}
\centering
\includegraphics{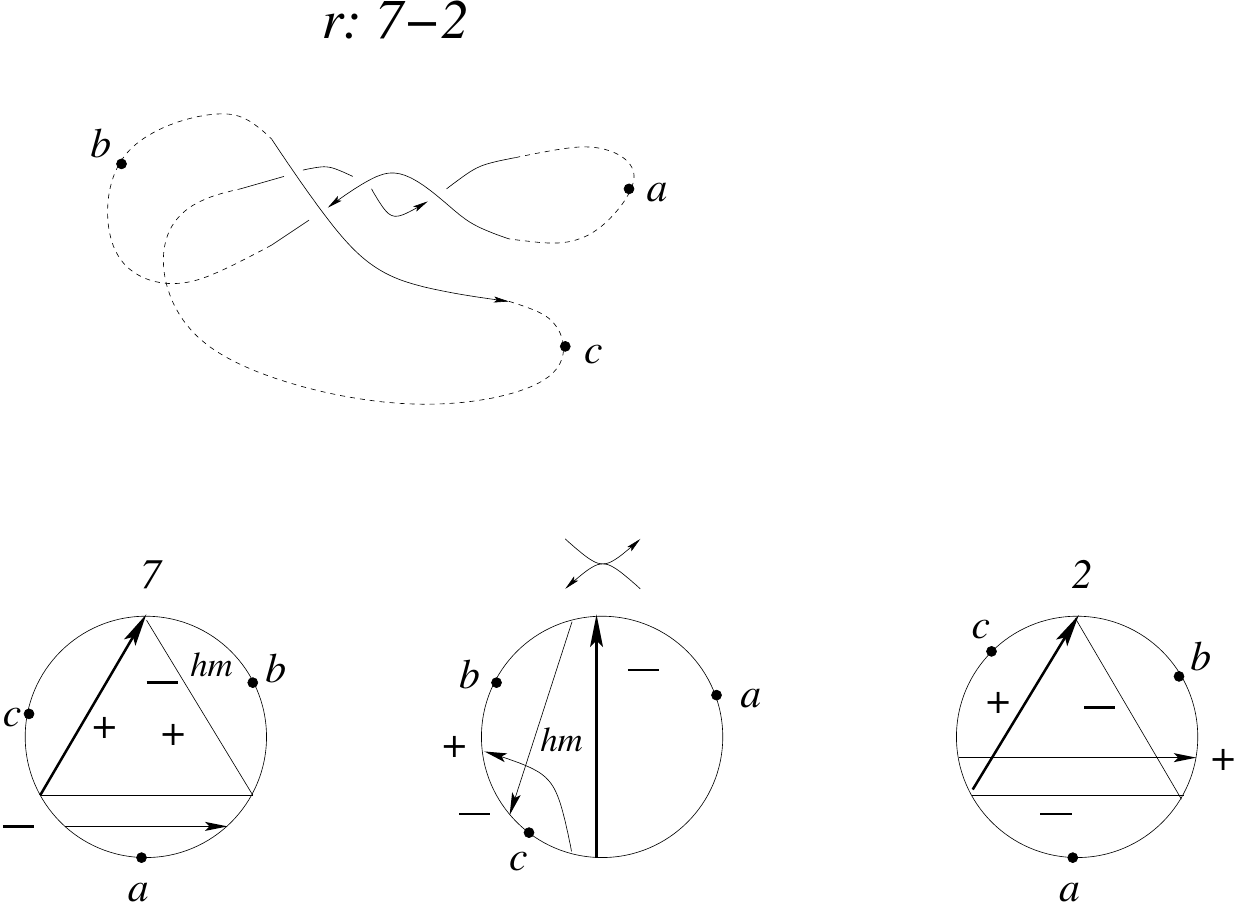}
\caption{\label{r7-2} r7-2}  
\end{figure}

\begin{figure}
\centering
\includegraphics{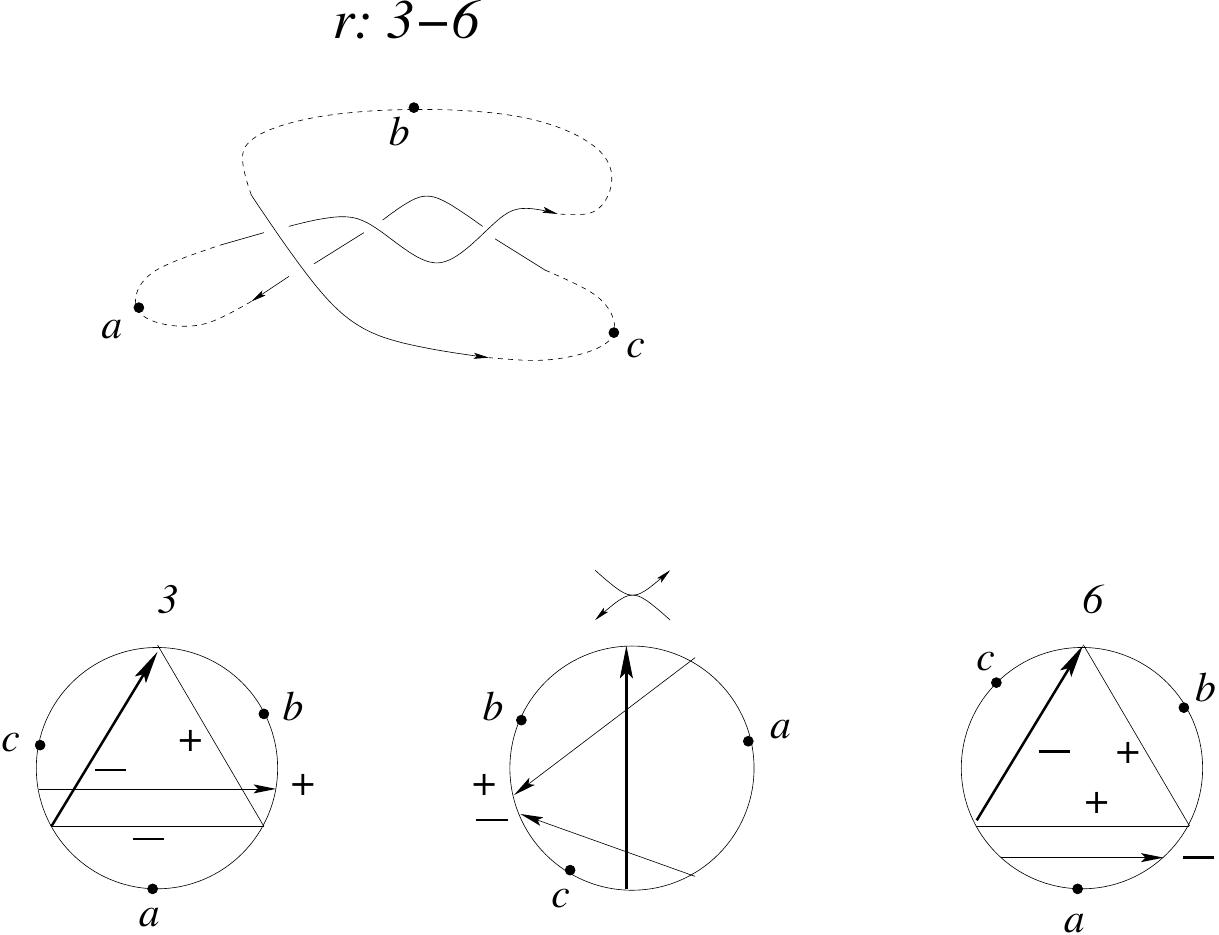}
\caption{\label{r3-6} r3-6}  
\end{figure}

\begin{figure}
\centering
\includegraphics{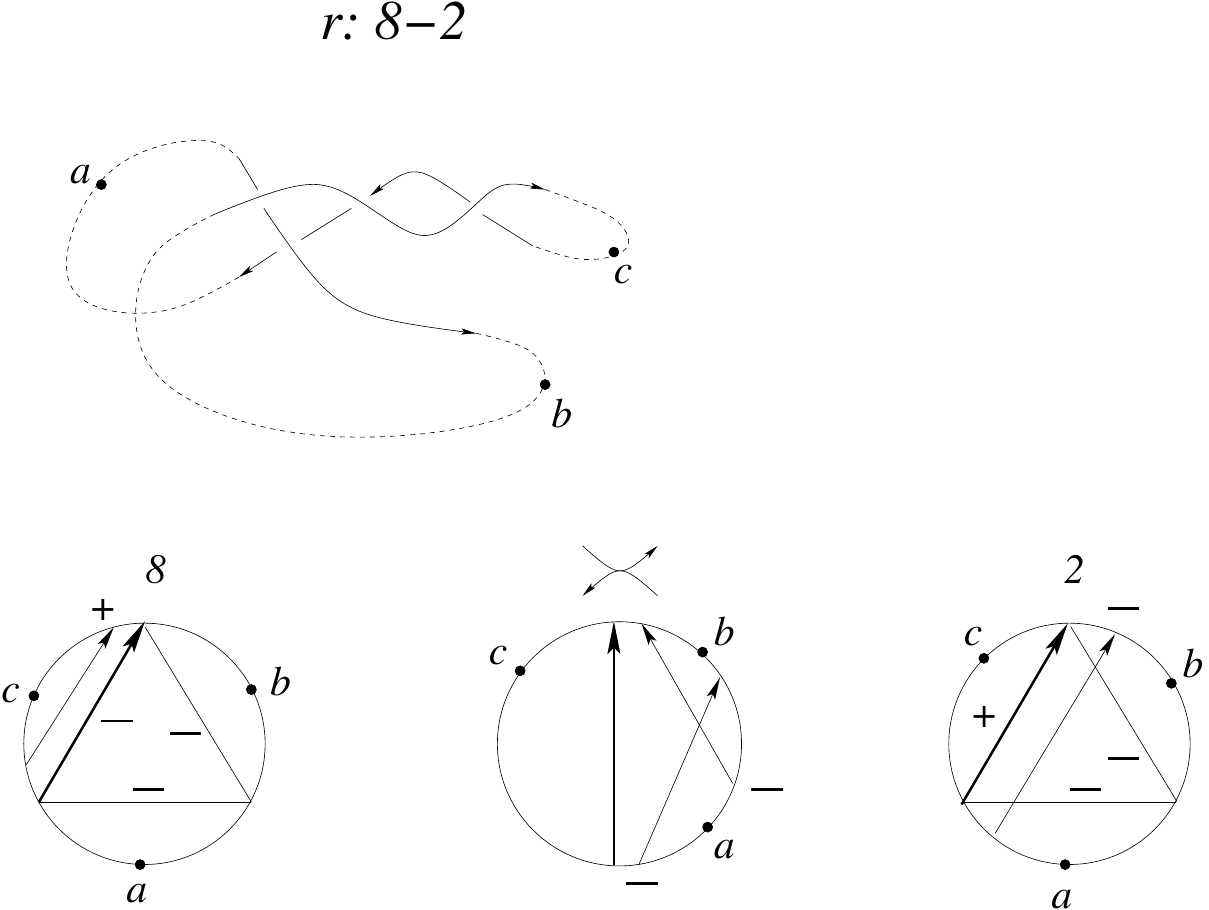}
\caption{\label{r8-2} r8-2}  
\end{figure}

\begin{figure}
\centering
\includegraphics{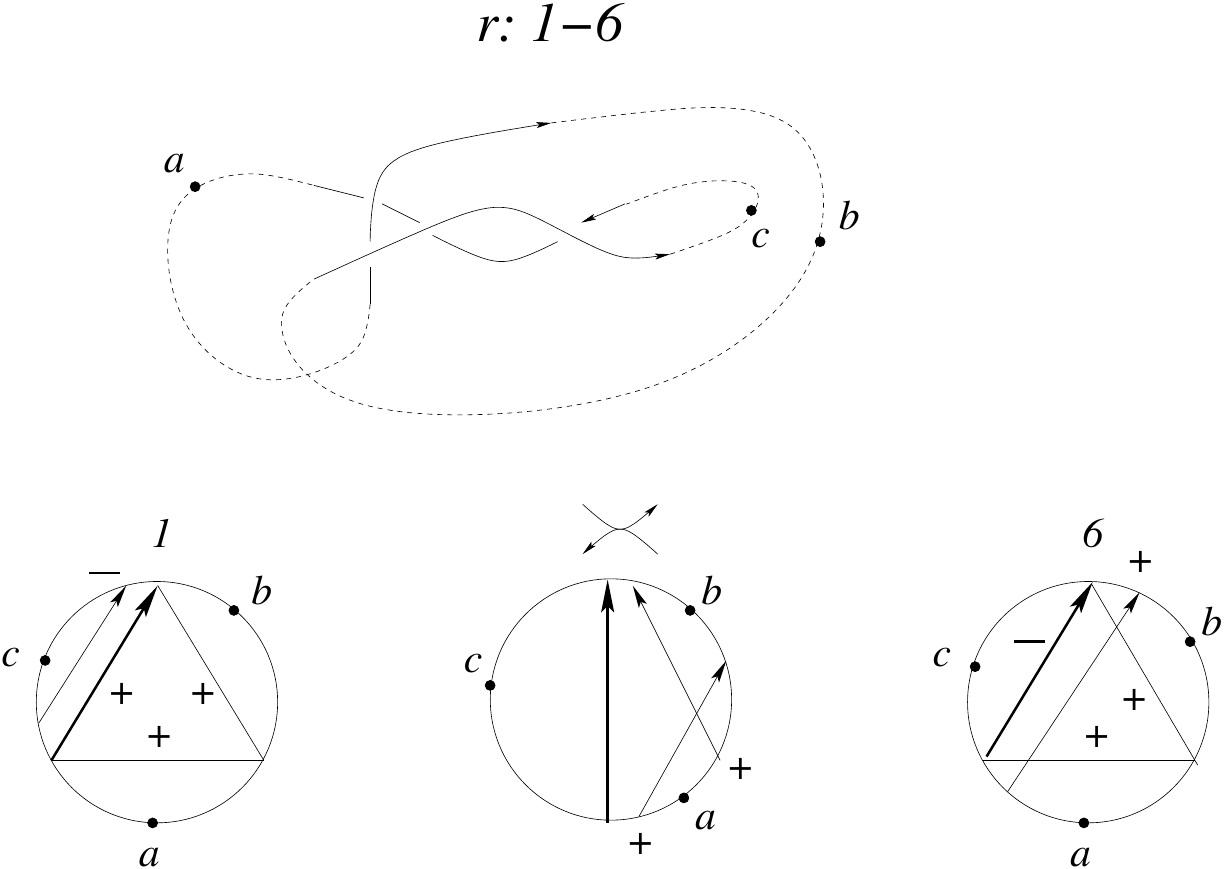}
\caption{\label{r1-6} r1-6}  
\end{figure}

It turns out that the new 1-cocycle has also the scan-property. We don't know whether it contains the same information as $R^{(1)}_{reg}(A)$. However, we will need the simpler partial smoothings given above in the construction of the 1-cocycle $\bar R^{(1)}$ (compare Lemma 10).
\end{remark}

\begin{figure}
\centering
\includegraphics{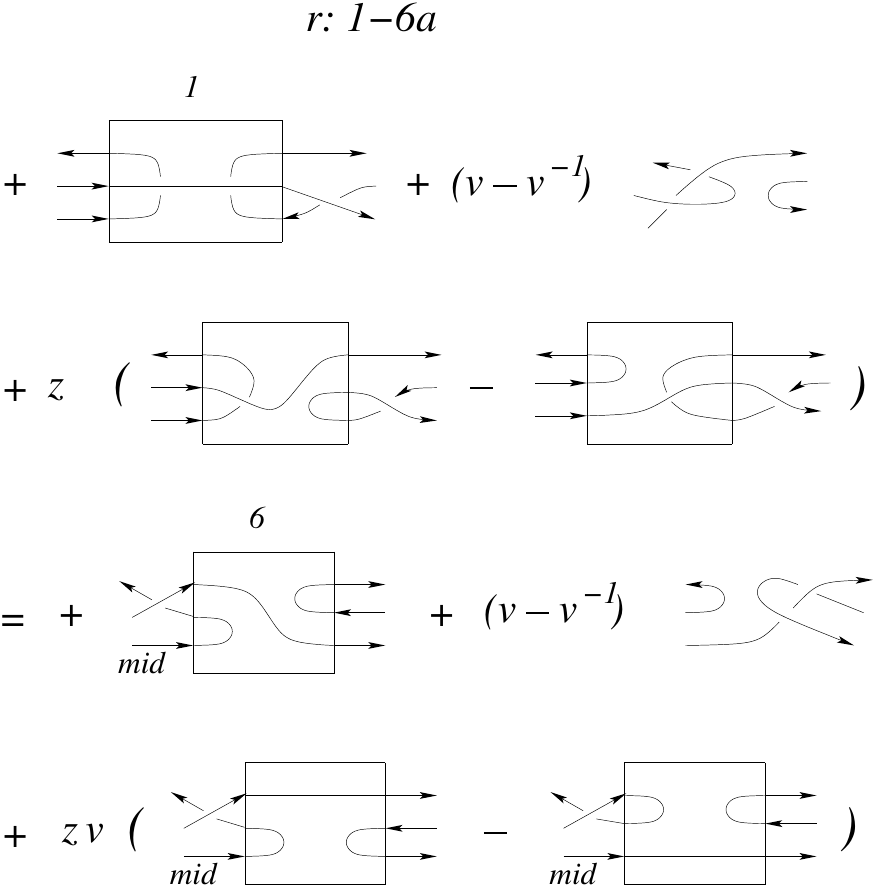}
\caption{\label{r1-6sa} smoothings for r1-6}  
\end{figure}

\begin{figure}
\centering
\includegraphics{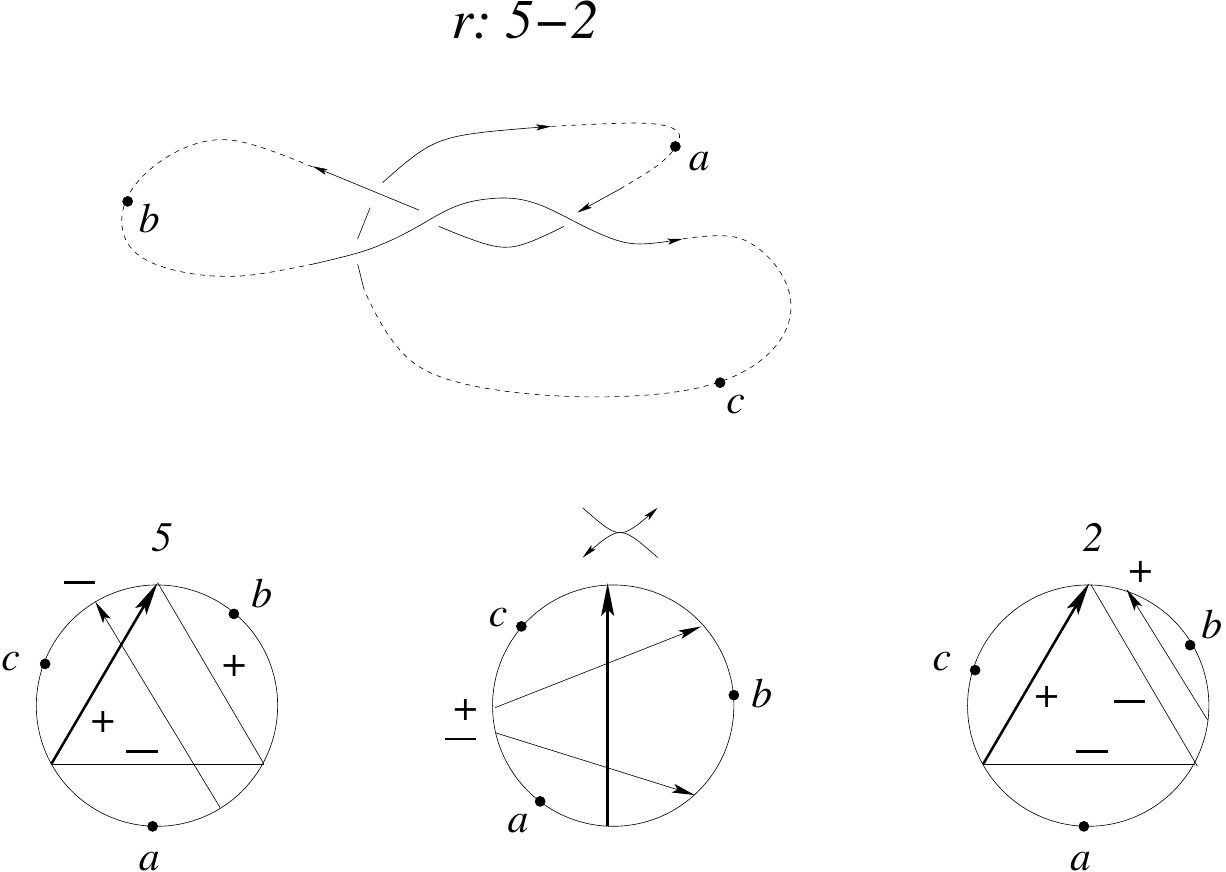}
\caption{\label{r5-2}  r5-2}  
\end{figure}

\begin{figure}
\centering
\includegraphics{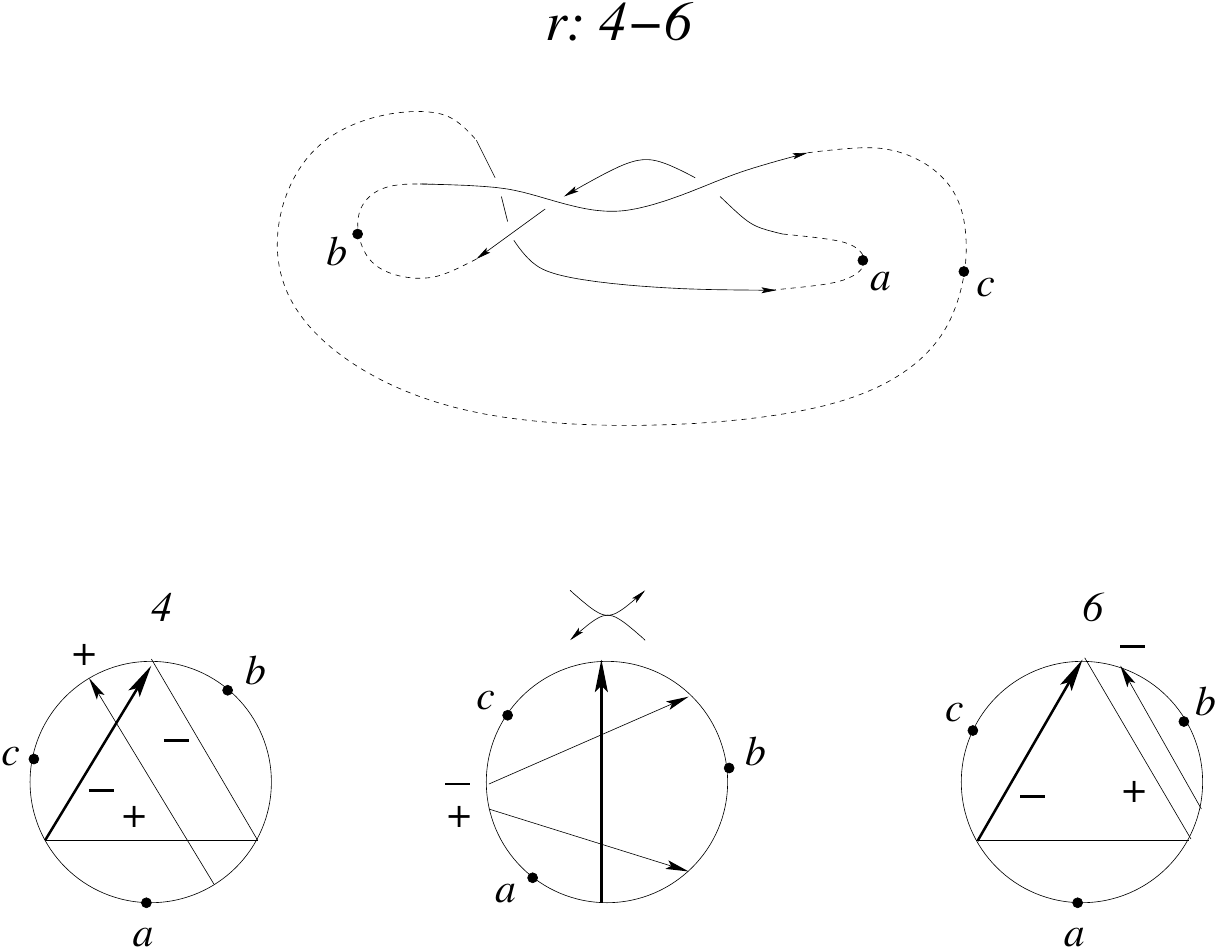}
\caption{\label{r4-6}  r4-6}  
\end{figure}

\begin{figure}
\centering
\includegraphics{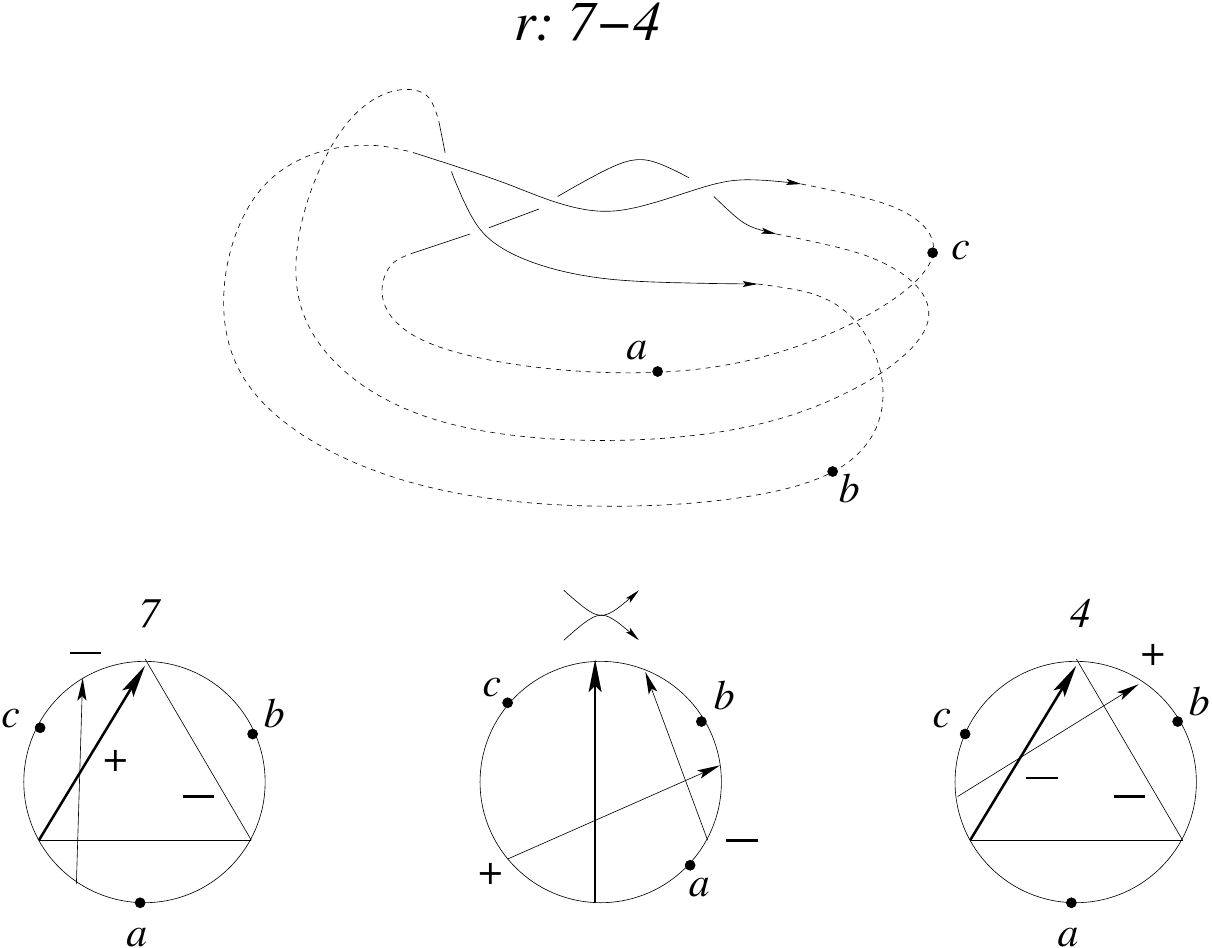}
\caption{\label{r7-4}  r7-4}  
\end{figure}

\section{Solution with quadratic weight of the cube equations}

We have to replace the constant weight $1$ by the linear weight $W_1(p)$ and the linear weight $W(p)$ by the quadratic weight $W_2(p)$. All partial smoothings are exactly the same as in the previous section besides an adjustment of the weights by a constant which depends on the global and on the local type.

The weights $W_1(p)$ and $W_2(p)$ for positive triple crossings (local type 1) were introduced in Definition 12 respectively Definition 11.

\begin{definition}
For the global types $l_b$, $l_c$ and $r_b$ the weights $W_1(p)$ and $W_2(p)$ do not depend on the local type.
For the global type $r_a$ the weights are defined as follows:

$W_1(type 1)=W_1(type 3)=W_1(type 7)=W_1(type 8)=W_1(p)$

$W_2(type 1)=W_2(type 3)=W_2(type 7)=W_2(type 8)=W_2(p)$

$W_1(type2)=W_1(type5)=W_1(p)-1$

$W_1(type4)=W_1(type6)=W_1(p)+1$

$W_2(type4)=W_2(type5)=W_2(p)-1$

$W_2(type2)=W_2(type6)=W_2(p)+1$
\end{definition}

\begin{definition}
Let $s$ be a generic oriented arc in $M^{reg}_T$ (with a fixed abstract closure $T \cup \sigma$ to an oriented circle).  The 
1-cochain $R^{(1)}$ is defined by \vspace{0,2 cm}

$ R^{(1)}(s)=\sum_{p \in s \cap l_c}sign(p)W_1(p)T_{l_c}(type)(p) + $ \\ \vspace{0,1 cm}

$ \sum_{p \in s \cap r_a}sign(p)[W_1(type)(p)T_{r_a}(type)(p)  + zW_2(type)(p)T(type)(p)] +$
\\ \vspace{0,2 cm}

$ \sum_{p \in s \cap r_b}sign(p)zW_2(p)T(type)(p) +$
\\ \vspace{0,2 cm}

$ \sum_{p \in s \cap l_b}sign(p)zW_2(p)T(type)(p) +$
\\ \vspace{0,2 cm}

$ \sum_{p \in s \cap II^+_0}sign(p)W_2(p)T_{II^+_0}(p)+$ 
\\ \vspace{0,2 cm}

$ \sum_{p \in s \cap II^-_0}sign(p)W_2(p)T_{II^-_0}(p)$ \vspace{0,2 cm}

Here all partial smoothings are exactly the same as in the previous section.
\end{definition}

\begin{proposition}
Let $m$ be a meridian of $\Sigma^{(2)}_{trans-self}$ or a loop in $\Gamma$. Then $R^{(1)}(m)=0$.
\end{proposition}

{\em Proof.} 
For the global type $l$ the proof is exactly the same as in the case of linear weight. Indeed, the fourth arrow which is not in the triangle is always almost identical with an arrow of the triangle. Consequently, in the case $l_c$ it can not be a r-crossing with respect to $hm$. In the case $l_b$ it 
could be a r-crossing with respect to some f-crossing. But the almost identical arrow in the triangle would be a r-crossing for the same f-crossing too. Their contributions cancel out, because they have different writhe.

The mutual position of the two arrows for a self-tangency does not change. Only the position with respect to the distinguished crossing $d$ changes. But this does not count because $d$ consists in fact of two crossings with opposite writhe. It remains to consider the edges "2-7" and "3-6" for $l_c$, where one of the two self-tangencies has a new f-crossing with respect to the other self-tangency. But one sees immediately from the figures that this new f-crossing is exactly the crossing $hm$ in the triple crossings. Consequently, all three summands in the first line of the figures are just multiplied by the same $W_1(p)$ and the equations are still satisfied.

Exactly the same arguments apply in the case $r_b$ too. In the case $r_c$ there are no contributions at all. The differences appear in the case $r_a$.

We examine the figures (don't forget the degenerate configurations):

edge "1-7": type 1 and type 7 share the same $W_1$ and $W_2$

edge "3-8": type 3 and type 8 share the same $W_1$ and $W_2$

edge "4-6":  same $W_1$ but $W_2+1$ for type 4 and $W_2-1$ for type 6

edge "5-2": same $W_1$ but $W_2+1$ for type 5 and $W_2-1$ for type 2

edge "1-5": $W_1-1$ and $W_2-1$ for type 1 with respect to type 5

edge "7-4": $W_1-1$ and $W_2+1$ for type 4 with respect to type 7

edge "5-3": $W_1+1$ and $W_2+1$ for type 5 with respect to type 3

edge "4-8": $W_1-1$ and $W_2+1$ for type 4 with respect to type 8

edge "7-2": $W_1+1$ and $W_2-1$ for type 2 with respect to type 7

edge "3-6": $W_1-1$ and $W_2-1$ for type 6 with respect to type 3

edge "8-2": $W_1-1$ and $W_2+1$ for type 8 with respect to type 2

edge "1-6": $W_1+1$ and $W_2+1$ for type 1 with respect to type 6

For both self-tangencies in "8-2" as well as in "1-6" we have the same $W_2$ as in type 8 respectively type 1.

It follows now easily that all cube equations are satisfied with the adjustments of the weights given in Definition 20.

$\Box$

\section{The 1-cocycle $R^{(1)}_{reg}$, the scan-property and first examples}

The 1-cochain $R^{(1)}_{reg}(A)$ was defined in Definition 19 and the 1-cochain $R^{(1)}$ in Definition 21. We have proven that they satisfy the global positive tetrahedron equation as well as the cube equations. Moreover, we have proven that  their value is 0 on the meridians of the strata of codimension 2 which correspond to a self-tangency in an ordinary flex as well as to the transverse intersections of two strata of codimension 1. These are the first four points from the list in Section 2. These are the only strata of codimension 2 which appear in generic homotopies of  generic regular isotopies in $M_T^{reg}$ (compare \cite{FK}). Consequently, we have proven the following theorem.

\begin{theorem}
The 1-cochains $R^{(1)}_{reg}(A)$ and $R^{(1)}$ are 1-cocycles in $M_T^{reg}$.
\end{theorem}

It is difficult to apply this theorem. The natural loops $rot$ and $hat$ in $M_T$ are not in $M_T^{reg}$. Of course we could approximate them by loops in $M_T^{reg}$ using Whitney tricks. However, this approximation is not unique and leads to very long calculations. But we can apply it to the loop which consists of dragging a knot through itself by a regular isotopy (compare the Introduction). We have calculated by hand that $R^{(1)}_{reg}(A)(drag 3_1^+)=0$. So, we don't know whether or not $R^{(1)}_{reg}(A)$ represents the trivial cohomology class in $M_T^{reg}$. Let $K$ and $K'$ be two long knots. It would be very interesting to calculate examples (with a computer program) of $R^{(1)}_{reg}(A)$ for the loop which consists of dragging $Cab_2(K)$ through $Cab_2(K')$ and then dragging $Cab_2(K')$ through $Cab_2(K)$ (compare the Introduction).

\begin{remark}
Even the specialization of  $R^{(1)}_{reg}(A)$ for $v=1$ can not be extended to a 1-cocycle for isotopies instead of regular isotopies.
\end{remark}

\begin{figure}
\centering
\includegraphics{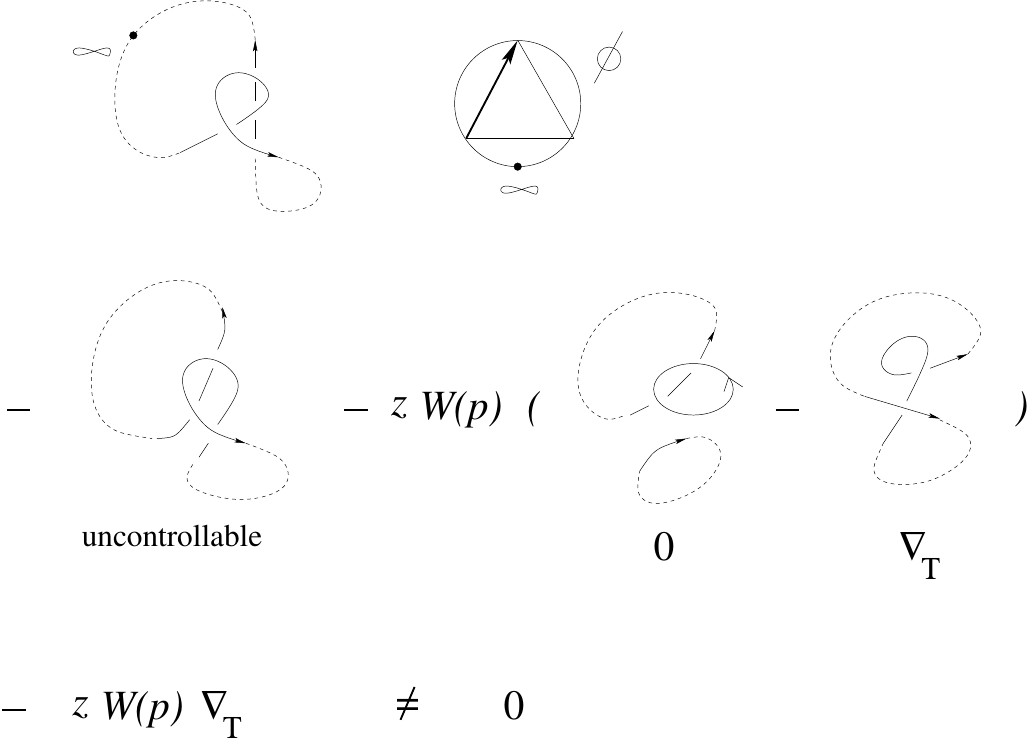}
\caption{\label{trans-cusp} the regular 1-cocycle from the Conway polynomial is
not controllable on meridians of $\Sigma_{trans-cusp}^{(2)}$}  
\end{figure}

Indeed, the intersection with a stratum $\Sigma^{(1)}_{cusp}$ does not change the calculation of the Conway invariants for the intersections with the other strata of $\Sigma^{(1)}$. There are two types of strata in $\Sigma^{(1)}_{cusp}$: the new crossing could be of type $0$ or of type $1$. Let's call them $\Sigma^{(1)}_{cusp 0}$ and $\Sigma^{(1)}_{cusp 1}$ respectively. The value of the 1-cocycle on a meridian of 
$\Sigma^{(1)} \cap \Sigma^{(1)}_{cusp 0}$ is 0 because we have set $v=1$. However, it is not controllable at all on a meridian of  $\Sigma^{(1)} \cap \Sigma^{(1)}_{cusp 1}$. Indeed, the new crossing of type $1$ could be a f-crossing for the stratum in $\Sigma^{(1)}$. Moreover, the 1-cocycle can behave uncontrollable on a meridian of $\Sigma^{(2)}_{trans-cusp}$ (compare Section 2). We show an example in  Fig.~\ref{trans-cusp}.

The good news are that $R^{(1)}_{reg}(A)$ and $R^{(1)}$ turn out to have the {\em scan-property} (compare the Introduction).

\begin{lemma}
Let $T$ be a diagram of a string link and let $t$ be a Reidemeister move of type II for $T$. Then the contribution of $t$ to $R^{(1)}_{reg}(A)$ respectively to $R^{(1)}$ does not change if a branch of $T$ is moved under $t$ from one side of $t$ to the other.
\end{lemma}
{\em Proof.} 
Besides $d$ of $t$ there are two crossings involved with the branch which goes under $t$. Moving the branch under $t$ from one side to the other slides the heads of the corresponding arrows over $d$, see e.g. Fig.~\ref{l1-7} and  Fig.~\ref{r1-7} as well as Fig.~\ref{l4-6} and Fig.~\ref{r4-6}. Consequently, the f-crossings and hence $W(p)$ do not change. Notice that the mutual position of the two arrows does not change. Consequently, there are no new r-crossings and $W_2(p)$ does not change neither. On the other hand it is clear that the partial smoothings of $t$ on both sides are regularly isotopic.

$\Box$

 The contribution of a Reidemeister I move will be defined in the next section (Definition 22). It does not depend on the place in the diagram and hence it has the scan-property.

We are now ready to prove Theorem 2 of the Introduction.

{\em Proof of Theorem 2.} Let $T$ be a diagram of a string link and let $s$ be a regular isotopy which connects $T$ with a diagram $T'$. We consider the loop 

-$s \circ $-$scan(T') \circ s \circ scan(T)$ in $M_T^{reg}$. This loop is contractible in $M_T^{reg}$ because $s$ and $scan$ commute and hence $R^{(1)}_{reg}(A)$ and $R^{(1)}$ vanish on this loop. Consequently, it suffices to prove that each contribution of a Reidemeister move $t$ in $s$ cancels out with the contribution of the same move $t$ in -$s$ (the signs are of course opposite). The difference for the two Reidemeister moves is in a branch which has moved under $t$. The partial smoothings are always regularly isotopic. Hence it suffices to study the weights. If $t$ is a positive triple crossing then the weights are the same just before the branch moves under $t$ and just after it has moved under $t$. This follows from the fact that for the positive global tetrahedron equation the contribution from the stratum -$P_2$ cancels out with that from the stratum $\bar P_2$ (compare the Sections 4 and 5). If we move the branch further away then the invariance follows from the  already proven fact that the values of the 1-cocycles do not change if the loop passes through a stratum of $\Sigma^{(1)} \cap \Sigma^{(1)}$. We use now again the graph $\Gamma$.
The meridian $m$ which corresponds to an arbitrary  edge of $\Gamma$ is a contractible loop in $M_T^{reg}$ no matter what is the position of the branch which moves under everything. Let's take an edge where one vertex is a triple crossing of type 1. We know from Lemma 6 that the contributions of the self-tangencies do not depend on the position of the moving branch. Consequently the contribution of the other vertex of the edge doesn't change neither because the contributions from all four Reidemeister moves together sum up to 0. Using the fact that the graph $\Gamma$ is connected we obtain the invariance with respect to the position of the moving branch for all (regular) Reidemeister moves $t$.

$\Box$

Notice that $R^{(1)}_{reg}(A)$ and $R^{(1)}$ do not have the scan-property for a branch which moves over everything else because the contributions of the strata $P_3$ and -$\bar P_3$ in the positive tetrahedron equation do not cancel out at all.

We start with very simple examples and we give all details of the calculations in order to make it easier for the reader to become familiar with the 1-cocycles $R^{(1)}_{reg}(A)$ and $R^{(1)}$.
In all examples we give names $x_1, x_2,...$ to the Reidemeister moves and sometimes we give names $c_1, c_2,...$ to the crossings too.

\begin{example}
Let $T= \sigma_1$. There is a unique closure to a circle and there are two choices for the point at infinity. We chose the small curl on the second branch like shown in Fig.~\ref{sig}. The Reidemeister move $x_1$ is a self-tangency with equal tangent direction. The Reidemeister move $x_2$ is a positive triple crossing (type 1) with $sign(x_2)=-1$. We see immediately from the Gauss diagrams that $W(x_2)=W_2(x_2)=0$ in the case $\infty_1$. Consequently, for the choice $\infty_1$ we obtain\vspace{0,2 cm} 

\begin{figure}
\centering
\includegraphics{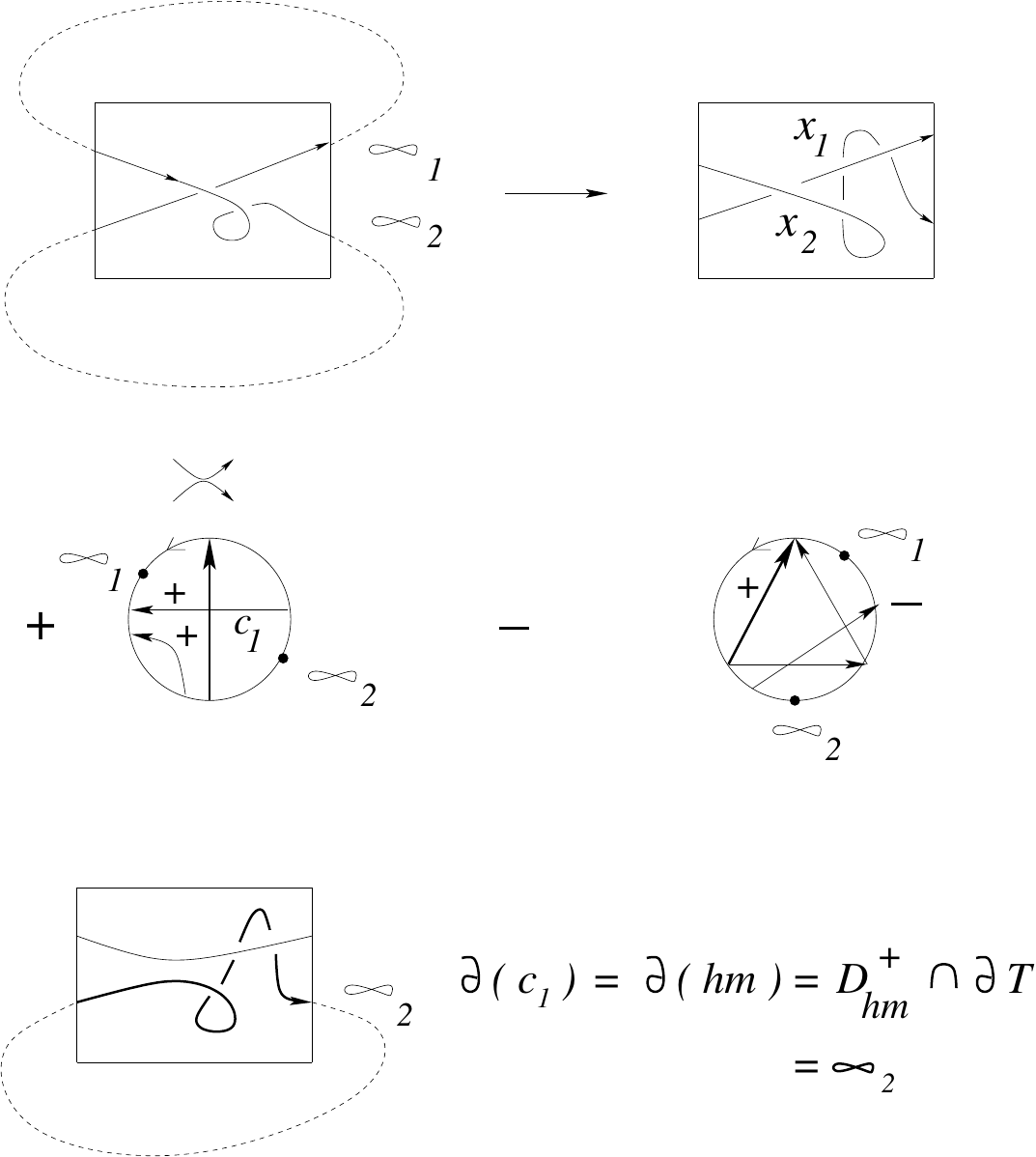}
\caption{\label{sig} the two Reidemeister moves in Example 1 for $scan_2$}  
\end{figure}

$R^{(1)}_{reg}(scan_2(\sigma_1))=R^{(1)}(scan_2(\sigma_1))=0$.\vspace{0,2 cm}

Let us consider the case $\infty_2$. The crossing $hm$ in the triangle is the only f-crossing. We obtain $\partial (hm)=  \infty_2$, $W(x_2)=1$, $W_2(x_2)=0$ (remember that the crossing $d$ is here a r-crossing too) and $W_1(x_2)=0$ too. $W(x_1)=1$ with the same grading as $hm$  and $W_2(x_1)=0$ because $d$ is a positive and a negative crossing. Consequently, \vspace{0,2 cm}

$R^{(1)}_{reg}(A= \partial \sigma_1)(scan_2(\sigma_1))=R^{(1)}(scan_2(\sigma_1))=0$.\vspace{0,2 cm}

The calculation of $R^{(1)}_{reg}(A= \infty_2)(scan_2(\sigma_1))$ is given in 
Fig.~\ref{sigcal} (by using the Definitions 6, 7, 14 and 15).
We obtain \vspace{0,2 cm}

$R^{(1)}_{reg}( \infty_2)(scan_2(\sigma_1))= -(2v-v^{-1}+vz^2).1$.\vspace{0,2 cm}

\begin{figure}
\centering
\includegraphics{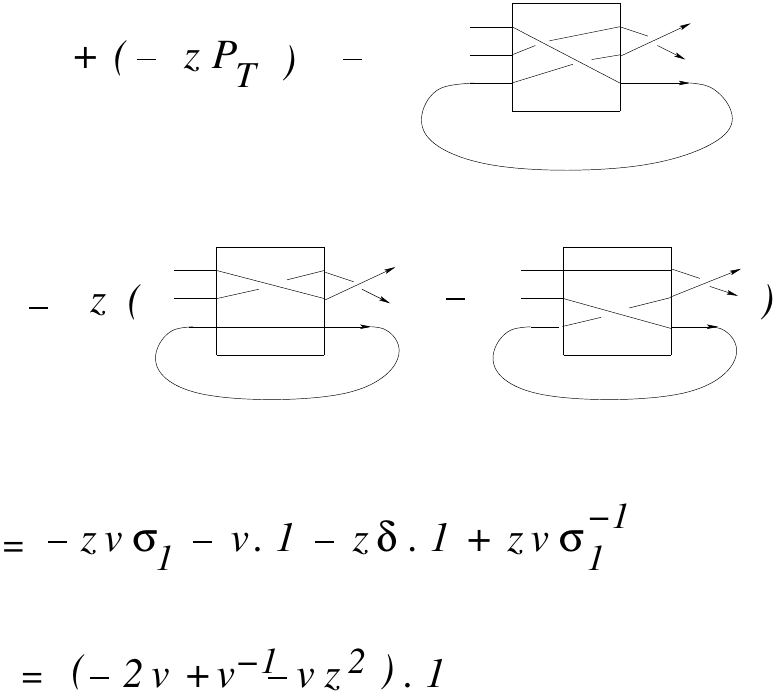}
\caption{\label{sigcal}  the calculation of $R^{(1)}_{reg} (scan_2)$ in Example 1}  
\end{figure}

Here $1$ and $\sigma_1$ are the generators of the skein module $S(\partial \sigma_1)$. It is amazing that the coefficient of the generator $1$ is just minus the HOMFLYPT polynomial of the right trefoil.\vspace{0,2 cm}

Let us consider the curl on the first branch. We show the corresponding regular isotopy in Fig.~\ref{sigco}. In the case $\infty_1$ we obtain $W(x_1)=W_2(x_1)=0$ and $W(x_2)=W_2(x_2)=0$. Consequently\vspace{0,2 cm}

\begin{figure}
\centering
\includegraphics{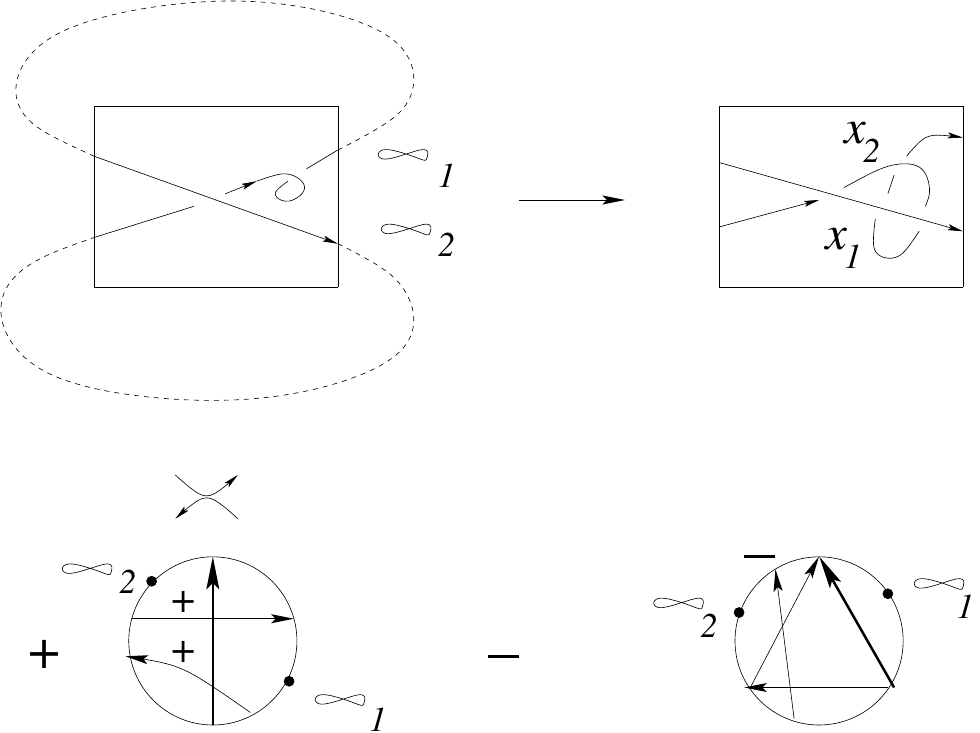}
\caption{\label{sigco}  the two Reidemeister moves in Example 1 for $scan_1$}  
\end{figure}

$R^{(1)}_{reg}(scan_1(\sigma_1))=R^{(1)}(scan_1(\sigma_1))=0$. \vspace{0,2 cm}  

In the case $\infty_2$ the move $x_1$ does not contribute because $d$ is of type 1. $\partial (hm)=  \infty_2$ for the move $x_2$ and $W_1(x_2)=0$. Consequently\vspace{0,2 cm}

$R^{(1)}(scan_1(\sigma_1))=0$.\vspace{0,2 cm}

The partial smoothing of the triple crossing $x_2$ is shown in Fig.~\ref{sigcocal}. It follows that\vspace{0,2 cm}

\begin{figure}
\centering
\includegraphics{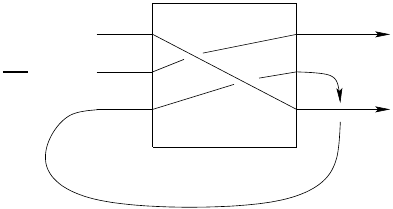}
\caption{\label{sigcocal} calculation of $R^{(1)}_{reg} (\infty_2)(scan_1)$}  
\end{figure}

$R^{(1)}_{reg}(A= \infty_2)(scan_1(\sigma_1))=-\delta \sigma_1$.\vspace{0,2 cm}

Let us consider now $T=\sigma_1^2$. We will consider only $scan_2(T)$. There is again a unique closure to a circle. We show the regular isotopy $scan(\sigma_1^2)$ and the corresponding Gauss diagrams in Fig.~\ref{sig2} where we give names to some crossings too.\vspace{0,2 cm}

\begin{figure}
\centering
\includegraphics{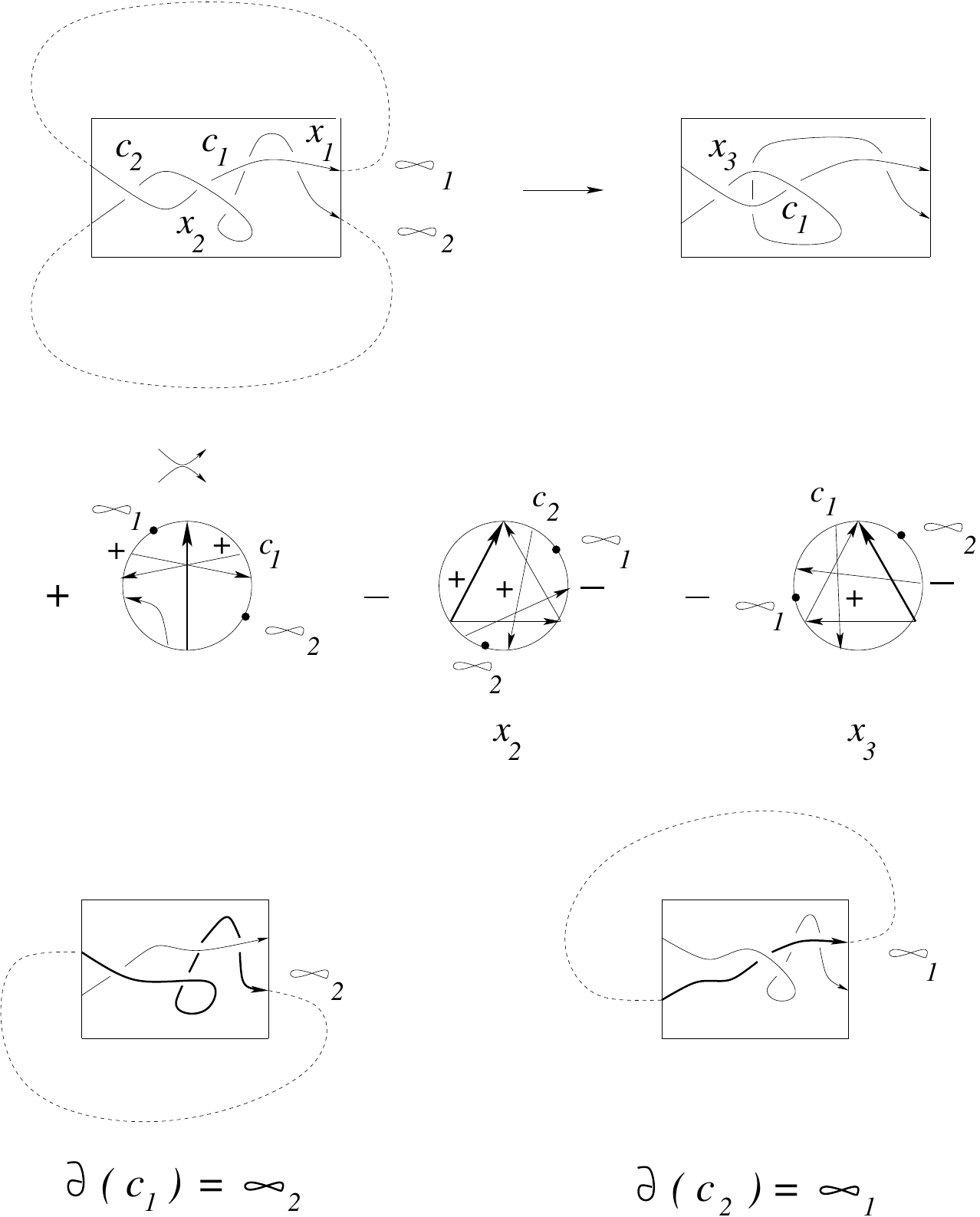}
\caption{\label{sig2} the Reidemeister moves for $scan(\sigma_1^2)$}  
\end{figure}

{\em The case $\infty_1$:}

The crossing $c_2$ is the only f-crossing for $x_2$ and the crossing $hm$ is the only r-crossing for $c_2$. We have $\partial c_2= \infty_1$, $W(x_2)=W_2(x_2)=1$. The Reidemeister move $x_3$ is also a positive triple crossing with $sign(x_3)=-1$. We have  $\partial (hm)= \partial (c_2)= \infty_1$ and $W_1(x_3)=1$ (because of the crossing $c_1$). It follows that
\vspace{0,2 cm} 

$R^{(1)}_{reg}(A= \partial \sigma_1^2)(scan(\sigma_1^2))=0$ \vspace{0,2 cm}

and that\vspace{0,2 cm}
 
$R^{(1)}_{reg}(A=  \infty_1)(scan(\sigma_1^2))=R^{(1)}(scan(\sigma_1^2))$\vspace{0,2 cm}

 in this case. The calculation of $R^{(1)}(scan(\sigma_1^2))$ is given in 
Fig.~\ref{sig2calc}. The result is\vspace{0,2 cm}

\begin{figure}
\centering
\includegraphics{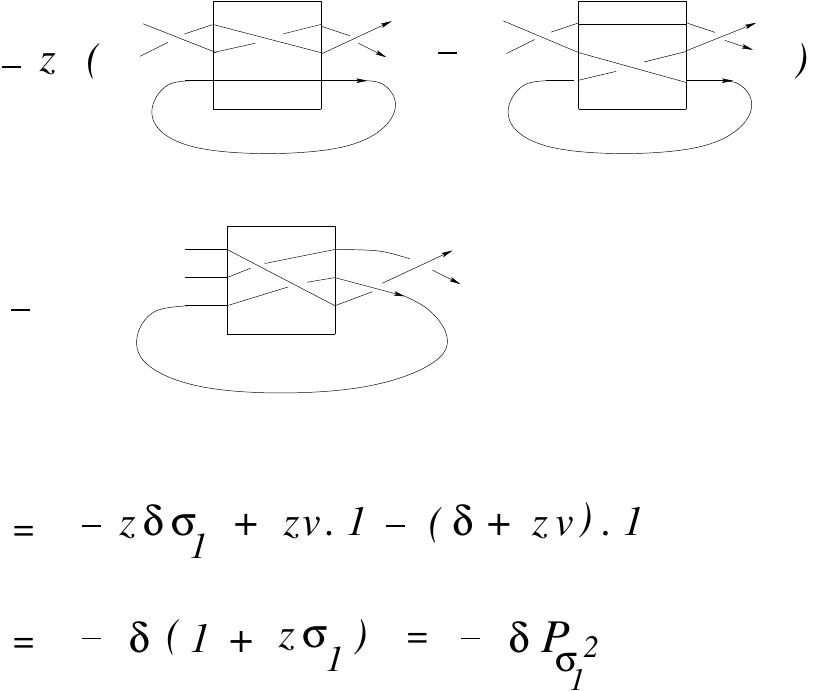}
\caption{\label{sig2calc}  calculation of $R^{(1)}(scan(\sigma_1^2))$ for $\infty_1$}  
\end{figure}

$R^{(1)}_{reg}(\infty_1)(scan(\sigma_1^2))=R^{(1)}(scan(\sigma_1^2))=-\delta (1+z.\sigma_1)= -\delta P_{\sigma_1^2}$.\vspace{0,2 cm}

Here as usual $\delta=(v-v^{-1})/z$ and $P_{\sigma_1^2}$ denotes the HOMFLYPT invariant of $\sigma_1^2$ in the skein module $S(\partial \sigma_1)$.\vspace{0,2 cm}

{\em The case $\infty_2$:}

The crossing $c_1$ is the only f-crossing for $x_1$ and there are no r-crossings. Consequently, $W(x_1)=1$ and $W_2(x_1)=0$. The crossing $hm$ is the only f-crossing for $x_2$. We have $\partial (hm)= \partial (c_1)=\infty_2$, $W(x_2)=1$ and $W_1(x_2)=W_2(x_2)=0$ (because there are two r-crossings, a negative one and the crossing $d$ which is positive). There is no r-crossing  in $x_3$.

Consequently,\vspace{0,2 cm}

$R^{(1)}_{reg}(A= \partial \sigma_1^2)(scan(\sigma_1^2))=R^{(1)}(scan(\sigma_1^2))=0$\vspace{0,2 cm}

A similar calculation gives now\vspace{0,2 cm}

$R^{(1)}_{reg}(A=  \infty_2)(scan(\sigma_1^2))=(-2v+v^{-1}-vz^2).\sigma_1$\vspace{0,2 cm}

For $T=\sigma_1$ and $T=\sigma_1^2$ we have $R^{(1)}(scan(T))=0$ for $\infty_2$. Therefore we consider also the next case $T=\sigma_1^3$ for $\infty_2$. The Reidemeister moves in $scan(\sigma_1^3)$ are shown in 
Fig.~\ref{sig3}. We have $W_2(x_1)=1$, $W_1(x_2)=W_2(x_2)=0$, $W_2(x_3)=1$ ($c_3$ is a f-crossing which contributes with the r-crossing $hm$), $W_1(x_4)=W_2(x_4)=1$ ($c_1$ is a r-crossing for $hm$, $c_2$ is a f-crossing but the corresponding r-crossings sum up to 0). The calculation of $R^{(1)}(scan(\sigma_1^3))$ is given in Fig.~\ref{calcsig3}.\vspace{0,2 cm}

$R^{(1)}(scan(\sigma_1^3))=(-3v+2v^{-1}-4vz^2+v^{-1}z^2-vz^4).1 $

\begin{figure}
\centering
\includegraphics{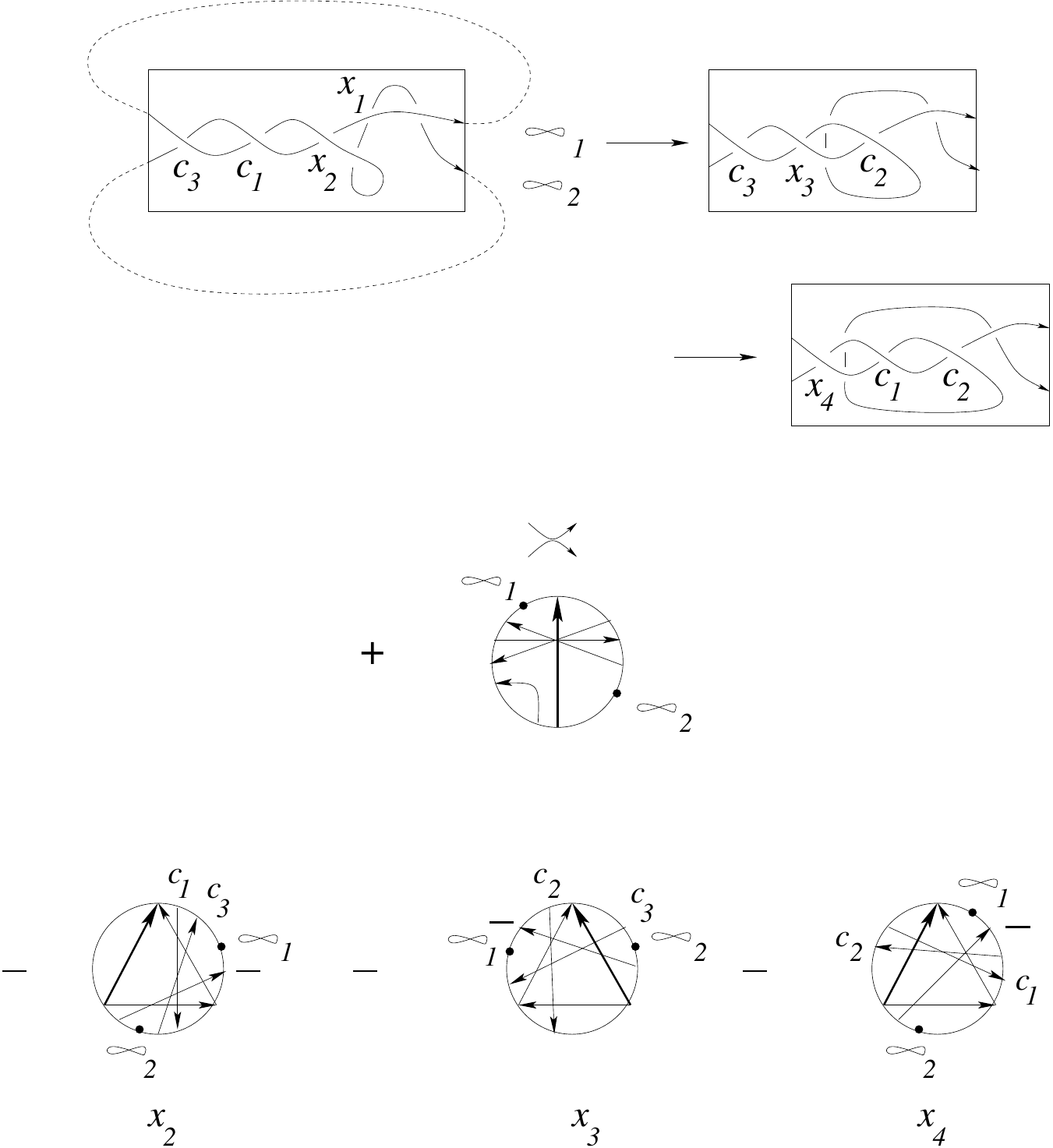}
\caption{\label{sig3} $scan(\sigma_1^3)$}  
\end{figure}

\begin{figure}
\centering
\includegraphics{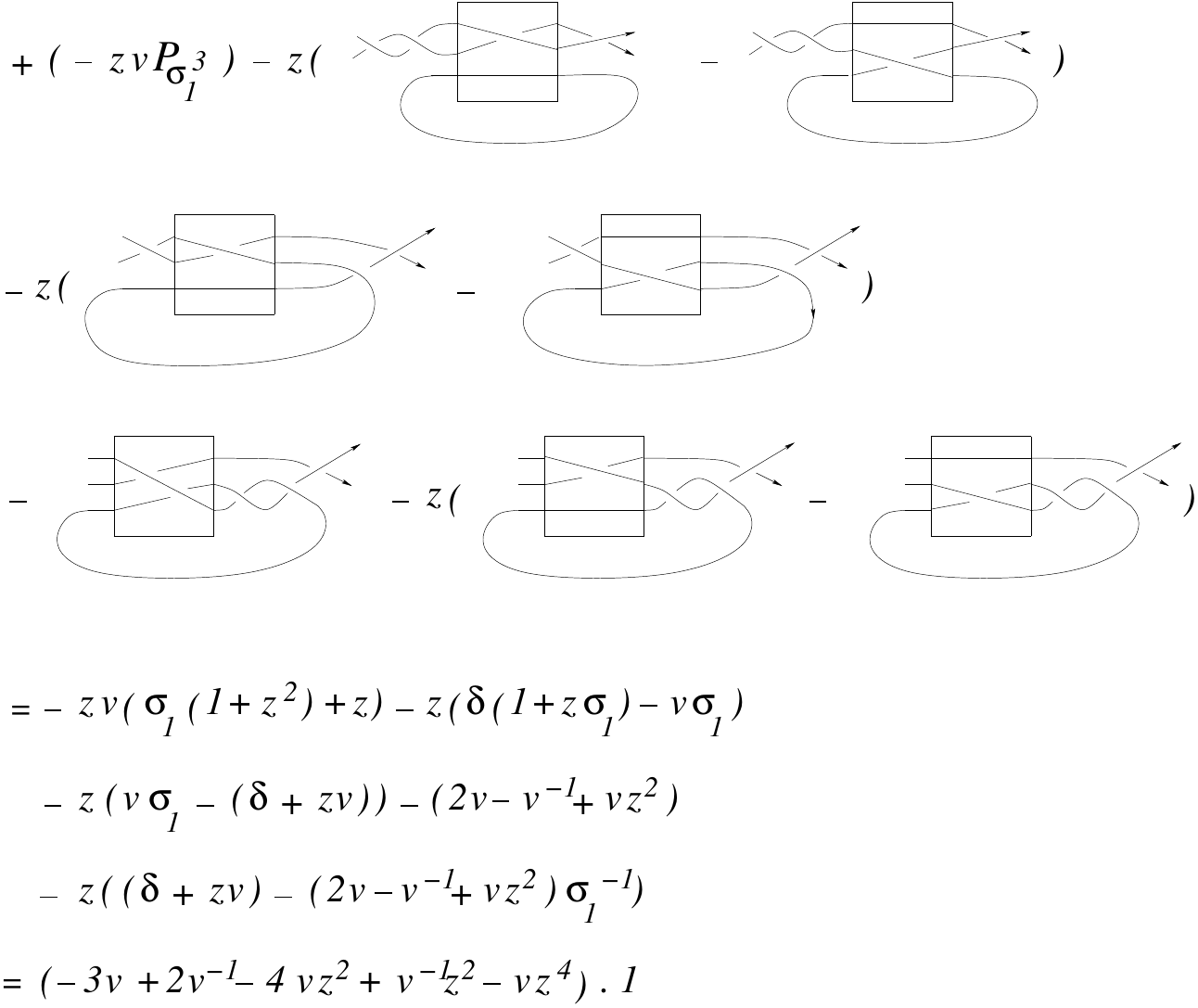}
\caption{\label{calcsig3}  calculation of $R^{(1)}(scan(\sigma_1^3)$ for $\infty_2$}  
\end{figure}

\end{example}

These simple examples show that $R^{(1)}_{reg}(A)$ depends on the grading $A$, that both $R^{(1)}_{reg}(A)$ and $R^{(1)}$ depend on the choice of the point at infinity and of the choice of the curl. Moreover, we see that the two 1-cocycles are independent and that both are {\em not} multiplicative, e.g.
$R^{(1)}(scan(\sigma_1^2)) \not= (R^{(1)}(scan(\sigma_1)))^2$ for $\infty_1$ and $R^{(1)}(scan(\sigma_1^3)) \not= (R^{(1)}(scan(\sigma_1)))^3$ for $\infty_2$. Moreover, we see that $R^{(1)}(scan(\sigma_1^3))$ is {\em not} a multiple of $\delta P_{\sigma_1^3}$.

\begin{figure}
\centering
\includegraphics{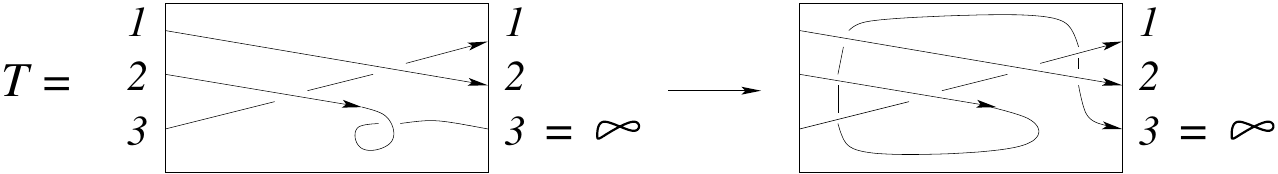}
\caption{\label{3tangle} the scan-arc for $T = \sigma_1\sigma_2$}  
\end{figure}

\begin{figure}
\centering
\includegraphics{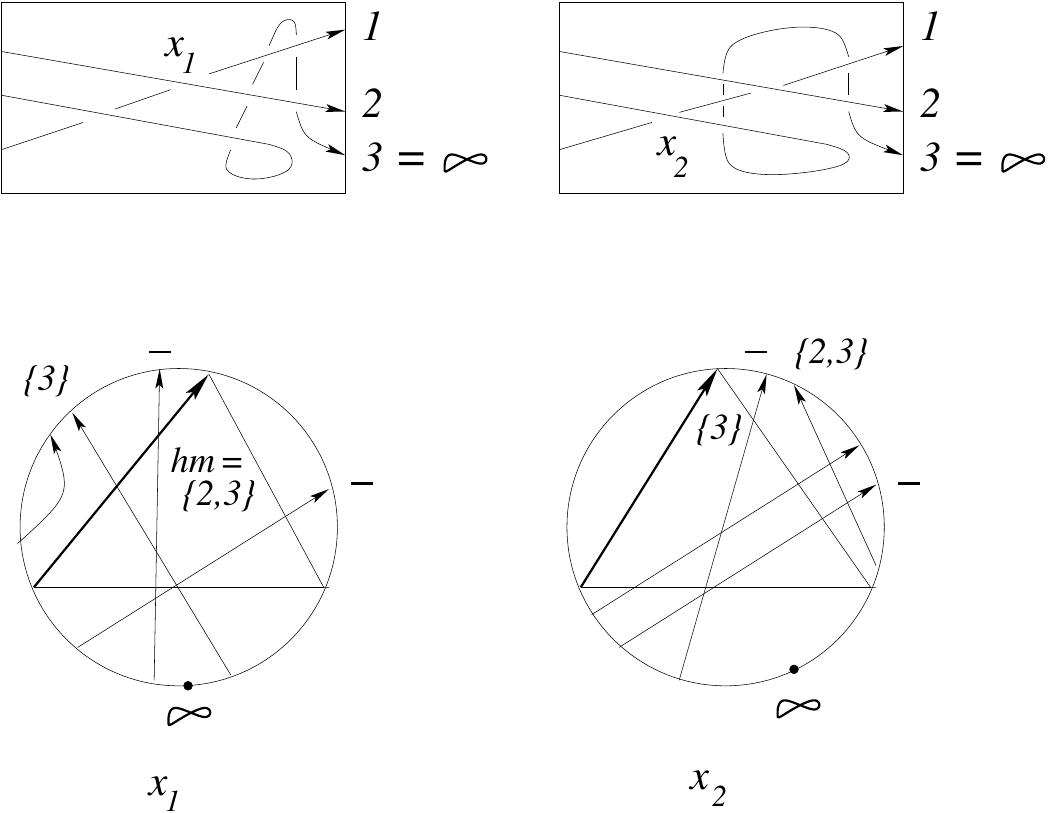}
\caption{\label{gauss3tangle} the Gauss diagrams in the scan-arc}  
\end{figure}

\begin{figure}
\centering
\includegraphics{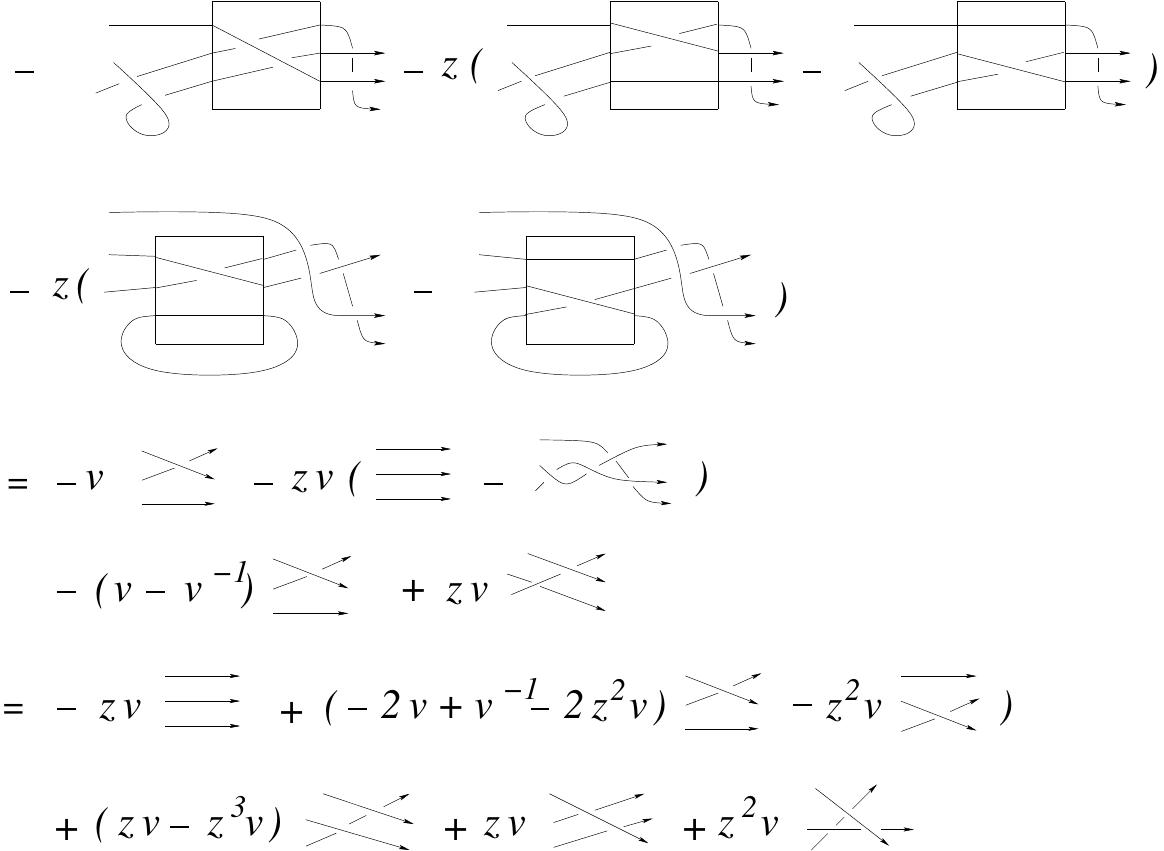}
\caption{\label{calc23} calculation of $R^{(1)}_{reg} (\{ 2, 3 \})(scanT)$}  
\end{figure}

\begin{example}
We consider the very simple 3-tangle $T=\sigma_1 \sigma_2$ with the standard braid closure  and the scan-arc shown in 
Fig.~\ref{3tangle}. Let us chose the point $3=\infty$. We calculate the version of $R^{(1)}_{reg}(A)$ and $R^{(1)}$ without contributions from $II^-_0$. Notice that the crossing $\sigma_1$ has marking $\{2,3\}$, the crossing $\sigma_2$ has marking $\{3\}$ and these are the only crossings with these markings. There are only two Reidemeister moves. They are positive Reidemeister III moves (i.e. type 1) of negative sign. Their Gauss diagrams are shown in Fig.~\ref{gauss3tangle}.
One easily sees that $R^{(1)}(scan(T))=0$, because all weights $W_1=W_2=0$. Moreover, $R^{(1)}_{reg}(\{2,3\})(scan(T))=R^{(1)}_{reg}(\{3\})(scan(T))$ because the partial smoothings $T_{r_a}$ for $x_1$ and $x_2$ lead to identical diagrams and the partial smoothings $T(type)$ for $x_1$ and $x_2$ enter into both invariants.

The calculation of $R^{(1)}_{reg}(\{2,3\})(scan(T))$ is given in Fig.~\ref{calc23}.

We see that even for this extremely simple tangle $T$ (just a basis element of the skein module of 3-braids) our new invariants are already rather complex, i.e. they contain lots of monomials and in particular all six basis elements of the skein module of 3-braids occur in the value of the invariant.
\end{example}

\section{The 1-cocycle $R^{(1)}$ in the complement of cusps with a transverse branch}

We know already that $R^{(1)}$ is a 1-cocycle in $M_T^{reg}$. First of all in $M_T$ we normalize the HOMFLYPT invariants as usual by $v^{-w(T)}$, where $w(T)$ is the writhe of $T$. Moreover we fix an (abstract) closure $\sigma$ of $T$ and a point at infinity $\infty$ in $\partial T$ (in the case of long knots there is only a canonical choice).

There are  several local types of Reidemeister I moves: the new crossing could be positive or negative. There are also several global types of Reidemeister I moves: the new crossing could be of type $0$ or $1$. 

\begin{definition}
The partial smoothing $T_I(p)$ of a Reidemeister I move $p$ of global type $0$ is defined by \vspace{0,2 cm} 

$T_I(p)=-1/2(v-v^{-1})\delta v_2(T)P_T$. \vspace{0,2 cm}

Here $\delta=(v-v^{-1})/z$ is the HOMFLYPT polynomial of the trivial link of two components and $P_T$ denotes the (normalized) HOMFLYPT invariant of $T$ in the skein module $S(\partial T)$. The integer $v_2(T)$ for a given (abstract) closure $\sigma$ of $T$ and a given choice of $\infty$ in $\partial T$ is defined by the first Polyak-Viro formula in Fig.~\ref{PV} with $\infty$ as the marked point.

The partial smoothing $T_I(p)$ of a Reidemeister I move $p$ of global type $1$ is $0$.

We add now $sign(p) T_I(p)$ to $R^{(1)}$ for each Reidemeister I move $p$.
\end{definition}

If $T$ is a long knot then $v_2(T)$ is just the Vassiliev invariant of degree 2 of $T$ (which is e.g. equal to the coefficient of $z^2$ in the Conway polynomial of $T$, see e.g. \cite{BN}).

So, Reidemeister I moves of global type $1$ do not contribute to the 1-cocycle $R^{(1)}$. Notice that the contribution of a Reidemeister I move of global type $0$ does not depend on the local type neither on its place in the diagram $T$. Hence in particular it has the scan-property.

Let $s$ be a loop in $M_T$. It is not difficult to see that the number of Reidemeister I moves of global type $0$ is always even (in fact the algebraic number of Reidemeister I moves of global type $0$ which are of positive local type is always equal to the algebraic number of Reidemeister I moves of global type $0$ which are of negative local type, because they come always in pairs in Reidemeister II moves). 

\begin{lemma}
The value of the 1-cocycle $R^{(1)}$ on a meridian of 
$\Sigma^{(1)} \cap \Sigma^{(1)}_{cusp}$ is zero.
\end{lemma}
 {\em Proof.}  The new crossing from $\Sigma^{(1)}_{cusp}$ could be a f-crossing for the Reidemeister move in $\Sigma^{(1)}$. However, it is an isolated crossing and in particular it has no r-crossings at all. The changing for the HOMFLYPT invariants for $\Sigma^{(1)}$, which comes from the new crossing in $\Sigma^{(1)}_{cusp}$, is compensated by the normalization. The contribution of $\Sigma^{(1)}_{cusp}$ doesn't change at all under isotopy of the tangle outside the cusp.

$\Box$

\begin{figure}
\centering
\includegraphics{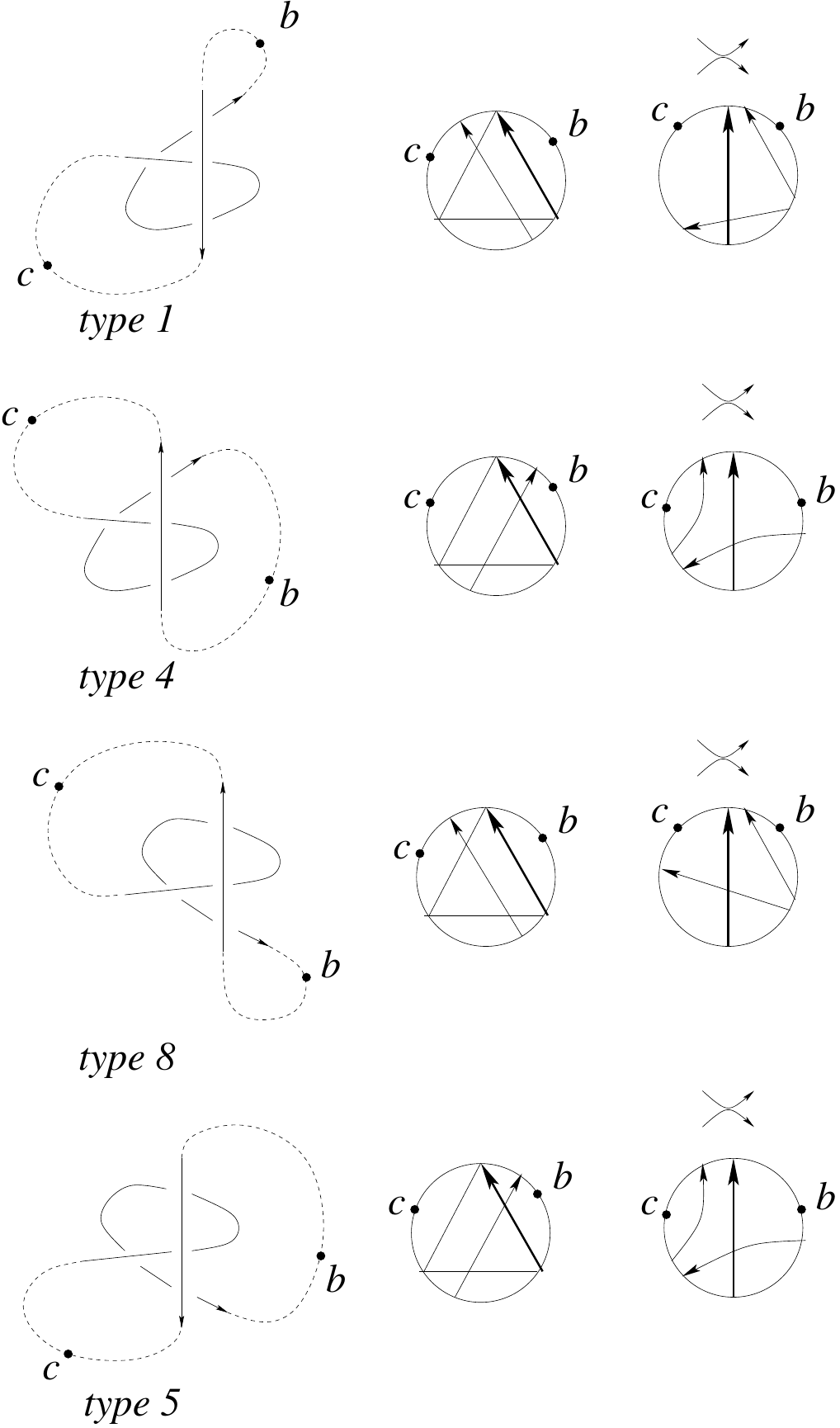}
\caption{\label{over-cusp0} the strata $\Sigma_{l_c}^{(2)}$ and 
$\Sigma_{l_b}^{(2)}$}  
\end{figure}

\begin{figure}
\centering
\includegraphics{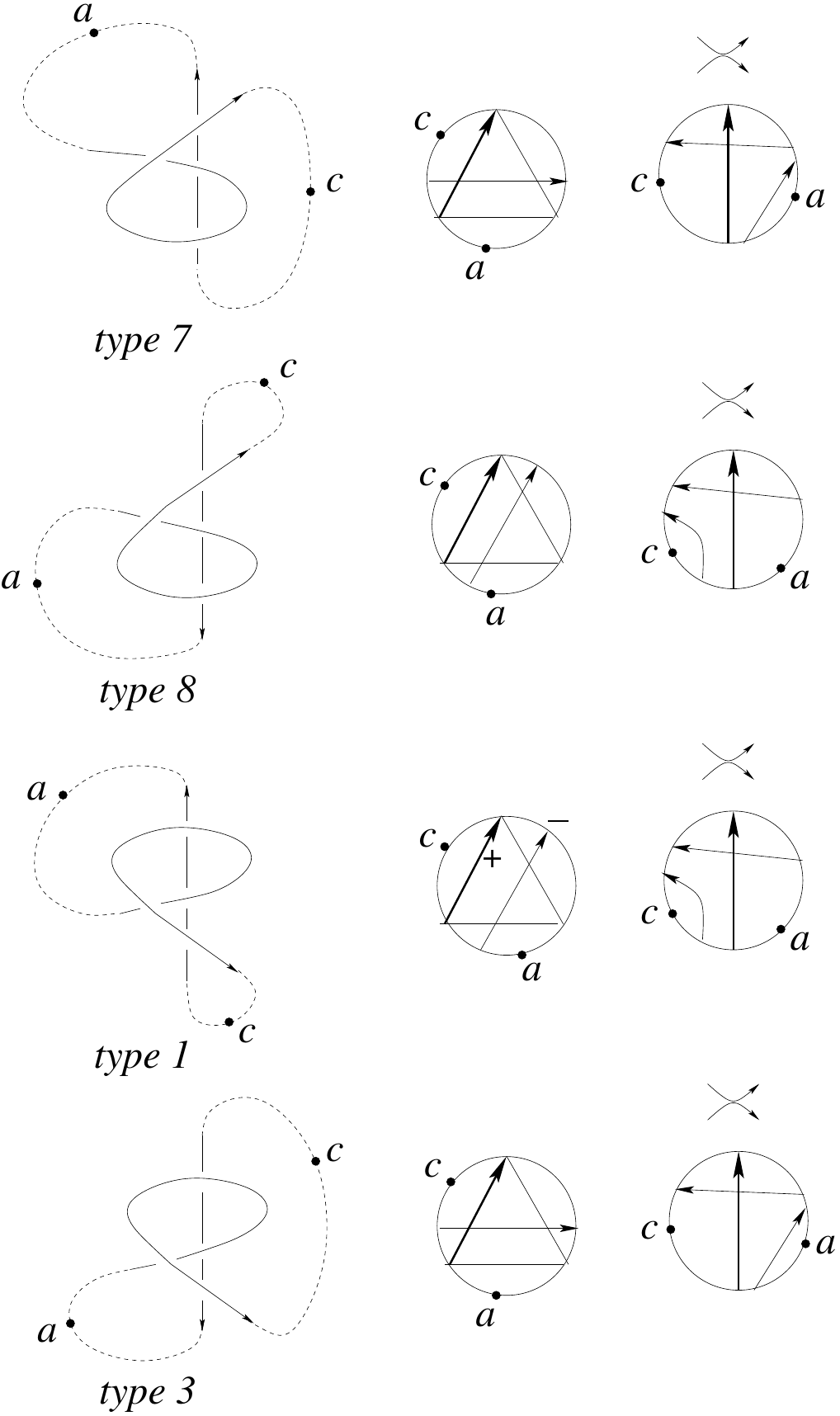}
\caption{\label{under-cusp1} the strata $\Sigma_{r_a}^{(2)}$ and    
$\Sigma_{r_c}^{(2)}$}  
\end{figure}

We have to deal now with the irreducible strata of codimension two which contain a diagram with a cusp. There are exactly sixteen different local types of strata corresponding to a cusp with a transverse branch. We list them in 
Fig.~\ref{over-cusp0}...Fig.~\ref{under-cusp1}, where we move the branch from the right to the left. For each local type we have exactly two global types. We give in the figure also the Gauss diagrams of the triple crossing and of the self-tangency with equal tangent direction. Notice that in each Gauss diagram of a triple crossing one of the three arcs is empty besides just one head or foot of an arrow. Let us denote each stratum of $\Sigma^{(2)}_{trans-cusp}$ simply by the global type of the corresponding triple crossing .

\begin{figure}
\centering
\includegraphics{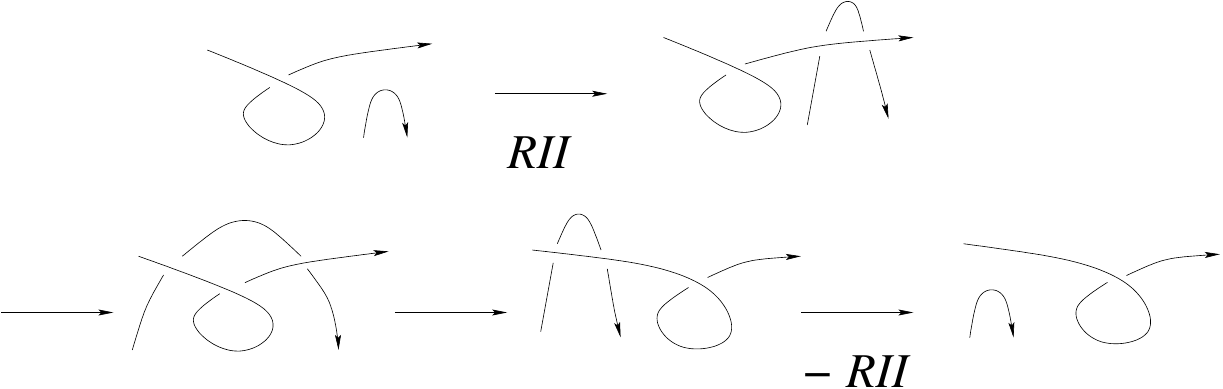}
\caption{\label{orbranch} two Reidemeister $II$ moves with cancelling 
contributions}  
\end{figure}

\begin{lemma}
The value of the 1-cocycle $R^{(1)}$ on a meridian is zero for all strata in $\Sigma^{(2)}_{trans-cusp}$
besides those in $\Sigma^{(2)}_{l_c}$. 

\end{lemma}
{\em Proof.}
It is clear that we have to consider only one orientation of the moving branch because the contributions of the two Reidemeister II moves in the 
Fig.~\ref{orbranch} always cancel out, no matter the orientations, the types and whether the branches move over or under the cusp.  We show first that $R^{(1)}$ vanishes on all Whitney tricks. 
It follows then from Fig.~\ref{over-cusp0}...Fig.~\ref{under-cusp1} that it suffices to prove the assertion in four cases:
$l_b$ of type 1 and the branch moves over the cusp, $l_b$  of type 1 and the branch moves under the cusp, $r_b$ of type 1 and $r_a$ of type 1.
Notice that the Reidemeister I move in the meridian does never contribute because it enters twice but with opposite signs.

\begin{figure}
\centering
\includegraphics{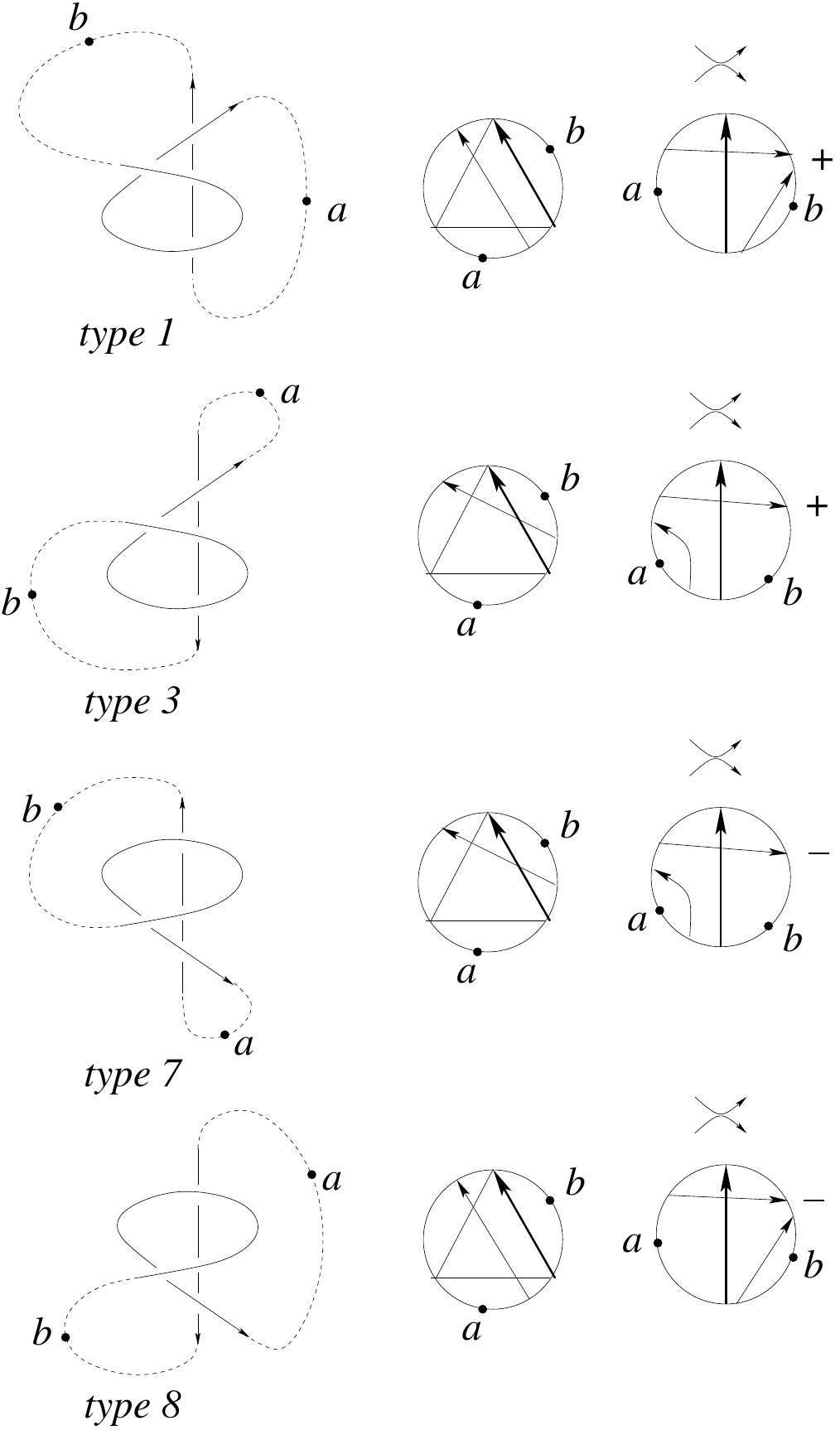}
\caption{\label{under-cusp0} the strata $\Sigma_{l_a}^{(2)}$ and  
$\Sigma_{l_b}^{(2)}$}  
\end{figure}

\begin{figure}
\centering
\includegraphics{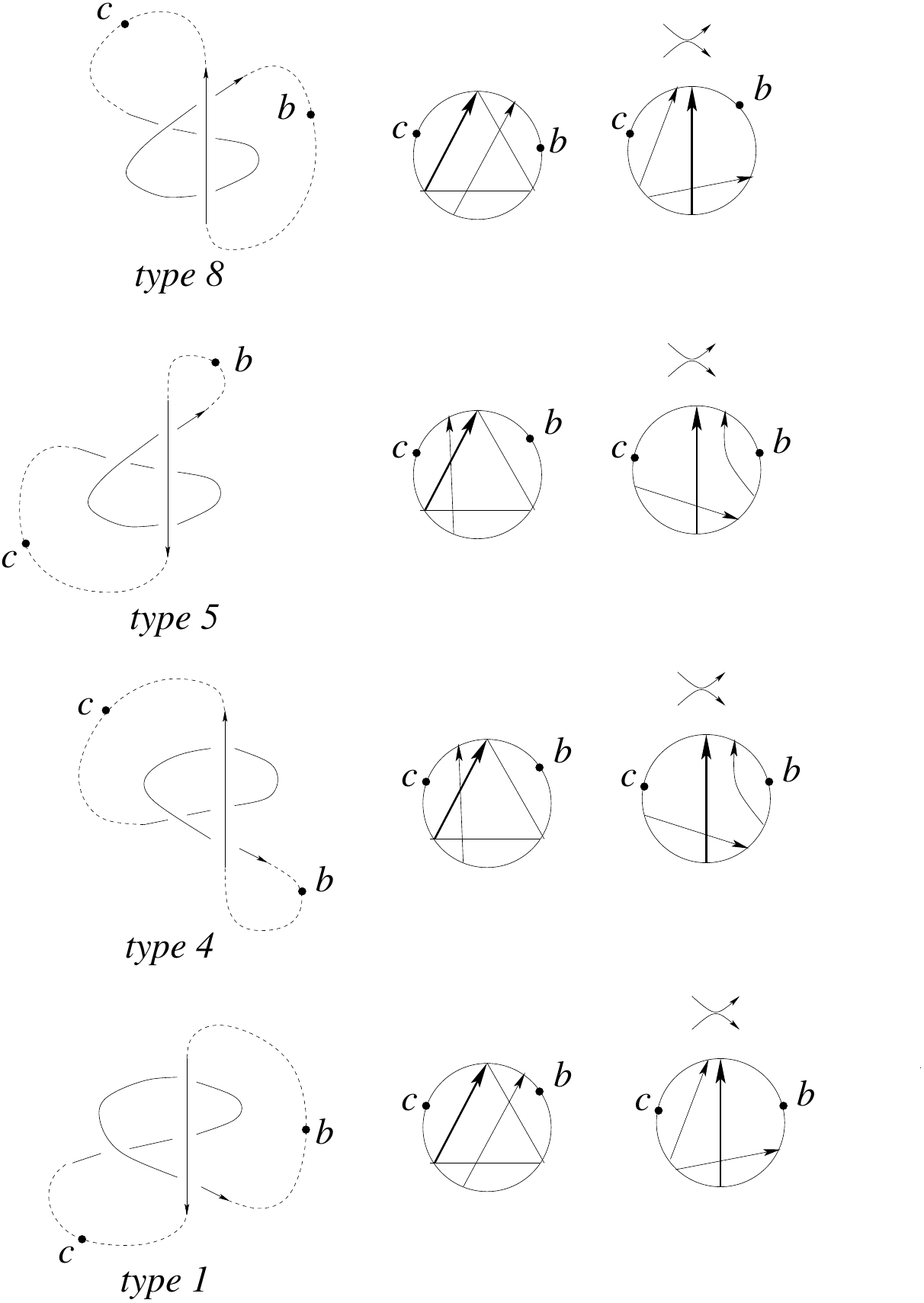}
\caption{\label{over-cusp1} the strata $\Sigma_{r_b}^{(2)}$ and  
$\Sigma_{r_c}^{(2)}$}  
\end{figure}

\begin{figure}
\centering
\includegraphics{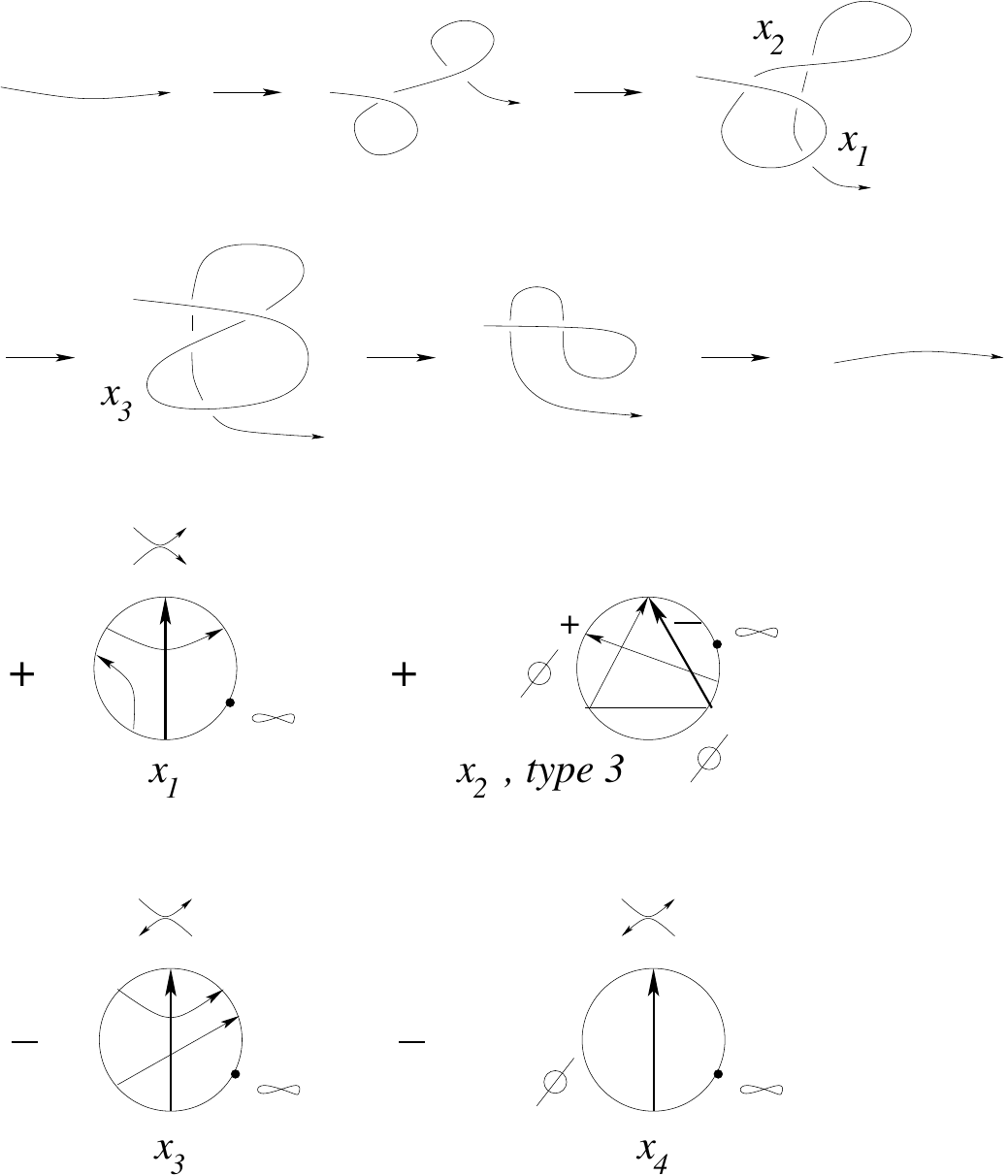}
\caption{\label{firstW} the first Whitney trick}  
\end{figure}

\begin{figure}
\centering
\includegraphics{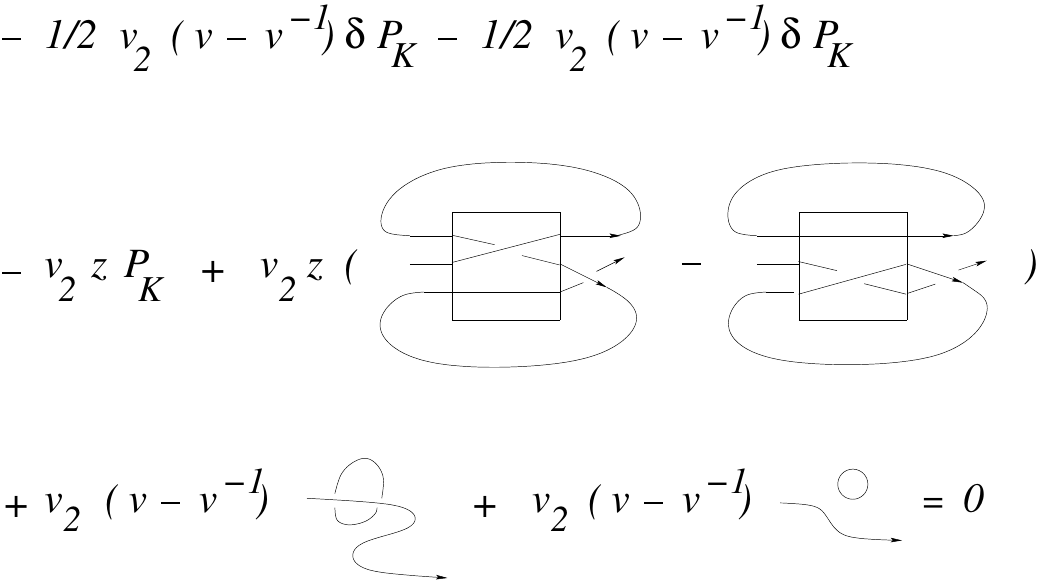}
\caption{\label{calfirstW} calculation of $R^{(1)}$ for the first Whitney trick}  
\end{figure}

\begin{figure}
\centering
\includegraphics{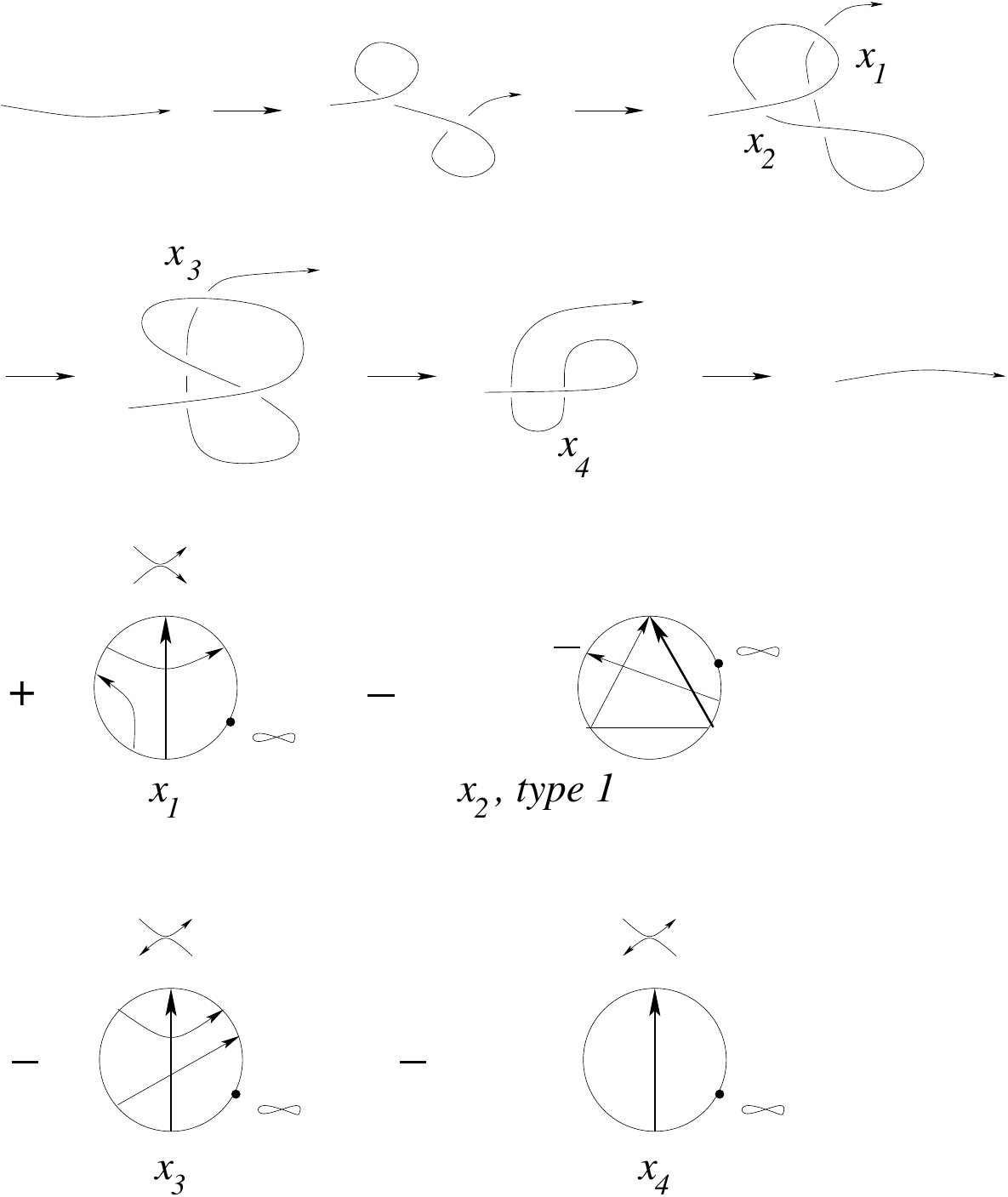}
\caption{\label{secondW} the second Whitney trick }  
\end{figure}

\begin{figure}
\centering
\includegraphics{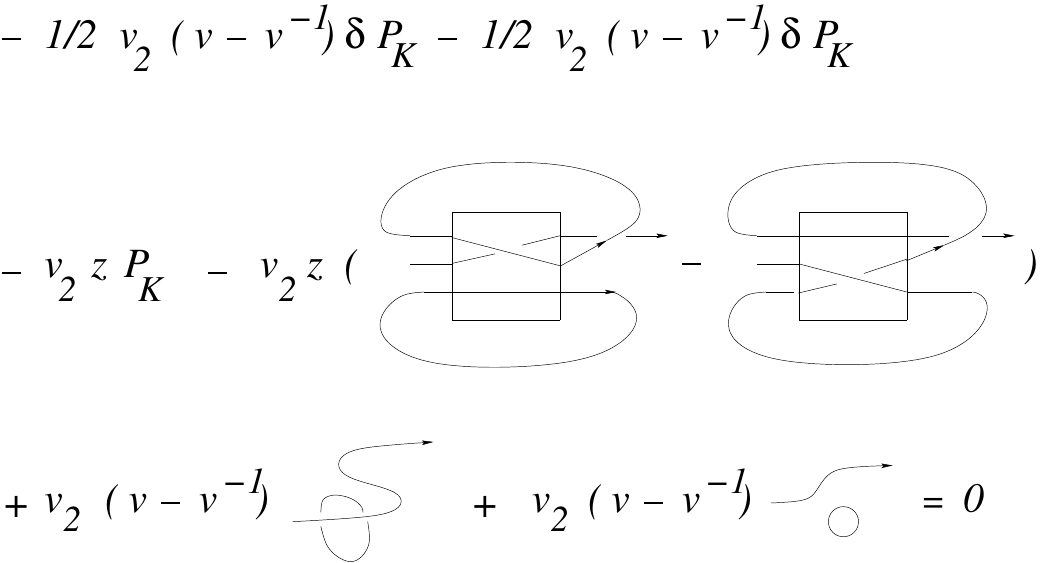}
\caption{\label{calsecondW} calculation of $R^{(1)}$ for the second Whitney trick }  
\end{figure}

\begin{figure}
\centering
\includegraphics{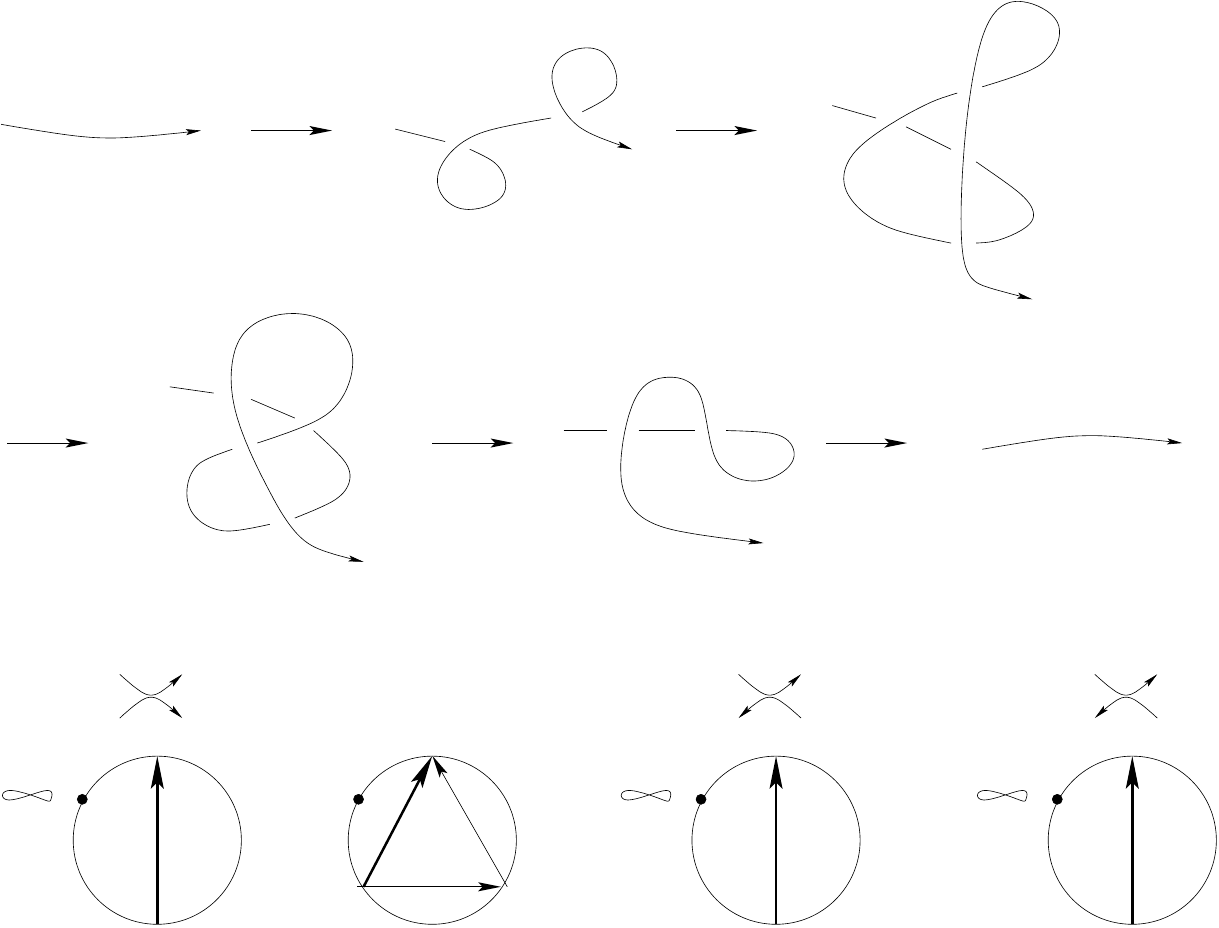}
\caption{\label{thirdW} the third Whitney trick}  
\end{figure}

\begin{figure}
\centering
\includegraphics{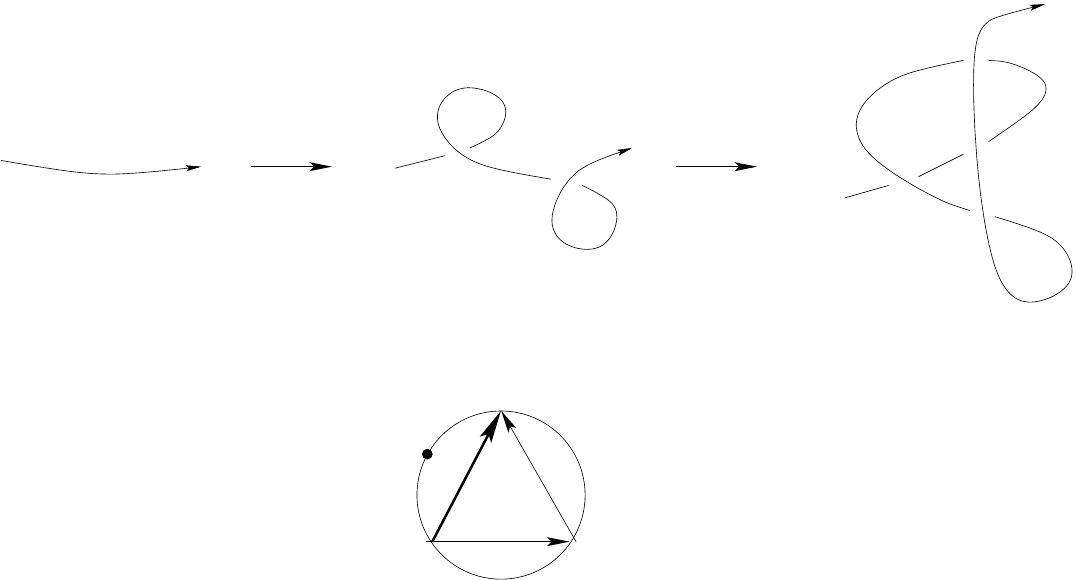}
\caption{\label{fourthW} the fourth Whitney trick}  
\end{figure}

The first Whitney trick is shown in Fig.~\ref{firstW} and the calculation in 
Fig.~\ref{calfirstW}. We write $K$ for $T$ and one easily sees that all weights $W_2$ are equal to $v_2$. The second Whitney trick is rather similar to the first, see 
Fig.~\ref{secondW} and Fig.~\ref{calsecondW}. In the third and fourth Whitney trick all distinguished crossings $d$ are of type 1. It suffices to establish that the triple crossing is not of type $l_c$. We do this in the Fig.~\ref{thirdW} and 
Fig.~\ref{fourthW}. Consequently, there are no contributions at all for these Whitney tricks.

\begin{figure}
\centering
\includegraphics{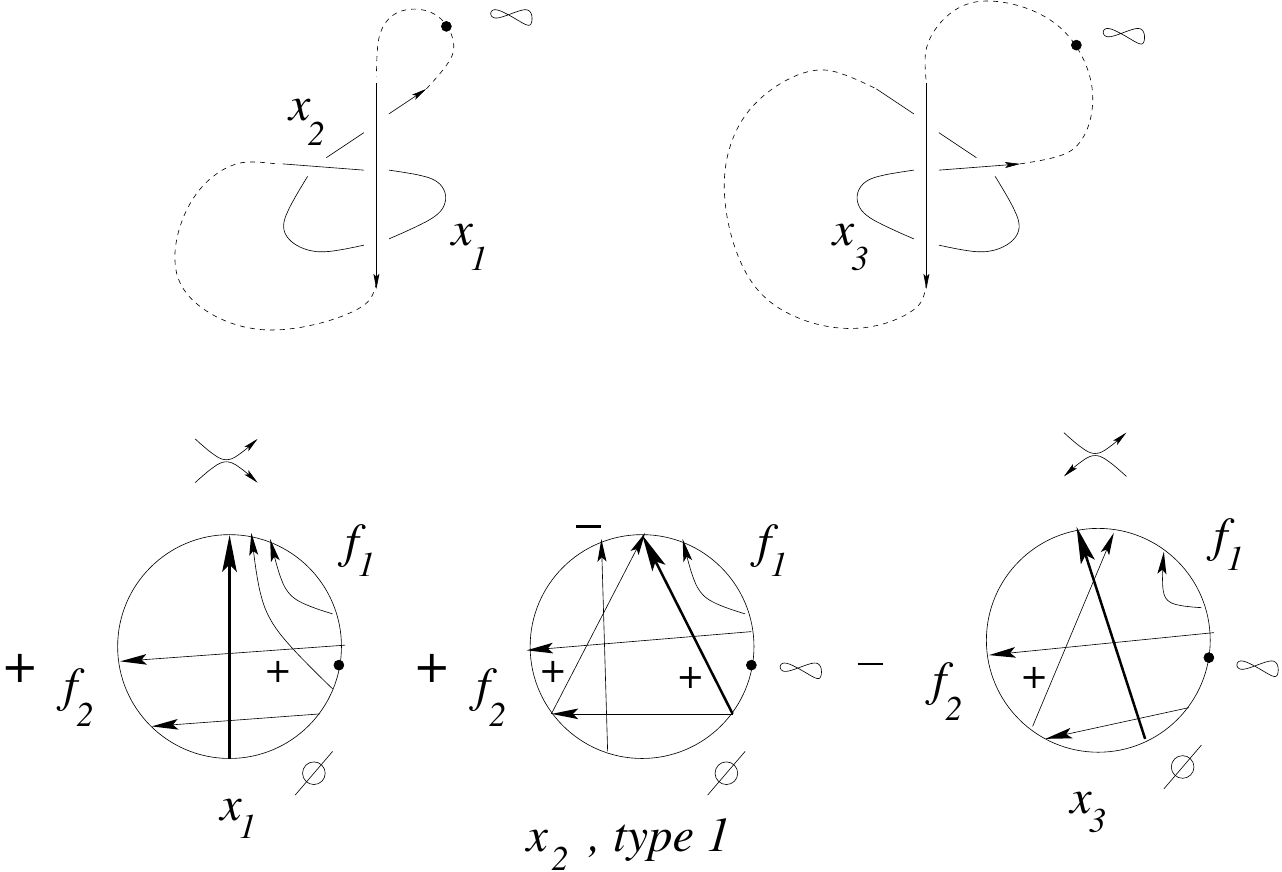}
\caption{\label{lbover}  the meridian for $\Sigma^{(2)}_{l_b}$ of local type $1$ with over branch}  
\end{figure}

\begin{figure}
\centering
\includegraphics{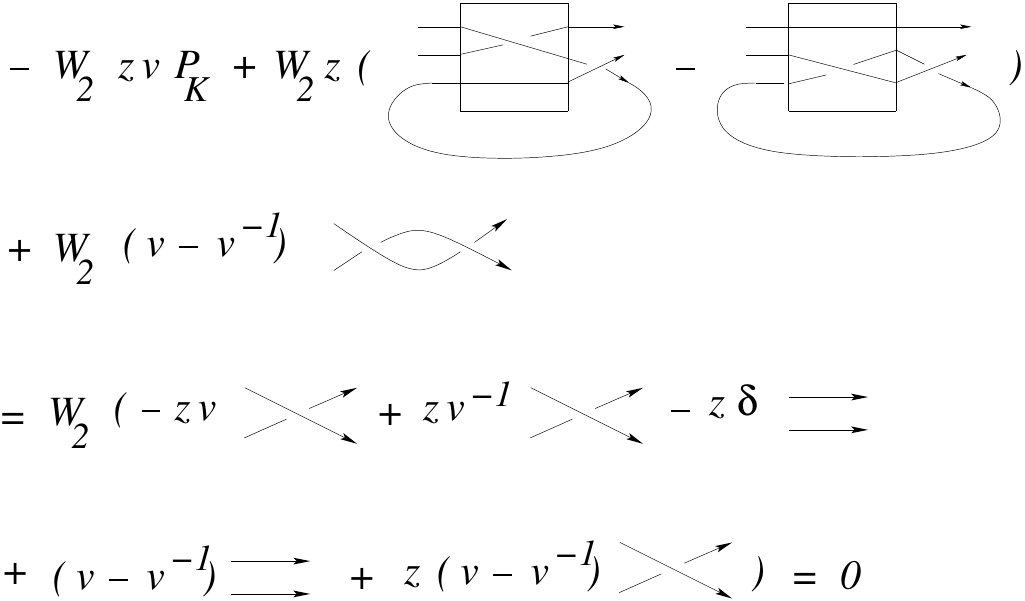}
\caption{\label{calbover} calculation for $\Sigma^{(2)}_{l_b}$ of local type $1$ with
over branch}  
\end{figure}

Let us consider now the above four cases. We study the Gauss diagrams in detail for the case $l_b$ with the branch moving over the cusp. The Reidemeister moves with the Gauss diagrams are given in Fig.~\ref{lbover}.  Notice that the self-tangency with opposite tangent direction and the triple crossing share the same distinguished crossing $d$. The distinguished crossing of the other self-tangency has almost the same head as $d$ but not the same foot. However their foots are connected by an empty arc in the circle. There are exactly two sorts of f-crossings, which we call $f_1$ and $f_2$ (compare Fig.~\ref{lbover}). In all three Gauss diagrams each corresponding $f_1$-crossing shares the same r-crossing. In $x_1$ each $f_2$-crossing has already exactly one positive r-crossing in the drawn part of the Gauss diagram. There are three r-crossings in $x_2$ but their sum is $+1$. In $x_3$ there is again a single positive r-crossing (but different from that in $x_1$). All other r-crossings (not drawn in the diagrams) are the same for each corresponding $f_2$-crossing. Consequently, we have shown that the weights $W_2$ are the same for all three Reidemeister moves. We give the calculation of $R^{(1)}$ now in 
Fig.~\ref{calbover}.

\begin{figure}
\centering
\includegraphics{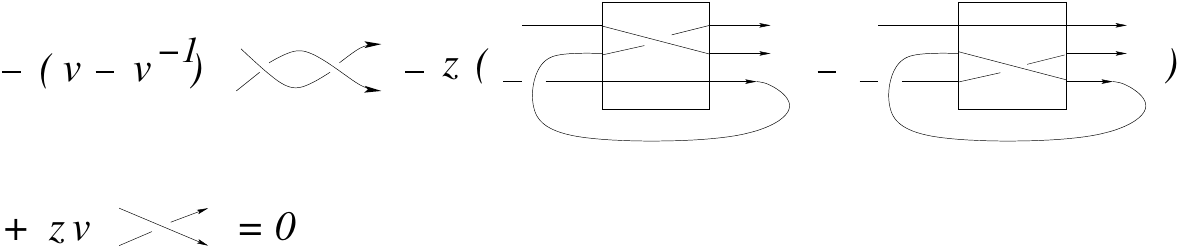}
\caption{\label{lbunder} calculation for $\Sigma^{(2)}_{l_b}$ of local type $1$ with
under branch}  
\end{figure}

\begin{figure}
\centering
\includegraphics{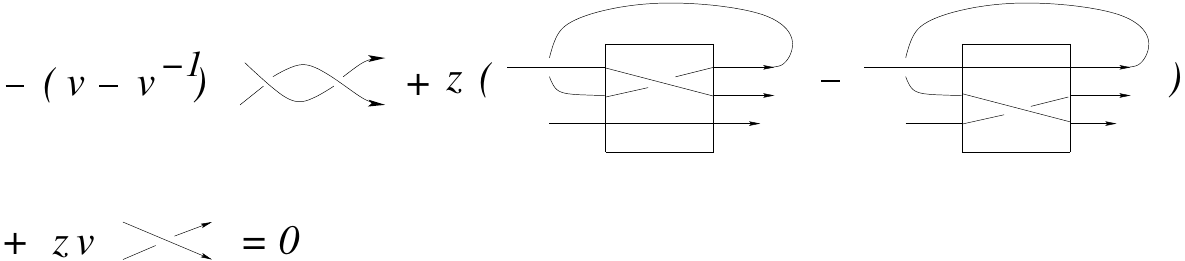}
\caption{\label{calrb} calculation for $\Sigma^{(2)}_{r_b}$ of local type $1$}  
\end{figure}

\begin{figure}
\centering
\includegraphics{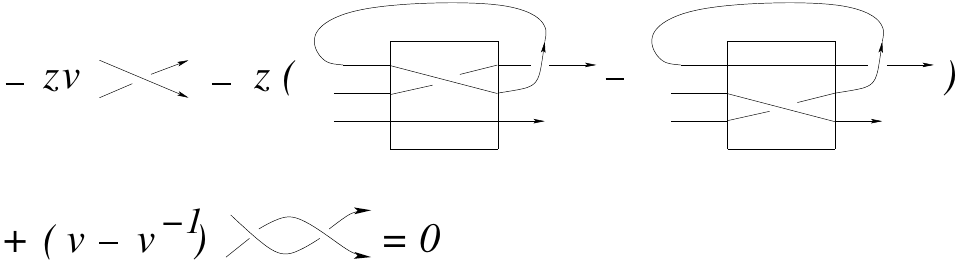}
\caption{\label{calra} calculation for $\Sigma^{(2)}_{r_a}$ of local type $1$}  
\end{figure}

In exactly the same way one shows that the weights $W_2$ of the three moves are the same in each of the remaining three cases. Therefore we give only the calculations of the sum of the partial smoothings in Fig.~\ref{lbunder}, 
Fig.~\ref{calrb} and Fig.~\ref{calra}. Notice that in the last case the triple crossing could also contribute with a partial smoothing associated to $W_1$. However, we see immediately from Fig.~\ref{under-cusp1} that $W_1=0$. Indeed, there are only exactly two crossings which contribute to $W_1$, because there is only exactly one head of an arrow in the arc in the circle from the middle crossing to the upper crossing. But the two crossings have opposite signs (which are indicated in the figure too).

$\Box$

\begin{lemma}
The value of the 1-cocycle  $R^{(1)}$  on a meridian of a degenerate cusp, locally given by $x^2=y^5$ and denoted by  $\Sigma^{(2)}_{cusp-deg}$, is zero.

\end{lemma}

\begin{figure}
\centering
\includegraphics{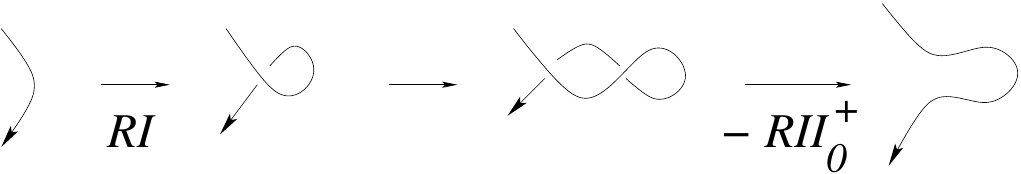}
\caption{\label{decusp0} meridian of a degenerated cusp of type $0$}  
\end{figure}

\begin{figure}
\centering
\includegraphics{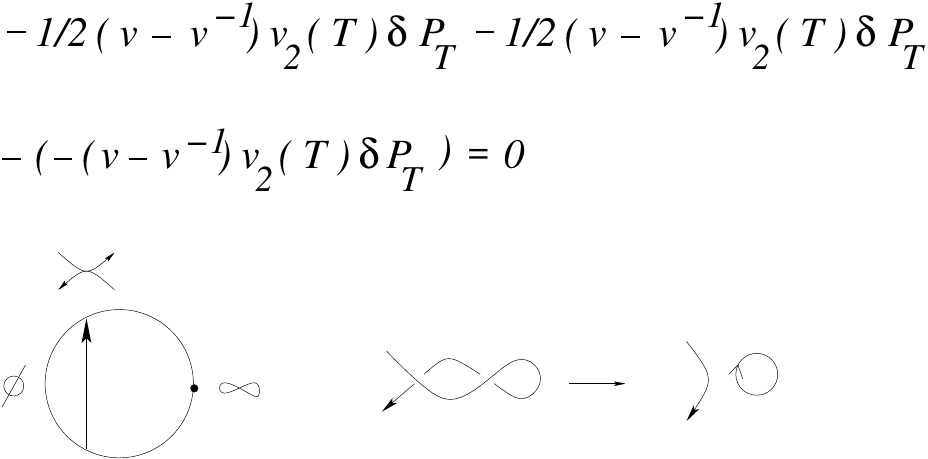}
\caption{\label{caldecusp0} calculation of $R^{(1)}$ on the meridian of a degenerated cusp of type $0$}  
\end{figure}

\begin{figure}
\centering
\includegraphics{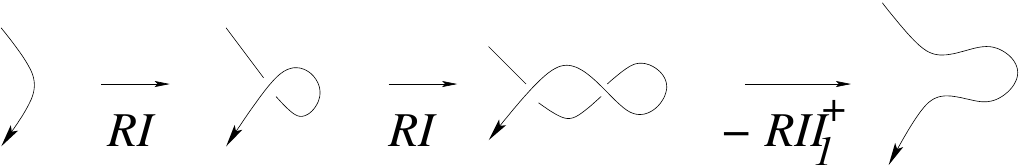}
\caption{\label{decusp1} meridian of a degenerated cusp of type $1$}  
\end{figure}

{\em Proof.} Again the cusp can be of type $0$ or of type $1$. We show the meridian of the degenerate cusp of type $0$ in Fig.~\ref{decusp0} and the calculation in Fig.~\ref{caldecusp0}. Notice that $W_2(p)=v_2(T)$ for the self-tangency $p$. The meridian of the degenerate cusp of type $1$ is shown in 
Fig.~\ref{decusp1}. In this case again no Reidemeister move at all contributes to $R^{(1)}$.

$\Box$

We have finally shown that $R^{(1)}(m)=0$ for the meridians $m$ of all six types of strata of codimension 2 (given in Section 2) besides for the four local types of strata from $\Sigma^{(2)}_{l_c}$. Theorem 2 implies that $R^{(1)}$ has the scan-property. Consequently we have proven the following proposition. 

\begin{proposition}
Let $s$ be a generic oriented arc in $M_T$ (with a fixed abstract closure $T \cup \sigma$ to an oriented circle and a fixed point at infinity in $\partial T$).  Then the 
1-cochain  \vspace{0,2 cm}

$ R^{(1)}(s)=\sum_{p \in s \cap l_c}sign(p)W_1(p)T_{l_c}(type)(p) $
 \\ \vspace{0,2 cm}

$+ \sum_{p \in s \cap r_a}sign(p)[W_1(type)(p)T_{r_a}(type)(p)  + zW_2(type)(p)T(type)(p)] $
\\ \vspace{0,2 cm}

$+ \sum_{p \in s \cap r_b}sign(p)zW_2(p)T(type)(p) +\sum_{p \in s \cap l_b}sign(p)zW_2(p)T(type)(p) $
\\ \vspace{0,2 cm}

$ + \sum_{p \in s \cap II^+_0}sign(p)W_2(p)T_{II}(p) +\sum_{p \in s \cap II^-_0}sign(p)W_2(p)zP_T$ 
\\ \vspace{0,2 cm}

$+\sum_{p \in s \cap I}sign(p)T_I(p)$ \vspace{0,2 cm}

is a 1-cocycle in  $M_T \setminus \Sigma^{(2)}_{l_c}$ and it has the scan-property for branches moving under the tangle.
 
\end{proposition}

\begin{figure}
\centering
\includegraphics{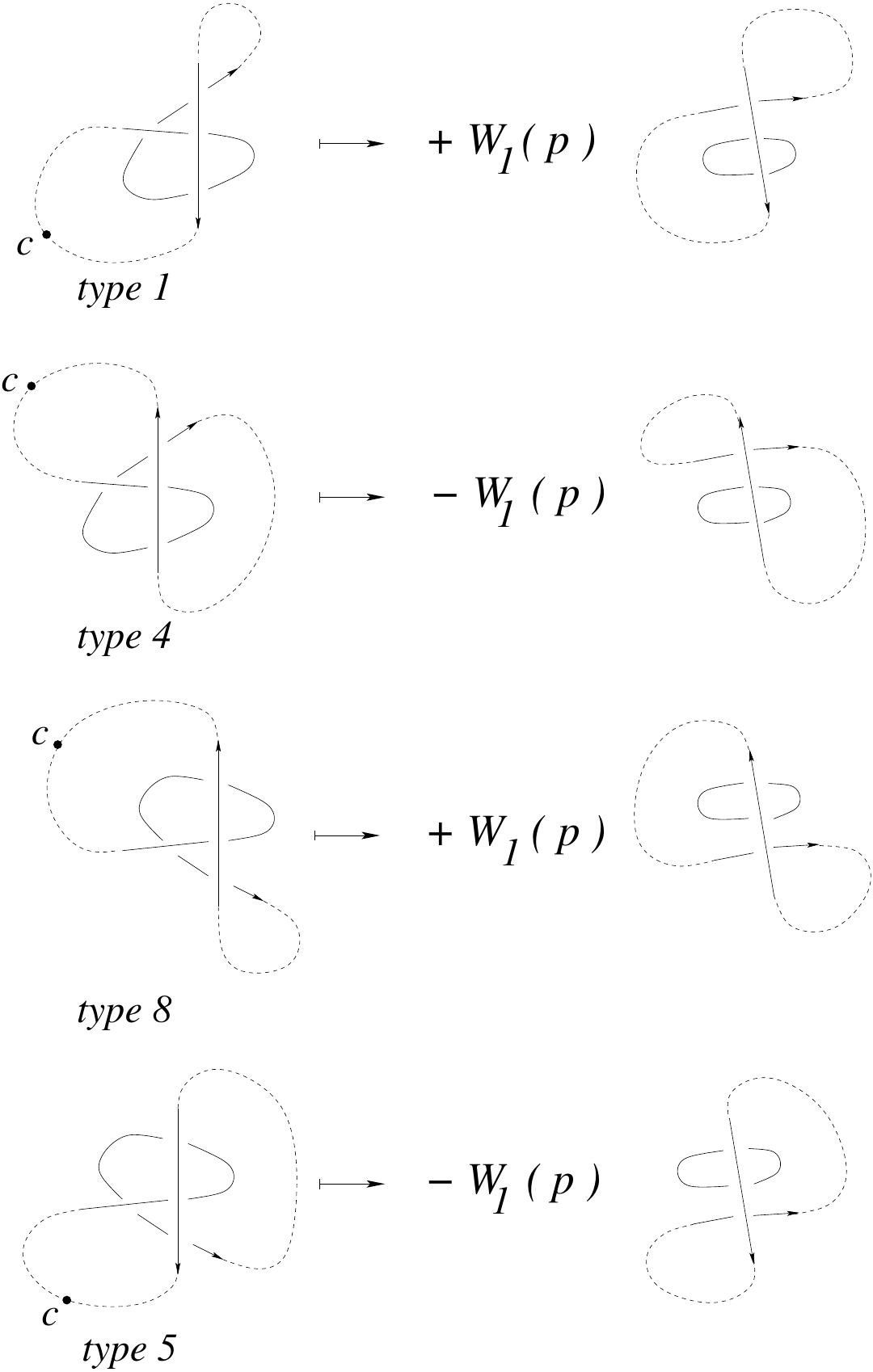}
\caption{\label{merover-cusp} the values of $R^{(1)}$ on the meridians of $\Sigma^{(2)}_{l_c}$}  
\end{figure}

\begin{lemma}
The value of $R^{(1)}$ on the meridian of $\Sigma^{(2)}_{l_c}$ is equal to $+W_1(p)\delta P_T$ or $-W_1(p)\delta P_T$ as shown in  Fig.~\ref{merover-cusp}.
\end{lemma}
{\em Proof.} We observe that the fourth involved crossing, which is not in the triangle, does not contribute to $W_1(p)$. The rest is a straightforward calculation. Notice that the local types 8 and 5 forces us to use the version of $R^{(1)}$ which associates 
$sign(p)W_2(p)zP_T$ to Reidemeister moves of type $II^-_0$ instead of the more complicated partial smoothings for the global types $r_a$ and $l_c$ in the case of the local types $2, 3, 5, 8$ and no contribution from $II^-_0$ at all (compare Remark 9).

$\Box$

\section{The finite type 1-cocycle $V^{(1)}$ in the complement of cusps with a transverse branch}

Let us summarize: the quantum part $R^{(1)}$ of the 1-cocycle uses the global types $r_a, r_b, l_b, l_c$ of Reidemeister III moves and the types $II_0^+$ and $II_0^-$ of Reidemeister II moves. Notice that in all cases the distinguished crossing $d$ is of type $0$ besides in the case $l_c$. The value of $R^{(1)}$ on the meridian of $\Sigma^{(2)}_{trans-cusp}$ is non zero exactly in the latter case. This happens because the distinguished crossing $d$ is of type $1$ for the Reidemeister II moves in the corresponding unfolding and it does not compensate the contribution of the triple crossing. This forces us to consider the moves of type $II^-_1$ too.

\begin{figure}
\centering
\includegraphics{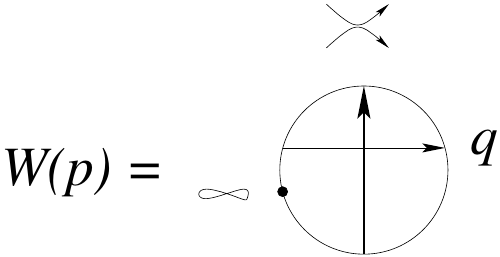}
\caption{\label{II1} the weight $W(p)$ for Reidemeister moves of type $II_1^-$}  
\end{figure}

\begin{definition}
The weight $W(p)$ of a Reidemeister II move of type $II^-_1$ is defined as the sum of the writhes $w(q)$ of all crossings $q$ of the configuration given in 
Fig.~\ref{II1}.
\end{definition}

Remember that we say simply that the weight is defined by the configuration in the figure (compare Section 3).

\begin{lemma}
The weight $W(p)$ of a Reidemeister II move of type $II^-_1$ satisfies the r-cube equations.
\end{lemma}

{\em Proof.} This follows immediately from the inspection of the figures of the r-cube equations in Section 6.

$\Box$

However, this weight doesn't satisfy the l-cube equations. Therefore we have to associate a weight to the remaining global types of Reidemeister III moves too.

\begin{figure}
\centering
\includegraphics{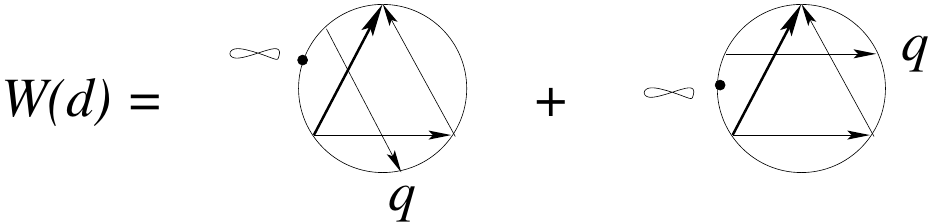}
\caption{\label{weightrc} the weight $W(d)$ of a Reidemeister III move of type $r_c$}  
\end{figure}

\begin{definition}
The weight $W(d)$ of a Reidemeister III move of type $r_c$ is defined by the configurations in Fig.~\ref{weightrc}.
\end{definition}

\begin{figure}
\centering
\includegraphics{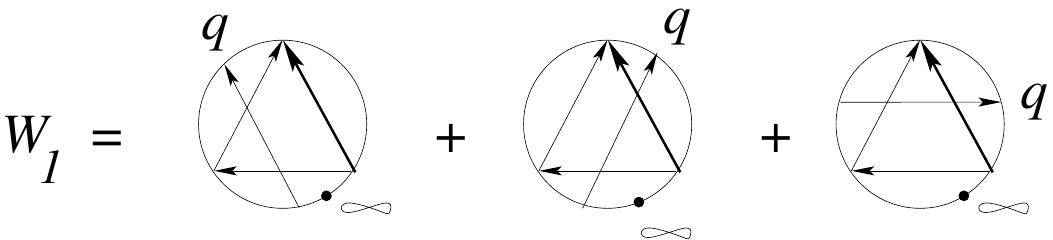}
\caption{\label{weightla} the preliminary weight $W_1$ of a Reidemeister III move of type $l_a$}  
\end{figure}

\begin{definition}
The preliminary weight $W_1$ of a Reidemeister III move of type $l_a$ is defined by the configurations in Fig.~\ref{weightla}.
The weight $W_1(type)$ depends moreover on the local type of the Reidemeister III move and is given below:

$W_1(type 1)=W_1(type 3)=W_1(type 6)=W_1+1$

$W_1(type 2)=W_1(type 7)=W_1(type 8)=W_1-1$

$W_1(type 4)=W_1(type 5)=W_1$
\end{definition}

We are now ready  to define the 1-cochain $V^{(1)}$ (the letter "V" stands for Vassiliev).

\begin{definition}
Let $s$ be a generic oriented arc in $M_T$ (with a fixed abstract closure $T \cup \sigma$ to an oriented circle and a fixed point at infinity in $\partial T$).  Then the 
1-cochain   $V^{(1)}$ is defined as\vspace{0,2 cm}

 $V^{(1)}(s)=-\sum_{p \in s \cap l_a}sign(p)W_1(type)+\sum_{p \in s \cap r_c}sign(p)W(d)$
\\ \vspace{0,2 cm}

$+\sum_{p \in s \cap II^-_1}sign(p)W(p)$ 
\end{definition}

\begin{lemma}
$V^{(1)}$ satisfies the positive global tetrahedron equation. Moreover, the contribution of $-P_2+\bar P_2$ is always zero.
\end{lemma}
{\em Proof.} First of all we observe that we can forget about the constant correction $W_1(type 1)=W_1+1$ because each stratum appears twice and with different signs. We inspect the figures in Section 3 and we sum up the contributions $V^{(1)}$ from $-P_1$ to $-\bar P_4$:

$I_1$: 0, $I_2$: 0, $I_3$: 0, $I_4$: -2+2.

$II_1$: 0, $II_2$: +1-1, $II_3$: 0, $II_4$: 0.

$III_1$: 0, $III_2$: 0, $III_3$: -1-1-1+2+1, $III_4$: 0.

$IV_1$: 0, $IV_2$: +1-1-1+1, $IV_3$: 0, $IV_4$: 0. 

$V_1$: 0, $V_2$: +1-1, $V_3$: 0, $V_4$: +2-2.

$VI_1$: 0, $VI_2$: 0, $VI_3$: +1-2-2+1+2, $VI_4$: 0.

Moreover, we see that the contributions of $-P_2$ and $+\bar P_2$ always cancel out.

$\Box$

\begin{lemma}
$V^{(1)}$ vanishes on the meridians of $\Sigma^{(1)} \cap \Sigma^{(1)}$, $\Sigma^{(2)}_{self-flex}$ and $\Sigma^{(2)}_{cusp-deg}$.
\end{lemma}
{\em Proof.} Completely obvious. But notice that the strata $\Sigma^{(2)}_{cusp-deg}$ force us to associate $sign(p)W(p)$ to Reidemeister moves of type $II^-_1$ instead of type $II^+_1$.

$\Box$

\begin{lemma}
$V^{(1)}$ is zero on all meridians of $\Sigma^{(2)}_{trans-cusp} \setminus \Sigma^{(2)}_{l_c}$. 

\end{lemma}
{\em Proof.} By inspection of Fig.~\ref{over-cusp0}...Fig.~\ref{under-cusp1} from the previous section. Notice that the four cases with a triple crossing of global type $l_a$ force us to define the above constant corrections for the local types 1,3,7 and 8.
Notice that this correction term is exactly $w(hm)$ for these four cases. The correction term for the remaining cases is forced by the cube equations.

$\Box$

\begin{lemma}
$V^{(1)}$ satisfies the cube equations.

\end{lemma}
{\em Proof.} We know already that the weight $W(p)$ of a Reidemeister II move of type $II^-_1$ satisfies the r-cube equations (Lemma 11). The weight $W(d)$ of a Reidemeister III move of type $r_c$ satisfies evidently the r-cube equations because the arrow $q$ can not slide onto an arrow of the triangle without sliding over $\infty$ (compare Observation 1 in Section 6). So, we are left with the l-cube equations for $\infty=a$. We inspect the figures from Section 6 (and of course the signs from Fig.~\ref{loc-trip} and from Definition 4). The contribution of the triple crossings is always in brackets (do not forget about the constant correction term).

"1-7": -(-1+1)=-(+1-1), "1-5": -(0+1)=-1-(0), "5-3": -(0)-1=-(0+1), "7-4": (0-1)+1=(0), "3-8": -(-1+1)=-(+1-1), 
"4-8": (0)=+1+(0-1).

We determine now the correction term for the star-like triple crossings. Reidemeister moves of type $II_1^-$ do no longer occur.

"1-6": -(0+1)=-(0+1), "7-2": -(0-1)=-(0-1), "3-6": (0+1)=(0+1), "4-6": (0)=(-1+1), "5-2": (0)=(+1-1), "8-2": -(0-1)=-(0-1).

$\Box$

The following remark is very important.
\begin{remark}
It is tempting to complete $V^{(1)}$ to a 1-cocycle in $M_T$ by adding $\sum_{p \in s \cap l_c}sign(p)W_1(p)$ (compare Definition 12). Then it would indeed be zero on all meridians of $\Sigma^{(2)}_{trans-cusp}$. However, 
$\sum_{p \in s \cap l_c}sign(p)W_1(p)$ does not satisfy the  global positive tetrahedron equation as one sees immediately from the figures in Section 3 for the global case $II$ and $\infty = 2$. Hence, the quantum part $R^{(1)}$ of the 1-cocycle is essential in order to obtain a 1-cocycle for the whole $M_T$.
\end{remark}

\begin{proposition}
$V^{(1)}$ is an integer valued 1-cocycle in $M_T \setminus \Sigma^{(2)}_{l_c}$. Its value on a meridian of $\Sigma^{(2)}_{l_c}$ is equal to $+W_1(p)$ or $-W_1(p)$ of the triple crossing. It has the scan-property for branches moving under the tangle $T$.

\end{proposition}
{\em Proof.} It follows from the lemmas in this section that $V^{(1)}$ vanishes on the meridians of five of the six types of strata of codimension two (compare Section 2). Moreover, according to Lemma 14 it vanishes also on the meridians of $\Sigma^{(2)}_{trans-cusp} \setminus \Sigma^{(2)}_{l_c}$. It follows that $V^{(1)}$ is a 1-cocycle in 
$M_T \setminus \Sigma^{(2)}_{l_c}$. It has the scan-property for positive triple crossings according to Lemma 12.
Inspecting the figures in Section 6 we see that it has the scan-property for Reidemeister moves of type $II^-_1$ (and hence for all Reidemeister moves of type II). It follows now that $V^{(1)}$ has the scan-property because the graph $\Gamma$ is connected. This is exactly the same argument like in the case of $R^{(1)}$. 

$\Box$

\begin{remark}
$V^{(1)}$ can not be defined as a 1-cocycle of finite type in Vassiliev's \cite{V1} sense because it is only defined in 
$M_T \setminus \Sigma^{(2)}_{l_c}$ (we just take out something of codimension two in the components of smooth non singular knots and we will not see anything dual  in the discriminant of singular knots). However, it is of finite type in the sense that it can be given by a Gauss diagram formula of finite type (i.e. a finite number of arrows in each configuration in the formula).  We define the {\em Gauss degree} of a Gauss diagram formula in the last section. According to this definition $V^{(1)}$ is of Gauss degree 3.

\end{remark}

\section{The 1-cocycle $\bar R^{(1)}=R^{(1)}-\delta P_T V^{(1)}$ represents a non trivial cohomology class in the topological moduli space}

\begin{definition}
The completion $\bar R^{(1)}$  is defined by $\bar R^{(1)}=R^{(1)}-\delta P_T V^{(1)}$.
\end{definition}

The following theorem is the most important result in this paper.

\begin{theorem}
Let $s$ be a generic oriented arc in $M_T$ (with a fixed abstract closure $T \cup \sigma$ to an oriented circle and a fixed point at infinity in $\partial T$). Then the 
1-cochain with values in $S(\partial T)$ \vspace{0,2 cm}

$\bar R^{(1)}(s)=R^{(1)}(s)-\delta P_T V^{(1)}(s)=$
 \\ \vspace{0,2 cm}

$+\sum_{p \in s \cap l_c}sign(p)W_1(p)T_{l_c}(type)(p) $
\\ \vspace{0,2 cm}

$ +\sum_{p \in s \cap r_a}sign(p)[W_1(type)(p)T_{r_a}(type)(p)  + zW_2(type)(p)T(type)(p)] $
\\ \vspace{0,2 cm}

$+ \sum_{p \in s \cap r_b}sign(p)zW_2(p)T(type)(p)  +\sum_{p \in s \cap l_b}sign(p)zW_2(p)T(type)(p) $
\\ \vspace{0,2 cm}

$ + \sum_{p \in s \cap II^+_0}sign(p)W_2(p)T_{II}(p) +\sum_{p \in s \cap II^-_0}sign(p)W_2(p)zP_T$
\\ \vspace{0,2 cm}
 
$+\sum_{p \in s \cap I}sign(p)T_I(p)$ 
\\ \vspace{0,2 cm}

$+\sum_{p \in s \cap l_a}sign(p)W_1(type)(p)\delta P_T- \sum_{p \in s \cap r_c}sign(p)W(p)\delta P_T $
\\ \vspace{0,2 cm}

$ - \sum_{p \in s \cap II^-_1}sign(p)W(p) \delta P_T$ 
\\ \vspace{0,2 cm}

is a 1-cocycle in  $M_T$. It represents a non trivial cohomology class and it has the scan-property for branches moving under the tangle T.

\end{theorem}

\begin{figure}
\centering
\includegraphics{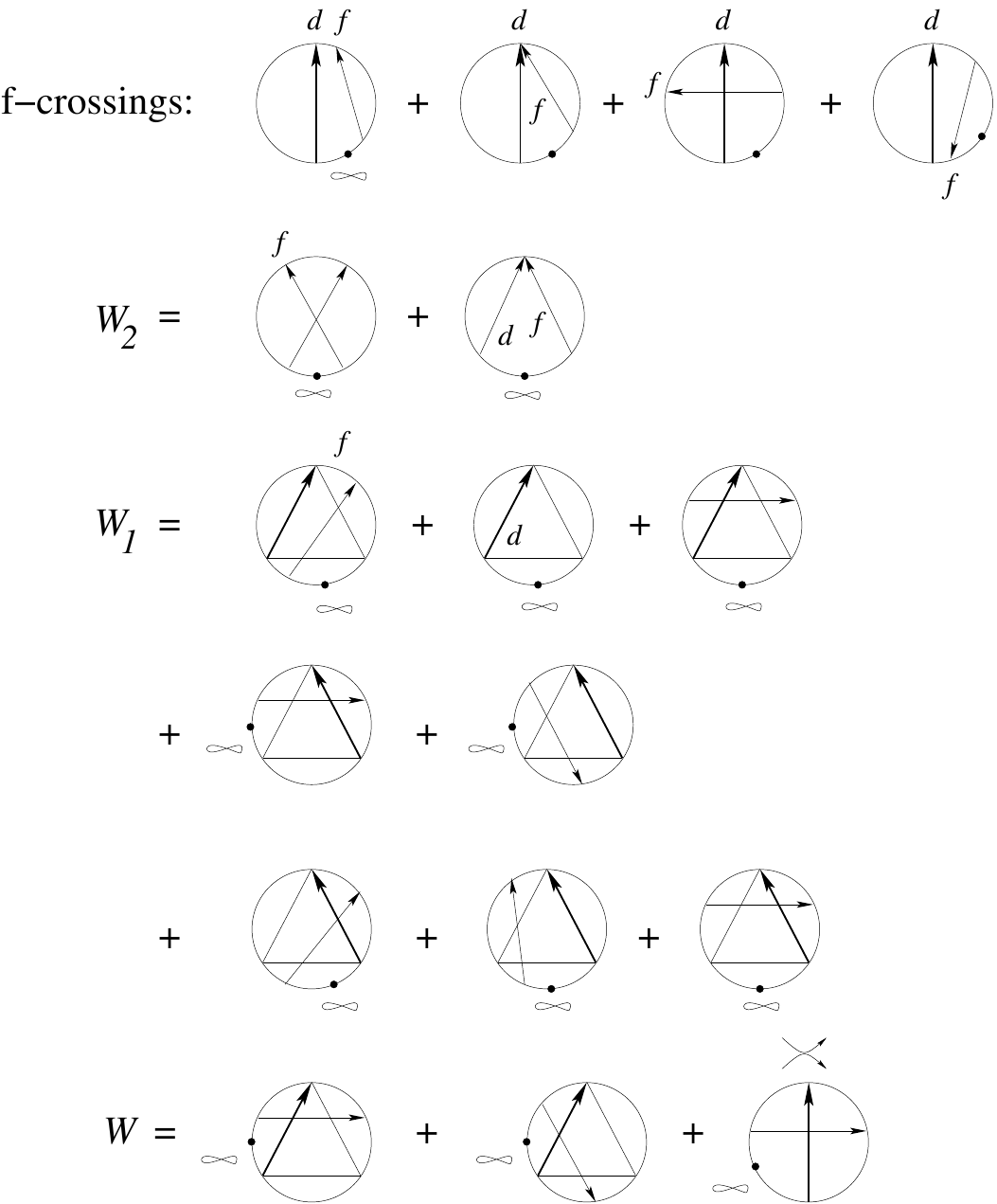}
\caption{\label{Rbarw} the weights for $\bar R^{(1)}$}  
\end{figure}

\begin{figure}
\centering
\includegraphics{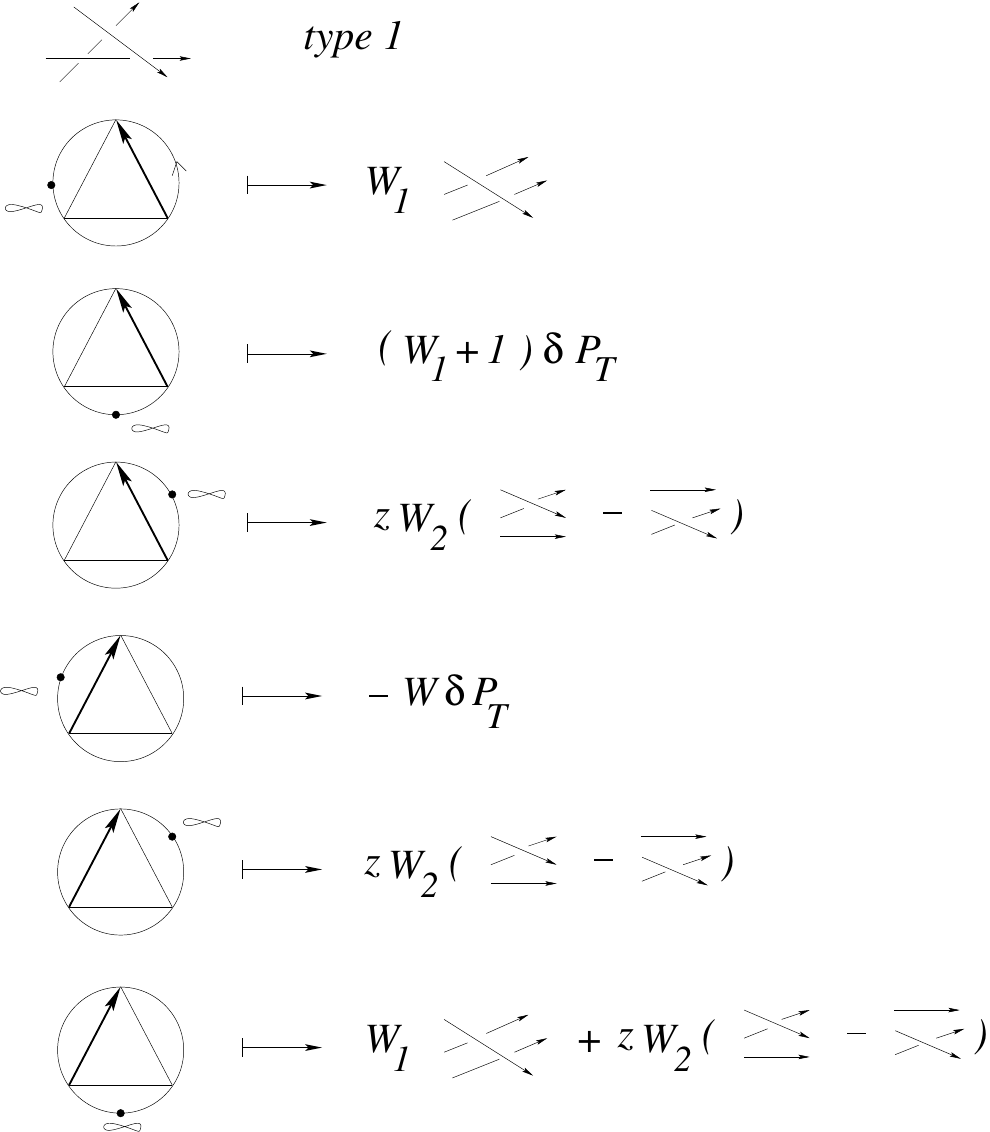}
\caption{\label{Rbar1} the partial smoothings for $RIII$ of type 1}  
\end{figure}

\begin{figure}
\centering
\includegraphics{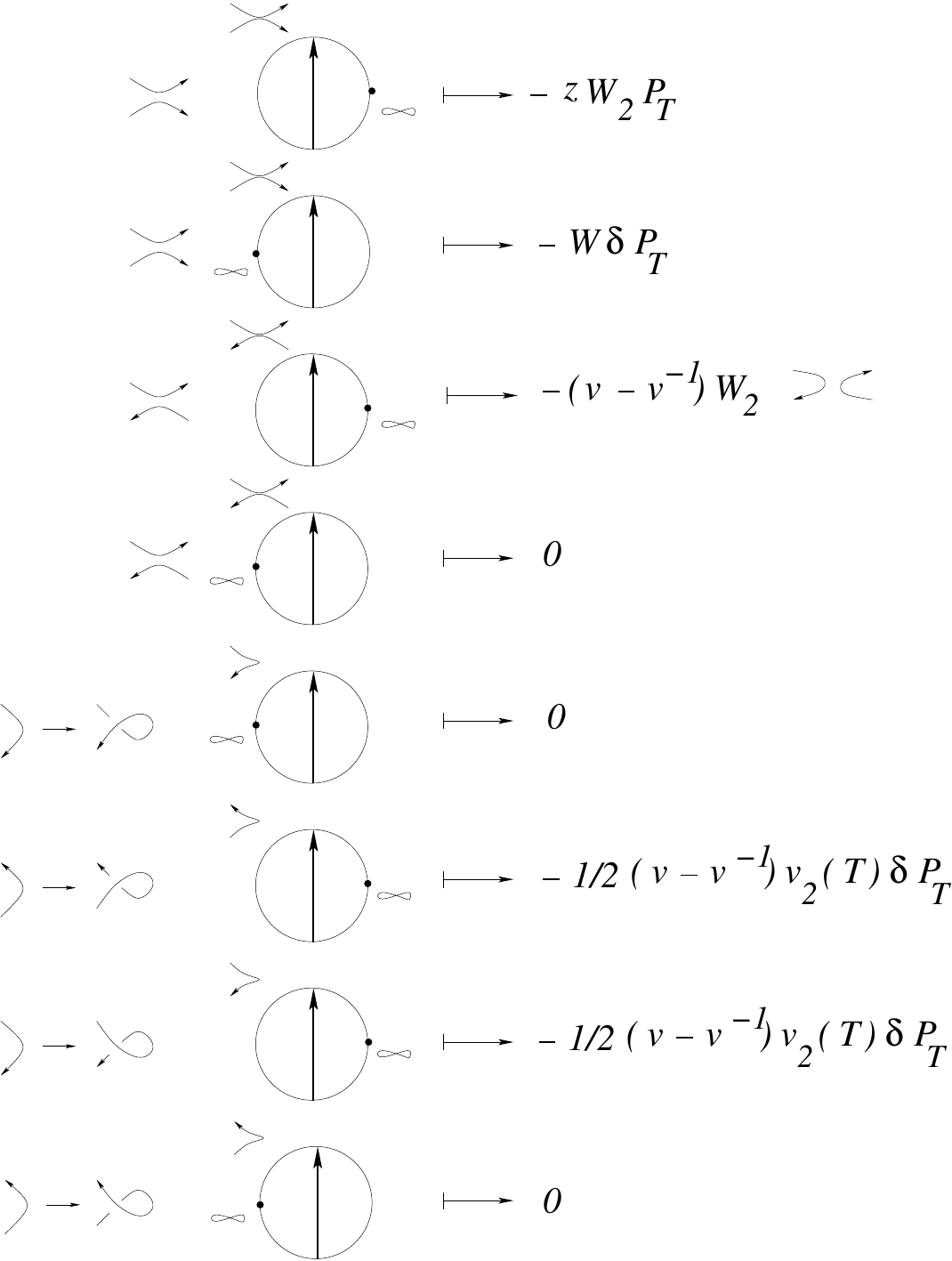}
\caption{\label{RbarII} the partial smoothings for Reidemeister moves of type II and I}  
\end{figure}

For the convenience of the reader all the weights are summarized in Fig.~\ref{Rbarw}. Here $d$ is the distinguished crossing of a R III move (i.e. the crossing between the highest and the lowest branch) or the two new crossings together of a R II move or the new crossing of a R I move. Remember our convention that the contribution of a configuration (i.e. a Gauss diagram without writhes) is the sum of the product of the writhe of the arrows (not in the triangle). Hence $W_1$ and $W$ are linear weights and $W_2$ is a quadratic weight. The partial smoothings are summarized in Fig.~\ref{Rbar1} up to Fig.~\ref{RbarII}. Here we show the triple crossing and the self-tangency on the negative side of the discriminant  $\Sigma^{(1)}$, i.e. the Reidemeister move to the other side of the discriminant enters with the sign $+1$. The distinguished crossing $d$ is always drawn with a thicker arrow. Remember also that $\delta$ is the HOMFLYPT polynomial of the trivial link of two components, that $v_2(T)$ is the invariant of degree two defined by the first Polyak-Viro formula in Fig.~\ref{PV} and that $P_T$ is the element in the HOMFLYPT skein module $S(\partial T)$ represented by $T$.

{\em Proof.} It follows from Lemma 6, Propositions 6 and 7 that $\bar R^{(1)}$ is a 1-cocycle in $M_T \setminus \Sigma^{(2)}_{l_c}$ and that it has the scan-property. It remains to show that it vanishes on the meridians of 
$\Sigma^{(2)}_{l_c}$. But we see this immediately from  Fig.~\ref{merover-cusp} and from the corresponding Gauss diagrams in Fig.~\ref{over-cusp0} (compare also the proof of Lemma 8).

$\Box$

\begin{figure}
\centering
\includegraphics{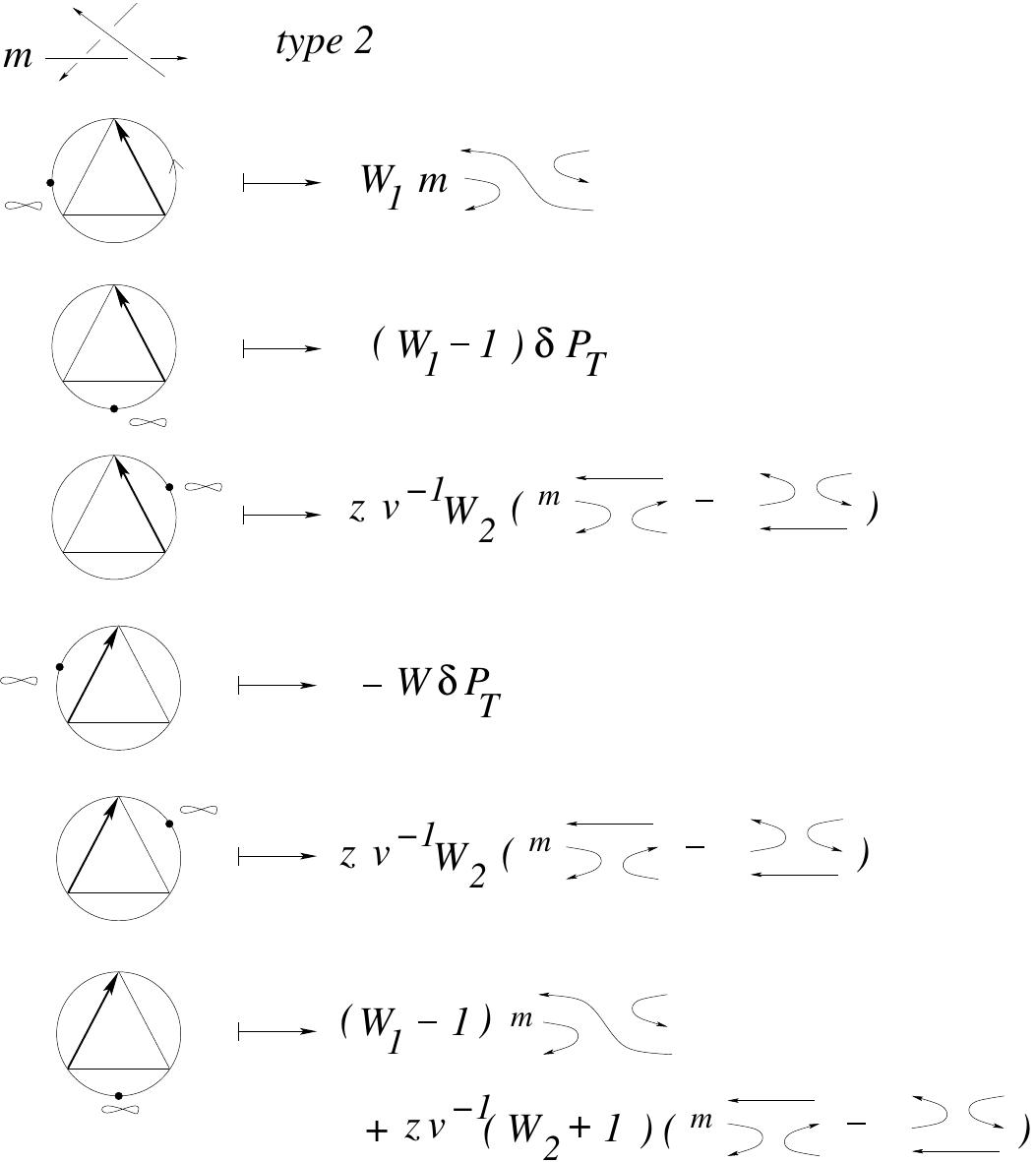}
\caption{\label{Rbar2} the partial smoothings for $RIII$ of type 2}  
\end{figure}

\begin{figure}
\centering
\includegraphics{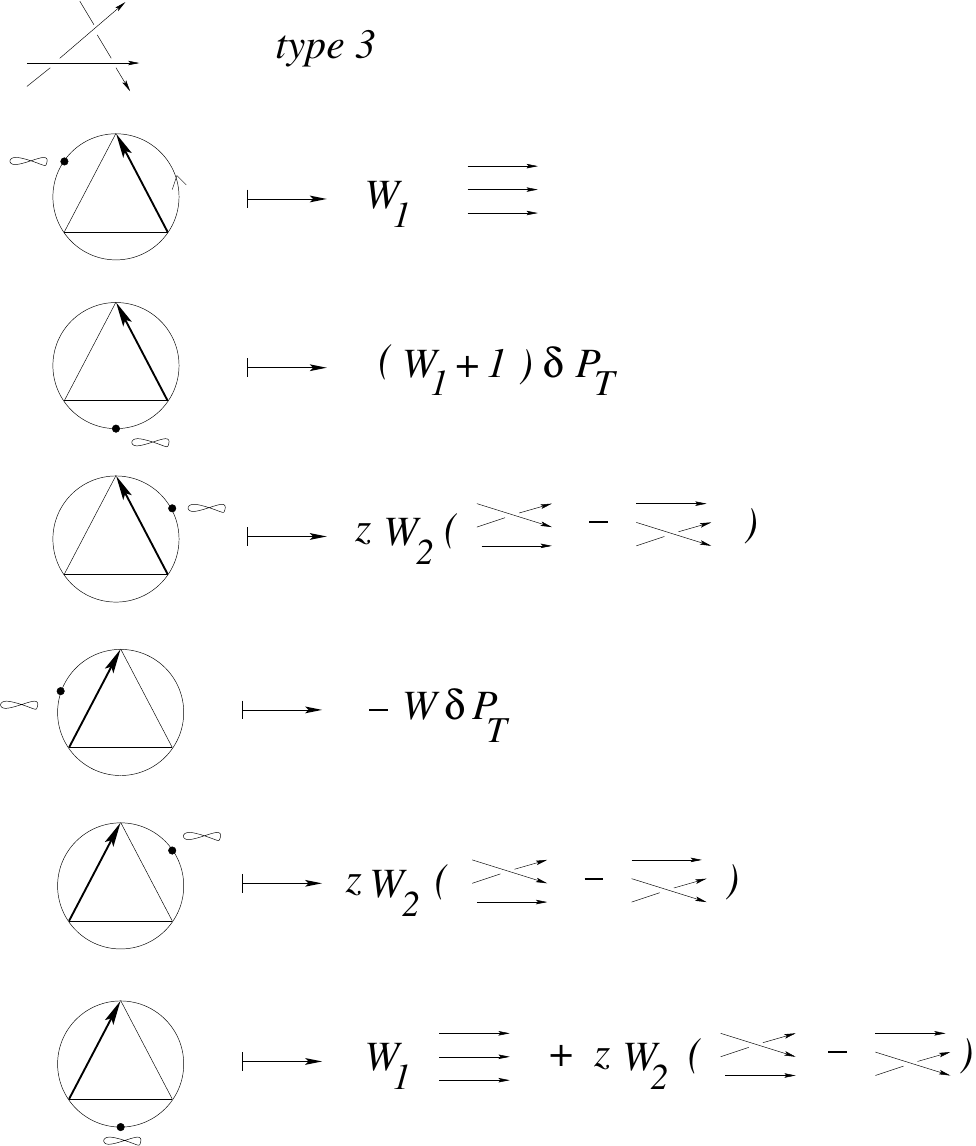}
\caption{\label{Rbar3} the partial smoothings for $RIII$ of type 3}  
\end{figure}

\begin{figure}
\centering
\includegraphics{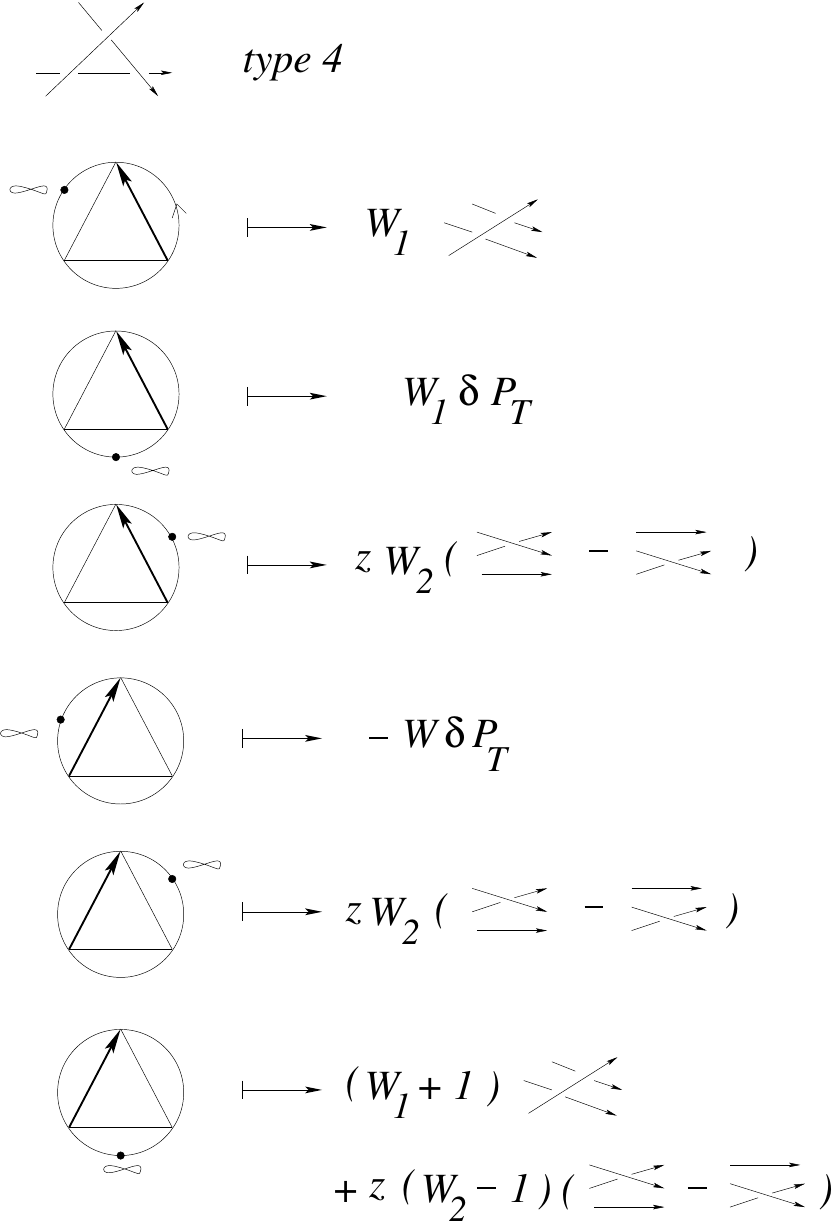}
\caption{\label{Rbar4} the partial smoothings for $RIII$ of type 4}  
\end{figure}

\begin{figure}
\centering
\includegraphics{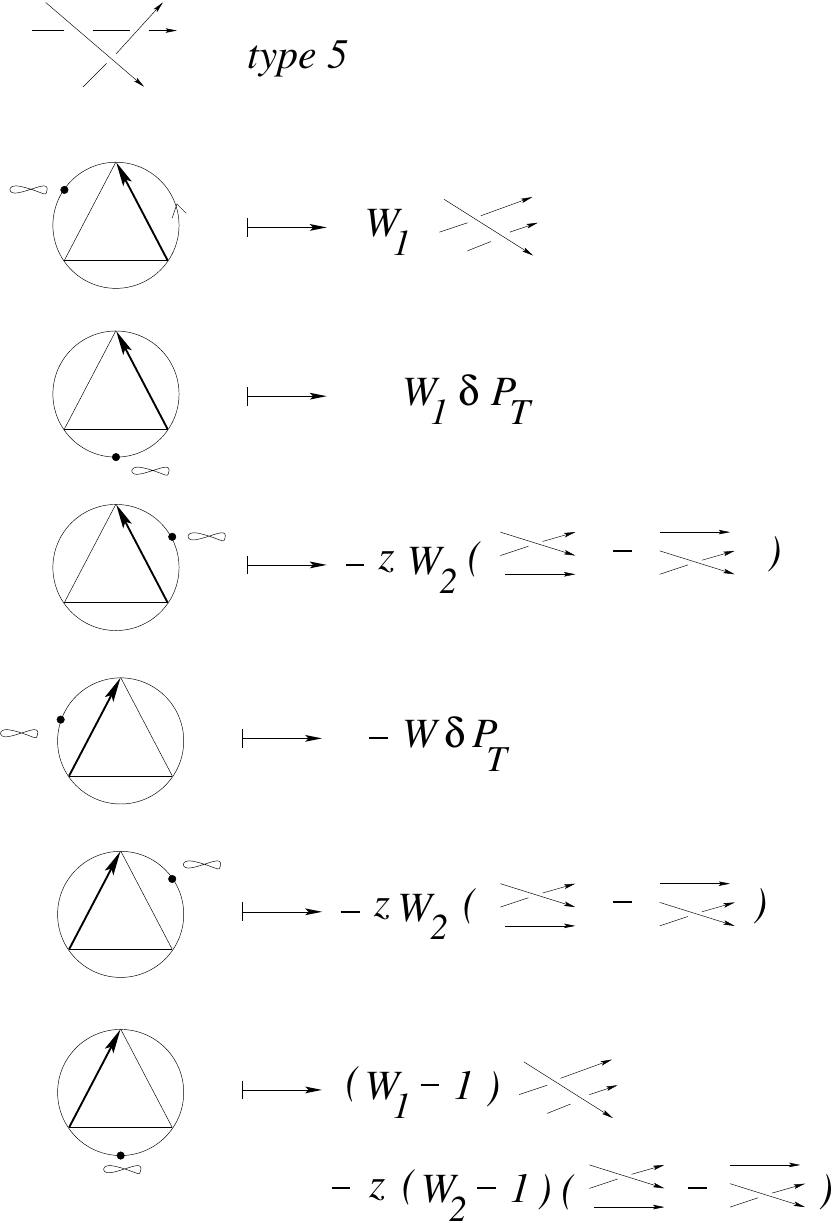}
\caption{\label{Rbar5} the partial smoothings for $RIII$ of type 5}  
\end{figure}

\begin{figure}
\centering
\includegraphics{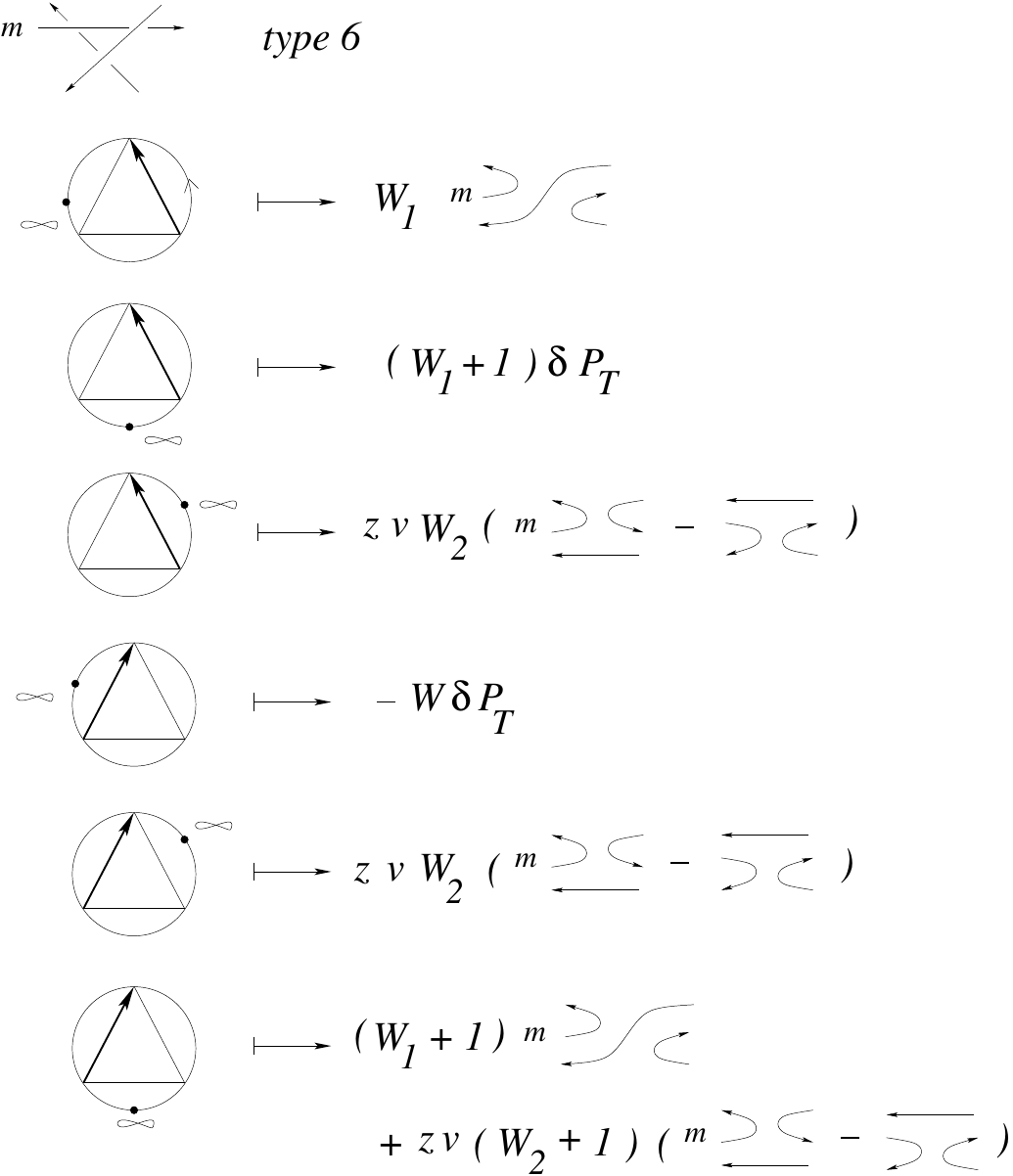}
\caption{\label{Rbar6} the partial smoothings for $RIII$ of type 6}  
\end{figure}

\begin{figure}
\centering
\includegraphics{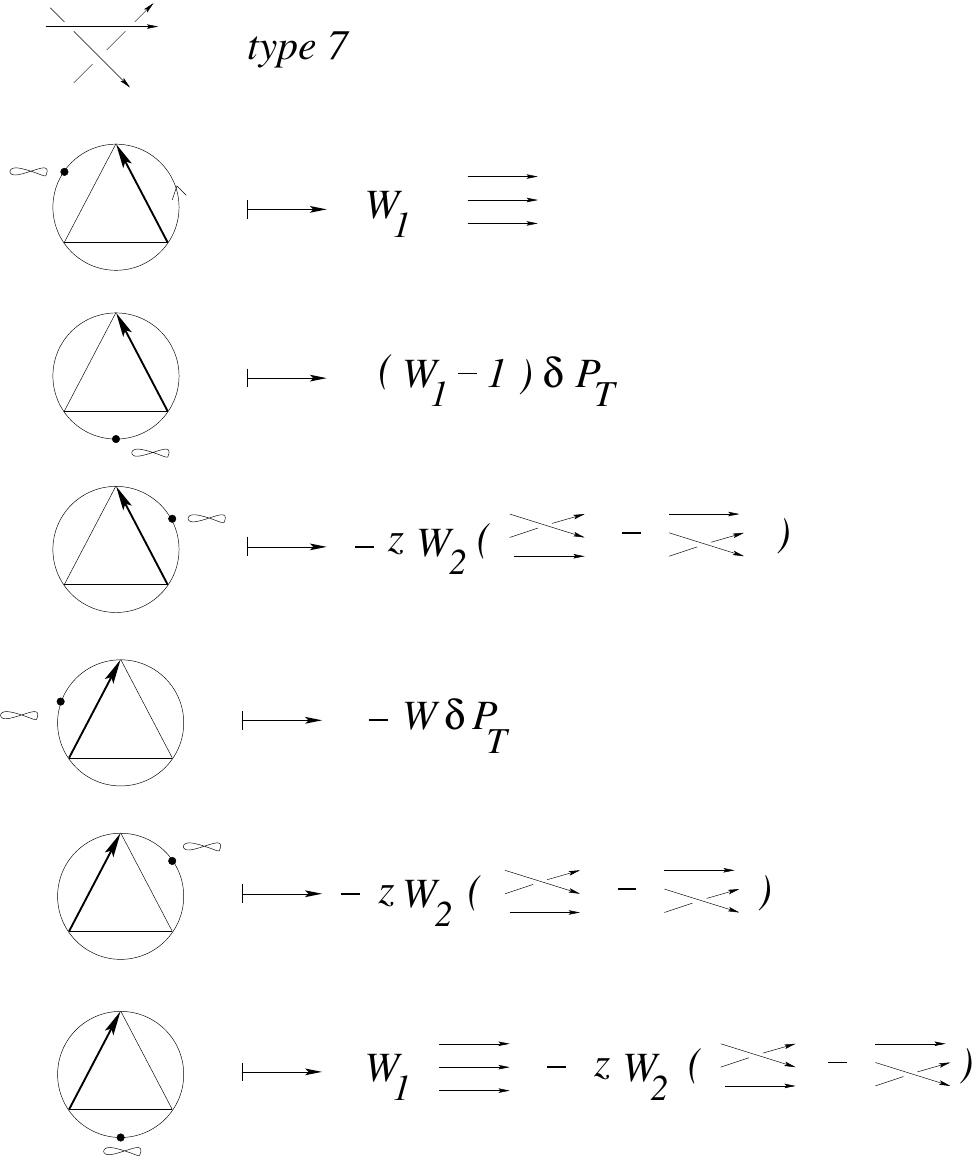}
\caption{\label{Rbar7} the partial smoothings for $RIII$ of type 7}  
\end{figure}

\begin{figure}
\centering
\includegraphics{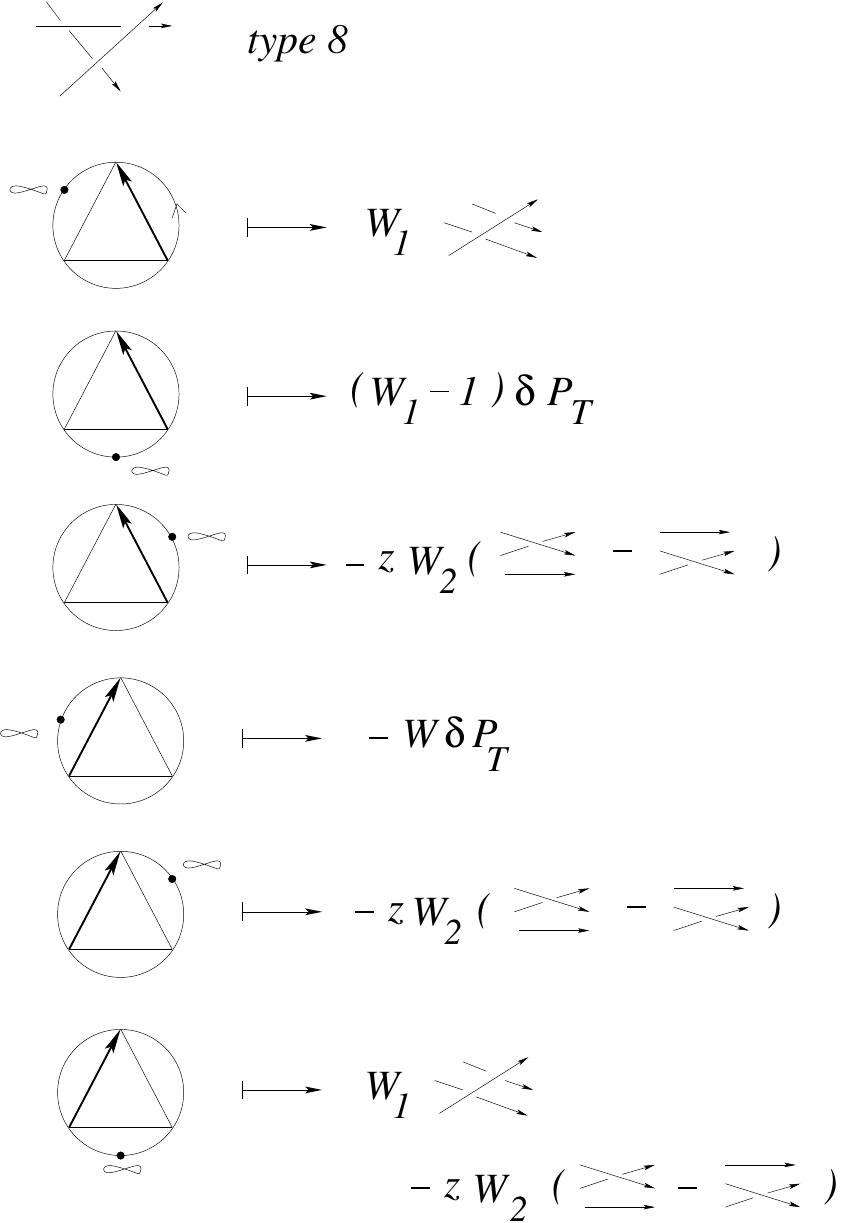}
\caption{\label{Rbar8} the partial smoothings for $RIII$ of type 8}  
\end{figure}

An important property of this 1-cocycle is the fact that it represents a non trivial cohomology class in particular for $\amalg_K M_K$, where $K$ runs over all isotopy classes of long knots.

\begin{figure}
\centering
\includegraphics{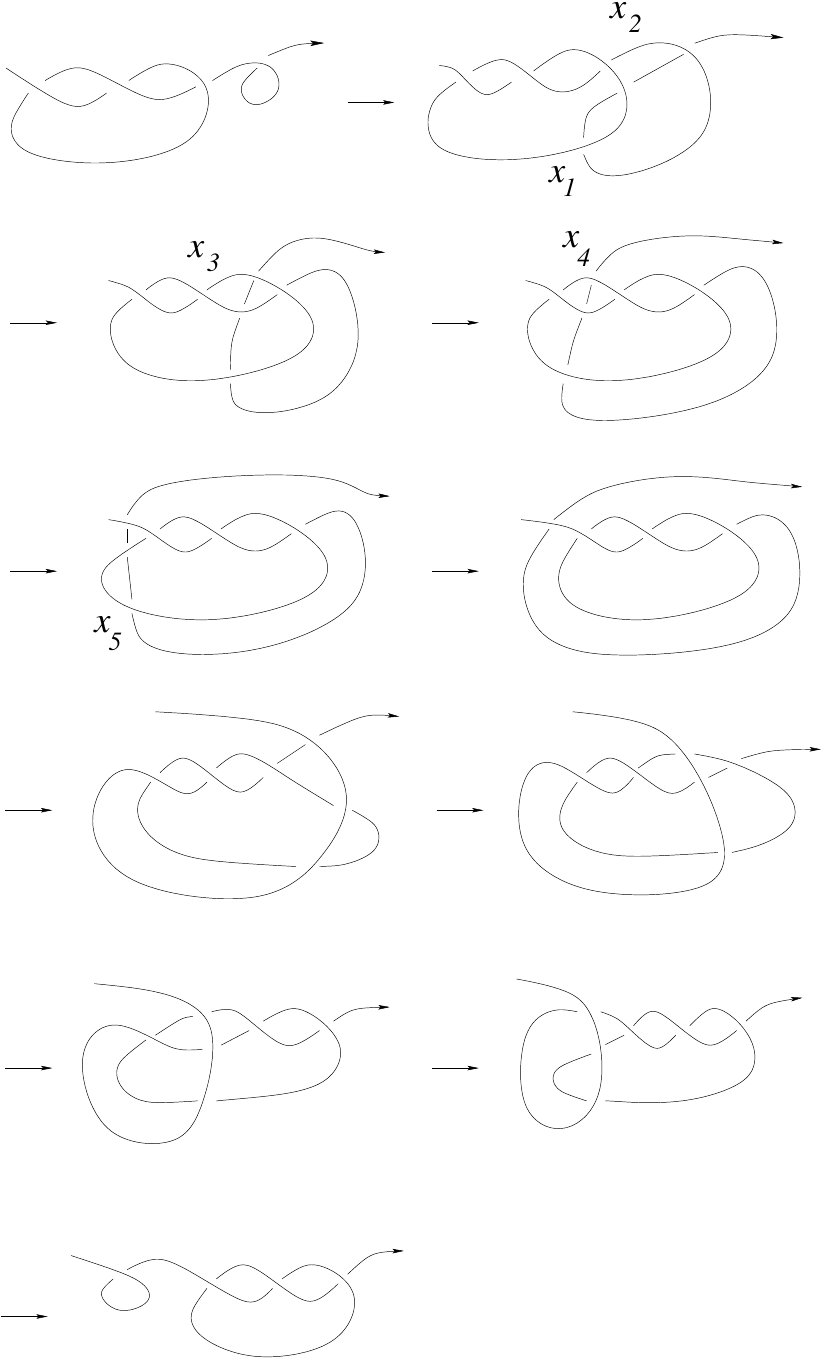}
\caption{\label{rotK} the regular part of the loop $rot(3_1^+)$}  
\end{figure}

\begin{figure}
\centering
\includegraphics{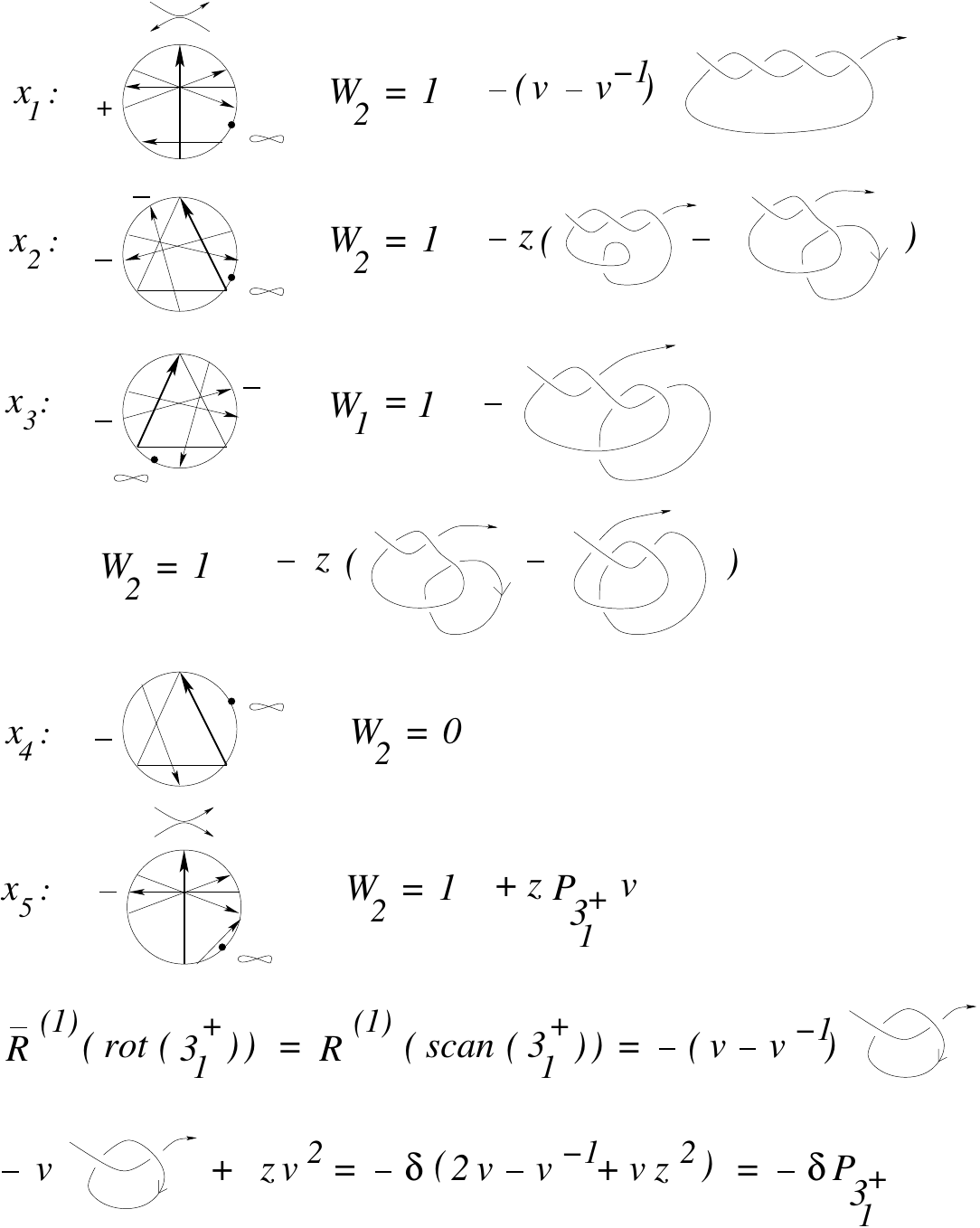}
\caption{\label{calrot} calculation of $\bar R^{(1)}(scan(3_1^+))$}  
\end{figure}

\begin{figure}
\centering
\includegraphics{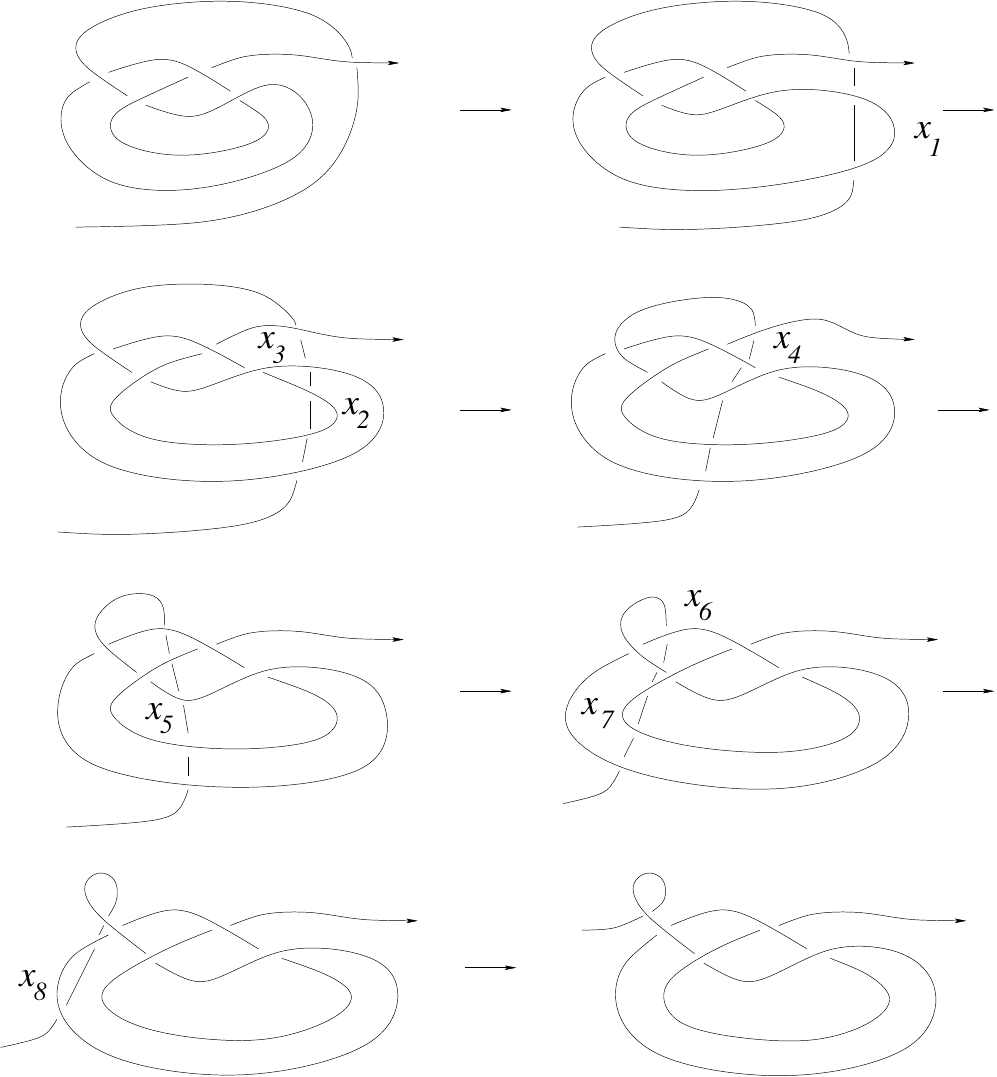}
\caption{\label{rot4} the second half of an alternative loop for $rot(4_1)$}  
\end{figure}

\begin{figure}
\centering
\includegraphics{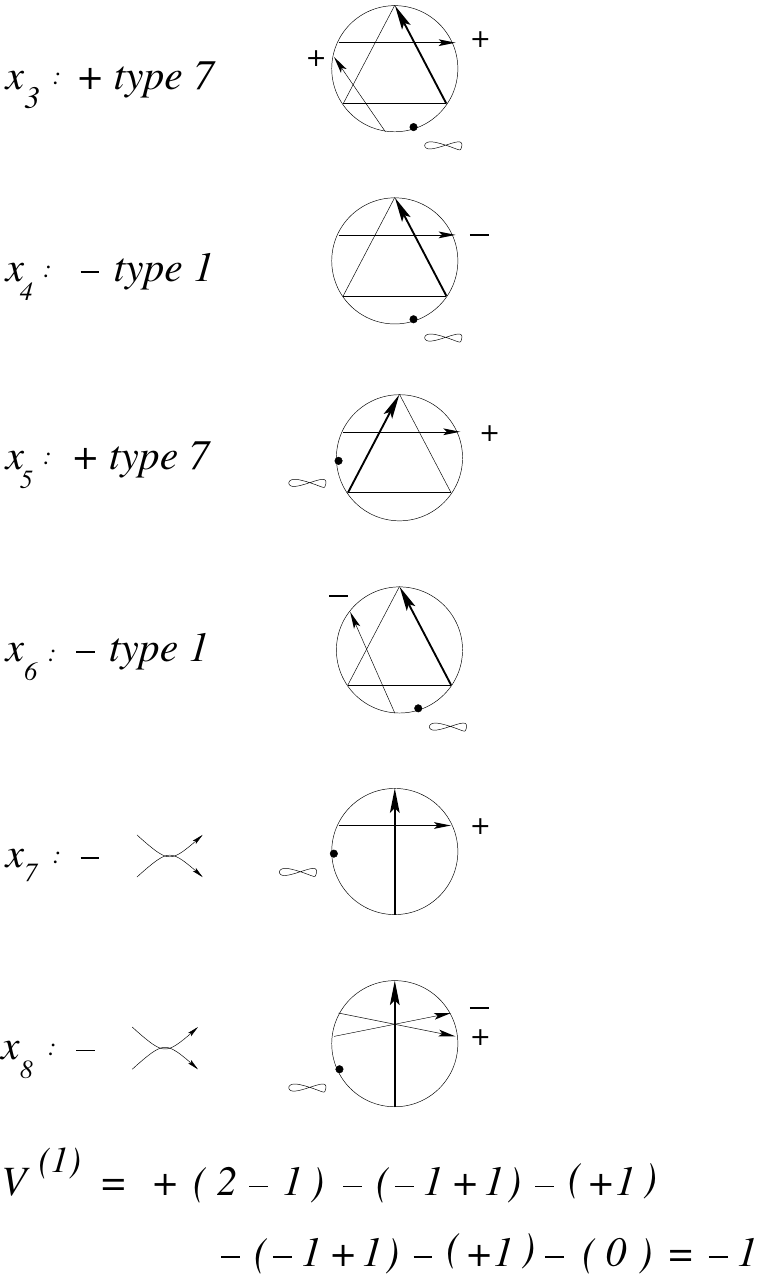}
\caption{\label{calrot4} calculation of $V^{(1)}(rot(4_1)$}  
\end{figure}

\begin{example}
We start with showing that $\bar R^{(1)}(rot K)=-\delta P_K$ for $K$ the right trefoil (compare the Introduction). A version of the loop $rot(K)$ as almost regular isotopy is shown in  Fig.~\ref{rotK} (compare \cite{T}) and the calculation of $\bar R^{(1)}$ is shown in  Fig.~\ref{calrot} (compare Remark 13 below).

We calculate $\bar R^{(1)}(rot K)=\delta P_K$ for $K$ the figure eight knot by using the alternative loop from Remark 14 below. The relevant part is contained in  Fig.~\ref{rot4} and the calculation of $V^{(1)}$ from the relevant part of the Gauss diagrams is given in Fig.~\ref{calrot4}.
\end{example}

\begin{figure}
\centering
\includegraphics{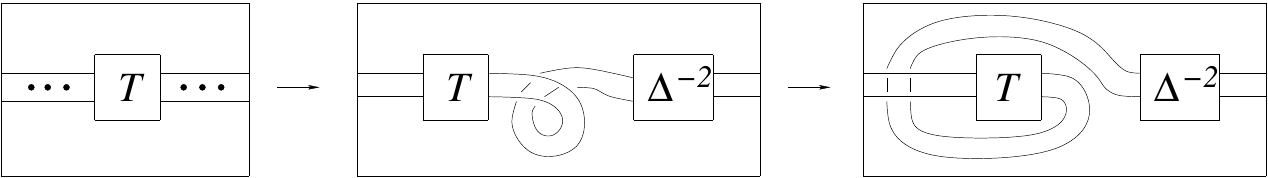}
\caption{\label{totscan} the scan-arc $rot(T)$}  
\end{figure}

\begin{remark}
The scan-arc $scan(K)$ is the first half of $rot(K)$. In the second half we move always from $\infty$ to the over cross of a distinguished crossing $d$. Consequently, strata of type $l_c$ and of type $r_a$ do not occur. Moreover, from $\infty$ we move only over the rest of the diagram up to the overcross of $d$. Hence there are no f-crossings at all for the strata of type $l_b$, $r_b$ and $II_0$. It follows that the second half does not contribute to $ R^{(1)}(rot (K))$. There aren't any strata of type $r_c, l_a, l_c, II_1$ at all in $rot(K)$ and the contributions of the two Reidemeister I moves cancel out. Consequently $V^{(1)}$ doesn't contribute and we obtain 

$\bar R^{(1)}(rot (K))=R^{(1)}(scan (K))$
 
for all long knots $K$. 

For a general tangle $T$  we could consider the scan-arc for all components simultaneously as shown in Fig.~\ref{totscan}. This scan-arc generalizes $rot(K)$ for long knots. Let us call it $rot(T)$.  Does $R^{(1)}(rot (T))$ also contain $P_T \in H_n(z,v)$ as a factor for each choice of a point at infinity in $\partial T$ and each abstract closure $\sigma$ of $T$ to a circle?
\end{remark}

\begin{example}
The loop $drag 3_1^+$ is a regular isotopy. Notice that many Reidemeister moves do not contribute to $R^{(1)}$ because they have not the right global type or they have the weight zero but they contribute to $V^{(1)}$. We leave the calculation of  $\bar R^{(1)}(drag 3_1^+)=-3\delta P^2_{(3_1^+)}$ to the reader as an exercise.
\end{example}

\begin{figure}
\centering
\includegraphics{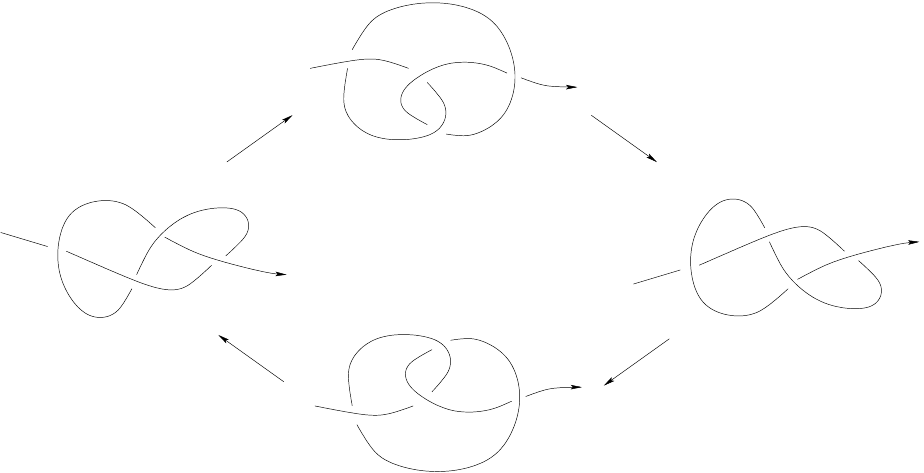}
\caption{\label{hatK} Hatcher's loop for the figure eight knot}  
\end{figure}

\begin{figure}
\centering
\includegraphics{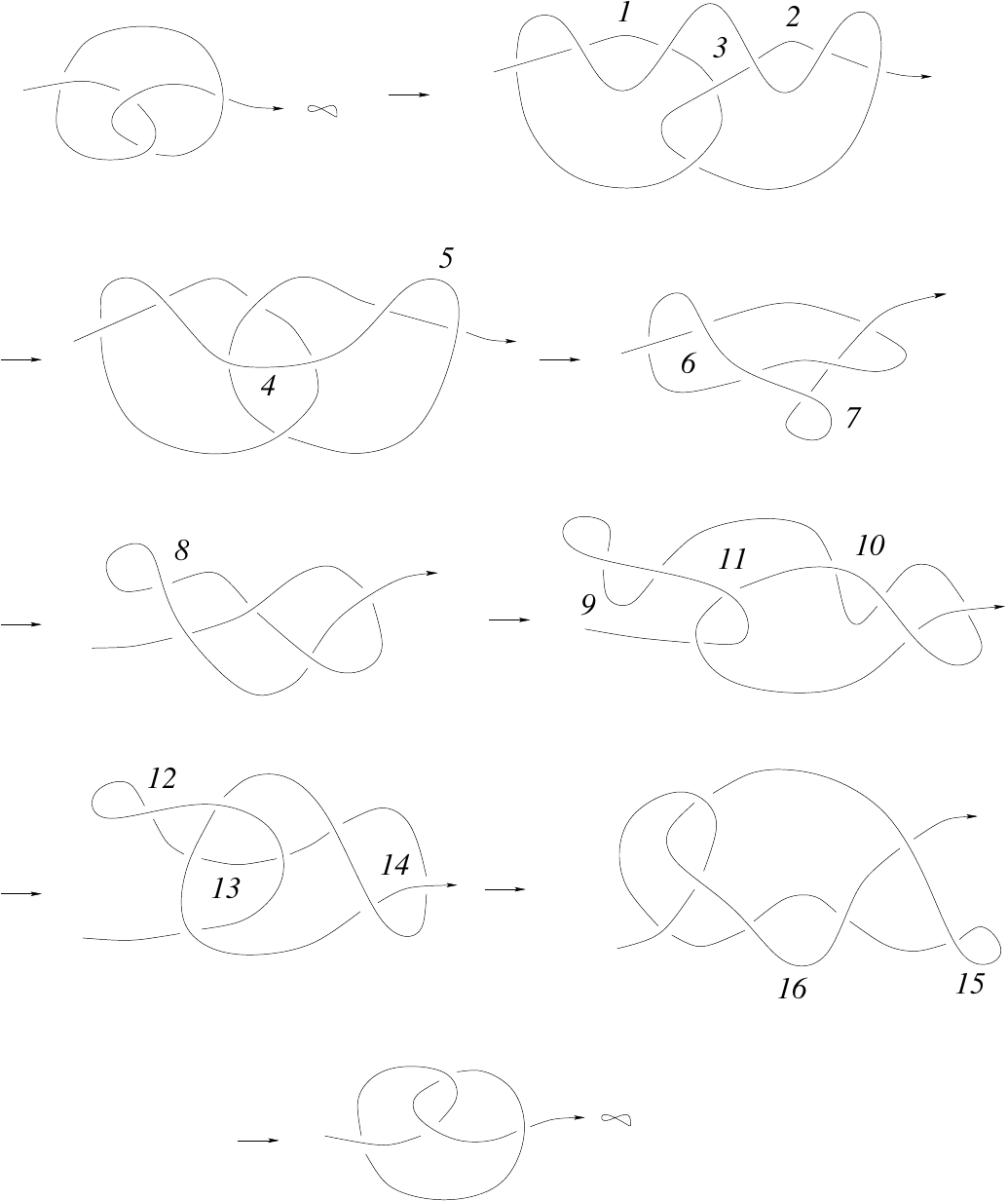}
\caption{\label{hatK1} first half of Hatcher's loop}  
\end{figure}

\begin{figure}
\centering
\includegraphics{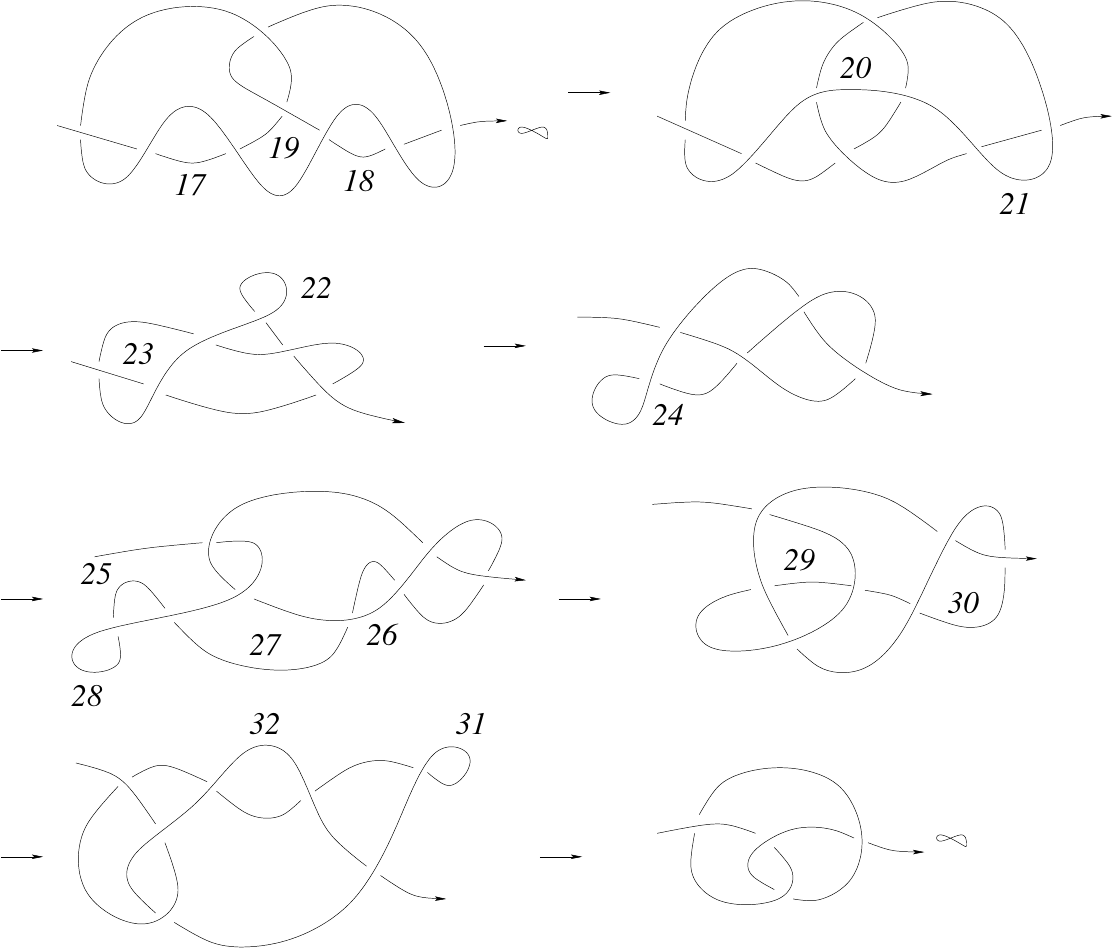}
\caption{\label{hatK2} second half of Hatcher's loop}  
\end{figure}

\begin{figure}
\centering
\includegraphics{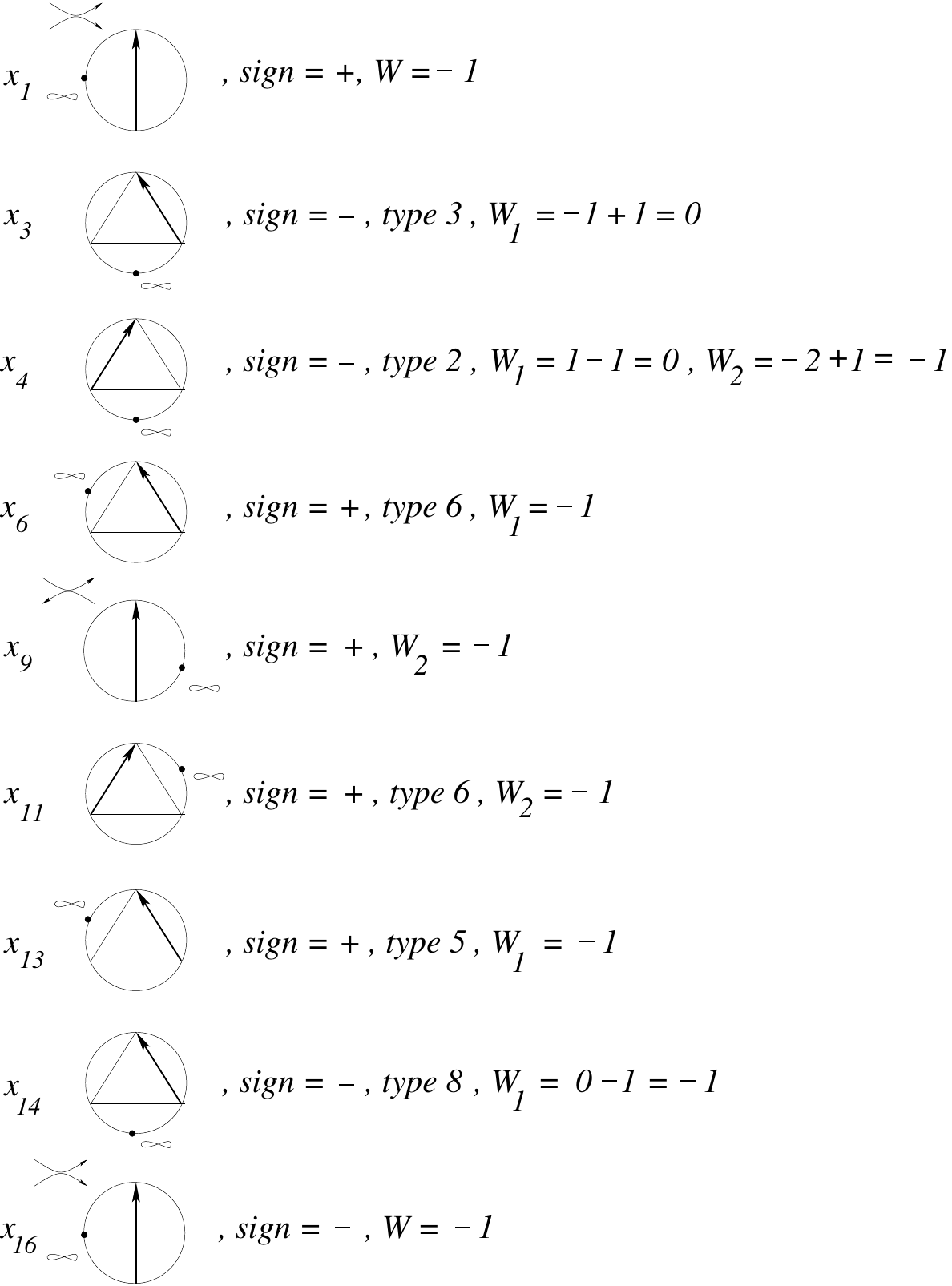}
\caption{\label{hatK3} contributing Reidemeister moves in $hat(K)$}  
\end{figure}

\begin{figure}
\centering
\includegraphics{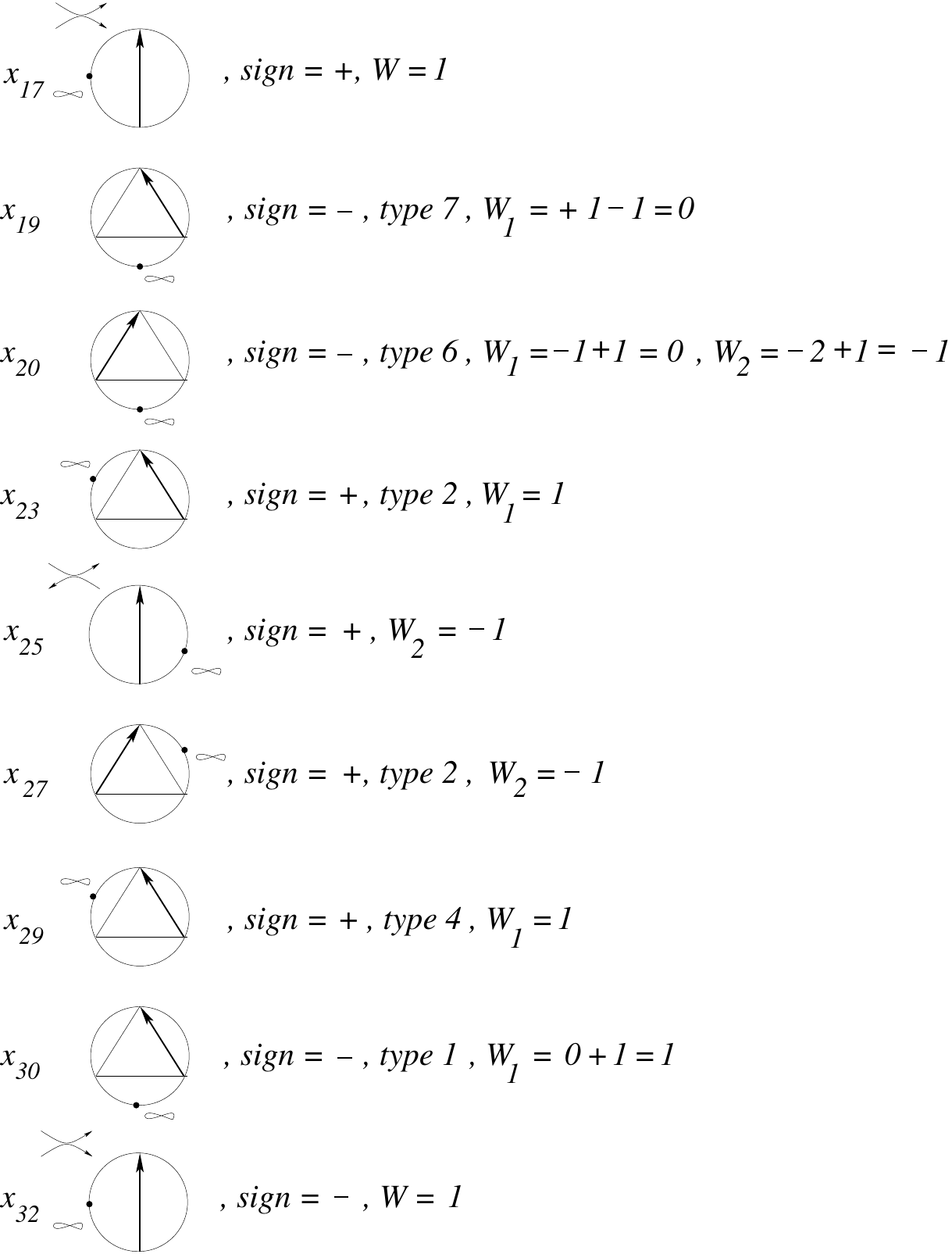}
\caption{\label{hatK4} the remaining Reidemeister moves that contribute in $hat(K)$}  
\end{figure}

\begin{example}
We show that $\bar R^{(1)}(hat (K))=0$ for  the figure eight knot $K$ with trivial framing. Because of its importance we give the calculation in some detail, so that the reader can easily check the calculation. The loop $hat K$ is shown in Fig.~\ref{hatK}, compare \cite{H}. We represent it by a loop which contains R I moves only of negative sign in Fig.~\ref{hatK1} and Fig.~\ref{hatK2}. (So, the reader could transform it into a regular loop by collecting the small curls at one place in the diagram and  eliminating them two by two using Whitney tricks. But attention, pushing a small curl along $K$ does not contribute to $\bar R^{(1)}$ but it does contribute to $R^{(1)}_{reg}$.) The contributions of the Reidemeister moves $x_2$ and $x_5$ as well as $x_{18}$ and $x_{21}$ cancel out. Moreover there are exactly four Reidemeister I moves of type $0$, namely $x_7$, $x_{12}$, $x_{22}$ and $x_{28}$. Each of them has $sign=-1$, and as well known $v_2(4_1)=-1$. Consequently, the R I moves contribute $-2(v-v^{-1}) \delta P_T$. The moves $x_{10}$ and $x_{26}$ do not contribute at all. We give the types of the remaining moves together with the sign and the corresponding weights in Fig.~\ref{hatK3} and Fig.~\ref{hatK4}. One easily calculates now that  $R^{(1)}(hat (K))=0$ and that $V^{(1)}(hat (K))=0$ too. It follows that $\bar R^{(1)}(hat (K))=0$.

Notice that the loops $hat(K)$ and $rot(K)$ are apriori of a different nature. The loop $hat(K)$ is in general not induced by a $SO(3)$-action on $\mathbb{R}^3$ and it mixes in general all six global types of R III moves. (Local types of R III moves can be replaced by each other but global types can not.)
\end{example}

\begin{figure}
\centering
\includegraphics{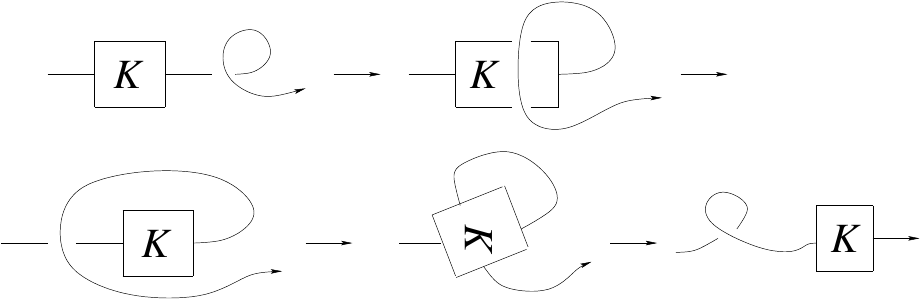}
\caption{\label{rotpos} a different loop representing $rot(K)$}  
\end{figure}

\begin{remark}
Let us consider a loop which represents $rot(K)$ but by using a positive curl with positive Whitney index as shown in  Fig.~\ref{rotpos}. We see immediately that no Reidemeister move at all contributes to  $R^{(1)}$. Hence 
$\bar R^{(1)}(rot K)=-V^{(1)}(rot(K))\delta P_K$. This means that Conjecture 1 (compare the Introduction) is not really difficult. It reduces to show that
 
$V^{(1)}(rot(K))=v_2(K)$ 

for this loop representing $rot(K)$.

In contrast to that Conjecture 2 seems to be rather difficult (if true) and needs some understanding of the relation of the hyperbolic volume of a knot $K$ with the value of $\bar R^{(1)}(rot(K) -hat(K))$.
\end{remark}

\section{Solutions of the positive global tetrahedron equation for the Kauffman polynomial}
The types, the signs, the weights, the grading in the Kauffman case are exactly the same as in the HOMFLYPT case. Only the partial smoothings are different.

\begin{proposition} The 1-cochain \vspace{0,2 cm}

$R_{F,reg}(m)(A)=\sum_{p \in m} sign(p) \sigma_2\sigma_1(p) + \sum_{p \in m}sing(p)W(p)(\sigma^2_2(p)-\sigma^2_1(p))$ \vspace{0,2 cm}

is a solution of the positive global tetrahedron equation. Here the first sum is over all triple crossings of the global types shown in Fig.~\ref{fglob}(i.e. the types $l_c$ and $r_a$) and such that $\partial  (hm)=A$.  The second sum is over all triple crossings $p$ which have a distinguished crossing $d$ of type 0 (i.e. the types $r_a$, $r_b$ and $l_b$).

\end{proposition}

\begin{figure}
\centering
\includegraphics{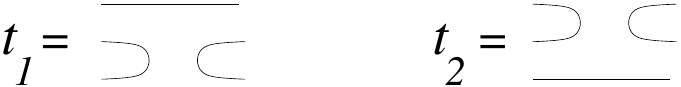}
\caption{\label{3kauf} new generators of the Kauffman skein module of $3$-braids}  
\end{figure}

\begin{figure}
\centering
\includegraphics{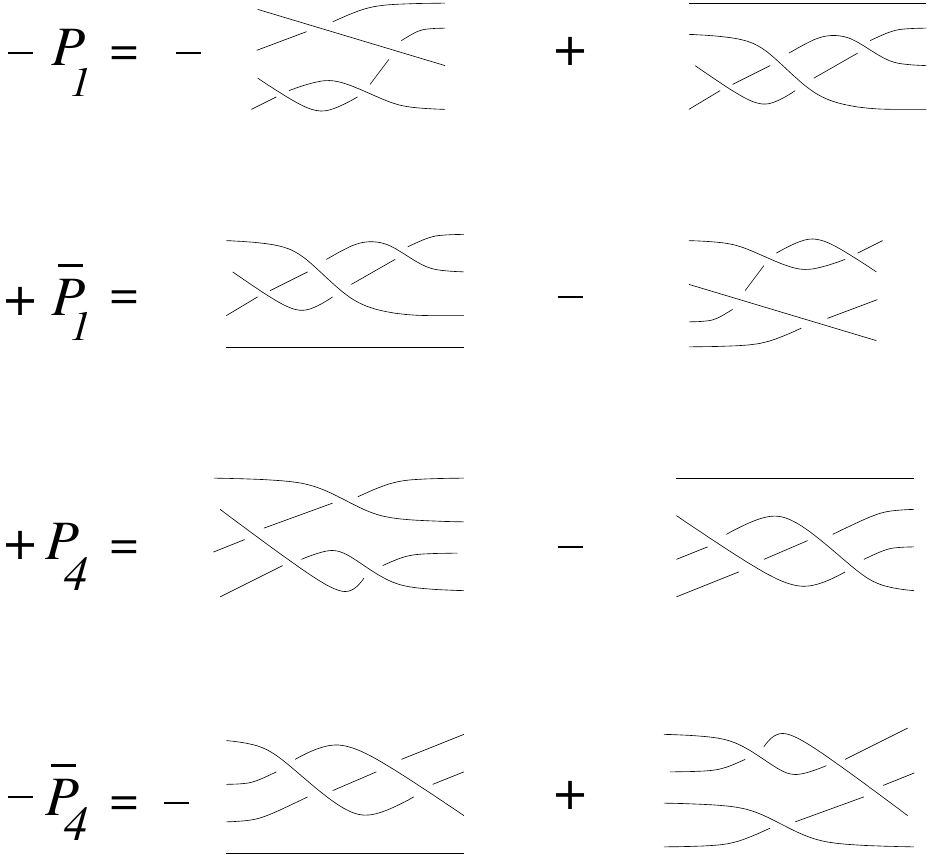}
\caption{\label{contsik} contributions of $\sigma_2(p) - \sigma_1(p)$}  
\end{figure}

\begin{figure}
\centering
\includegraphics{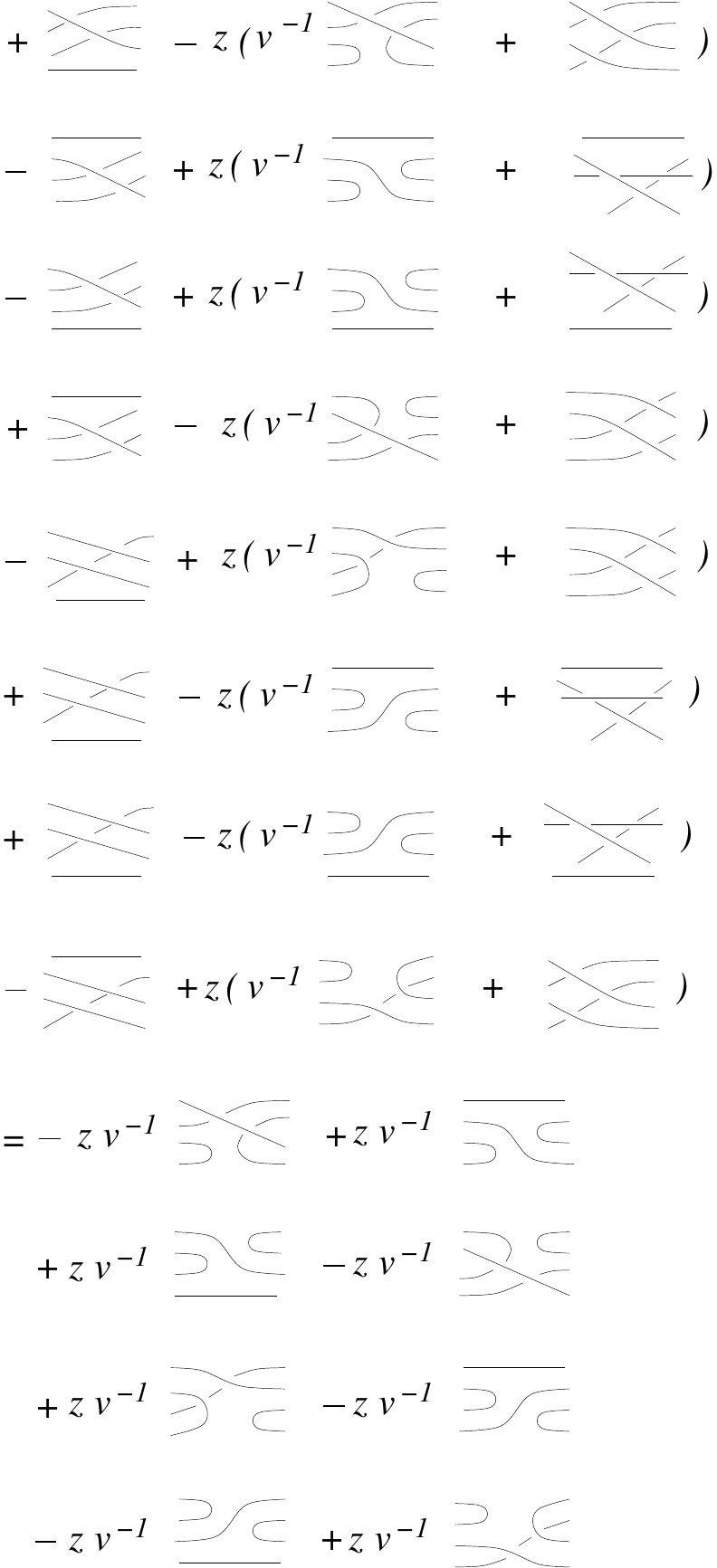}
\caption{\label{calsig} calculation of $\sigma_2(p) - \sigma_1(p)$ for $- P_1 +
\bar P_1 + P_4 - \bar P_4$}  
\end{figure}

\begin{figure}
\centering
\includegraphics{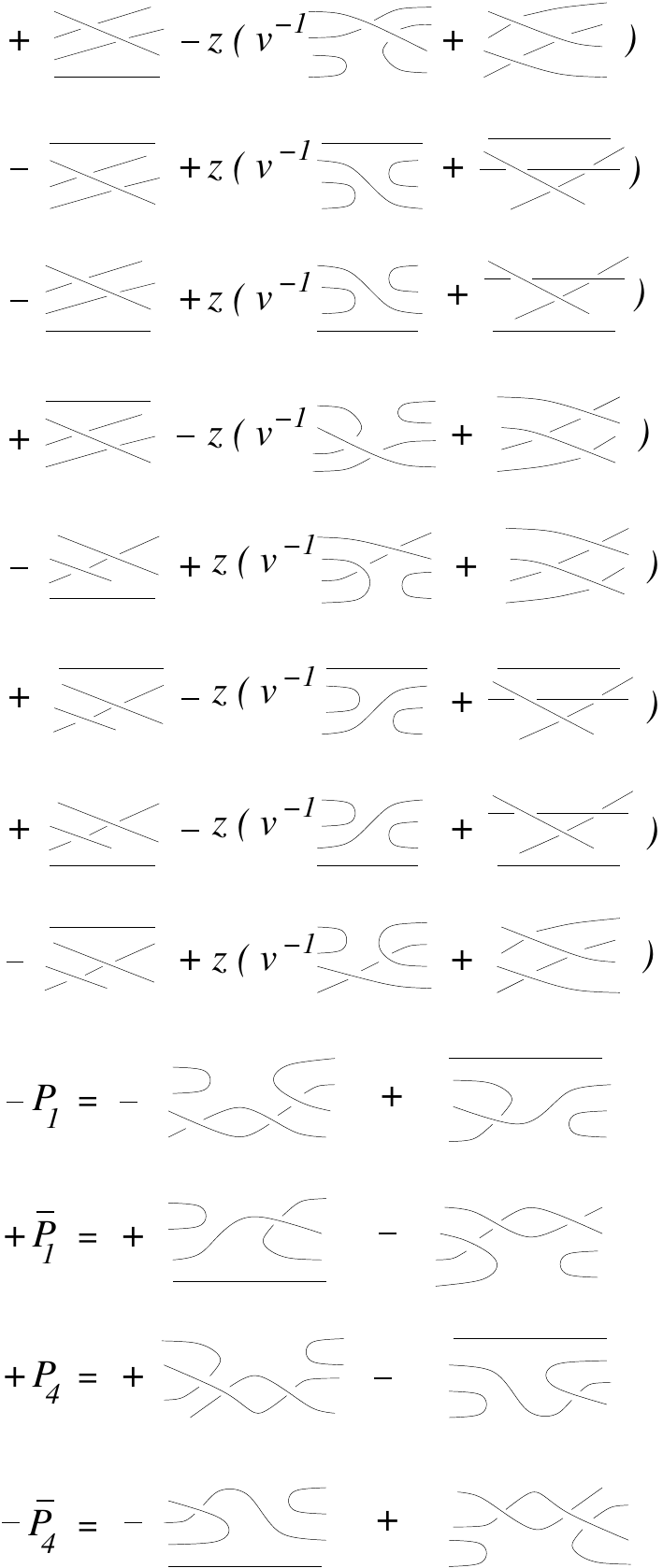}
\caption{\label{conttt} contributions of $t_2(p) - t_1(p)$}  
\end{figure}

\begin{figure}
\centering
\includegraphics{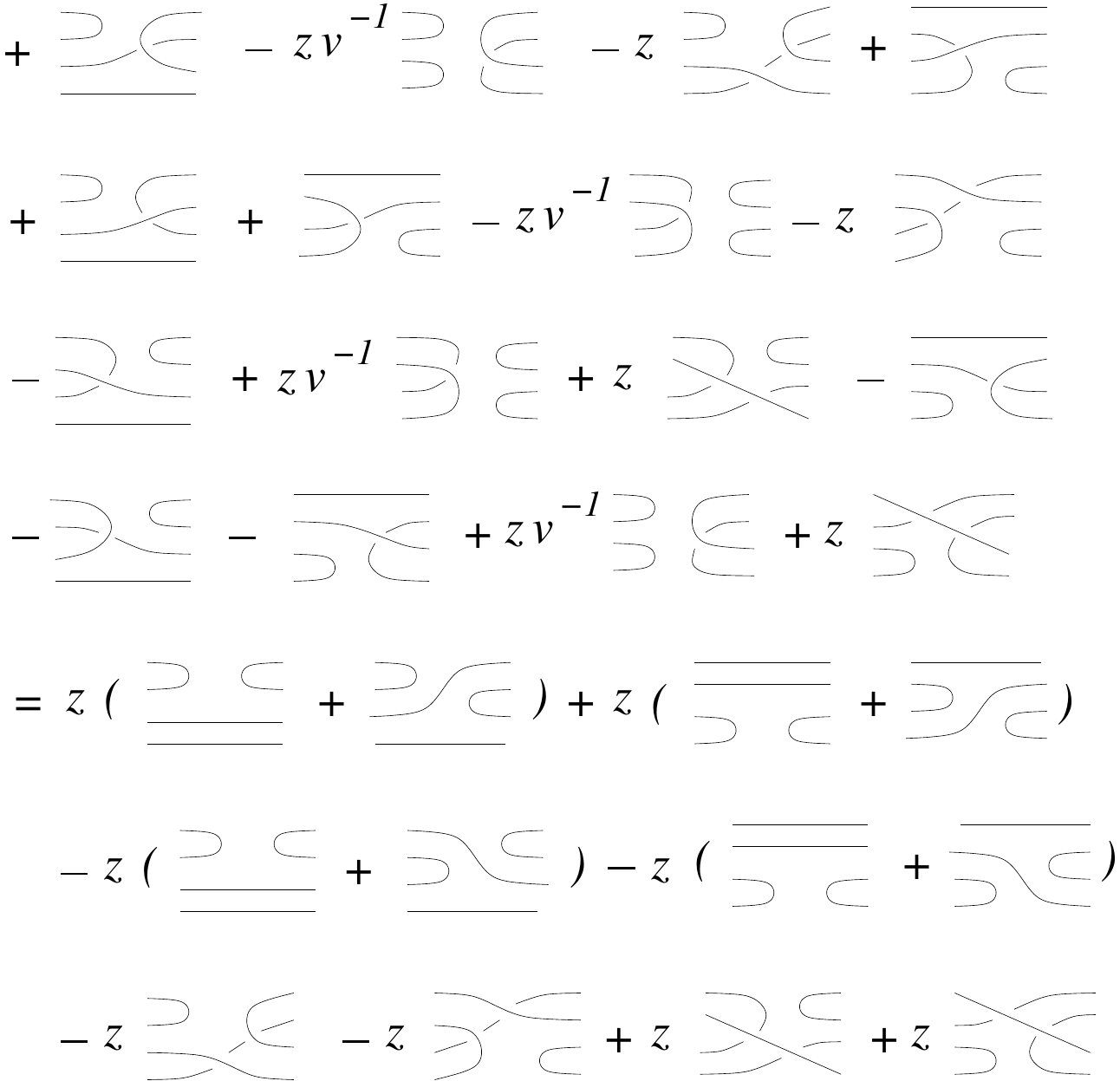}
\caption{\label{calt} calculation of of $t_2(p) - t_1(p)$ for $- P_1 +
\bar P_1 + P_4 - \bar P_4$}  
\end{figure}

\begin{lemma}
$(\sigma^2_2(p)-\sigma^2_1(p))$ is a solution with constant weight of the positive tetrahedron equation.
\end{lemma}
{\em Proof.} Let $t_2$ and $t_1$ be the new generators of the Kauffman skein module of 3-braids as shown in Fig.~\ref{3kauf}. Using the Kauffman skein relations we see that $(\sigma^2_2(p)-\sigma^2_1(p))$ is equal to 
$z(\sigma_2(p)-\sigma_1(p))+zv^{-1}(t_2(p)-t_1(p))$. Like in the HOMFLYPT case the strata $P_2, \bar P_2, P_3, \bar P_3$ do not contribute to the solution with constant weight.
We give the calculation of $(\sigma_2(p)-\sigma_1(p))$ in Fig.~\ref{contsik} and 
Fig.~\ref{calsig}.
and the calculation of $(t_2(p)-t_1(p))$ in Fig.~\ref{conttt} and  Fig.~\ref{calt}. We see that they cancel out.

$\Box$

\begin{lemma}
$\sigma_2\sigma_1(p)$ cancels out in $P_4-\bar P_4$.

\end{lemma}
{\em Proof.} The proof in the HOMFLYPT case doesn't use the skein relations and hence it is still valid in the Kauffman case.

$\Box$

\begin{figure}
\centering
\includegraphics{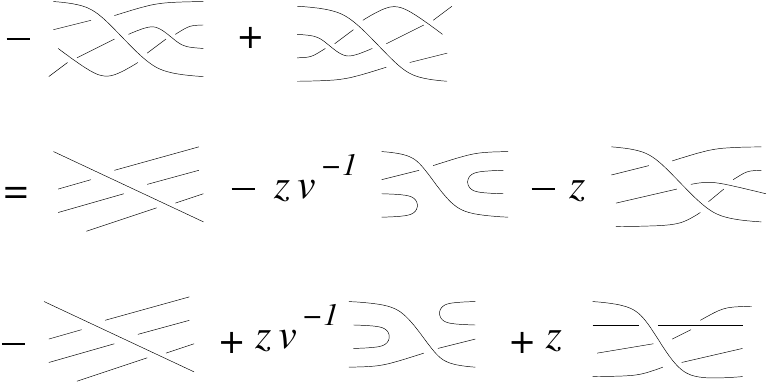}
\caption{\label{sig211} calculation of of $\sigma_2\sigma_1(p)$ for $- P_1 +
\bar P_1$}  
\end{figure}

The calculation of $\sigma_2\sigma_1(p)$ for $-P_1+\bar P_1$ is given in 
Fig.~\ref{sig211}.

\begin{figure}
\centering
\includegraphics{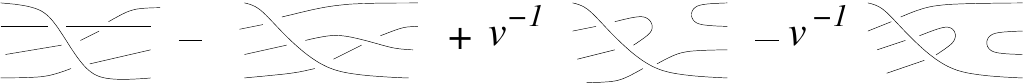}
\caption{\label{calR} calculation of of $(\sigma_2(p) - \sigma_1(p)) + v^{-1}(t_2(p) - t_1(p))$ for $P_3$}  
\end{figure}

{\em Proof of Proposition 8.}
It is completely analogous to the proof of Proposition 2 besides the calculations for $-P_1+\bar P_1+P_3-\bar P_3$. We give the remaining calculation of the contribution of $(\sigma_2(p)-\sigma_1(p))+v^{-1}(t_2(p)-t_1(p))$ from $P_3$ in Fig.~\ref{calR}.

$\Box$

We haven't solved the cube equations in this case. It seems to us that one should perhaps associate to a self-tangency with equal tangent direction (i.e. of type $II^-_0$) a partial smoothing which consists of some combination of the four 2-tangles involved in the Kauffman skein relation. But we haven't carried this out.

\begin{proposition}
The 1-cochain  \vspace{0,2 cm}

$R_F(m)=\sum_{p \in m} sign(p)W_1(p) \sigma_2\sigma_1(p) + \sum_{p \in m}sing(p)W_2(p)(\sigma^2_2(p)-\sigma^2_1(p))$ \vspace{0,2 cm}

is a solution of the positive global tetrahedron equation. Here the first sum is over all triple crossings of the global types shown in Fig.~\ref{fglob} (i.e. the types $l_c$ and $r_a$).  The second sum is over all triple crossings $p$ which have a distinguished crossing $d$ of type 0 (i.e. the types $r_a$, $r_b$ and $l_b$).

\end{proposition}
{\em Proof.} This follows immediately from the proofs of Proposition 8 and Proposition 3.

$\Box$

Surprisingly there is a second solution of the global positive tetrahedron equation in Kauffman's case.

\begin{proposition} The 1-cochain \vspace{0,2 cm}

$R^{(1)}_{F,reg}(m)(A)=\sum_{p \in m} sign(p) t_1t_2(p) + \sum_{p \in m}sing(p)W(p)z(-\sigma_1^{-1}t_2-\sigma_2t_1+t_1\sigma_2^{-1}+t_2\sigma_1)$ 
\vspace{0,2 cm}

(compare Fig.~\ref{solkauf}, where we give it in shorter form)
is also a solution of the positive global tetrahedron equation. Here the first sum is over all triple crossings of the global types shown in Fig.~\ref{fglob} (i.e. the types $l_c$ and $r_a$) and such that $\partial  (hm)=A$.  The second sum is over all triple crossings $p$ which have a distinguished crossing $d$ of type 0 (i.e. the types $r_a$, $r_b$ and $l_b$).

\end{proposition}

\begin{figure}
\centering
\includegraphics{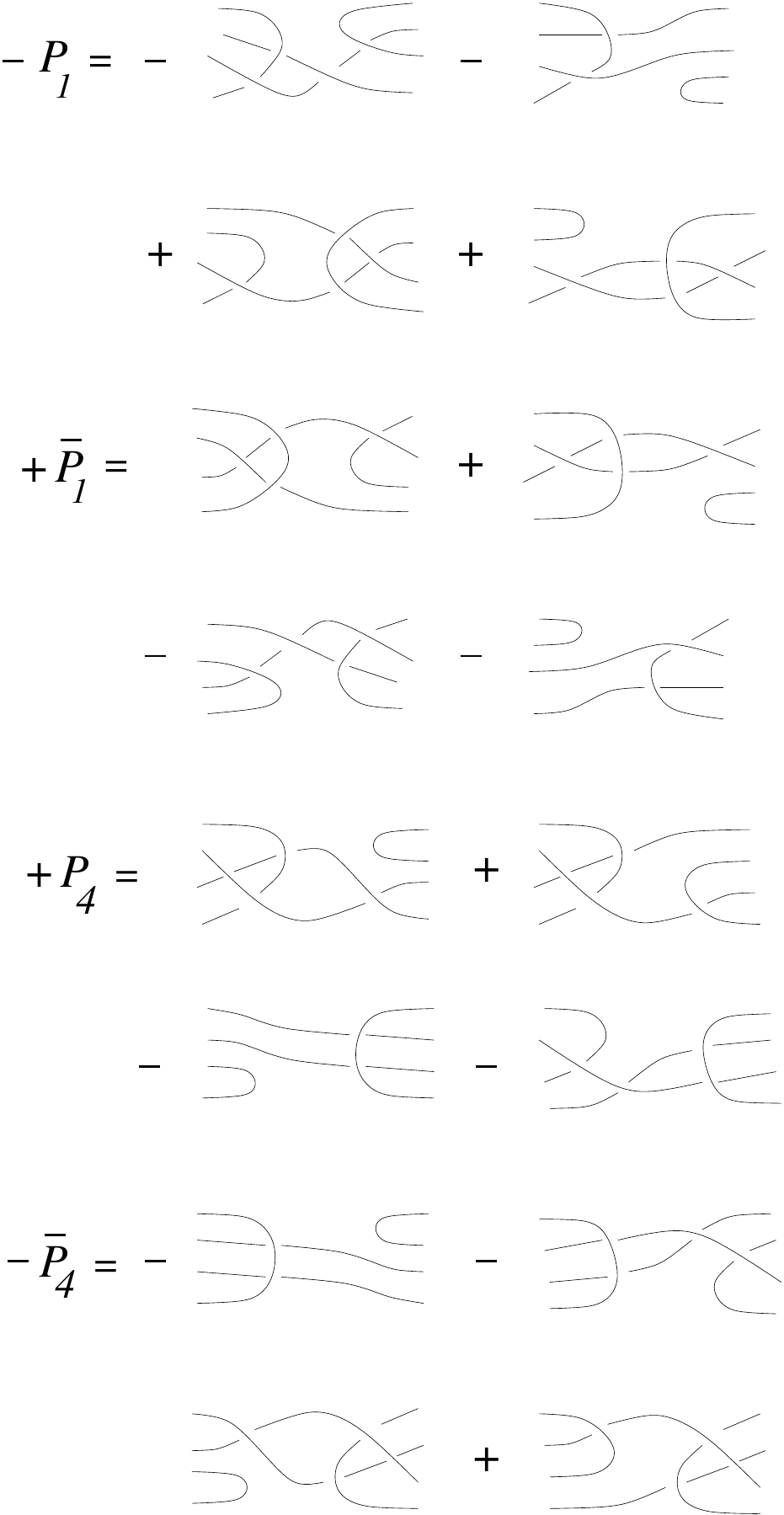}
\caption{\label{contri} contribution of $\sigma_1^{-1}t_2 + \sigma_2t_1 - 
t_1\sigma_2^{-1} - t_2\sigma_1$}  
\end{figure}

\begin{lemma}
$(\sigma_1^{-1}t_2+\sigma_2t_1-t_1\sigma_2^{-1}-t_2\sigma_1)$ is a solution with constant weight of the positive tetrahedron equation.
\end{lemma}
{\em Proof.} The calculations are contained in the figures Fig.~\ref{contri}, Fig.~\ref{calP1}, Fig.~\ref{calP4}.

$\Box$
%%%%%%%%%%%%%%%%%%%%
%%%%%%%%%%%%%%%%%%%%
\begin{figure}
\centering
\includegraphics{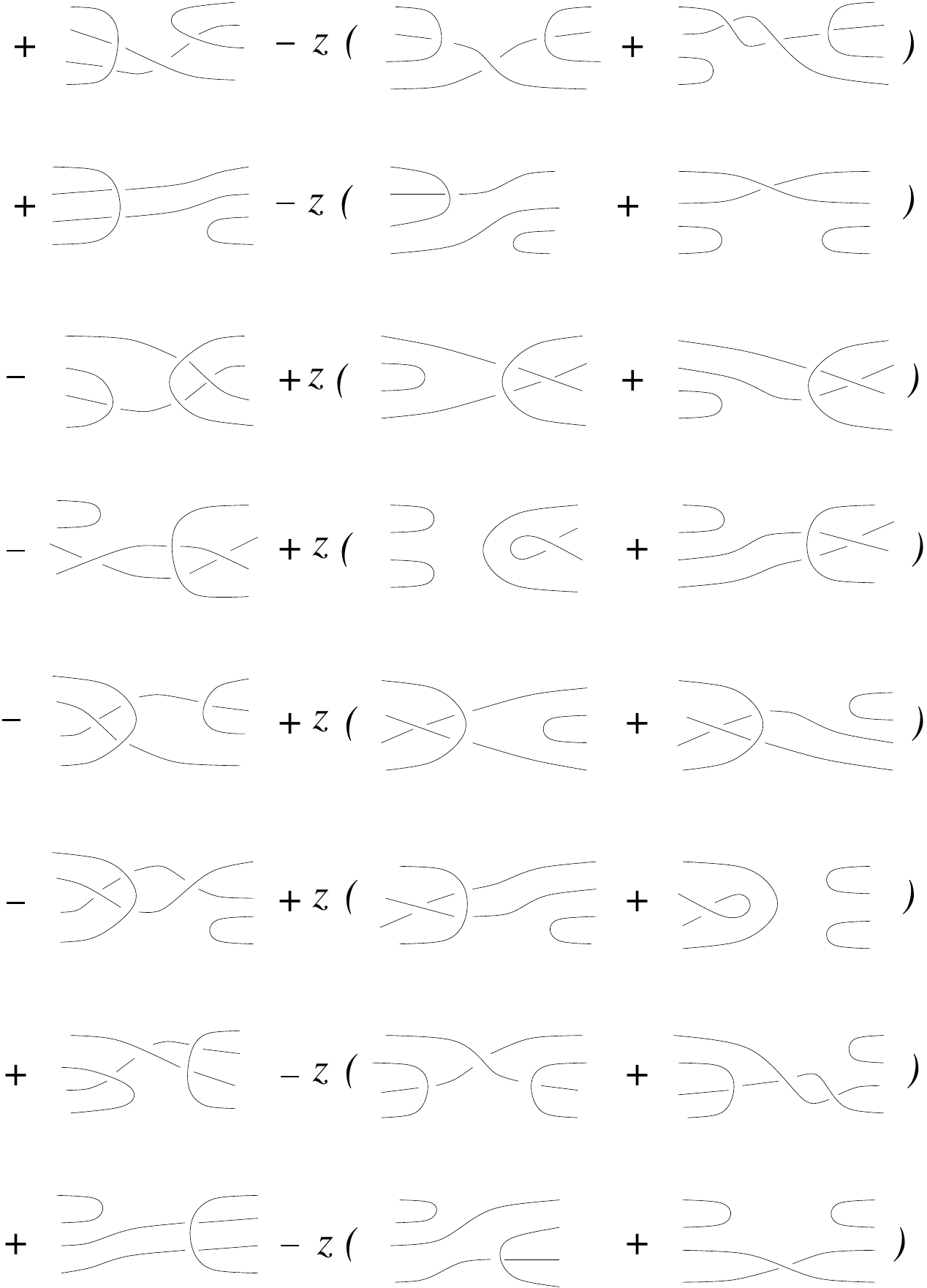}
\caption{\label{calp1} calculation of $\sigma_1^{-1}t_2 + \sigma_2t_1 - 
t_1\sigma_2^{-1} - t_2\sigma_1$ for $-P_1 + \bar P_1$}  
\end{figure}
%%%%%%%%%%%%%%%%%%%%
\begin{figure}
\centering
\includegraphics{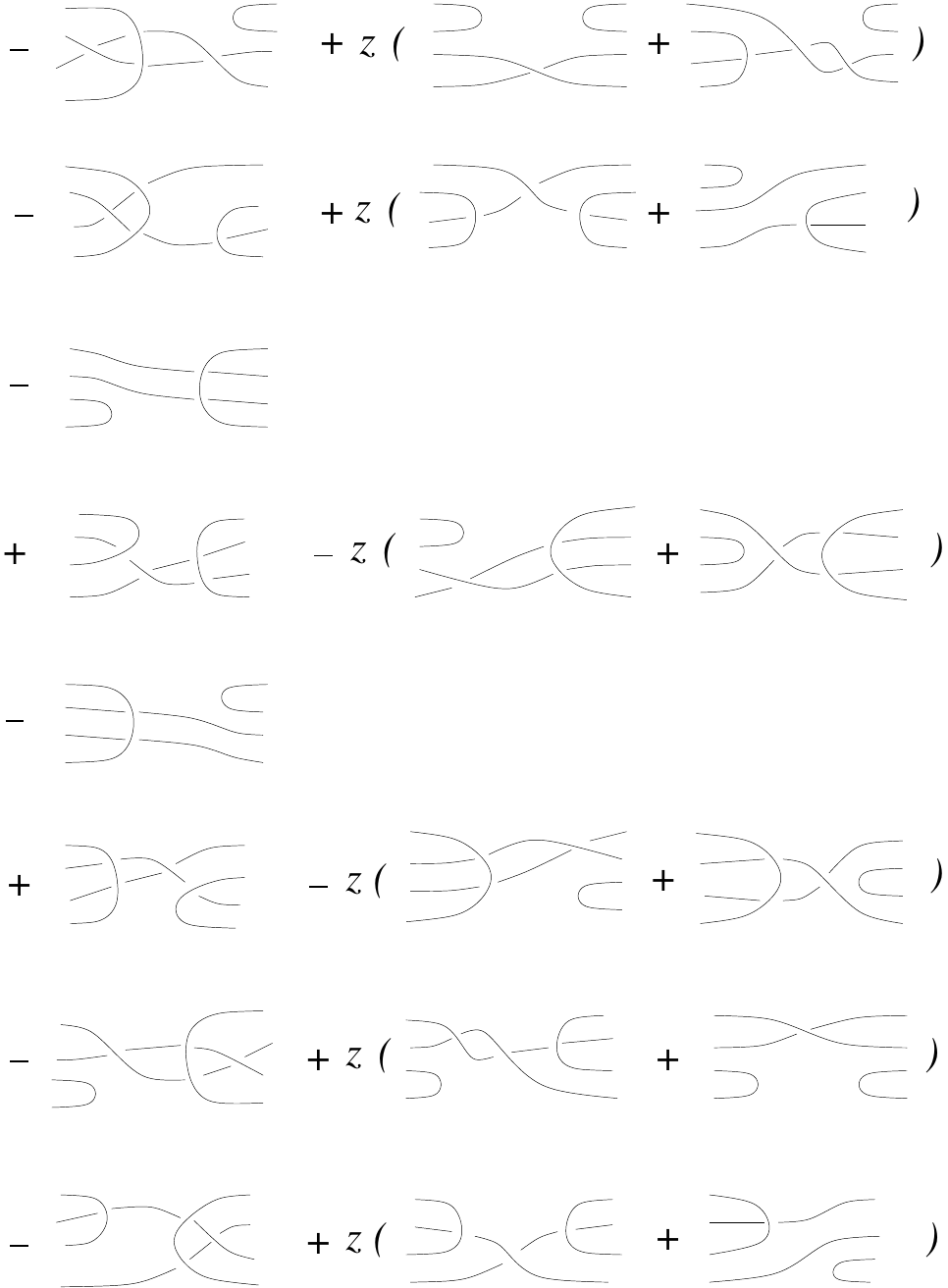}
\caption{\label{calp1} calculation of $\sigma_1^{-1}t_2 + \sigma_2t_1 - 
t_1\sigma_2^{-1} - t_2\sigma_1$ for $P_4 - \bar P_4$}  
\end{figure}

\begin{figure}
\centering
\includegraphics{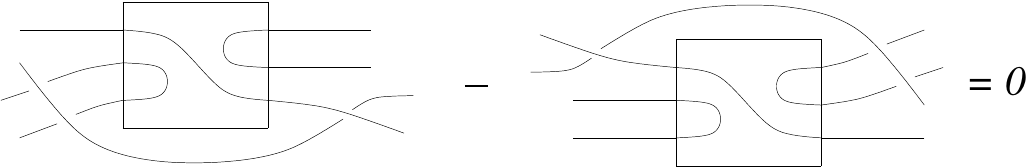}
\caption{\label{KP4} $t_1t_2(p)$ in $P_4 - \bar P_4$}  
\end{figure}

\begin{lemma}
$t_1t_2(p)$ cancels out in $P_4-\bar P_4$.

\end{lemma}
{\em Proof.} The proof is in Fig.~\ref{KP4}.

$\Box$

\begin{figure}
\centering
\includegraphics{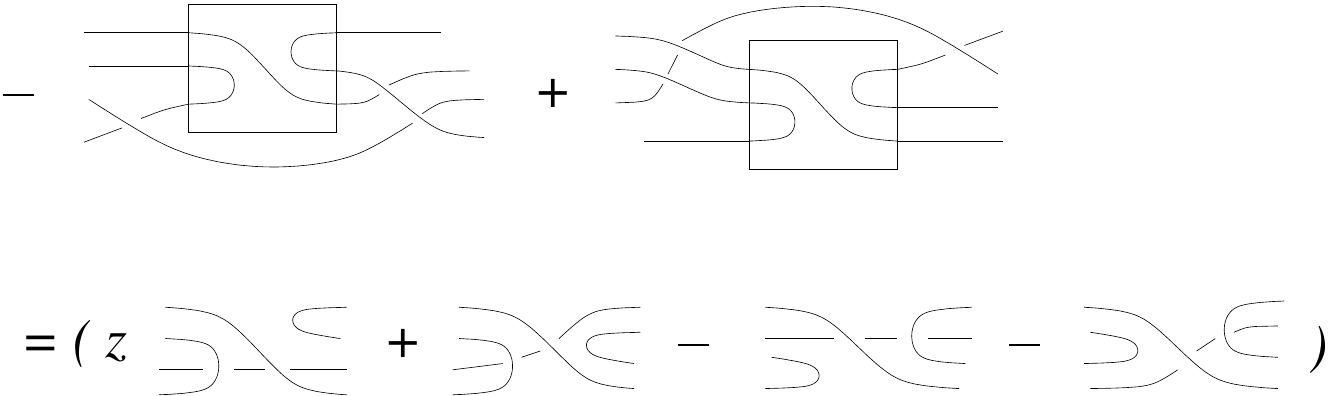}
\caption{\label{KP1} $t_1t_2(p)$ in $-P_1 + \bar P_1$}  
\end{figure}

\begin{figure}
\centering
\includegraphics{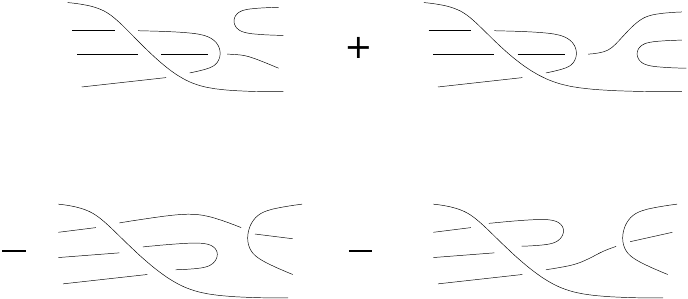}
\caption{\label{KP3} calculation of $\sigma_1^{-1}t_2 + \sigma_2t_1 - 
t_1\sigma_2^{-1} - t_2\sigma_1$ in $P_3$}  
\end{figure}

{\em Proof of Proposition 10.}
The contribution of $-P_1+\bar P_1$ is shown in Fig.~\ref{KP1} and the remaining contribution of  $(\sigma_1^{-1}t_2+\sigma_2t_1-t_1\sigma_2^{-1}-t_2\sigma_1)$ from $P_3$ is shown in Fig.~\ref{KP3}. We see that their combination in $R^{(1)}_{F,reg}(m)(A)$  cancels out.

$\Box$
%%%%%%%%%%%%%%%%%%%%
%%%%%%%%%%%%%%%%%%%%

\begin{remark}
The Kauffman polynomial as an invariant of regular isotopy depends only on the unoriented link in contrast to our 1-cocycle $R^{(1)}_{F,reg}(m)(A)$. Indeed, almost all our definitions use in an essential way the orientation of the circle $T \cup \sigma$. Besides the signs of the Reidemeister moves they are not even invariant under "$flip$" (compare the Introduction).
\end{remark}

\begin{proposition}
The 1-cochain \vspace{0,2 cm}

$R^{(1)}_F(m)=\sum_{p \in m} sign(p)W_1(p) t_1t_2(p) + \sum_{p \in m}sing(p)W_2(p)z(-\sigma_1^{-1}t_2-\sigma_2t_1+t_1\sigma_2^{-1}+t_2\sigma_1)$ 
\vspace{0,2 cm}

is a solution of the positive global tetrahedron equation. Here the first sum is over all triple crossings of the global types shown in Fig.~\ref{fglob} (i.e. the types $l_c$ and $r_a$).  The second sum is over all triple crossings $p$ which have a distinguished crossing $d$ of type 0 (i.e. the types $r_a$, $r_b$ and $l_b$).
\end{proposition}

{\em Proof.} This follows immediately from the proofs of Proposition 10 and Proposition 3.

$\Box$

\section{The 1-cocycles $R^{(1)}_{F, reg}$, $R^{(1)}_F$ and $\bar R^{(1)}_F$}

We proceed in a way completely analogous to the case of the HOMFLYPT invariant. However, we will solve the cube equations only with coefficients in $\mathbb{Z}/2\mathbb{Z}$ (again, one has probably to associate more general partial smoothings to self-tangencies in order to get a solution with integer coefficients). 

\begin{figure}
\centering
\includegraphics{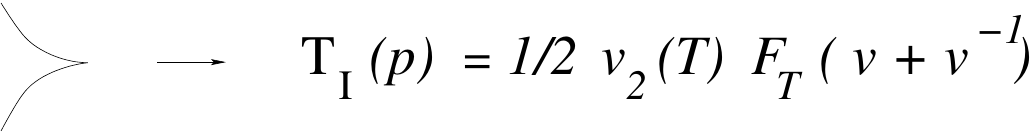}
\caption{\label{KRI} partial smoothing of a $RI$ move}  
\end{figure}

\begin{definition}
For $R^{(1)}_{F, reg}$:

The partial smoothing $T_{II^+_0}(p)$ of a self-tangency with opposite tangent direction and $d$ of type 0 is defined in Fig.~\ref{Ksts}.

The partial smoothing $T_{II^-_0}(p)$ of a self-tangency with equal tangent direction and $d$ of type 0 is defined in Fig.~\ref{Keq}.

\begin{figure}
\centering
\includegraphics{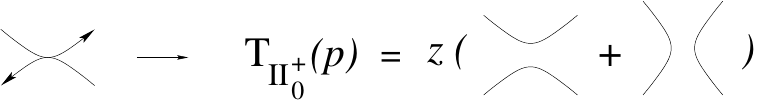}
\caption{\label{Ksts} partial smoothing of a $RII$ move with opposite tangent direction}  
\end{figure}

\begin{figure}
\centering
\includegraphics{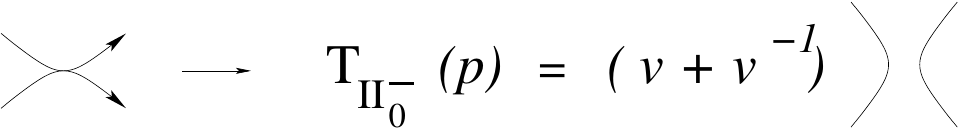}
\caption{\label{Keq} partial smoothing of a $RII$ move with same tangent direction}  
\end{figure}

For $R^{(1)}_F$:

We replace in the above definitions $W(p)$ by $W_2(p)$ and we normalize as usual $F_T$ by $v^{-w(T)}$.

The partial smoothing $T_I(p)$ of a  Reidemeister I move with $d$ of type 0 is defined in Fig.~\ref{KRI}.
\end{definition}

\begin{definition}
Let $p$ be a self-tangency with opposite tangent direction. Its contribution to $R^{(1)}_{F, reg}$ is defined by

 $R^{(1)}_{F, reg}=sign(p) W(p)T_{II^+_0}(p)$.\vspace{0,2 cm}

Its contribution to $R^{(1)}_F$ is defined by

 $R^{(1)}_F=sign(p) W_2(p)T_{II^+_0}(p)$.\vspace{0,2 cm}

Let $p$ be a self-tangency with equal tangent direction. Its contribution to $R^{(1)}_{F, reg}$ is defined by

 $R^{(1)}_{F, reg}=sign(p) W(p)T_{II^-_0}(p)$.\vspace{0,2 cm}

Its contribution to $R^{(1)}_F$ is defined by

 $R^{(1)}_F=sign(p) W_2(p)T_{II^-_0}(p)$.\vspace{0,2 cm}

The contribution of a  Reidemeister I move with $d$ of type 0 is defined by

 $R^{(1)}_F=sign(p) T_I(p)$.

\end{definition}

\begin{lemma}
Let $m$ be the meridian of $\Sigma^{(2)}_{self-flex}$ (compare Section 2). Then $R^{(1)}_{F, reg}(m)=0$ and 
$R^{(1)}_F(m)=0$.
\end{lemma}

\begin{figure}
\centering
\includegraphics{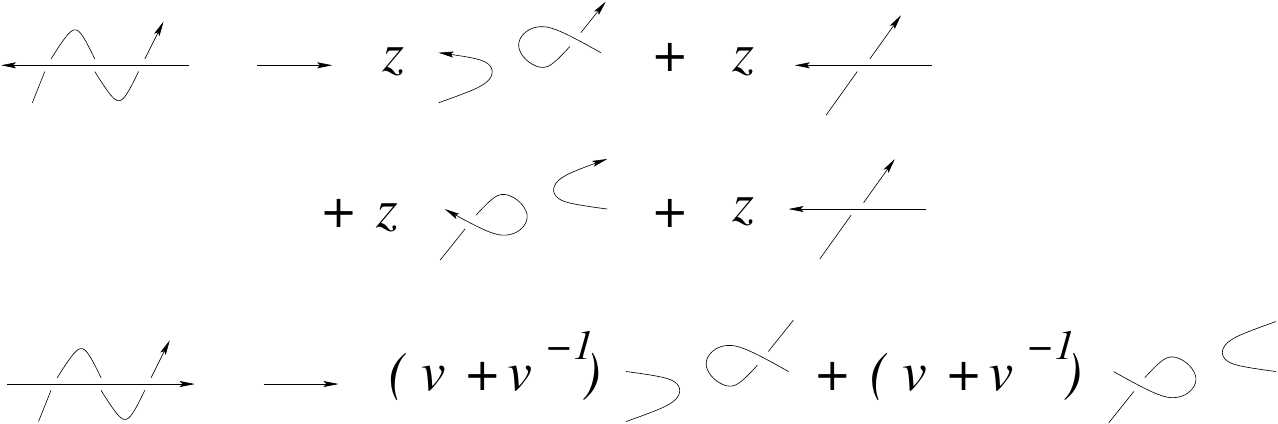}
\caption{\label{Ksf} the meridians of a self-tangency in a flex}  
\end{figure}

{\em Proof.} The weights are exactly the same as in the HOMFLYPT case. We calculate the values on the meridians in Fig.~\ref{Ksf}.

$\Box$

\begin{lemma}
The value of the 1-cocycle  $R^{(1)}_F$  on a meridian of a degenerate cusp, locally given by $x^2=y^5$ and denoted by  $\Sigma^{(2)}_{cusp-deg}$, is zero.

\end{lemma}

\begin{figure}
\centering
\includegraphics{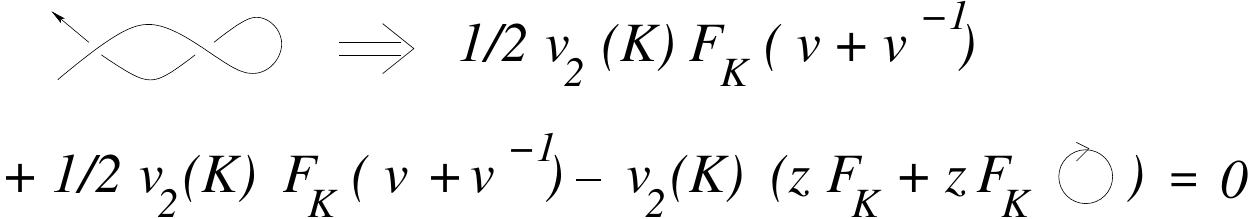}
\caption{\label{Kdeg} $R^{(1)}_F$ for the meridians of a degenerated cusp}  
\end{figure}

{\em Proof.} Again the cusp can be of type $0$ or of type $1$. Only the case of type $0$ is interesting and we give the calculation in Fig.~\ref{Kdeg}.

$\Box$

We will solve now the cube equations.

\begin{figure}
\centering
\includegraphics{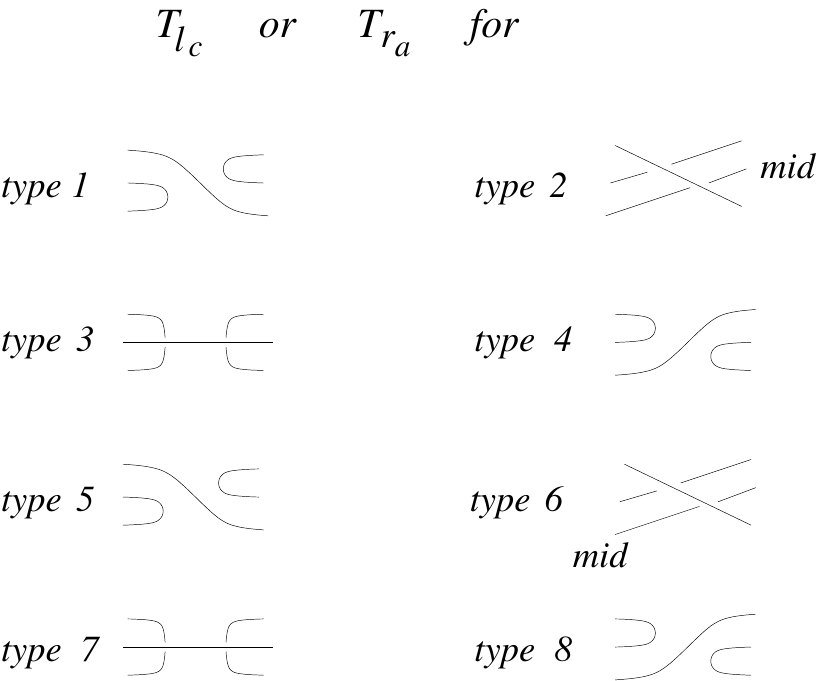}
\caption{\label{Ktypesmooth} the partial smoothings $T_{l_c}$ and $T_{r_a}$ in the
Kauffman case}  
\end{figure}

\begin{figure}
\centering
\includegraphics{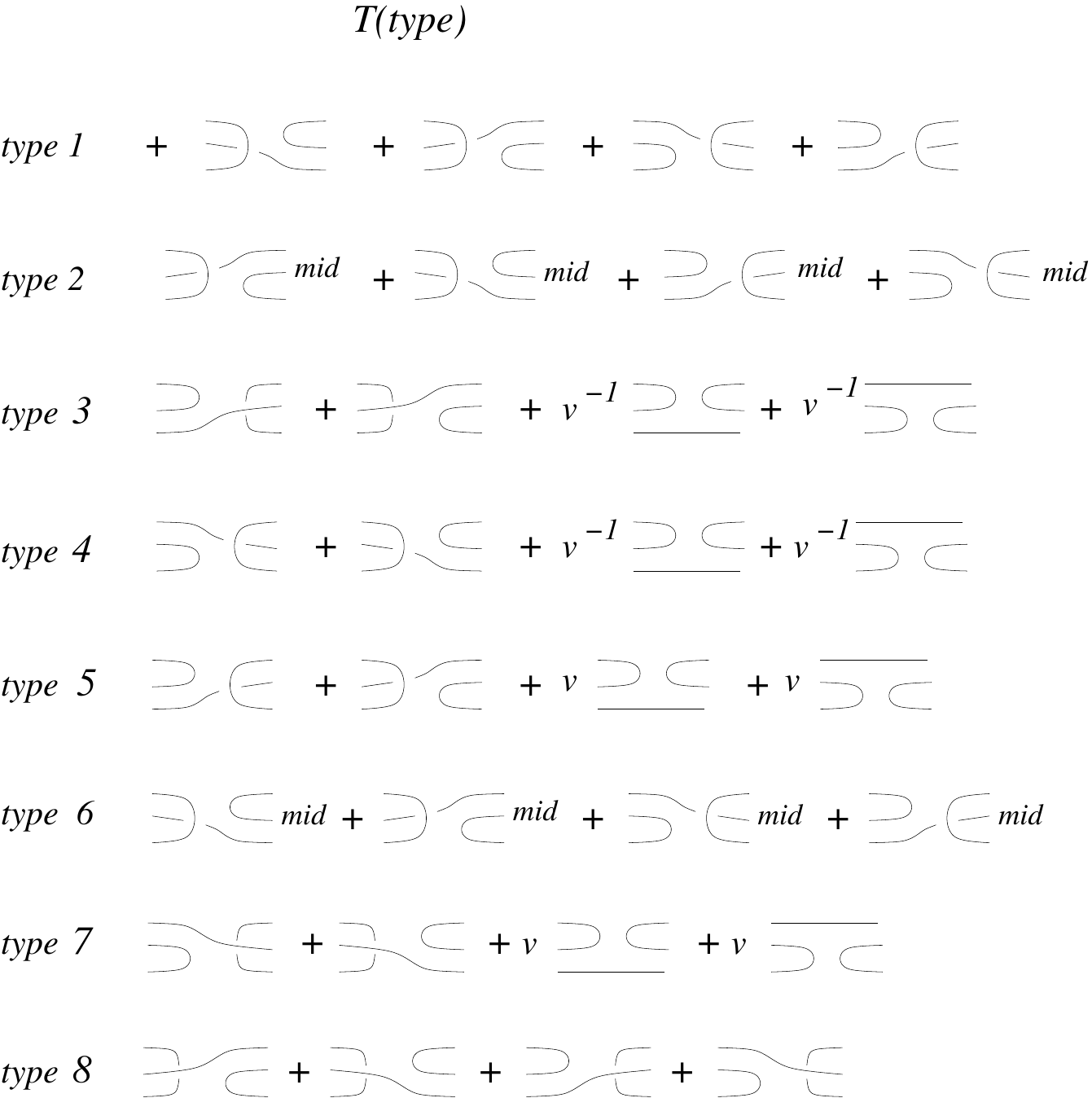}
\caption{\label{Ksmooth} the partial smoothings $T(type)$ in the
Kauffman case}  
\end{figure}

\begin{definition}
The partial smoothings for the local and global types of triple crossings in Kauffman's case are given in Fig.~\ref{Ktypesmooth} and  Fig.~\ref{Ksmooth} (remember that {\em mid} denotes the ingoing middle branch for a star-like triple crossing). 

\end{definition}

\begin{definition}
Let $s$ be a generic oriented arc in $M^{reg}_T$ (with a fixed abstract closure $T \cup \sigma$ to an oriented circle). Let $A \subset \partial T$ be a given grading.  The 
1-cochain $R^{(1)}_{F, reg}(A)$ is defined by \vspace{0,2 cm}

$R^{(1)}_{F,reg}(A)(s ) = \sum_{p \in s \cap l_c} sign(p)T_{l_c} (type)(p) +\sum_{p \in s \cap r_a }sign(p)T_{r_a}(type)(p) +$\vspace{0,2 cm}

 $\sum_{p \in s \cap r_a }sign(p)zW(p)T(type)(p) +$\vspace{0,2 cm}

 $\sum_{p \in s \cap r_b }sign(p)zW(p)T(type)(p) +$\vspace{0,2 cm}

 $\sum_{p \in s \cap l_b }sign(p)zW(p)T(type)(p) +$\vspace{0,2 cm}

 $\sum_{p \in s \cap II^+_0 } sign(p)W(p)T_{II^+_0}(p)$+\vspace{0,2 cm} 

 $\sum_{p \in s \cap II^-_0 } sign(p)W(p)T_{II^-_0}(p)$\vspace{0,2 cm}

Here all weights $W(p)$ are defined only over the f-crossings $f$ with $\partial f= A$ and in the first two sums (i.e. for $T_{l_c}$ and $T_{r_a}$) we require that $\partial (hm)=A$ for the triple crossings.

\end{definition}

\begin{definition}
Let $s$ be a generic oriented arc in $M^{reg}_T$ (with a fixed abstract closure $T \cup \sigma$ to an oriented circle).  The 
1-cochain $R^{(1)}_F$ is defined by \vspace{0,2 cm}

$ R^{(1)}_F(s)=\sum_{p \in s \cap l_c}sign(p)W_1(p)T_{l_c}(type)(p) + $ \\ \vspace{0,1 cm}

$ \sum_{p \in s \cap r_a}sign(p)[W_1(type)(p)T_{r_a}(type)(p)  + zW_2(type)(p)T(type)(p)] +$
\\ \vspace{0,2 cm}

$ \sum_{p \in s \cap r_b}sign(p)zW_2(p)T(type)(p) +$
\\ \vspace{0,2 cm}

$ \sum_{p \in s \cap l_b}sign(p)zW_2(p)T(type)(p) +$
\\ \vspace{0,2 cm}

$  \sum_{p \in s \cap II^+_0}sign(p)W_2(p)T_{II^+_0}(p)+ \sum_{p \in s \cap II^-_0}sign(p)W_2(p)T_{II^-_0}(p)+$ 
\\ \vspace{0,2 cm}

$ \sum_{p \in s \cap I}sign(p)T_I(p)$ \vspace{0,2 cm}

\end{definition}

\begin{proposition}
The 1-cochains $R^{(1)}_{F, reg}(A)$ and $R^{(1)}_F$ with the adjustments from Definition 20 (in Section 7) satisfy the cube equations.

\end{proposition}
{\em Proof.} The weights were already studied for the corresponding proofs in the case of the HOMFLYPT invariant.
We just have to check that our partial smoothings satisfy the cube equations. This is completely analogous to the HOMFLYPT case and we left the verification to the reader (using the corresponding figures in Section 6). But remember that we consider only coefficients in $\mathbb{Z}/2\mathbb{Z}$).

$\Box$

Exactly the same arguments as in the HOMFLYPT case imply that $R^{(1)}_{F, reg}(A)$ and $R^{(1)}_F$ vanish on the meridians of $\Sigma^{(1)} \cap \Sigma^{(1)}$.

\begin{figure}
\centering
\includegraphics{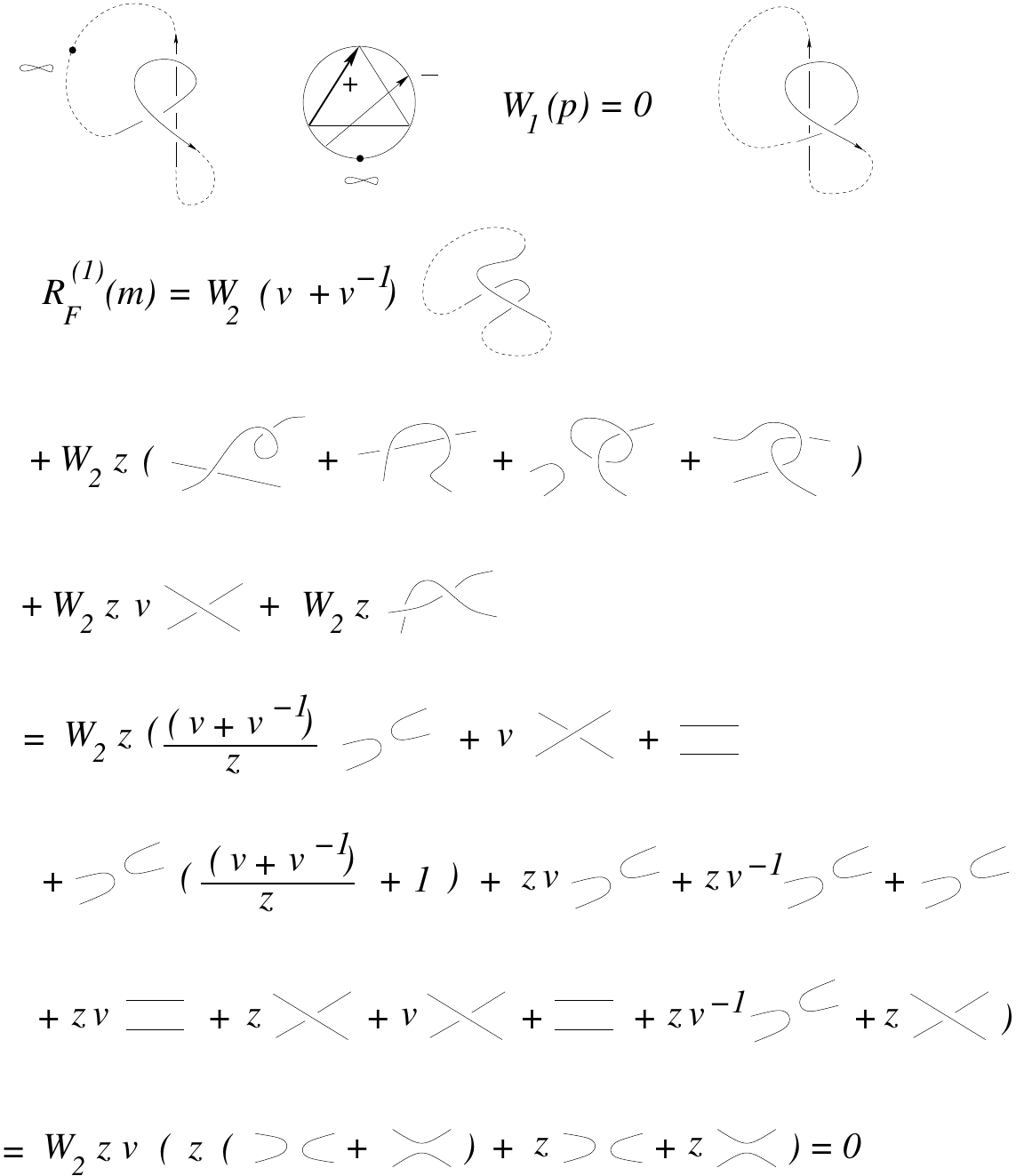}
\caption{\label{Ktranscusp} $R^{(1)}_F$ on a meridian of a cusp with a transverse branch}  
\end{figure}

\begin{lemma}
The value of the 1-cocycle $R^{(1)}_F$ on a meridian of $\Sigma^{(2)}_{trans-cusp}$ is zero besides for $\Sigma^{(2)}_{l_c}$.

\end{lemma}
{\em Proof.} We consider just one case shown in Fig.~\ref{Ktranscusp}. All other cases are similar and we left the verification to the reader.

$\Box$

$R^{(1)}_{F, reg}(A)$ and $R^{(1)}_F$ have both the scan-property for branches moving under the tangle $T$. The proof is exactly the same as in the HOMFLYPT case. Consequently, we have proven the following proposition.

\begin{proposition}
The 1-cochains $R^{(1)}_{F, reg}(A)$ and $R^{(1)}_F$ are 1-cocycles in $M^{reg}_T$ respectively $M_T \setminus \Sigma^{(2)}_{l_c}$. They have both the scan-property.

\end{proposition}

\begin{definition}
The completion $\bar R^{(1)}_F$ $mod$ 2 is defined by $\bar R^{(1)}_F=R^{(1)}_F+F_T V^{(1)}$.
Here $ V^{(1)}$ is the same as in the HOMFLYPT case. 
\end{definition}

\begin{figure}
\centering
\includegraphics{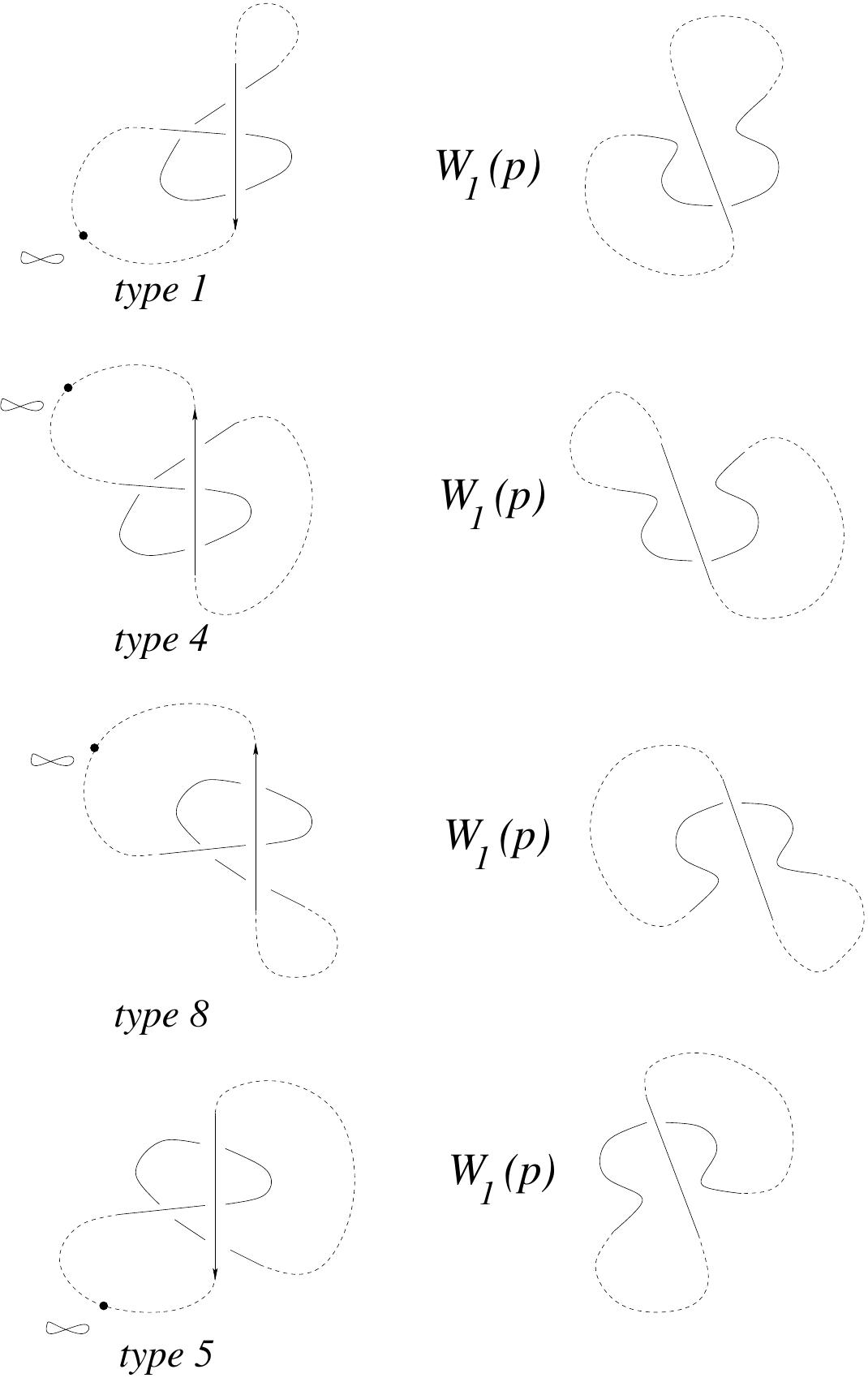}
\caption{\label{F-over-cusp} $R^{(1)}_F$ on the meridians of $\Sigma^{(2)}_{l_c}$}  
\end{figure}

\begin{theorem}
$\bar R^{(1)}_F$ is a 1-cocycle in $M_T$ which represents a non trivial cohomology class and which has the scan-property for branches moving under the tangle $T$.
\end{theorem}
{\em Proof.}
The value of $R^{(1)}_F$ on the meridians of $\Sigma^{(2)}_{l_c}$ is equal to $W_1(p)F_T$ as shown in Fig.~\ref{F-over-cusp}. Consequently, $R^{(1)}_F+F_T V^{(1)}$ is already a 1-cocycle for the whole $M_T$, exactly like in the HOMFLYPT case. $\bar R^{(1)}_F$ has the scan-property because both $R^{(1)}_F$ and $V^{(1)}$ have this property.

$\Box$

We consider just the easiest example in order to show that $\bar R^{(1)}_F$ represents a non trivial cohomology class.

\begin{figure}
\centering
\includegraphics{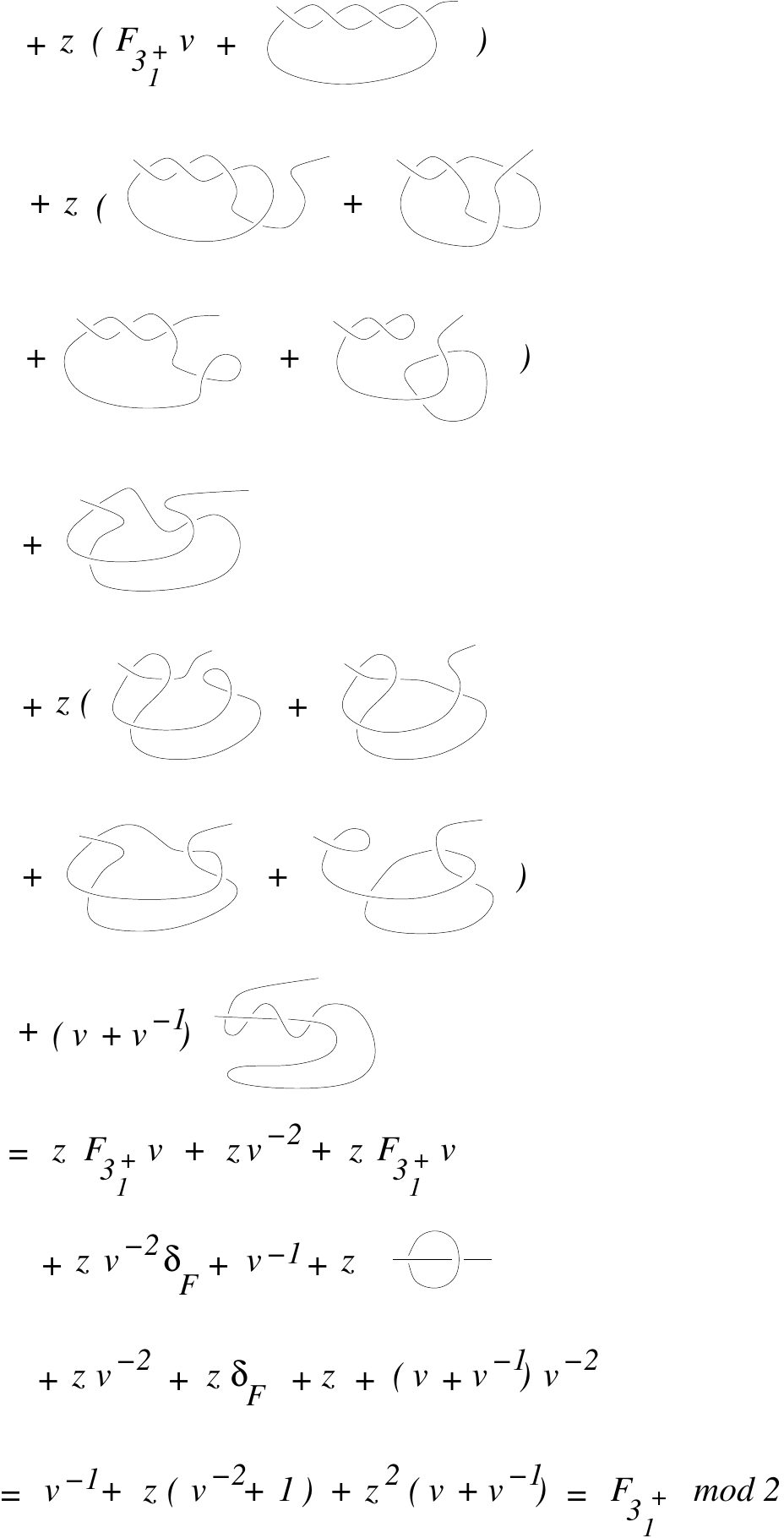}
\caption{\label{Frottrefoil} calculation of $R^{(1)}_F(rot(3^+_1))$}  
\end{figure}

\begin{example}
$\bar R^{(1)}_F(rot(3^+_1))=F_{3^+_1}$. We have only to calculate $\bar R^{(1)}_F(scan(3^+_1)$ because the second half of the rotation does not contribute at all like in the HOMFLYPT case. Moreover, the types and the weights are exactly the same as in the HOMFLYPT case. Only the partial smoothings are different. We give the calculation in Fig.~\ref{Frottrefoil}.
\end{example}

\begin{question}
Do the cube equations have a solution with integer coefficients in Kauffman's case?
\end{question}

\begin{question}
As well known the Jones polynomial is a specialization of both the HOMFLYPT and the Kauffman polynomial (see e.g. \cite{K}).
We can consider the corresponding specializations (with coefficients in $\mathbb{Z}/2\mathbb{Z}$) of our 1-cocycles $\bar R^{(1)}$ and $\bar R^{(1)}_F$. Do they represent the same cohomology class for $M_T$? Do they contain the same information when they are applied to scan-arcs?
\end{question}

\section{Digression: finite type 1-cocycles of degree three for long knots}

It seems that the Teiblum-Turchin 1-cocycle $v^1_3$ and its lift to $\mathbb{R}$ by Sakai are the only other known 1-cocycles for long knots which represent a non trivial cohomology class. The Teiblum-Turchin 1-cocycle is an integer valued 1-cocycle of degree 3 in the sense of Vassiliev's theory. Its reduction mod 2 has a combinatorial description and can be calculated (see \cite{V1}, \cite{V} and \cite{T}). Sakai has defined a $\mathbb{R}$ valued version of the Teiblum-Tourtchine 1-cocycle via configuration space integrals (see \cite{S}). Moreover, he has shown that the value of his 1-cocycle on $rot(K)$ is equal to $v_2(K)$.

Vassiliev has obtained his formula for the reduction mod 2 of the Teiblum-Turchin 1-cocycle $v^1_3$ by combining  the simplicial resolution of his discriminant with the beautiful idea of a discrete moving frame. For the convenience of the reader we recall Vassiliev's formula here.
 
Let $x$ be a constant vector field on the plane which is transverse near infinity to the projection $pr(K)$ of the long knot $K$ into the plane. Let us consider a moving frame along $K$ (the vectors need not to be orthogonal but they form just an oriented basis of the 3-space). The first vector is oriented tangential to the knot, the second vector is constant orthogonal to the plane and the third vector is determined by the orientation of the ambient space. There are exactly two new types of strata of codimension one: namely diagrams with an ordinary crossing for which the projection in the plan of one of the oriented tangent directions of the knot coincides with $x$, and diagrams for which the projection is oriented tangential to $x$ in an ordinary flex.
Moreover, the moving frame degenerates exactly for Reidemeister I moves.  

\begin{figure}
\centering
\includegraphics{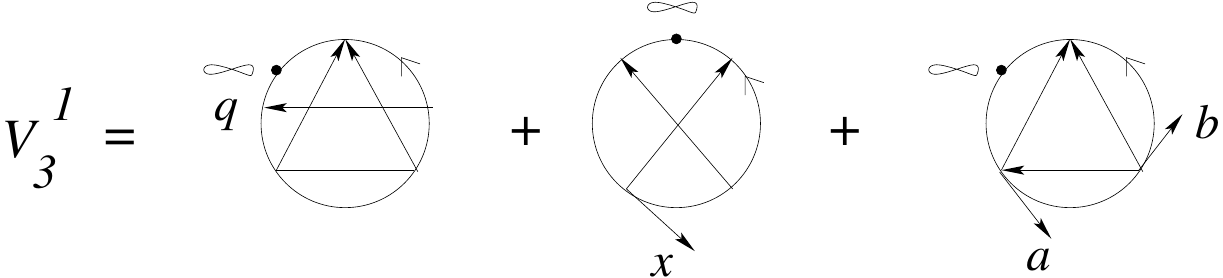}
\caption{\label{vasv3} Vassiliev's formula for the Teilblum-Turchin 1-cocycle mod 2}  
\end{figure}

Using our language of Gauss diagram formulas Vassiliev's formula is given in 
Fig.~\ref{vasv3}. The first term is of Gauss degree 3 (see the definition below). In the second term the projection of the undercross has to be oriented tangential to $x$. In the third term it is required that the sector spanned by the projections of the tangent vectors $a$ and $b$ does not contain the vector $x$.

Using our figures from Section 3 we see immediately that the first term in Vassiliev's formula satisfies the positive global tetrahedron equation but only mod 2. Moreover, it has the scan-property for branches moving over everything, i.e. the contributions of $P_3$ and $\bar P_3$ cancel always out (there is of course a dual formula for which it has the usual scan-property). The third part has obviously the scan-property. It turns out that the second part has also the scan-property for branches moving over everything. Indeed, Reidemeister II moves don't appear in the formula and hence it is sufficient to prove it only for positive Reidemeister III moves (i.e. of type 1). Consequently, we have only to study the changing of the second part in Vassiliev's formula when the undercross of crossing 3 in Fig.~\ref{III1typer} and in Fig.~\ref{III2typel} is tangential to $x$. The only interesting case is $\infty=a$ in 
Fig.~\ref{III1typer}. Indeed, we see that it doesn't change. 

{\em It follows that Vassiliev's formula has the scan-property for branches moving over the knot.}

Turchin \cite{T} has shown that $v^1_3(rot(K))=v_2(K)$ mod 2 for each long knot $K$ and he has conjectured that the equality is still true over the integers. In view of Conjecture 1 it seems likely that for $rot(K)$ the value of the Teiblum-Turchin 1-cocycle of finite type coincides with the value of our quantum 1-cocycle by sending $P_K$ to $-1/\delta$.

It is well known that each finite type knot invariant can be represented by a Gauss diagram formula (see \cite{GPV}). This is not known for finite type 1-cocycles. It is therefore natural to distinguish the degree in Vassiliev's sense of a 1-cocycle from the degree of a 1-cocycle given by a Gauss diagram formula.

\begin{definition}

The {\em Gauss degree} of a 1-cocycle in $M_K$, $M^{reg}_K$ or in $M_K \setminus \bar \Sigma^{(2)}_{trans-cusp}$ is defined as the minimum over all Gauss diagram formulas for it of the the maximal number of arrows minus one amongst all configurations in the Gauss diagram formula. 
\end{definition}

Vassiliev's formula for the Teiblum-Turchin 1-cocycle 
$v^1_3$ $mod$ 2 is  of Gauss degree 3. Does the Gauss degree always coincide with the degree in Vassiliev's sense, when the latter is well defined ? (See \cite{F2} for examples of 1-cocycles for closed braids of arbitrary Gauss degree.) \vspace{0,2 cm}

We use the occasion (the techniques developed in this paper)  to introduce four integer valued 1-cocycles of Gauss degree 3 in 
$M_K \setminus \bar \Sigma^{(2)}_{trans-cusp}$.

\begin{figure}
\centering
\includegraphics{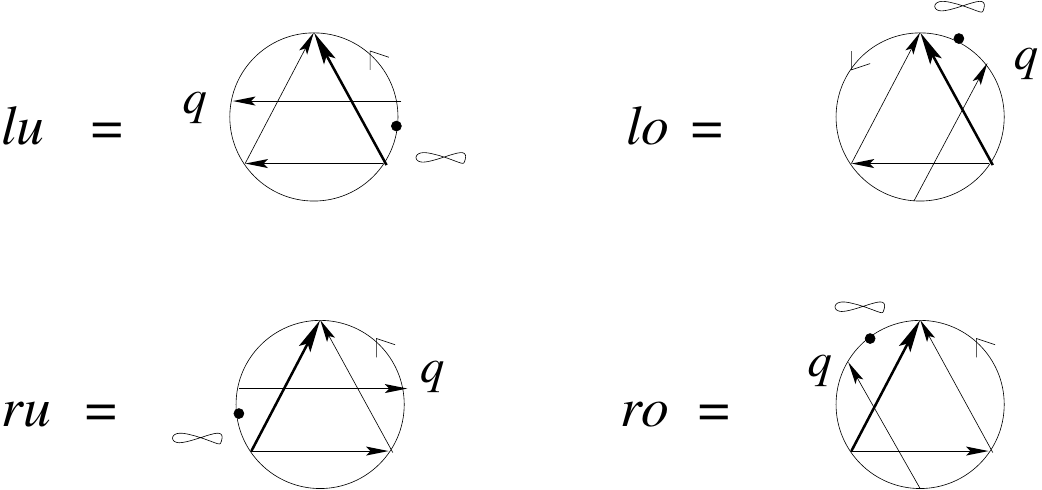}
\caption{\label{vtri} 1-cocycles of Gauss degree $3$}  
\end{figure}

\begin{definition}
Let $s$ be a generic oriented loop in $M_K \setminus \bar \Sigma^{(2)}_{trans-cusp}$. The four 1-cochaines $v^{(1)}_{lu}$, $v^{(1)}_{lo}$, $v^{(1)}_{ru}$, $v^{(1)}_{ro}$ are defined by \vspace{0,2 cm}

$\sum_{p \in s \cap \Sigma^{(1)}_{tri}}sign(p) \sum w(q)$ \vspace{0,2 cm}

where the second sum is only over the corresponding configurations shown in 
Fig.~\ref{vtri}. The sign is defined by using the global coorientation from  
Fig.~\ref{homoco}.

\end{definition}

(The letter "l" stands for the arrow $ml$ goes to the left and the letter "u" stands for $\infty$ is under the arrow $q$.)

\begin{proposition}
Each of the 1-cochains $v^{(1)}_{lu}$, $v^{(1)}_{lo}$, $v^{(1)}_{ru}$, $v^{(1)}_{ro}$ is a 1-cocycle of Gauss  degree 3 in $M_K \setminus \bar\Sigma^{(2)}_{trans-cusp}$.
\end{proposition}
{\em Proof.} Obviously each of the 1-cochains vanishes on the meridians of 

$\Sigma^{(1)} \cap \Sigma^{(1)}$ (the crossing $q$ can not be the new crossing from a Reidemeister I move) and on the meridians of the degenerate cusp and of the self-tangency in a flex.

Inspecting our figures in Section 3 we see that each of the 1-cochains satisfies the global positive tetrahedron equation (but non of them has a scan-property). 

For all four 1-cochains the crossing $q$ is {\em rigid}, i.e. the arrow can not slide into an arrow of the triangle without moving over the point at infinity (compare \cite{F2} for the case of closed braids). This implies immediately that the 1-cochains satisfy the cube equations.
It follows that the 1-cochains are 1-cocycles in the complement (of the closure) of the remaining strata of codimension 2 (compare Section 2), i.e. all strata which correspond to a cusp with a transverse branch.

Each of the 1-cocycles is defined by a unique configuration and which contains four arrows. Hence they are of Gauss  degree 3. It is not difficult to prove that we can not represent them by configurations with fewer arrows. But we left this to the reader.

$\Box$

It is easy to see that each of the four 1-cocycles is non zero on some meridian of $\Sigma^{(2)}_{trans-cusp}$. Therefore the values of the 1-cocycles on $rot(K)$ depend now on the choice of the loop for $rot(K)$ in the complement of $\bar \Sigma^{(2)}_{trans-cusp}$.

Let $rot(K)$ be represented by the loop which was introduced in 

Example 3.

\begin{proposition}
For any long knot $K$ we have\vspace{0,2 cm}

$v^{(1)}_{lu}(rot (K))=-v_2(K)$, $v^{(1)}_{lo}(rot (K))=v_2(K)$, \vspace{0,2 cm}

$v^{(1)}_{ru}(rot (K)) =v^{(1)}_{ro}(rot (K))=0$.\vspace{0,2 cm}

If we replace the positive curl by a positive curl with opposite Whitney index then $v^{(1)}_{l}(rot (K))$ interchanges with $v^{(1)}_{r}(rot (K))$.

\end{proposition}

The proof is based on a careful analysis of the changing of the values of the $v^{(1)}$ under a crossing change for $K$. We leave the details to the reader.

The proposition implies for example that the loop $rot(K)$ defined by using a curl with writhe $+1$ and Whitney index $-1$ and the loop $rot(K)$ defined by using a curl with writhe $+1$ and Whitney index $+1$ are not homologous in 
$M_K \setminus \bar \Sigma^{(2)}_{trans-cusp}$ if $v_2(K)$ is non zero (but they are of course always homologous in $M_K$).

\begin{question}
Is it possible to complete any of the above four 1-cocycles $v^{(1)}$ to a formula for the integer Teiblum-Turchin 1-cocycle by using 
Vassiliev's methode of a discrete moving frame?

\end{question}

Institute de Math\'ematiques de Toulouse

Universit\'e Paul Sabatier

118, route de Narbonne 

31062 Toulouse Cedex 09, France

fiedler@math.univ-toulouse.fr

\end{document}